\documentclass[11pt]{amsart}
\usepackage[
letterpaper,
left=0.8in,
right=0.8in,
top=1in,
bottom=0.75in
]{geometry}
\usepackage{amsmath,amssymb,amsthm,mathtools}
\usepackage{tikz,tikz-cd}
\usetikzlibrary{arrows.meta,decorations.markings,decorations.pathmorphing,decorations.pathreplacing,calc,positioning,shapes.geometric,fit,backgrounds,patterns}
\usepackage{enumitem,booktabs,array,longtable}
\usepackage{xcolor}
\usepackage[colorlinks=true,linkcolor=blue!60!black,citecolor=blue!60!black,urlcolor=blue!60!black]{hyperref}

\newtheorem{theorem}{Theorem}[section]
\newtheorem{proposition}[theorem]{Proposition}
\newtheorem{lemma}[theorem]{Lemma}
\newtheorem{corollary}[theorem]{Corollary}
\theoremstyle{definition}
\newtheorem{definition}[theorem]{Definition}

\newtheorem{example}[theorem]{Example}

\theoremstyle{remark}
\newtheorem{remark}[theorem]{Remark}
\numberwithin{equation}{section}
\newcommand{\Z}{\mathbb{Z}}
\newcommand{\N}{\mathbb{N}}
\newcommand{\Q}{\mathbb{Q}}
\newcommand{\R}{\mathbb{R}}
\newcommand{\C}{\mathbb{C}}
\newcommand{\T}{\mathbb{T}}
\newcommand{\Sph}{\mathbb{S}}
\newcommand{\EE}{\mathbb{E}}
\newcommand{\JS}{\mathcal{Z}}
\newcommand{\Nor}{\mathcal{N}}
\newcommand{\Kt}{\widetilde{K}}
\newcommand{\KU}{KU}

\newcommand{\Aff}{\operatorname{Aff}}
\newcommand{\Tor}{\operatorname{Tor}}
\newcommand{\tc}{\operatorname{tc}}
\newcommand{\diag}{\operatorname{diag}}
\newcommand{\Ad}{\operatorname{Ad}}
\newcommand{\End}{\operatorname{End}}
\newcommand{\ev}{\operatorname{ev}}
\newcommand{\id}{\mathrm{id}}
\newcommand{\ib}{\bar\imath}
\newcommand{\gt}{\widetilde\gamma}
\newcommand{\cc}{\mathfrak{c}}
\newcommand{\lcm}{\operatorname{lcm}}
\newcommand{\lims}{\varinjlim}
\newcommand{\limp}{\varprojlim}
\newcommand{\perm}{\mathrm{perm}}
\newcommand{\pt}{\mathrm{pt}}
\newcommand{\tr}{\mathrm{tr}}
\newcommand{\Tr}{\operatorname{Tr}}
\newcommand{\rank}{\operatorname{rank}}
\newcommand{\supp}{\operatorname{supp}}
\DeclareMathOperator{\hofib}{hofib}
\DeclareMathOperator{\Hom}{Hom}

\newcommand{\K}{\mathcal{K}}

\newcommand{\Ped}{\operatorname{Ped}}

\newcommand{\Ell}{\operatorname{Ell}}
\newcommand{\dd}{\mathfrak{d}}
\newcommand{\Tt}{\widetilde{T}}
\newcommand{\bb}{\mathbf{b}}
\newcommand{\cv}{\mathbf{c}}
\newcommand{\fw}{\mathrm{fw}}
\newcommand{\bw}{\mathrm{bw}}
\newcommand{\var}{\mathrm{var}}
\newcommand{\her}{\operatorname{her}}
\newcommand{\ret}{\operatorname{ret}}

\title[Untwisted groupoid models for classifiable $C^*$-algebras]{Every classifiable stably finite C$^*$-algebra has an untwisted \'etale groupoid model}
\date{\today}


\begin{document}
	
	\author{Nikolay A. Ivanov}
	\thanks{The first author is partially supported by the Bulgarian National Science Fund under Grant No.~KP-06-H92/6 (December 08, 2025) and by Grant No 80-10-31/01.04.2026 of the Research Fund of Sofia University.}
	\address  {Nikolay A. Ivanov, Faculty of Mathematics and Informatics, Sofia University ``St.~Kliment Ohridski'', 5 James Bourchier Blvd., Sofia, Bulgaria} \email{nivanov@math.uni-sofia.bg}\

	\author{Ali I. Raad}
	\address  {Ali I. Raad, Mathematics and Science Department, American University in Bulgaria, Office 307 BAC, ul. Svoboda Bacharova 8, 2700 Blagoevgrad, Bulgaria} \email{araad@aubg.edu}\

	\begin{abstract}
		We show that every classifiable stably finite $C^*$-algebra is isomorphic to the
		$C^*$-algebra of a second countable, locally compact Hausdorff, principal étale groupoid. Equivalently, every such $C^*$-algebra has a $C^*$-diagonal whose Weyl twist is trivial. Together with Evington and Sibbel's recent result for all UCT
		Kirchberg algebras, this provides untwisted principal groupoid models for every classifiable $C^*$-algebra.
	\end{abstract}
	
	\maketitle
	\tableofcontents
	
	%======================================================================
	\section{Introduction}\label{sec:intro}
	
	\subsection{Context and motivation}
	An interesting and important question in the theory of $C^*$-algebras is: when can a $C^*$-algebra $A$ be realized as a groupoid $C^*$-algebra? A topological groupoid \(G\) is a small category in which every morphism is invertible, the sets of objects and morphisms are equipped with topologies, and all the category's structure maps (source, range, composition, and inversion) are continuous. They capture various interesting topological structures, such as dynamical systems, group actions, topological equivalence relations, foliations, and various structures in differential geometry, to name a few. From a topological groupoid $G$ with certain properties one can build a (reduced) groupoid $C^*$-algebra $C^*_r(G)$. Various important classes of $C^*$-algebras are groupoid $C^*$-algebras, for instance graph $C^*$-algebras, Cuntz-Krieger algebras, AF-algebras, crossed products $C_0(X)\rtimes_r\Gamma\cong C^*_r(\Gamma\ltimes X)$ of group actions on locally compact spaces, to name a few. Whenever a $C^*$-algebra $A$ is isomorphic to a groupoid $C^*$-algebra $C^*_r(G)$ we say that $A$ has a \emph{groupoid model}. Groupoid models make the dynamical and geometric structure of a $C^*$-algebra visible. Various structural properties of a groupoid $C^*$-algebra can be understood at the groupoid level, such as the structure of ideals, nuclearity, and separability, among others. 
	
	Groupoid $C^*$-algebras are closely tied to Cartan subalgebras:
	by Renault's theorem \cite[Theorem~5.9]{Re08} (building on Kumjian \cite{Ku}), a separable $C^*$-algebra $A$ with a Cartan subalgebra $B$ is
	canonically of the form $C^*_r(G,\Sigma)$, where $G$ is an effective étale groupoid, the \emph{Weyl groupoid} of $(A,B)$, and $\Sigma$ is a \emph{twist}
	over $G$, the \emph{Weyl twist}. Roughly speaking, a twist over a topological groupoid $G$ can be interpreted as an obstruction to making a globally consistent choice of phases compatible with groupoid multiplication. More precisely, a twist $\Sigma$ is a central extension of $G$ by the circle group $\mathbb{T}$, $\mathbb{T}\times G^0\hookrightarrow \Sigma \twoheadrightarrow G$, and its nontriviality reflects the obstruction to choosing a continuous multiplicative section of this extension. When the twist admits a continuous global section, it can be described by a continuous $\mathbb{T}$-valued (2)-cocycle on $G$, which measures the failure of the chosen section to preserve multiplication. Changing the section modifies the cocycle by a coboundary, so its cohomology class captures the obstruction to removing this phase ambiguity. In particular, the vanishing of this class is equivalent to the existence of a continuous multiplicative section, and hence to the triviality of the twist as a groupoid extension.
	
	Li proved that every classifiable $C^*$-algebra has a Cartan subalgebra \cite{Li}. More precisely, he constructed, for every possible
	value of the classifying invariant, in the unital and in the stably projectionless case, a model which is an inductive limit of building blocks
	whose Cartan pairs are preserved by the connecting maps \cite[Theorems~1.2, 1.3, Corollaries~1.6, 1.7]{Li}; the limit Cartan pair then has a
	principal Weyl groupoid, but possibly a nontrivial twist. In Li's construction the twist is only needed to accommodate torsion in $K_0$: when
	$K_0$ is torsion-free, his stably finite models have trivial twist \cite[Corollary~1.8(i)]{Li}. Whether the twist can always be removed, i.e.\
	whether every classifiable $C^*$-algebra is the $C^*$-algebra of an \emph{untwisted}  étale groupoid, is one of the problems
	collected by Schafhauser, Tikuisis and White \cite[Problem~XLIII,~\S13]{STW}. Untwisted models were known for the Jiang--Su algebra $\JS$ (Deeley--Putnam--Strung
	\cite{DPS18}, Austin--Mitra \cite{AM}, Li \cite[Corollary~1.9]{Li}) and for the Razak--Jacelon algebra $\mathcal W$ \cite{AM}, \cite[Corollary~1.9]{Li};
	for classes of unital $C^*$-algebras with finite-dimensional trace simplex arising from minimal $\Z$-actions, with quite a lot of torsion allowed in $K_0$, and in the
	real rank zero case \cite{DPS24}; Ala, Liao and Tikuisis constructed étale equivalence relations with certain prescribed torsion in their zeroth homology
	\cite{ALT}. Evington and Sibbel have now settled the classifiable purely infinite case \cite[Theorem~A]{ES}. In this paper we treat all stably finite classifiable
	$C^*$-algebras; together these results settle the problem for every classifiable $C^*$-algebra.

	Recall that a $C^*$-algebra $A$ is \emph{stably projectionless} if $A\otimes\K$ has no nonzero projection. Throughout, a $C^*$-algebra is called
	\emph{classifiable} if it is separable, simple, nuclear, $\JS$-stable and satisfies the UCT; here $\JS$ is the Jiang--Su algebra \cite{JS}. For
	separable, simple, nuclear $C^*$-algebras which are not elementary (that is, not isomorphic to the algebra of compact operators on a finite- or
	infinite-dimensional Hilbert space), $\JS$-stability is equivalent to finite nuclear dimension \cite{WZ}, \cite{Ti}, \cite{CE}, \cite{CETWW}.
	Elementary $C^*$-algebras have to be excluded here: $M_k$ and $\K$ have finite nuclear dimension but are not $\JS$-stable.
	
	\begin{theorem}[Main theorem]\label{thm:main}
		Let $A$ be a classifiable, stably finite C$^*$-algebra. There is a second countable, locally compact Hausdorff, principal étale groupoid $G$ such that
		\[
		A\cong C^*(G).
		\]
		In particular, $A$ has a $C^*$-diagonal whose Weyl twist is trivial.
	\end{theorem}
	
	A classifiable stably finite $C^*$-algebra is either unital, or stably projectionless, or non-unital with a nonzero projection in $A\otimes\K$. The unital
	case is Theorem~\ref{thm:unital}, proved in Section~\ref{sec:proof}; the stably projectionless case is Theorem~\ref{thm:SPA}, proved in
	Section~\ref{sec:SPscale}; and the third case is Theorem~\ref{thm:proj}, proved in Section~\ref{sec:proj}. Finite-dimensional simple $C^*$-algebras $M_k$ have
	the obvious untwisted model, the pair groupoid on $k$ points (Example~\ref{ex:pairgroupoid}).
	
	In the stably projectionless case the heart of the proof is an untwisted version of the stably projectionless range theorem. For a classifiable
	stably projectionless $C^*$-algebra with \emph{continuous scale} (Definition~\ref{def:cs}) the Elliott invariant reduces to a quadruple $(K_0,T,\rho,K_1)$,
	where $T$ is the (compact) tracial state space and $\rho\colon K_0\to\Aff(T)$ the pairing; the quadruples which occur are exactly those for which
	$\rho(K_0)$ contains no strictly positive function (Li \cite[Theorem~1.3]{Li}, Gong--Lin \cite[Theorem~14.8]{GL4}; see Theorem~\ref{thm:P} for the
	necessity).
	
	\begin{theorem}[Untwisted range theorem, stably projectionless case]\label{thm:SPB}
		Let $G_0$ and $G_1$ be countable abelian groups, $T$ a nonempty metrizable Choquet simplex and $\rho\colon G_0\to\Aff(T)$ a homomorphism such that
		no $\rho(g)$ is strictly positive on $T$. There is a second countable, locally compact Hausdorff, principal étale groupoid $G_c$ such that
		$C=C^*(G_c)$ is a classifiable stably projectionless $C^*$-algebra with continuous scale and
		\[
		\big(K_0(C),T(C),\rho_C,K_1(C)\big)\cong(G_0,T,\rho,G_1).
		\]
		Moreover every densely finite lower semicontinuous trace on $C$ is bounded.
	\end{theorem}
	
	Combining the recent result of Evington and Sibbel that every UCT Kirchberg algebra has an untwisted principal étale groupoid model \cite[Theorem~A]{ES} with Theorem~\ref{thm:main} gives the conclusion for all classifiable $C^*$-algebras.
	
	\begin{corollary}[All classifiable $C^*$-algebras]\label{cor:allclassifiable}
		Let $A$ be a classifiable $C^*$-algebra. There is a second countable, locally compact Hausdorff, principal étale groupoid $G$ such that
		$A\cong C^*(G)$. In particular, $A$ has a $C^*$-diagonal whose Weyl twist is trivial.
	\end{corollary}

	\subsection{The method of proof}
	In Li's stably finite constructions the torsion subgroup of $K_0$ is carried by the reduced $K$-theory of finite complexes, and the building blocks
	contain homogeneous algebras of the form $P M_N(C(Z))P$, where $P$ is a sum of projections representing line bundles. The natural diagonal of such an
	algebra consists of the corners $p_i M_N(C(Z))p_i$, and its Weyl twist is assembled from the circle bundles of the line bundles $\Hom(p_i,p_j)$, which
	are in general nontrivial (Remark~\ref{ex:hopf}) leading to nontrivial twists.
	
	The idea of this paper is to replace these vector bundles by \emph{permutation bundles}. A finite covering $\pi\colon Y\to X$ gives a vector bundle
	$V_\pi=\pi_*\C$ whose fibre over $x$ is the vector space with the canonical basis $\pi^{-1}(x)$. Its transition functions are permutation matrices, and
	consequently
	\[
	\Gamma(\End(\pi_*\C))\cong C^*(Y\times_XY)
	\]
	is the $C^*$-algebra of an equivalence relation, with its natural diagonal $C(Y)$; it carries no twist (Proposition~\ref{prop:coveringalgebra}). The
	difficulty is to arrange enough permutation bundles, maps between the underlying spaces, and sufficiently large matrix multiplicities in one
	induction. This requires three new ingredients.
	\begin{enumerate}[label=(\arabic*)]
		\item \emph{Finite complexes realizing torsion by coverings} (Section~\ref{sec:top}). We construct finite CW complexes $T_n$ with $\Kt^0(T_n)\cong E_n$
		and $K^1(T_n)=0$ for prescribed finite abelian groups $E_n$, on which every element of $\Kt^0(T_n)$ is the reduced class $[\pi_*\C]-\deg\pi$ of a finite
		covering $\pi$, together with based maps $g_n\colon T_{n+1}\to T_n$ realizing prescribed homomorphisms $E_n\to E_{n+1}$. This is not possible on
		low-dimensional complexes (Remark~\ref{rem:dim3}); our complexes are built from Moore spaces whose dimensions are multiples of Adams' periods, refined
		by acyclic maps so as to acquire large perfect fundamental groups which carry the required coverings.
		\item \emph{A delayed algebraic presentation} (Section~\ref{sec:alg}). A covering realizing a torsion correction may have a large degree, which must be
		absorbed by matrix multiplicities. In the usual presentation of $K_0$ the multiplicities are chosen before the torsion correction (and hence the covering)
		is known. We delay every torsion correction by one stage, so that the covering degree is known before the multiplicity which has to absorb it is chosen.
		\item \emph{A criterion for a trivial limit twist} (Proposition~\ref{prop:twist}). Rather than describing groupoid models of all connecting maps, we track
		normalizers with nonnegative coefficients in canonical bases, and show that a coherent family of such normalizers trivializes the twist of the limit. This
		allows us to modify the connecting maps freely, as long as nonnegative monomial matrices are preserved.
	\end{enumerate}
	
	The construction depends on an \emph{auxiliary order unit} $u$ of the ordered group
	\[
	H=G_0\oplus G',\qquad
	H^+=\{0\}\cup\bigl\{(g,f)\in H:\ \rho(g)(\tau)+f(\tau)>0\ \text{for every }\tau\in T\bigr\},
	\]
	which presents the invariant; here $G'$ is a countable dense subgroup of $\Aff(T)$ containing the constant functions $1$ and $\sqrt2$
	(Section~\ref{ss:enlarge}).
	\begin{itemize}
		\item In the \emph{unital case} $u=([1_A],0)$ is the class of the unit of $A$, which lies in $G_0\subset H$. The endpoint maps of the interval algebras of
		the building blocks are unital, and so are the building blocks and the connecting maps. The endpoint values of the interval algebras are determined by
		the evaluation of the homogeneous part of the building block at a basepoint, and the natural connecting maps pass on to the matrix and interval summands
		only these basepoint values. The limit of these maps is therefore never simple (Section~\ref{ss:whymodify}). To make the limit simple we add, at every
		stage, finitely many point evaluations away from the basepoint: we move basepoint evaluations along paths of representations, which changes the connecting
		maps only within their homotopy classes, so that $K$-theory is unchanged; the modified maps still preserve nonnegative monomial elements
		(Section~\ref{sec:simple}). The traces are identified for the unmodified maps, which intertwine exactly with an AF system whose trace simplex is $T$, and
		the modification is controlled by a summable perturbation argument (Section~\ref{sec:traces}). The first limit $A_\infty=\lims(A_n,\xi_n)$ is not known to
		be $\JS$-stable, so we tensor with Li's untwisted model of $\JS$. The order on $K_0$, which the construction does not track, is then determined by the
		pairing, by strict comparison (Proposition~\ref{prop:cone}). Finally, the classification theorem of \cite{CGSTW} completes the proof
		(Section~\ref{sec:proof}).
		\item In the \emph{stably projectionless case} there is no unit, and we take $u=(0,1)\notin G_0$. The vector of differences of the endpoint ranks of the
		interval algebras of the building blocks is then nonzero, so the two endpoint maps have different ranks in at least one interval summand. The endpoint
		maps therefore have zero blocks, and the building blocks are non-unital and stably projectionless (Lemma~\ref{lem:An}). The device used for simplicity in
		the unital case is not available, because the endpoint values of interval representations are not basepoint blocks. Instead, the detecting
		representations are inserted directly into one matrix summand, with an explicit $K_0$-correction (Section~\ref{sec:SPmaps}). The trace simplex is
		identified by a quantitative comparison with an AF algebra, all densely finite traces are shown to be bounded, and the limit, tensored with $\JS$, has
		continuous scale (Section~\ref{sec:SPtraces}). Finally, an arbitrary scale is realized by an element of the canonical diagonal of the stabilization, whose
		hereditary subalgebra is the $C^*$-algebra of a restricted groupoid. The Gong--Lin classification theorem \cite{GL4} completes the proof
		(Section~\ref{sec:SPscale}).
		\item In the \emph{remaining stably finite case}, $A$ is non-unital but $A\otimes\K$ contains a projection $p$. Then $A$ is a hereditary subalgebra of
		$p(A\otimes\K)p\otimes\K$, which has an untwisted model by the unital case. Since $A$ is not unital, the spectrum of a generator $b$ of this hereditary
		subalgebra accumulates at $0$, and we show that such elements are determined up to Cuntz equivalence by their rank functions, and that Cuntz equivalent
		elements generate isomorphic hereditary subalgebras. Realizing the rank function of $b$ in the canonical diagonal, by the method of the stably
		projectionless case, gives the model (Section~\ref{sec:proj}). No classification theorem for non-unital $C^*$-algebras is needed.
	\end{itemize}
	
	\subsection{Organization}
	The paper has four parts. \emph{Part~I} (Sections~\ref{sec:diag}--\ref{sec:blocks}) is common to the unital and stably projectionless case. Section~\ref{sec:diag} recalls étale
	groupoids and their $C^*$-algebras, introduces covering algebras and the basic homomorphisms between them, recalls Cartan pairs and Li's inductive limit
	theorem, and proves the twist criterion. Section~\ref{sec:top} constructs the complexes $T_n$ and the maps $g_n$ (Theorem~\ref{thm:realize}).
	Section~\ref{sec:alg} constructs the delayed algebraic presentation, with an auxiliary order unit, and discusses the two cases. Section~\ref{sec:blocks}
	constructs the homogeneous part of the building blocks, the building blocks themselves (with possibly non-unital endpoint maps), their groupoid models and
	their normalizers. \emph{Part~II} (Sections~\ref{sec:unitalmaps}--\ref{sec:proof}) treats the unital case: exact connecting maps, simplicity, traces, and the
	proof. \emph{Part~III} (Sections~\ref{sec:SPprelim}--\ref{sec:SPscale}) treats the stably projectionless case: traces and scales, the reduction to continuous
	scale, the connecting maps, traces and the range theorem, and the realization of the scale. \emph{Part~IV} (Section~\ref{sec:proj}) treats the
	non-unital $C^*$-algebras with projections in the stabilization, completing the stably finite case. Appendix~\ref{app:P} proves the necessity of the
	pairing condition.
	
	\subsection{Conventions}\label{ss:conventions}
	All $C^*$-algebras are complex. $M_k=M_k(\C)$ is the algebra of $k\times k$ matrices, with matrix units $e_{ij}$, and $[k]=\{1,\dots,k\}$. $\K$ denotes the
	compact operators on a separable infinite-dimensional Hilbert space. For a compact space $X$, $C(X)$ is the algebra of continuous functions, and
	$K^j(X)=K_j(C(X))$. For a pointed connected space $X$, $\Kt^j(X)$ denotes reduced $K$-theory (Section~\ref{ss:Ktheory}). On $M_k$ we write $\Tr_k$ for the
	unnormalized trace, $\Tr_k(1_k)=k$, and $\tr_k=\frac1k\Tr_k$ for the normalized one; $\Tr$ is the usual trace on $\K$. For a trace $\tau$ on a
	$C^*$-algebra $A$ and a projection $p\in M_k(A)$ we put $\tau_*[p]=(\tau\otimes\Tr_k)(p)$. This depends only on the class of $p$ in $K_0(A)$. For a unital
	$C^*$-algebra, a \emph{trace} means a tracial state, $T(A)$ is the set of tracial states with the weak$^*$ topology, and $\Aff T(A)$ is the real Banach
	space of continuous affine functions $T(A)\to\R$ with the supremum norm. Traces on non-unital $C^*$-algebras are discussed in Section~\ref{ss:SPtraces}.
	Inductive limits $\lims(A_n,\phi_n)$ of $C^*$-algebras are taken in the category of $C^*$-algebras; $\phi_{n,m}=\phi_{m-1}\circ\dots\circ\phi_n$ for
	$n\le m$, and $\phi_{n,\infty}\colon A_n\to\lims A_n$ is the canonical map. 
	
	\subsection{Statement on the use of Artificial Intelligence} This project is the culmination of an extended period of both independent research by the authors together with extensive interactions with ChatGPT (5.6 Sol, 6.0 Astra, and 6.1 Sol) and Claude (Opus 5, Opus 5.5, and Fable 5.1). The initial AI-generated constructions and arguments were humanly incomprehensible, and so the authors devoted considerable effort to making the mathematics humanly accessible. The authors intended this document to be understandable by the general operator algebraist without having specialized background in algebraic topology which many of the arguments are based on. This motivated the extensive algebraic topology preliminaries. Finally, the recent preprint by Evington and Sibbel allowed the authors to omit their treatment of the purely infinite case, shortening the manuscript. The authors take full responsibility for the mathematical correctness of the results, arguments, and constructions presented in this manuscript.
	%======================================================================
	\part*{Part I: Common foundations}
	
	%======================================================================
	\section{Groupoids, coverings and nonnegative normalizers}\label{sec:diag}
	
	This section collects the operator-algebraic preliminaries. Sections~\ref{ss:groupoids}--\ref{ss:basic} are elementary:
	étale groupoids, finite coverings, the covering algebras $\Gamma(\End\pi_*\C)$ and the homomorphisms between them that we
	use throughout. Section~\ref{ss:cartan} recalls Renault's theory of Cartan pairs, Section~\ref{ss:Li} Li's inductive limit
	theorem, and Section~\ref{ss:twistcrit} proves the criterion (Proposition~\ref{prop:twist}) by which the twist of our
	limit will be shown to be trivial.
	
	\subsection{Étale groupoids and their reduced \texorpdfstring{$C^*$}{C*}-algebras}\label{ss:groupoids}
	
	\begin{definition}[Groupoids]\label{def:groupoid}
		A \emph{groupoid} is a set $G$ together with a subset $G^{(0)}\subset G$ (the \emph{units}), maps $r,s\colon G\to G^{(0)}$
		(\emph{range} and \emph{source}), an inversion $\gamma\mapsto\gamma^{-1}$ and a multiplication $(\gamma,\eta)\mapsto\gamma\eta$
		defined on the set of \emph{composable pairs} $G^{(2)}=\{(\gamma,\eta)\in G\times G:s(\gamma)=r(\eta)\}$, such that: the
		multiplication is associative whenever defined; $r(\gamma\eta)=r(\gamma)$ and $s(\gamma\eta)=s(\eta)$; $r(\gamma)\gamma=\gamma
		=\gamma s(\gamma)$; and $\gamma^{-1}\gamma=s(\gamma)$, $\gamma\gamma^{-1}=r(\gamma)$. The elements of $G$ are called
		\emph{arrows}: one thinks of $\gamma$ as an arrow from the unit $s(\gamma)$ to the unit $r(\gamma)$, and of $\gamma\eta$ as
		the composite ``first $\eta$, then $\gamma$''. For a unit $x$ we write $G_x=s^{-1}(x)$, $G^x=r^{-1}(x)$, and $G^x_x=G_x\cap G^x$
		for the \emph{isotropy group} at $x$.
	\end{definition}
	
	\begin{definition}[Topological and étale groupoids]
		A \emph{topological groupoid} is a groupoid with a topology in which multiplication (on $G^{(2)}$ with the relative product
		topology) and inversion are continuous. An open set $U\subset G$ is an \emph{open bisection} if $r|_U$ and $s|_U$ are
		homeomorphisms onto open subsets of $G^{(0)}$. A topological groupoid is \emph{étale} if its topology has a basis of open
		bisections; equivalently, $r$ (and hence $s$) is a local homeomorphism. In an étale groupoid $G^{(0)}$ is open in $G$, and
		the fibres $G_x$, $G^x$ are discrete. A topological groupoid is \emph{principal} if every isotropy group is trivial, i.e.\ the
		map $(r,s)\colon G\to G^{(0)}\times G^{(0)}$ is injective; it is \emph{topologically principal} if the
		set of units with trivial isotropy is dense in $G^{(0)}$, and \emph{effective} if the interior of the isotropy bundle
		$\{\gamma\in G:r(\gamma)=s(\gamma)\}$ is $G^{(0)}$. For locally compact Hausdorff étale groupoids, topological principality implies effectiveness; if the groupoid is also
		second countable, the two notions coincide \cite[\S3]{Re08}. (Without such topological hypotheses the two notions can differ.)
	\end{definition}
	
	A principal groupoid is the same thing as an equivalence relation $R=\{(r(\gamma),s(\gamma))\}\subset G^{(0)}\times
	G^{(0)}$, with arrows the pairs $(x,y)$ of equivalent points, $r(x,y)=x$, $s(x,y)=y$, $(x,y)(y,z)=(x,z)$ and
	$(x,y)^{-1}=(y,x)$; the topology, however, may be finer than the relative product topology.
	
	\begin{example}[Pair groupoids]\label{ex:pairgroupoid}
		For a set $Z$, the \emph{pair groupoid} is $Z\times Z$ with the structure just described; its units are the pairs $(z,z)$,
		identified with $Z$. For the finite set $[k]$ with the discrete topology, the pair groupoid $[k]^2$ is a principal étale
		groupoid, and its arrows $(i,j)$ correspond to the matrix units $e_{ij}$ (see Example~\ref{ex:Mk}).
	\end{example}
	
	\begin{definition}[Fibred products]\label{def:fibred}
		Let $\pi\colon Y\to X$ be a continuous map. The \emph{fibred product}
		\[
		Y\times_XY=\{(y,y')\in Y\times Y:\ \pi(y)=\pi(y')\}
		\]
		carries the relative topology of $Y\times Y$. It is the equivalence relation ``lying in the same fibre of $\pi$'', so it is a
		principal groupoid with $r(y,y')=y$, $s(y,y')=y'$, $(y,y')(y',y'')=(y,y'')$, $(y,y')^{-1}=(y',y)$, and unit space the diagonal
		$\{(y,y)\}\cong Y$. An arrow $(y,y')$ is an arrow from $y'$ to $y$ inside one fibre.
	\end{definition}
	
	\begin{definition}[Reduced groupoid $C^*$-algebras {\cite{Re80}, \cite{Sims}}]\label{def:Cred}
		Let $G$ be a locally compact Hausdorff étale groupoid. The space $C_c(G)$ of compactly supported continuous functions is a
		$*$-algebra under
		\[
		(f*g)(\gamma)=\sum_{\alpha\beta=\gamma}f(\alpha)g(\beta)=\sum_{\eta\in G_{s(\gamma)}}f(\gamma\eta^{-1})g(\eta),\qquad
		f^*(\gamma)=\overline{f(\gamma^{-1})}
		\]
		(the sums are finite since $f$ has compact support and the fibres are discrete). For each unit $x$ the \emph{regular
			representation} $\lambda_x\colon C_c(G)\to B(\ell^2(G_x))$ is $(\lambda_x(f)\xi)(\gamma)=\sum_{\eta\in G_x}f(\gamma\eta^{-1})\xi(\eta)$.
		The \emph{reduced norm} is $\|f\|_r=\sup_x\|\lambda_x(f)\|$, and the \emph{reduced $C^*$-algebra} $C^*_r(G)$ is the
		completion of $C_c(G)$ in this norm. Restriction of functions to the open and closed set $G^{(0)}$ extends to a faithful
		conditional expectation $\EE\colon C^*_r(G)\to C_0(G^{(0)})$, and $C_0(G^{(0)})\subset C^*_r(G)$ as the functions supported on
		the unit space. The formula $j(a)(\gamma)=\langle\lambda_{s(\gamma)}(a)\delta_{s(\gamma)},\delta_\gamma\rangle$ defines an injective
		contractive linear map $j\colon C^*_r(G)\to C_0(G)$ which is the identity on $C_c(G)$; thus elements of $C^*_r(G)$ ``are''
		functions on $G$ \cite[Proposition~II.4.2]{Re80}, \cite{Sims}.
	\end{definition}
	
	If every fibre $G_x$ has at most $N$ points and $G$ is compact, then $\|f\|_\infty\le\|f\|_r\le N\|f\|_\infty$ for $f\in C(G)$, so
	$C_c(G)=C(G)$ is already complete and $C^*_r(G)=C(G)$ with the convolution product. The same holds with $C_0(G)$ when $G$ is
	not compact but the fibres are uniformly bounded. This is the situation for the covering groupoids of Section~\ref{ss:coverings} and for the
	groupoids $\mathcal G_n$ of the building blocks (Lemma~\ref{lem:groupoid}). It is not the situation for the limit groupoids. If all source fibres
	have at most $N$ points, the regular representations have dimension at most $N$, so $C^*_r(G)$ is $N$-subhomogeneous, and a simple subhomogeneous
	$C^*$-algebra is a matrix algebra; hence a groupoid modelling an infinite-dimensional simple $C^*$-algebra has source fibres of unbounded size. In fact, if $G$ is
	principal and $C^*_r(G)$ is simple and infinite-dimensional, then every source fibre is infinite. Indeed, simplicity forces $G$ to be minimal, since a
	nonempty proper open invariant subset $U$ of $G^{(0)}$ would give the proper nonzero ideal generated by $C_0(U)$. As $G$ is principal, $r$ maps the source
	fibre $G_x$ bijectively onto the orbit of $x$. A finite orbit is closed (the unit space is Hausdorff) and invariant, hence by minimality equal to
	$G^{(0)}$; then $G$ is a principal transitive groupoid on finitely many units, i.e.\ a finite pair groupoid, and $C^*_r(G)$ is a matrix algebra. The pair
	groupoid of $\N$ used for stabilization also has infinite source fibres.
	
	\begin{example}[Matrix algebras]\label{ex:Mk}
		For $G=[k]^2$, a function $f\colon[k]^2\to\C$ is a matrix $(f(i,j))_{i,j}$, the convolution $(f*g)(i,k)=\sum_jf(i,j)g(j,k)$
		is matrix multiplication, and $f^*$ is the conjugate transpose. So $C^*_r([k]^2)=M_k$, and $C([k])=C(G^{(0)})$ is the
		diagonal subalgebra $D_k$. The expectation $\EE$ is the diagonal compression $x\mapsto\diag(x)$. More generally, for a compact
		space $Z$, $C^*_r(Z\times[k]^2)=M_k(C(Z))$ with diagonal $C(Z\times[k])$.
	\end{example}
	
	\subsection{Finite coverings and covering algebras}\label{ss:coverings}
	
	\begin{definition}[Finite coverings]\label{def:covering}
		Let $X$ be a Hausdorff space. A \emph{finite covering} of $X$ is a continuous map $\pi\colon Y\to X$ from a Hausdorff space
		$Y$ such that every point of $X$ has an open neighbourhood $U$ which is \emph{evenly covered}: there are a finite set $F$
		and a homeomorphism $\pi^{-1}(U)\cong U\times F$ over $U$. The cardinality of $\pi^{-1}(x)$ is then locally constant, and we
		always assume it is constant, equal to the \emph{degree} $\deg\pi$ (this is automatic if $X$ is connected), and that $\deg\pi\ge1$, i.e.\
		that $\pi$ is surjective. The \emph{trivial
			covering} of degree $d$ is $X\times[d]\to X$. The \emph{disjoint union} of coverings $Y_1\to X$, $Y_2\to X$ is $Y_1\sqcup
		Y_2\to X$, and for a continuous map $u\colon A\to X$ the \emph{pullback} $u^*\pi\colon A\times_XY\to A$ is again a finite
		covering of the same degree; its fibre over $a$ is canonically $\pi^{-1}(u(a))$.
	\end{definition}
	
	Over an evenly covered set $U$, $\pi^{-1}(U)$ is a disjoint union of $d$ open \emph{sheets} $U_1,\dots,U_d$, each mapped
	homeomorphically onto $U$ by $\pi$. If $U$ is connected the sheets are the components of $\pi^{-1}(U)$ and are therefore
	canonical; in general they depend on the choice of the homeomorphism $\pi^{-1}(U)\cong U\times F$, but we never need to choose them globally.
	
	\begin{lemma}\label{lem:Rpietale}
		Let $\pi\colon Y\to X$ be a finite covering of degree $d$ of a compact Hausdorff space $X$. Then $R_\pi=Y\times_XY$ is a compact
		Hausdorff principal étale groupoid, all of whose source fibres have exactly $d$ points.
	\end{lemma}
	
	\begin{proof}
		$Y$ is compact, being locally a product of an open subset of $X$ with a finite set, and $R_\pi$ is closed in $Y\times Y$
		because $X$ is Hausdorff; so $R_\pi$ is compact Hausdorff. The groupoid operations are restrictions of continuous maps, and
		$R_\pi$ is principal by Definition~\ref{def:fibred}. For étaleness, let $(y,y')\in R_\pi$ and let $U$ be an evenly covered open
		neighbourhood of $x=\pi(y)$ with sheets $U_1\ni y$ and $U_k\ni y'$. The set
		\[
		U_{1k}=\{(z,z')\in R_\pi:\ z\in U_1,\ z'\in U_k\}=\big\{\big(z,(\pi|_{U_k})^{-1}\pi(z)\big):z\in U_1\big\}
		\]
		is open in $R_\pi$, contains $(y,y')$, and is the graph of the homeomorphism $(\pi|_{U_k})^{-1}\circ\pi|_{U_1}\colon U_1\to U_k$.
		Hence $r$ maps $U_{1k}$ homeomorphically onto $U_1$ and $s$ maps it homeomorphically onto $U_k$: $U_{1k}$ is an open bisection.
		The source fibre of $y'$ is $\pi^{-1}(\pi(y'))\times\{y'\}$, which has $d$ points.
	\end{proof}
	
	\begin{definition}[Permutation bundles and covering algebras]\label{def:Vpi}
		Let $\pi\colon Y\to X$ be a finite covering of degree $d$ of a compact Hausdorff space. The \emph{permutation bundle}
		$V_\pi=\pi_*\C$ is the complex vector bundle over $X$ whose fibre over $x$ is the space $\C[\pi^{-1}(x)]$ of functions on the
		finite set $\pi^{-1}(x)$, with the orthonormal basis $\{\delta_y:y\in\pi^{-1}(x)\}$, the \emph{sheet basis}. Over an evenly
		covered set $U$ with sheets $U_1,\dots,U_d$, the map $U\times\C^d\to V_\pi|_U$, $(x,(c_k))\mapsto\sum_kc_k\delta_{y_k(x)}$ with
		$y_k(x)=(\pi|_{U_k})^{-1}(x)$, is a local trivialization; two such trivializations differ by a locally constant permutation
		matrix. Thus $V_\pi$ is a Hermitian vector bundle of rank $d$ whose transition functions are permutation matrices. Put
		\[
		F_\pi=\Gamma(\End V_\pi),\qquad D_\pi=C(Y),\qquad R_\pi=Y\times_XY .
		\]
		An element $a\in F_\pi$ is a continuous choice of endomorphisms $a(x)\in\End(\C[\pi^{-1}(x)])$; its \emph{matrix entries}
		are $a(x)_{y,y'}=\langle a(x)\delta_{y'},\delta_y\rangle$ for $y,y'\in\pi^{-1}(x)$. A function $f\in C(Y)$ acts as the diagonal
		endomorphism $\delta_y\mapsto f(y)\delta_y$; this embeds $D_\pi$ into $F_\pi$ as the \emph{sheet-diagonal} elements.
	\end{definition}
	
	Since the sheet basis of each fibre is canonical (it is a set, not an ordered set), it makes sense to say that $a(x)$ is
	\emph{diagonal}, \emph{monomial} or \emph{nonnegative monomial} in the sheet basis (Definition~\ref{def:monomial}); the
	ordering of the basis plays no role.
	
	\begin{proposition}[The covering algebra is an untwisted groupoid $C^*$-algebra]\label{prop:coveringalgebra}
		The map $a\mapsto k_a$, $k_a(y,y')=a(\pi(y))_{y,y'}$, is a $*$-isomorphism $F_\pi\cong C^*_r(R_\pi)=C(R_\pi)$ (with convolution).
		It carries $D_\pi$ onto $C(R_\pi^{(0)})=C(Y)$, and it carries the fibrewise diagonal compression $a\mapsto(x\mapsto\diag(a(x)))$ onto
		the canonical expectation $\EE$ of Definition~\ref{def:Cred}. The full and reduced $C^*$-algebras of $R_\pi$ coincide.
	\end{proposition}
	
	\begin{proof}
		\emph{$k_a$ is continuous.} Near an arrow $(y,y')$ choose an evenly covered neighbourhood $U$ of $x=\pi(y)$ with a decomposition into
		sheets, the sheets $U_j\ni y$, $U_k\ni y'$, and the corresponding local trivialization of Definition~\ref{def:Vpi}. In it, $a|_U$ is a continuous matrix-valued
		function $(a_{il}(x))$, and on the open bisection $U_{jk}$ of Lemma~\ref{lem:Rpietale} we have $k_a(z,z')=a_{jk}(\pi(z))$, which
		is continuous.
		
		\emph{Every continuous function on $R_\pi$ arises.} Conversely, a continuous $k\colon R_\pi\to\C$ defines endomorphisms $a(x)$ by
		$a(x)_{y,y'}=k(y,y')$, and in local trivializations their matrix entries are $x\mapsto k(y_i(x),y_l(x))$, which are continuous. So
		$a\mapsto k_a$ is a linear bijection $F_\pi\to C(R_\pi)$.
		
		\emph{Products and adjoints.} For $(y,y'')\in R_\pi$ with $x=\pi(y)$,
		\[
		(k_a*k_b)(y,y'')=\sum_{y'\in\pi^{-1}(x)}k_a(y,y')k_b(y',y'')=\sum_{y'}a(x)_{y,y'}b(x)_{y',y''}=(a(x)b(x))_{y,y''}=k_{ab}(y,y''),
		\]
		because the arrows of $R_\pi$ composing to $(y,y'')$ are exactly the pairs $(y,y'),(y',y'')$ with $y'$ in the same fibre. Similarly
		$k_a^*(y,y')=\overline{k_a(y',y)}=\overline{a(x)_{y',y}}=(a(x)^*)_{y,y'}$. So $a\mapsto k_a$ is a $*$-isomorphism of $*$-algebras.
		
		\emph{Norms.} For a unit $y'$, the source fibre is $\{(y,y'):y\in\pi^{-1}(x)\}$ ($x=\pi(y')$), and under the identification
		$\ell^2(G_{y'})=\C[\pi^{-1}(x)]$, $\delta_{(y,y')}\mapsto\delta_y$, the regular representation $\lambda_{y'}(k_a)$ is the matrix
		$a(x)$. Hence $\|k_a\|_r=\sup_x\|a(x)\|=\|a\|$, and the isomorphism is isometric. Since $R_\pi$ is compact with finite fibres,
		$C^*_r(R_\pi)=C(R_\pi)$ (remark after Definition~\ref{def:Cred}). Every $*$-representation $\sigma$ of $C(R_\pi)$ satisfies
		$\|\sigma(k)\|\le\|k\|_r$ because $C(R_\pi)\cong F_\pi$ is already a $C^*$-algebra in the norm $\|\cdot\|_r$, and a $*$-homomorphism
		from a $C^*$-algebra is contractive; so the full and reduced norms agree.
		
		\emph{Diagonal and expectation.} $k_a$ is supported on the unit space $\{(y,y)\}$ if and only if every $a(x)$ is diagonal, and
		then $k_a(y,y)=a(\pi(y))_{y,y}$ is the function on $Y$ defining $a$; so $D_\pi$ corresponds to $C(Y)$. Restricting $k_a$ to the
		unit space keeps exactly the diagonal entries, i.e.\ it is the diagonal compression.
	\end{proof}
	
	\begin{figure}[ht]
		\centering
		\begin{tikzpicture}[scale=1.0,>={Stealth[length=2mm]}]
			\draw[thick] plot[domain=0:360,samples=80] ({2.2*cos(\x)},{0.55*sin(\x)});
			\node at (0,-0.95) {$X=S^1$};
			\foreach \k/\h/\col in {0/2.0/blue!70!black,1/2.7/red!70!black,2/3.4/green!50!black}{
				\draw[thick,\col] plot[domain=330:660,samples=90] ({2.2*cos(\x)},{0.55*sin(\x)+\h});}
			\draw[thick,blue!70!black] ({2.2*cos(300)},{0.55*sin(300)+2.0}) -- ({2.2*cos(330)},{0.55*sin(330)+2.7});
			\draw[thick,red!70!black] ({2.2*cos(300)},{0.55*sin(300)+2.7}) -- ({2.2*cos(330)},{0.55*sin(330)+3.4});
			\draw[thick,green!50!black,dashed] ({2.2*cos(300)},{0.55*sin(300)+3.4}) .. controls (3.3,4.0) and (3.6,1.3) .. ({2.2*cos(330)},{0.55*sin(330)+2.0});
			\coordinate (x) at ({2.2*cos(240)},{0.55*sin(240)});
			\fill (x) circle (1.6pt) node[below]{$x$};
			\foreach \h/\col/\n in {2.0/blue!70!black/1,2.7/red!70!black/2,3.4/green!50!black/3}{
				\fill[\col] ({2.2*cos(240)},{0.55*sin(240)+\h}) circle (1.8pt) node[left=2pt,font=\small,fill=white,inner sep=1pt]{$y_\n$};}
			\draw[densely dotted] (x) -- ({2.2*cos(240)},{0.55*sin(240)+3.4});
			\draw[->] (-3.3,2.6) -- node[left]{$\pi$} (-3.3,0.4);
			\node[align=left,font=\small] at (7.0,1.9) {sheet basis of $V_\pi$ at $x$: $\{\delta_{y_1},\delta_{y_2},\delta_{y_3}\}$\\[3pt]
				going once around $X$ permutes\\ the sheets: $y_1\mapsto y_2\mapsto y_3\mapsto y_1$\\[3pt]
				arrows of $R_\pi$ over $x$: the nine pairs $(y_i,y_j)$;\\
				$a(x)$ is a $3\times3$ matrix in this basis};
		\end{tikzpicture}
		\caption{The $3$-sheeted covering $z\mapsto z^3$ of the circle. The sheet basis is canonical; going around the circle permutes
			it, so the transition functions of $V_\pi$ are permutation matrices and carry no phases.}\label{fig:cover}
	\end{figure}
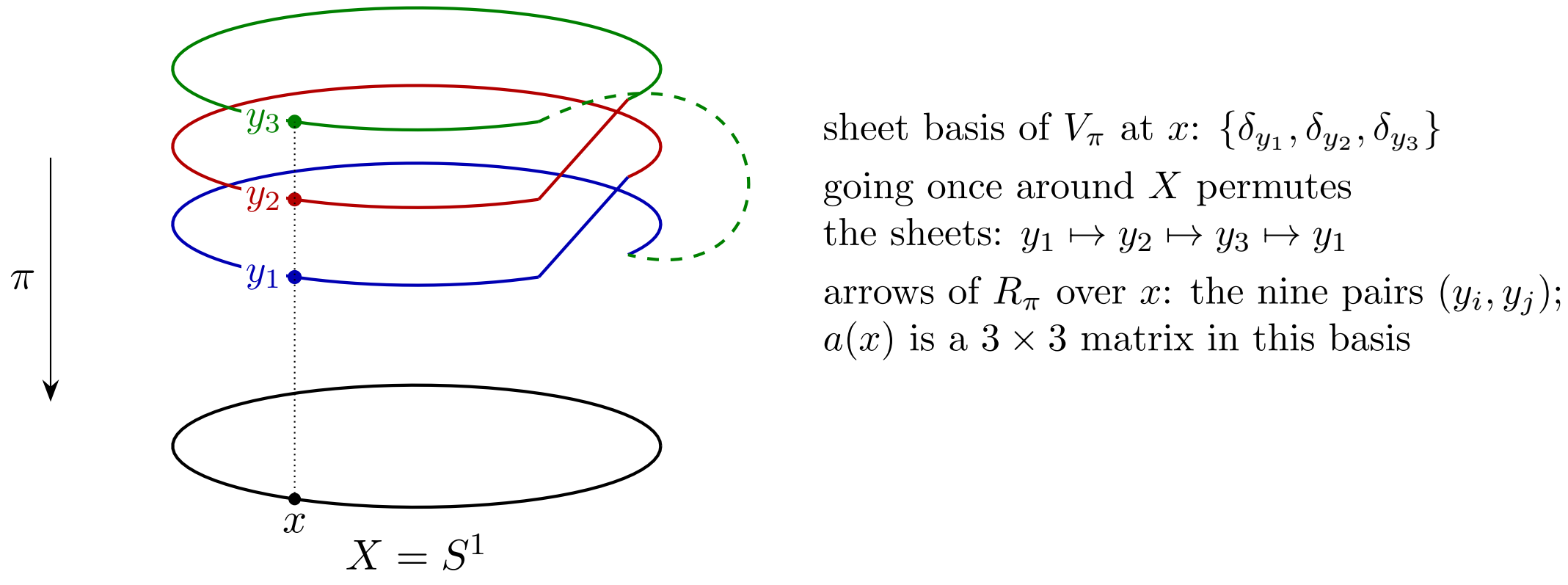
	
	\begin{remark}[Comparison with line bundles]\label{ex:hopf}
		Compare with a homogeneous algebra $F=\Gamma(\End(L_1\oplus L_2))$ for two line bundles $L_1,L_2$ over $X$, with the diagonal
		$\Gamma(\End L_1)\oplus\Gamma(\End L_2)\cong C(X)\oplus C(X)$. The off-diagonal corner is $\Gamma(\Hom(L_2,L_1))$, and its
		normalizers are nowhere vanishing only where $\Hom(L_2,L_1)$ admits nonvanishing sections. The Weyl groupoid is $X\times[2]^2$,
		but the Weyl twist over the component $X\times\{(1,2)\}$ is the unit circle bundle of $\Hom(L_2,L_1)$. For $X=S^2$, $L_1$
		trivial and $L_2$ the tautological line bundle, this is the Hopf fibration $S^3\to S^2$, which is nontrivial; so this Cartan pair
		is twisted. Homogeneous blocks of this kind occur in Li's construction \cite[\S6.1]{Li}. For a permutation bundle the analogous
		circle bundles are trivialized by the sheet bases, and Proposition~\ref{prop:coveringalgebra} says that no twist occurs.
	\end{remark}
	
	\subsection{\texorpdfstring{$K$}{K}-theory of covering algebras}\label{ss:Kcov}
	
	Recall that for a compact space $X$, $K^0(X)=K_0(C(X))$ is the Grothendieck group of isomorphism classes of complex vector
	bundles over $X$ (Serre--Swan \cite{Swan}): a bundle $V$ corresponds to the projective module $\Gamma(V)$, and choosing a
	complement $V\oplus V'\cong X\times\C^N$ realizes $\Gamma(V)$ as $pC(X)^N$ for a projection $p\in M_N(C(X))$, whose class in
	$K_0(C(X))$ is $[V]$. If $X$ is connected, the rank of a bundle defines a homomorphism $\rank\colon K^0(X)\to\Z$, and its
	kernel is the \emph{reduced} group $\Kt^0(X)$ (Section~\ref{ss:Ktheory}). In this subsection we compute the $K$-theory of a covering
	algebra, isolate the part of it which is not determined by the degree (the reduced covering class), and work out in detail which
	reduced classes occur on small examples. The examples explain why the topology of Section~\ref{sec:top} is needed.
	
	\begin{lemma}[$K$-theory of $F_\pi$]\label{lem:KFpi}
		Let $\pi\colon Y\to X$ be a finite covering of degree $d$ of a compact space. There is an isomorphism $K_*(F_\pi)\cong K^*(X)$
		under which a projection $e\in M_k(F_\pi)$ corresponds to the vector bundle $e(V_\pi\otimes\C^k)$; in particular
		\begin{equation}\label{eq:unitclass}
			[1_{F_\pi}]=[V_\pi]=d+\tc(\pi),\qquad\tc(\pi):=[V_\pi]-d\in\Kt^0(X)\quad(X\text{ connected}).
		\end{equation}
	\end{lemma}
	
	\begin{proof}
		Choose a bundle $V'$ with $V_\pi\oplus V'\cong X\times\C^N$ and let $p\in M_N(C(X))$ be the projection onto $V_\pi$. Then
		$F_\pi=\Gamma(\End V_\pi)\cong pM_N(C(X))p$. This corner is \emph{full}: $p$ has rank $d\ge1$ at every point, so the ideal it
		generates contains a projection of constant rank one locally, hence everything. The inclusion of a full corner induces an
		isomorphism on $K$-theory (Morita invariance; see \cite{Bl}, \cite{RLL}), and $M_N(C(X))$ has the same
		$K$-theory as $C(X)$. A projection $e\in M_k(pM_N(C(X))p)$ is a projection over $C(X)$ whose range bundle is
		$e(V_\pi\otimes\C^k)$, and its class is the class of that bundle. For $e=1$ we get $[V_\pi]$.
	\end{proof}
	
	We call $\tc(\pi)$ the \emph{reduced covering class} of $\pi$. It measures how far the permutation bundle is from being trivial, in
	$K$-theory. Two properties are immediate.
	\begin{enumerate}[label=(C\arabic*)]
		\item Degree and reduced covering class are additive under disjoint union, and trivial coverings have $\tc=0$; so adding
		trivial sheets changes only the degree. Indeed, the fibre of $V_{\pi_1\sqcup\pi_2}$ over $x$ is $\C[\pi_1^{-1}(x)\sqcup\pi_2^{-1}(x)]=
		\C[\pi_1^{-1}(x)]\oplus\C[\pi_2^{-1}(x)]$, so $V_{\pi_1\sqcup\pi_2}=V_{\pi_1}\oplus V_{\pi_2}$; and the permutation bundle of $X\times[d]\to X$ is
		$X\times\C^d$.
		\item For $u\colon A\to X$, $V_{u^*\pi}=u^*V_\pi$, hence $\tc(u^*\pi)=u^*\tc(\pi)$: the fibre of $u^*\pi$ over $a$ is $\pi^{-1}(u(a))$, so the fibre of
		$V_{u^*\pi}$ over $a$ is $\C[\pi^{-1}(u(a))]$, the fibre of $u^*V_\pi$, and the local trivializations correspond.
	\end{enumerate}
	A finite covering of a connected, locally path-connected space is the disjoint union of its restrictions to the connected components of the total space,
	and each of these is again a finite covering \cite[\S1.3]{Hat}. By (C1), the reduced covering classes of all finite coverings are therefore the sums of the
	classes of \emph{connected} coverings, and to determine them it suffices to know the connected coverings. Over a connected CW complex these correspond to
	conjugacy classes of subgroups of finite index of the fundamental group \cite[\S1.3]{Hat}.
	
	\subsubsection*{Detecting reduced classes in low dimensions}
	On complexes of dimension at most $3$, reduced $K$-theory is detected by the first Chern class of the determinant bundle. This gives a convenient way to
	compute reduced covering classes, and at the same time shows that such complexes cannot carry interesting covering classes.
	
	\begin{lemma}\label{lem:lowdim}
		Let $X$ be a connected finite CW complex of dimension at most $3$.
		\begin{enumerate}[label=\textup{(\alph*)}]
			\item The map $\delta\colon\Kt^0(X)\to H^2(X;\Z)$, $[V]-\rank V\mapsto c_1(\det V)$, is a well defined isomorphism of groups.
			\item For every finite covering $\pi$ of $X$, $2\,c_1(\det V_\pi)=0$; hence $2\,\tc(\pi)=0$.
		\end{enumerate}
	\end{lemma}
	
	\begin{proof}
		(a) Since $\det(V\oplus W)\cong\det V\otimes\det W$ and $c_1$ is additive on tensor products of line bundles, $V\mapsto c_1(\det V)$ is additive and
		vanishes on trivial bundles; so $\delta$ is a well defined homomorphism. Under the representability isomorphisms $\Kt^0(X)\cong[X,BU]$ and
		$H^2(X;\Z)\cong[X,K(\Z,2)]=[X,BU(1)]$ (Sections~\ref{ss:Ktheory} and~\ref{ss:cov3}; both targets are simply connected, so based and free homotopy classes agree), $\delta$
		is induced by the determinant map $\det\colon BU\to BU(1)$. On homotopy groups, $\pi_k(BU)=\pi_{k-1}(U)$ is $0,\Z,0,\Z$ for $k=1,2,3,4$, and
		$\pi_k(BU(1))$ is $\Z$ for $k=2$ and $0$ otherwise. The determinant $U(N)\to U(1)$ induces an isomorphism on $\pi_1$ for every $N$ (the inclusion
		$U(1)\to U(N)$ as the upper left corner induces an isomorphism on $\pi_1$, and its composite with $\det$ is the identity). Hence $\det\colon BU\to BU(1)$ is
		an isomorphism on $\pi_k$ for $k\le3$ and onto for $k=4$, i.e.\ it is $4$-connected. By Theorem~\ref{thm:standardhtpy}(d), for a CW complex $X$ of
		dimension at most $3$ it induces a bijection $[X,BU]\to[X,BU(1)]$.
		
		(b) Let $N=\deg\pi$. By (C4) of Section~\ref{ss:cov3}, $V_\pi\cong c^*\C^N_\perm$ for a map $c\colon X\to B\Sigma_N$, where $\C^N_\perm$ is the bundle of the
		permutation representation of $\Sigma_N$ on $\C^N$. The determinant of a permutation matrix is the sign of the permutation, so $\det V_\pi\cong c^*\varepsilon$,
		where $\varepsilon$ is the line bundle of the sign character. Since the square of the sign character is trivial, $(\det V_\pi)^{\otimes2}$ is trivial, and
		$2c_1(\det V_\pi)=0$. By (a), $2\tc(\pi)=0$.
	\end{proof}
	
	For the examples we also need the Chern classes of flat line bundles. A homomorphism $\chi\colon\pi_1(X)\to U(1)$ defines the \emph{flat line bundle}
	$L_\chi=\widetilde X\times_{\pi_1(X)}\C$, where $\widetilde X$ is the universal covering and $\pi_1(X)$ acts on $\C$ through $\chi$. Its transition functions
	can be chosen locally constant with values in $U(1)$, so $L_\chi$ defines a class in $H^1(X;U(1)_\delta)=\Hom(\pi_1(X),U(1))$, namely $\chi$; and its first Chern
	class is the image of $\chi$ under the Bockstein homomorphism $\partial\colon H^1(X;U(1)_\delta)\to H^2(X;\Z)$ of the exact sequence of coefficient groups
	$0\to\Z\to\R\to U(1)_\delta\to0$ (the exponential sequence with discrete coefficients) \cite[Chapter~2]{Hus}. Moreover $L_{\chi\chi'}=L_\chi\otimes L_{\chi'}$.
	
	\begin{lemma}[Flat line bundles on $P^2(m)$]\label{lem:flatP2}
		Let $X=P^2(m)=S^1\cup_{z^m}e^2$, $m\ge2$ (Section~\ref{ss:moore}). Then $\pi_1(X)\cong\Z/m$, generated by the class $g$ of the circle,
		$H^1(X;\R)=H^2(X;\R)=0$, and $H^2(X;\Z)\cong\Z/m$. The map $\chi\mapsto c_1(L_\chi)$ is an isomorphism $\Hom(\Z/m,U(1))\to H^2(X;\Z)$. In particular, if
		$\chi(g)=e^{2\pi i/m}$ and $L=L_\chi$, then $c_1(L^{\otimes k})=k\,c_1(L)$, the bundles $L^{\otimes k}$, $0\le k<m$, are pairwise non-isomorphic, and $c_1(L)$
		generates $H^2(X;\Z)$.
	\end{lemma}
	
	\begin{proof}
		By the Seifert--van Kampen theorem (Theorem~\ref{thm:vankampen}(b)), attaching the $2$-cell along $z\mapsto z^m$ adds the relation $g^m=1$ to $\pi_1(S^1)=\Z$; so $\pi_1(X)=\Z/m$. The cellular cochain
		complex of $X$ with coefficients in an abelian group $M$ is $M\xrightarrow{0}M\xrightarrow{m}M$ in degrees $0,1,2$ (the $1$-cell is a loop, so the first
		coboundary vanishes, and the $2$-cell is attached by a map of degree $m$). Hence $H^1(X;M)=\ker(m\colon M\to M)$ and $H^2(X;M)=M/mM$; for $M=\R$ both vanish,
		and for $M=\Z$ we get $H^2(X;\Z)=\Z/m$. In the long exact cohomology sequence of $0\to\Z\to\R\to U(1)_\delta\to0$, the terms $H^1(X;\R)$ and $H^2(X;\R)$
		vanish, so $\partial\colon H^1(X;U(1)_\delta)\to H^2(X;\Z)$ is an isomorphism; and $H^1(X;U(1)_\delta)=\Hom(\Z/m,U(1))$. Combined with
		$c_1(L_\chi)=\partial\chi$ this gives the first assertion, and the rest follows because $\chi\mapsto\chi^k$ corresponds to multiplication by $k$.
	\end{proof}
	
	\subsubsection*{Examples}
	
	\begin{example}[Coverings of the circle]\label{ex:S1cov}
		The connected $k$-sheeted covering of $S^1$ is $z\mapsto z^k$ (Figure~\ref{fig:cover}); it corresponds to the subgroup $k\Z$ of $\pi_1(S^1)=\Z$. Since
		$\Kt^0(S^1)=0$ (Section~\ref{ss:Ktheory}), every reduced covering class over $S^1$ vanishes, although the coverings are far from trivial: going once around
		the circle permutes the sheets cyclically. The permutation bundle is nevertheless trivial, because every vector bundle over $S^1$ is trivial ($U(k)$ is
		connected). The point of this example is that \emph{nontrivial monodromy does not by itself produce nonzero classes}; one needs $K$-theory in which the
		permutation bundle is detected.
	\end{example}
	
	\begin{example}[$\R P^2$]\label{ex:RP2}
		Let $\pi\colon S^2\to\R P^2$ be the double covering, $y\mapsto\{y,-y\}$. We show that $\tc(\pi)$ generates $\Kt^0(\R P^2)\cong\Z/2$.
		
		\emph{(a) The permutation bundle splits.} The fibre of $V_\pi$ over $[y]=\{y,-y\}$ is $\C[\{y,-y\}]$, with sheet basis $\delta_y,\delta_{-y}$. The vector
		$\delta_y+\delta_{-y}$ does not depend on the choice of the representative $y$, so $[y]\mapsto\delta_y+\delta_{-y}$ is a nowhere vanishing continuous section; it
		spans a trivial line subbundle $V^+$ (the \emph{even} functions on the fibre). Its orthogonal complement $V^-$ is spanned at $[y]$ by $\delta_y-\delta_{-y}$,
		which changes sign if $y$ is replaced by $-y$ (the \emph{odd} functions). So $V_\pi=V^+\oplus V^-\cong\C\oplus V^-$.
		
		\emph{(b) The odd part is the complexified tautological bundle.} Let $\gamma\to\R P^2$ be the tautological real line bundle, whose fibre over $[y]$ is the
		line $\R y\subset\R^3$. The map $\R y\to V^-_{[y]}$, $ty\mapsto t(\delta_y-\delta_{-y})$, is well defined (replacing $y$ by $-y$ replaces $t$ by $-t$ and
		$\delta_y-\delta_{-y}$ by its negative), linear, injective and continuous in $[y]$; so it is an isomorphism of real line bundles from $\gamma$ onto a real form of
		$V^-$, and $V^-\cong\gamma\otimes_\R\C=:L$. Equivalently, $L=L_\chi$ is the flat line bundle of the nontrivial character $\chi\colon\pi_1(\R P^2)=\Z/2\to\{\pm1\}$,
		since $S^2$ is the universal covering and the deck transformation acts on the odd functions by $-1$.
		
		\emph{(c) The class.} By Lemma~\ref{lem:lowdim}(a) ($\dim\R P^2=2$) and Lemma~\ref{lem:flatP2} ($\R P^2=P^2(2)$), $\delta\colon\Kt^0(\R P^2)\to
		H^2(\R P^2;\Z)=\Z/2$ is an isomorphism and $\delta([L]-1)=c_1(L)\neq0$. Hence
		\[
		\tc(\pi)=[V_\pi]-2=[L]-1
		\]
		is the nonzero element of $\Kt^0(\R P^2)\cong\Z/2$.
		
		\emph{(d) All classes are covering classes.} The connected coverings of $\R P^2$ are the identity and $\pi$ (the subgroups of $\Z/2$). By (C1) every reduced
		covering class is a multiple of $\tc(\pi)$, and both elements of $\Kt^0(\R P^2)$ occur.
		
		\emph{(e) Adams operations.} The operation $\psi^2$ of Section~\ref{ss:adams} satisfies $\psi^2[L]=[L^{\otimes2}]=1$ (as $\chi^2=1$), so $\psi^2(\tc(\pi))=0\neq
		\tc(\pi)$. Covering classes need not be fixed by Adams operations $\psi^k$ with $k$ dividing the order of the monodromy group; we come back to this in
		Section~\ref{ss:topoverview}.
	\end{example}
	
	\begin{example}[$P^2(p)$, $p$ odd]\label{ex:P2p}
		Let $X=P^2(p)$ for an odd prime $p$. In contrast to $\R P^2=P^2(2)$, \emph{no} covering of $X$ has a nonzero reduced class, although $\Kt^0(X)\cong\Z/p$.
		
		\emph{(a) The group.} By Lemmas~\ref{lem:lowdim}(a) and~\ref{lem:flatP2}, $\delta\colon\Kt^0(X)\to H^2(X;\Z)=\Z/p$ is an isomorphism. Let $L=L_\chi$ with
		$\chi(g)=e^{2\pi i/p}$. Since $\delta([L^{\otimes k}]-1)=c_1(L^{\otimes k})=k\,c_1(L)=\delta\big(k([L]-1)\big)$ and $\delta$ is injective,
		\[
		[L^{\otimes k}]-1=k\,\lambda\qquad(k\in\Z),\qquad\lambda:=[L]-1,
		\]
		and $\lambda$ generates $\Kt^0(X)\cong\Z/p$.
		
		\emph{(b) The connected coverings.} The subgroups of $\pi_1(X)=\Z/p$ are $\{1\}$ and $\Z/p$, so the connected coverings are the identity and the universal
		covering $\pi\colon\widetilde X\to X$, of degree $p$. The latter is a principal $\Z/p$-bundle, and its permutation bundle is the bundle associated to the regular
		representation: $V_\pi=\widetilde X\times_{\Z/p}\C[\Z/p]$. The regular representation of the cyclic group $\Z/p$ is the direct sum of its characters
		$\chi^k$, $0\le k<p$ (with eigenvectors $\sum_{j}\chi(g)^{-jk}\delta_{g^j}$). Hence
		\[
		V_\pi\cong\bigoplus_{k=0}^{p-1}L^{\otimes k}.
		\]
		
		\emph{(c) The class vanishes.} By (a) and (b),
		\[
		\tc(\pi)=\sum_{k=0}^{p-1}\big([L^{\otimes k}]-1\big)=\Big(\sum_{k=0}^{p-1}k\Big)\lambda=\frac{p(p-1)}2\,\lambda=0,
		\]
		because $p$ is odd, so $\frac{p-1}2$ is an integer and $p$ divides $\frac{p(p-1)}2$. (Alternatively: $2\tc(\pi)=0$ by Lemma~\ref{lem:lowdim}(b), and
		$\Kt^0(X)$ has odd order.) By (C1) and the remark preceding Lemma~\ref{lem:lowdim}, every finite covering of $X$ has reduced class $0$.
		
		So the torsion of $\Kt^0(P^2(p))$ is there, but permutation bundles cannot see it. For $p=2$ the same computation gives $\tc(\pi)=\frac{2\cdot1}2\lambda=\lambda\ne0$,
		which is Example~\ref{ex:RP2}.
	\end{example}
	
	\begin{remark}[The need for complexes of high dimension]\label{rem:dim3}
		By Lemma~\ref{lem:lowdim}(b), on a connected finite CW complex of dimension at most $3$ every reduced covering class has order at most $2$. Hence torsion of
		order other than $2$, such as the $\Z/p$ of Example~\ref{ex:P2p}, can never be carried by coverings over complexes of dimension at most $3$. In dimension
		$4$ the determinant no longer detects $\Kt^0$ ($\pi_4(BU)=\Z$ is killed by $\det$), and higher-dimensional complexes may carry more covering classes. This is why
		Section~\ref{sec:top} uses Moore spaces $P^d(p^r)$ whose dimension $d$ is a multiple of an Adams period ($d\ge4$, and $d\ge8$ at $p=2$), refined to complexes
		of unbounded dimension; and it explains why the fundamental group of these complexes has to be large (Remark~\ref{rem:simplyconnected}).
	\end{remark}
	
	\subsection{Monomial matrices and the basic homomorphism}\label{ss:basic}
	
	\begin{definition}[Monomial matrices]\label{def:monomial}
		A complex square matrix is \emph{monomial} if each row and each column contains at most one nonzero entry, and
		\emph{nonnegative monomial} if moreover all its entries are $\ge0$. For a square matrix $x$, $\diag(x)$ denotes its diagonal
		compression (the diagonal matrix with the same diagonal entries).
	\end{definition}
	
	\begin{lemma}\label{lem:monomialprops}
		\begin{enumerate}[label=\textup{(\alph*)}]
			\item Products and adjoints of (nonnegative) monomial matrices are (nonnegative) monomial; limits of monomial matrices are
			monomial.
			\item If $m$ is monomial and $b$ is diagonal, then $mbm^*$ and $m^*bm$ are diagonal.
			\item Conversely, if $me_{kk}m^*$ and $m^*e_{kk}m$ are diagonal for every matrix unit $e_{kk}$, then $m$ is monomial.
			\item Block-diagonal sums $m_1\oplus\dots\oplus m_k$, and conjugates $\sigma m\sigma^{-1}$ by permutation matrices $\sigma$, of
			diagonal, monomial, or nonnegative monomial matrices have the same property, and $\diag(m_1\oplus\dots\oplus
			m_k)=\diag(m_1)\oplus\dots\oplus\diag(m_k)$, $\diag(\sigma m\sigma^{-1})=\sigma\diag(m)\sigma^{-1}$.
		\end{enumerate}
	\end{lemma}
	
	\begin{proof}
		(a) A monomial matrix $m$ is a partial permutation matrix times an invertible diagonal matrix: there are an injective partial map
		$\sigma$ and nonzero scalars $c_j$ with $me_j=c_je_{\sigma(j)}$ for $j$ in the domain of $\sigma$ and $me_j=0$ otherwise. Composites
		of such data are of the same form, and $m^*e_{\sigma(j)}=\bar c_je_j$. If $m_k\to m$ and $m$ had two nonzero entries in a row,
		these entries would be nonzero in $m_k$ for large $k$. Nonnegativity of entries is preserved by products, adjoints and limits.
		(b) $mbm^*e_{\sigma(j)}=|c_j|^2b_{jj}e_{\sigma(j)}$, and $mbm^*$ kills the vectors not in the range of $\sigma$. (c)
		$me_{kk}m^*=(me_k)(me_k)^*$ is diagonal iff the column $me_k$ has at most one nonzero entry; similarly $m^*e_{kk}m$ is diagonal
		iff the $k$th row has at most one nonzero entry. (d) is clear.
	\end{proof}
	
	We now describe the homomorphisms between covering algebras used in this paper. Let $\pi\colon Y\to X$ be a finite covering of
	degree $r$, let $\nu\colon W\to X'$ be a finite covering of degree $d$, and let $g\colon W\to X$ be continuous. The pullback
	$g^*\pi\colon U=W\times_XY\to W$ is an $r$-sheeted covering, and the composite
	\[
	\nu\circ g^*\pi\colon U=\{(w,y)\in W\times Y:\ g(w)=\pi(y)\}\longrightarrow X',\qquad(w,y)\mapsto\nu(w),
	\]
	is an $rd$-sheeted covering of $X'$, whose fibre over $x'$ is $\bigsqcup_{w\in\nu^{-1}(x')}\{w\}\times\pi^{-1}(g(w))$.
	
	\begin{definition}[The basic homomorphism]\label{def:basic}
		With these data, define $\Phi=\Phi_{\nu,g}\colon F_\pi\to F_{U\to X'}$ by
		\begin{equation}\label{eq:basic}
			\Phi(a)(x')=\bigoplus_{w\in\nu^{-1}(x')}a(g(w)),\qquad\text{i.e.}\qquad
			\Phi(a)(x')_{(w,y),(w',y')}=\delta_{w,w'}\,a(g(w))_{y,y'} .
		\end{equation}
		In words: pull $a$ back along $g$, and then add over the sheets of $\nu$.
	\end{definition}
	
	\begin{figure}[ht]
		\centering
		\begin{tikzcd}[column sep=large,row sep=large]
			U=W\times_XY \arrow[r] \arrow[d,"g^*\pi"'] & Y \arrow[d,"\pi"] \\
			W \arrow[r,"g"] \arrow[d,"\nu"'] & X \\
			X' &
		\end{tikzcd}
		\qquad
		\begin{tikzpicture}[baseline=(c.base),scale=0.8]
			\node (c) at (0,0) {};
			\draw[thick] (0,-2.1) rectangle (4.2,2.1);
			\fill[blue!20] (0,2.1) rectangle (1.4,0.7);
			\fill[red!20] (1.4,0.7) rectangle (2.8,-0.7);
			\fill[green!20] (2.8,-0.7) rectangle (4.2,-2.1);
			\node[font=\scriptsize] at (0.7,1.4) {$a(g(w_1))$};
			\node[font=\scriptsize] at (2.1,0) {$a(g(w_2))$};
			\node[font=\scriptsize] at (3.5,-1.4) {$a(g(w_3))$};
			\node[font=\small,align=center] at (2.1,-2.8) {$\Phi(a)(x')$, with\\ $\nu^{-1}(x')=\{w_1,w_2,w_3\}$};
		\end{tikzpicture}
		\caption{The basic homomorphism \eqref{eq:basic}: pull back along $g$, then add over the sheets of $\nu$.}\label{fig:basic}
	\end{figure}
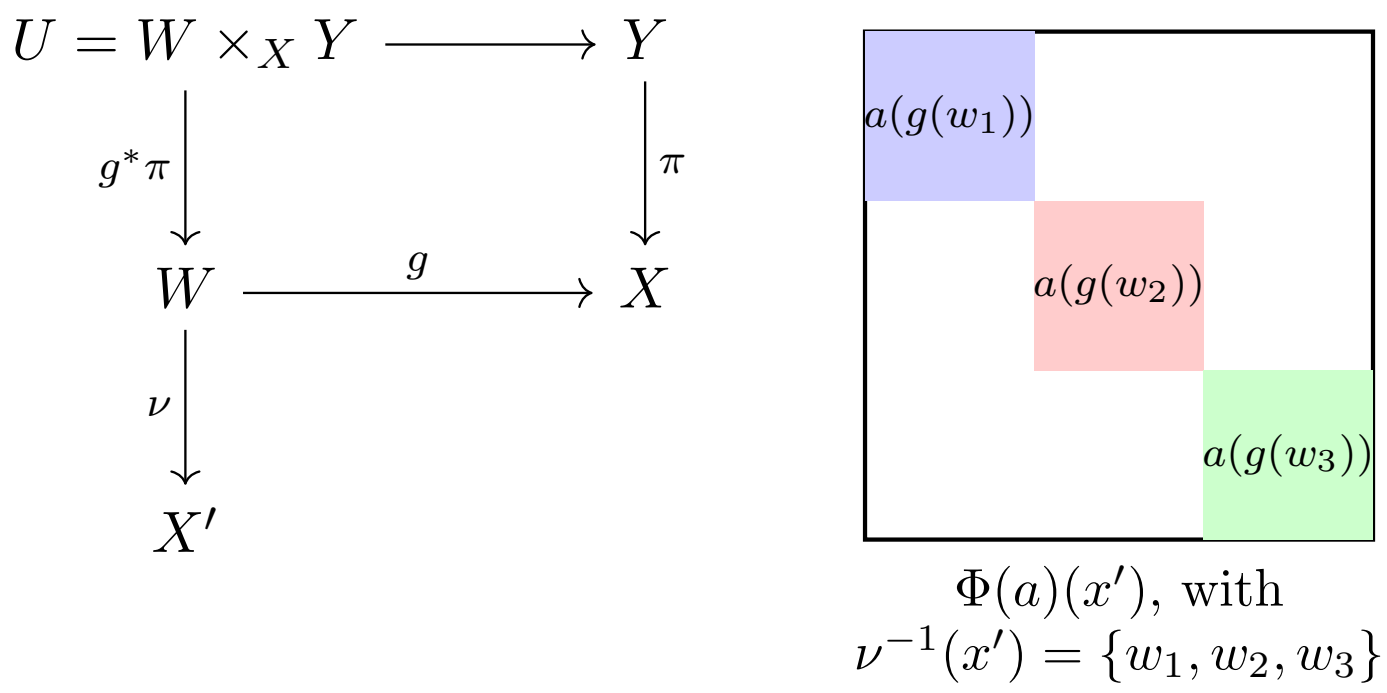
	
	\begin{example}[Blocks]\label{ex:blocks}
		Let $\pi$ have degree $2$ and $\nu$ degree $3$, and fix $x'$ with $\nu^{-1}(x')=\{w_1,w_2,w_3\}$, and
		$\pi^{-1}(g(w_i))=\{y_i,y_i'\}$. Order the fibre of $U$ over $x'$ as
		$(w_1,y_1)$, $(w_1,y_1')$, $(w_2,y_2)$, $(w_2,y_2')$, $(w_3,y_3)$, $(w_3,y_3')$. Then
		\[
		\Phi(a)(x')=\begin{pmatrix}
			a_{11}^{(1)}&a_{12}^{(1)}&0&0&0&0\\
			a_{21}^{(1)}&a_{22}^{(1)}&0&0&0&0\\
			0&0&a_{11}^{(2)}&a_{12}^{(2)}&0&0\\
			0&0&a_{21}^{(2)}&a_{22}^{(2)}&0&0\\
			0&0&0&0&a_{11}^{(3)}&a_{12}^{(3)}\\
			0&0&0&0&a_{21}^{(3)}&a_{22}^{(3)}
		\end{pmatrix},\qquad a^{(i)}=\big(a(g(w_i))\big)\text{ in the basis }(\delta_{y_i},\delta_{y_i'}).
		\]
		If, say, $a(g(w_1))$ is the monomial matrix $\begin{psmallmatrix}0&2\\3&0\end{psmallmatrix}$ and $a(g(w_2)),a(g(w_3))$ are
		diagonal, the image is a $6\times6$ monomial matrix: each row and column still has at most one nonzero entry. When $x'$ moves
		around a loop, the three blocks may be permuted (monodromy of $\nu$), and inside each block the two basis vectors may be
		interchanged (monodromy of $\pi$); but the image is always block-diagonal in the sheet basis of $U$, and the blocks are moved
		around only by permutations.
	\end{example}
	
	Two special cases are important.
	\begin{itemize}
		\item \emph{Trivial $\nu$.} If $W=X'\times[d]$ and $g(x',k)=g_k(x')$, then
		\[
		\Phi(a)(x')=\diag\big(a(g_1(x')),\dots,a(g_d(x'))\big),
		\]
		a block-diagonal sum of ordinary \emph{eigenvalue maps} $a\mapsto a\circ g_k$.
		\item \emph{Matrix source.} If $X$ is a point and $Y=[r]$, then $F_\pi=M_r$, $g$ is constant, and $\Phi(a)(x')=\bigoplus_{w}a$, a
		direct sum of copies of the constant matrix $a$, one for each sheet of $\nu$; as $x'$ moves around a loop the copies are
		permuted. So the target covering algebra is $M_r\otimes F_\nu$, and $\Phi$ is $a\mapsto a\otimes1$; its image is the constant copy
		$M_r\otimes1\cong M_r$, not all of $M_r\otimes F_\nu$. (The unit of this image is the unit of $M_r\otimes F_\nu$, whose class has a nonzero reduced
		part $r\,\tc(\nu)$ in general; compare Lemma~\ref{lem:basicK}.)
	\end{itemize}
	
	\begin{lemma}\label{lem:basicK}
		$\Phi$ is a unital $*$-homomorphism. On $K$-theory, under the identifications of Lemma~\ref{lem:KFpi}, it induces
		$\nu_!\circ g^*\colon K^*(X)\to K^*(X')$, where $g^*$ is pullback and $\nu_!$ is the direct image $(\nu_!E)_{x'}=\bigoplus_{w\in
			\nu^{-1}(x')}E_w$. In particular, for a matrix source $M_r$ the class $[e_{11}]$ of a minimal projection is mapped to
		$[V_\nu]=\deg\nu+\tc(\nu)$, and $[1]$ to $r(\deg\nu+\tc(\nu))$.
	\end{lemma}
	
	\begin{proof}
		Continuity of $\Phi(a)$: near $x'$ choose an evenly covered neighbourhood $U'$ for $\nu$ with sheets $W_1,\dots,W_d$, and for each $k$
		an evenly covered neighbourhood of $g(w_k(x'))$ for $\pi$; in the resulting local trivializations $\Phi(a)$ is a block-diagonal
		matrix whose blocks are $a\circ g\circ w_k$ in local frames, which are continuous. Being a fibrewise direct sum of
		$*$-homomorphisms, $\Phi$ is a unital $*$-homomorphism. If $e\in M_k(F_\pi)$ is a projection with range bundle
		$E=e(V_\pi\otimes\C^k)$, then $\Phi(e)$ is a projection whose range at $x'$ is $\bigoplus_{w\in\nu^{-1}(x')}E_{g(w)}=(\nu_!g^*E)_{x'}$.
		By Lemma~\ref{lem:KFpi} this proves the statement on $K_0$; for $K_1$ apply the same argument to the suspensions (or to unitaries).
		For a matrix source, $e_{11}$ corresponds to the trivial line bundle over a point, whose image is $\nu_!(W\times\C)=V_\nu$.
	\end{proof}
	
	\begin{lemma}\label{lem:basic}
		The maps \eqref{eq:basic} send diagonal, monomial and nonnegative monomial elements to elements of the same kind, and they commute
		with fibrewise diagonal compression: $\diag(\Phi(a)(x'))=\Phi(\diag a)(x')$.
	\end{lemma}
	
	\begin{proof}
		At each $x'$ the matrix of $\Phi(a)(x')$ in the sheet basis of $U$ is, by \eqref{eq:basic}, the block-diagonal sum of the matrices of
		$a(g(w))$ in the sheet bases of $\pi^{-1}(g(w))$ (Example~\ref{ex:blocks}). Lemma~\ref{lem:monomialprops}(d) gives all assertions.
	\end{proof}
	
	\subsection{Cartan pairs, \texorpdfstring{$C^*$}{C*}-diagonals and twists}\label{ss:cartan}
	
	\begin{definition}[Cartan pairs {\cite[Definition~5.1]{Re08}}; $C^*$-diagonals {\cite{Ku}}]\label{def:cartan}
		Let $B\subset A$ be a $C^*$-subalgebra. The \emph{normalizers} of $B$ are the elements of
		\[
		\Nor(B)=\{a\in A:\ aBa^*\subset B\ \text{and}\ a^*Ba\subset B\}.
		\]
		$B$ is \emph{regular} in $A$ if $\Nor(B)$ generates $A$ as a $C^*$-algebra. A \emph{conditional expectation} $\EE\colon A\to B$ is a contractive positive
		linear idempotent map onto $B$ with $\EE(b_1ab_2)=b_1\EE(a)b_2$; it is \emph{faithful} if $\EE(a^*a)=0$ implies $a=0$. The pair $(A,B)$ is a \emph{Cartan pair}
		(and $B$ a \emph{Cartan subalgebra}) if $B$ is a maximal abelian $C^*$-subalgebra of $A$ which is regular in $A$, and there is a faithful conditional expectation
		$A\to B$. A Cartan pair is a \emph{$C^*$-diagonal} if moreover $B$ has the \emph{unique extension property}: every pure state of $B$ extends uniquely to a
		state of $A$.
	\end{definition}
	
	\begin{remark}[Approximate units]\label{rem:pitts}
		Renault's original definition \cite[Definition~5.1]{Re08} also requires that $B$ contain an approximate unit of $A$. This requirement is redundant: Pitts
		proved that if $B$ is a regular maximal abelian $C^*$-subalgebra of $A$, then every approximate unit of $B$ is an approximate unit of $A$ \cite[Theorem~2.6]{Pi}.
		The key observation is that for every normalizer $n$ the elements $n^*n$ and $nn^*$ commute with $B$, hence lie in $B$ by maximality \cite[Proposition~2.1]{Pi};
		if $(u_k)$ is an approximate unit of $B$, then $\|n-u_kn\|^2=\|(1-u_k)nn^*(1-u_k)\|\to0$, and since the normalizers generate $A$ it follows that $(u_k)$ is a left,
		hence (by symmetry) a two-sided, approximate unit of $A$. So Definition~\ref{def:cartan} agrees with Renault's. We shall use the consequence repeatedly:
		\emph{a Cartan subalgebra contains an approximate unit of the ambient $C^*$-algebra}; in particular, a Cartan pair $(A,B)$ has $A$ unital if and only if $B$ is unital.
	\end{remark}
	
	The structure of Cartan pairs is described by twisted groupoids.
	
	\begin{definition}[Twists {\cite[\S4]{Re08}}]\label{def:twist}
		Let $G$ be a locally compact Hausdorff étale groupoid. A \emph{twist} over $G$ is a locally compact Hausdorff groupoid $\Sigma$
		with an extension of topological groupoids
		\[
		G^{(0)}\times\T\hookrightarrow\Sigma\xrightarrow{\ p\ }G
		\]
		in which $G^{(0)}\times\T$ is central and $p$ is locally trivial: every $\gamma\in G$ has a neighbourhood $U$ with a continuous
		section $U\to\Sigma$ of $p$. The group $\T$ acts freely on $\Sigma$, and $p$ identifies $\Sigma/\T$ with $G$. The associated
		Hermitian line bundle is $L=\Sigma\times_\T\C\to G$. The twisted groupoid $C^*$-algebra $C^*_r(G,\Sigma)$ is the completion of the space
		of compactly supported continuous sections of $L$ (equivalently, of $\T$-equivariant functions $f(z\sigma)=\bar zf(\sigma)$ on
		$\Sigma$) in the reduced norm, with convolution and involution defined as in Definition~\ref{def:Cred} using the multiplication of
		$L$ \cite[\S4]{Re08}. The twist is \emph{trivial} if $\Sigma\cong G^{(0)}\times\T\times_{G^{(0)}}G=G\times\T$ as extensions; this
		holds if and only if $p$ has a continuous section which is a groupoid homomorphism. Then $C^*_r(G,\Sigma)\cong C^*_r(G)$.
	\end{definition}
	
	The map $j$ of Definition~\ref{def:Cred} extends to twisted groupoids: it is an injective contractive linear map
	$j\colon C^*_r(G,\Sigma)\to C_0(G;L)$ into the continuous sections of $L$ vanishing at infinity, which is the identity on compactly
	supported sections \cite[Proposition~4.1 and the subsequent discussion]{Re08}, \cite[Proposition~II.4.2]{Re80}. The canonical expectation is $\EE(a)=j(a)|_{G^{(0)}}$. The convolution formulas for $j(ab)$ and for
	$\EE(a^*b)$, which are central in Section~\ref{ss:twistcrit}, are proved there (Lemma~\ref{lem:jformulas}).
	
	\begin{theorem}[{\cite{Re08}}, {\cite{Ku}}]\label{thm:renault}
		\begin{enumerate}[label=\textup{(\alph*)}]
			\item \cite[Theorem~5.9]{Re08} Let $(A,B)$ be a Cartan pair with $A$ separable. There are a second countable, locally compact
			Hausdorff, topologically principal étale groupoid $G$ (the \emph{Weyl groupoid}) and a twist $\Sigma$ over $G$ (the \emph{Weyl
				twist}), unique up to isomorphism, and an isomorphism $A\cong C^*_r(G,\Sigma)$ carrying $B$ onto $C_0(G^{(0)})$.
			\item \cite[Theorem~5.2]{Re08} Conversely, for every such twisted groupoid, $C_0(G^{(0)})$ is a Cartan subalgebra of
			$C^*_r(G,\Sigma)$, and \cite[Propositions~4.13(ii), 4.15]{Re08} the Weyl groupoid and Weyl twist of this Cartan pair are $G$ and
			$\Sigma$.
			\item \cite[Corollary~5.10]{Re08} There is a unique conditional expectation onto a Cartan subalgebra.
			\item \cite[Proposition~4.8]{Re08} If $G$ is topologically principal, an element $a\in C^*_r(G,\Sigma)$ is a normalizer of
			$C_0(G^{(0)})$ if and only if its \emph{open support} $\{\gamma:j(a)(\gamma)\ne0\}$ is an open bisection.
			\item \cite[Proposition~5.11]{Re08}, \cite{Ku} A Cartan pair is a $C^*$-diagonal if and only if its Weyl groupoid is principal.
		\end{enumerate}
	\end{theorem}
	
	\begin{remark}[Separability is not needed]\label{rem:raad}
		Raad has shown that part (a) holds without the separability assumption \cite[Theorem~1.2]{Raad}: every Cartan pair $(A,B)$, separable or not, is isomorphic to
		$(C^*_r(G,\Sigma),C_0(G^{(0)}))$ for a twist $\Sigma$ over a locally compact Hausdorff étale groupoid $G$ which is \emph{effective} (the interior of the isotropy is
		the unit space). Without separability, $G$ need not be second countable, and topological principality has to be replaced by effectiveness; for second countable
		groupoids the two notions agree. Conversely, $C_0(G^{(0)})$ is a Cartan subalgebra of $C^*_r(G,\Sigma)$ for every such twisted effective groupoid
		\cite[Theorem~1.2]{Raad}. All $C^*$-algebras in this paper are separable, so we only use Renault's theorem in the form (a).
	\end{remark}
	
	We explain where the twist comes from, since this is what the criterion of Section~\ref{ss:twistcrit} exploits. Let $(A,B)$
	be a Cartan pair with $B=C_0(X)$. A normalizer $n$ induces a partial homeomorphism $\alpha_n$ of $X$, from the open set
	$\{x:(n^*n)(x)>0\}$ onto $\{x:(nn^*)(x)>0\}$, characterized by $n^*bn=(b\circ\alpha_n)\,n^*n$ for $b\in B$. The Weyl groupoid
	consists of the germs $[\alpha_n(x),\alpha_n,x]$ of these partial homeomorphisms; it records how normalizers move points, but not
	their phases. The Weyl twist consists of the classes $[\alpha_n(x),n,x]$, where two normalizers $n,n'$ define the same element
	at $x$ if $nb=n'b'$ for some $b,b'\in B$ with $b(x),b'(x)>0$. Multiplying $n$ by a function $z\in C(X,\T)$ changes $[\alpha_n(x),n,x]$ by
	the phase $z(x)$; so $\Sigma\to G$ is a circle bundle, and it is trivial exactly when phases can be chosen consistently and
	multiplicatively over all of $G$ \cite[\S4--5]{Re08}. Remark~\ref{ex:hopf} shows an example where this is impossible.
	
	\begin{example}[Untwisted examples]\label{ex:untwisted}
		For a principal étale groupoid $G$, $(C^*_r(G),C_0(G^{(0)}))$ is a $C^*$-diagonal with trivial twist. Examples:
		$(M_k,D_k)$ with $G=[k]^2$; $(M_k(C(Z)),D_k(C(Z)))$ with $G=Z\times[k]^2$; and, by Proposition~\ref{prop:coveringalgebra},
		$(F_\pi,D_\pi)$ with $G=R_\pi$.
	\end{example}
	
	\subsection{Li's inductive limit theorem}\label{ss:Li}
	
	\begin{theorem}[{\cite[Theorem~1.10]{Li}}]\label{thm:Li}
		Let $(A_n,B_n)$ be Cartan pairs with conditional expectations $\EE_n\colon A_n\to B_n$, and let $\varphi_n\colon A_n\to A_{n+1}$
		be injective $*$-homomorphisms such that
		\begin{enumerate}[label=\textup{(L\arabic*)}]
			\item $\varphi_n(B_n)\subset B_{n+1}$;
			\item $\varphi_n(\Nor(B_n))\subset\Nor(B_{n+1})$;
			\item $\EE_{n+1}\circ\varphi_n=\varphi_n\circ\EE_n$.
		\end{enumerate}
		Then $B=\lims(B_n,\varphi_n|_{B_n})$ is a Cartan subalgebra of $A=\lims(A_n,\varphi_n)$. If all $B_n$ are $C^*$-diagonals, then $B$
		is a $C^*$-diagonal in $A$.
	\end{theorem}
	
	The theorem does not require the $C^*$-algebras or the maps to be unital; this matters in the stably projectionless case.
	
	The hypotheses have the following roles. (L1) makes $B$ a subalgebra of $A$. (L2) ensures that the normalizers of the $B_n$
	become normalizers of $B$ in the limit (an element normalizing a dense subset of $B$ normalizes $B$), so that $\Nor(B)$ generates
	$A$. (L3) makes the expectations compatible, so that they define a faithful conditional expectation $\EE\colon A\to B$ (faithfulness
	in the limit is part of Li's argument). The difficult part of Li's theorem is that $B$ is maximal abelian in $A$, and the unique
	extension property in the $C^*$-diagonal case. Theorem~\ref{thm:Li} says nothing about the Weyl twist of the limit. In
	\cite{Li} the twist is controlled by describing groupoid models of the connecting maps, which requires a specific list of allowed
	connecting maps; when $K_0$ has torsion, Li's twists are in general nontrivial. We replace this by the following criterion,
	which needs much less information about the connecting maps.
	
	\subsection{A criterion for a trivial limit twist}\label{ss:twistcrit}
	
	This subsection contains the operator-algebraic mechanism by which the twist is removed, and it is one of the central points of the paper. We first explain
	the problem, then prove two lemmas about the coefficient functions $j(a)$, state and prove the criterion, and finally discuss its scope with examples. To our knowledge the criterion is new; it is a consequence of Renault's theory and
	of Li's theorem.
	
	\subsubsection*{The problem}
	Our $C^*$-algebras will be inductive limits $A=\lims(A_n,\varphi_n)$ of $C^*$-diagonals $(A_n,B_n)$, each of which is given as an \emph{untwisted} groupoid $C^*$-algebra
	$(C^*_r(\mathcal G_n),C_0(\mathcal G_n^{(0)}))$ (Section~\ref{ss:groupoid}), and the connecting maps satisfy (L1)--(L3). By Li's theorem the limit $(A,B)$ is a
	$C^*$-diagonal, and by Renault's theorem $(A,B)\cong(C^*_r(G,\Sigma),C_0(G^{(0)}))$ for a principal groupoid $G$ and a twist $\Sigma$. We must show that $\Sigma$ is
	trivial. This does not follow from the triviality of the twists of the $(A_n,B_n)$: a trivialization of the twist of $(A_n,B_n)$ is an additional structure (a
	consistent choice of phases of normalizers, see Remark~\ref{rem:modeldep}), the connecting maps need not respect these structures, and an inverse system of
	trivializations which are not compatible with each other gives no trivialization of the limit. What is needed is that the trivializations at the finite stages
	are compatible with the connecting maps. The criterion expresses this compatibility in terms of normalizers alone, so that it can be checked without computing
	the groupoid models of the connecting maps.
	
	The idea is the following. In an untwisted groupoid $C^*$-algebra $C^*_r(G)$ the fibres of the line bundle $L=G\times\C$ have a distinguished positive half-line
	$\R_{>0}\subset\C$, and the elements $a$ with $j(a)\ge0$ form a ``positive cone of coefficients''. Conversely, a continuous and multiplicative choice of a positive
	half-line in every fibre of the line bundle $L$ of a twist is exactly a trivialization of the twist (a unit vector in each fibre is a point of $\Sigma$). The
	criterion asks for a family of normalizers whose coefficients are ``nonnegative relative to each other'' (condition (S3)), which is closed under products
	(condition (S1)), large enough to see every arrow (condition (S2)), and preserved by the connecting maps (condition (S4)); from such a family we construct the
	positive half-lines in the limit.
	
	\subsubsection*{Two lemmas on coefficient functions}
	We use the notation of Definition~\ref{def:twist}: $G$ is a locally compact Hausdorff étale groupoid, $\Sigma$ a twist over $G$ with line bundle $L$ and Hermitian
	metric $\langle\ ,\ \rangle$ (conjugate-linear in the first variable), and $j\colon C^*_r(G,\Sigma)\to C_0(G;L)$ the coefficient map. The multiplication of $\Sigma$
	induces isometric bilinear maps $L_\alpha\times L_\beta\to L_{\alpha\beta}$ for $(\alpha,\beta)\in G^{(2)}$, written $(u,v)\mapsto uv$, and the inversion induces
	conjugate-linear isometries $L_\alpha\to L_{\alpha^{-1}}$, $u\mapsto\bar u$. In the untwisted case $L=G\times\C$, $uv$ is the product of complex numbers, $\bar u$
	is complex conjugation and $\langle u,v\rangle=\bar uv$.
	
	\begin{lemma}[Coefficient formulas]\label{lem:jformulas}
		Let $a,b\in C^*_r(G,\Sigma)$ and $\gamma\in G$.
		\begin{enumerate}[label=\textup{(\alph*)}]
			\item $j(a^*)(\gamma)=\overline{j(a)(\gamma^{-1})}$, and $\sum_{\beta\in G_x}\|j(a)(\beta)\|^2\le\|a\|^2$ for every unit $x$.
			\item The series $\sum_{\alpha\in G^{r(\gamma)}}j(a)(\alpha)\,j(b)(\alpha^{-1}\gamma)$ converges absolutely in $L_\gamma$, and
			\begin{equation}\label{eq:jprod}
				j(ab)(\gamma)=\sum_{\alpha\in G^{r(\gamma)}}j(a)(\alpha)\,j(b)(\alpha^{-1}\gamma).
			\end{equation}
			\item For every unit $x$,
			\begin{equation}\label{eq:Ests}
				\EE(a^*b)(x)=\sum_{\gamma\in G_x}\langle j(a)(\gamma),j(b)(\gamma)\rangle .
			\end{equation}
		\end{enumerate}
	\end{lemma}
	
	\begin{proof}
		We give the proof in the untwisted case; in the twisted case the same computation applies verbatim, with products and conjugates taken in $L$, because these
		operations are isometric \cite[\S4]{Re08}.
		
		(a) For $f\in C_c(G)$ the first formula is the definition of $f^*$, and it extends to $C^*_r(G)$ because $j$ and the involution are continuous. For the
		second, recall that for $\beta\in G_x$ the number $j(a)(\beta)$ is the $\beta$-coordinate of the vector $\lambda_x(a)\delta_x\in\ell^2(G_x)$
		(Definition~\ref{def:Cred}; for $f\in C_c(G)$ indeed $(\lambda_x(f)\delta_x)(\beta)=f(\beta)$), so
		$\sum_{\beta\in G_x}|j(a)(\beta)|^2=\|\lambda_x(a)\delta_x\|^2\le\|a\|^2$ by Parseval's identity in $\ell^2(G_x)$.
		
		(b) As $\alpha$ runs through $G^{r(\gamma)}$, $\alpha^{-1}$ runs through $G_{r(\gamma)}$ and $\alpha^{-1}\gamma$ runs through $G_{s(\gamma)}$, each bijectively. By (a)
		and the Cauchy--Schwarz inequality,
		\[
		\sum_{\alpha\in G^{r(\gamma)}}|j(a)(\alpha)|\,|j(b)(\alpha^{-1}\gamma)|\le\Big(\sum_{\alpha'\in G_{r(\gamma)}}|j(a^*)(\alpha')|^2\Big)^{1/2}
		\Big(\sum_{\beta\in G_{s(\gamma)}}|j(b)(\beta)|^2\Big)^{1/2}\le\|a\|\,\|b\|.
		\]
		So the right side of \eqref{eq:jprod} converges absolutely and defines, for fixed $\gamma$, a bilinear form $\Phi_\gamma(a,b)$ with $|\Phi_\gamma(a,b)|\le\|a\|\|b\|$. The
		left side $j(ab)(\gamma)$ is also bilinear with $|j(ab)(\gamma)|\le\|ab\|\le\|a\|\|b\|$. For $a,b\in C_c(G)$ the two sides agree: this is the convolution formula
		of Definition~\ref{def:Cred}, written with $\alpha=\gamma\eta^{-1}$. Since $C_c(G)$ is dense and both sides are jointly continuous, they agree for all $a,b$.
		
		(c) Apply (b) to $a^*$ and $b$ at the unit $x$: $\EE(a^*b)(x)=j(a^*b)(x)=\sum_{\alpha\in G^x}j(a^*)(\alpha)j(b)(\alpha^{-1})$. Substituting $\gamma=\alpha^{-1}\in G_x$ and
		using $j(a^*)(\gamma^{-1})=\overline{j(a)(\gamma)}$ from (a) gives \eqref{eq:Ests}.
	\end{proof}
	
	\begin{lemma}[Normalizers with nonnegative overlaps]\label{lem:overlap}
		Let $G$ be principal, and let $s,t\in C^*_r(G,\Sigma)$ be normalizers of $C_0(G^{(0)})$ with open supports $U_s=\{j(s)\neq0\}$ and $U_t=\{j(t)\ne0\}$. Then $U_s$ and
		$U_t$ are open bisections, and for $\gamma\in U_s$ with source $x$,
		\begin{equation}\label{eq:overlap}
			\EE(s^*t)(x)=\langle j(s)(\gamma),j(t)(\gamma)\rangle .
		\end{equation}
		In particular, if $\EE(s^*t)\ge0$ and $\gamma\in U_s\cap U_t$, then $j(t)(\gamma)$ is a positive multiple of $j(s)(\gamma)$.
	\end{lemma}
	
	\begin{proof}
		$U_s$ and $U_t$ are open because $j(s),j(t)$ are continuous sections, and they are bisections by Theorem~\ref{thm:renault}(d) ($G$ is principal, in particular
		topologically principal). In \eqref{eq:Ests} a term $\langle j(s)(\eta),j(t)(\eta)\rangle$ with $\eta\in G_x$ can be nonzero only if $\eta\in U_s$. Since $U_s$ is a
		bisection, the source map is injective on $U_s$, so $\gamma$ is the only arrow of $U_s$ with source $x$, and the sum reduces to the single term
		\eqref{eq:overlap}. If $\gamma\in U_s\cap U_t$, then $j(s)(\gamma)$ and $j(t)(\gamma)$ are nonzero vectors in the one-dimensional space $L_\gamma$, so
		$j(t)(\gamma)=c\,j(s)(\gamma)$ with $c\in\C\setminus\{0\}$ and $\langle j(s)(\gamma),j(t)(\gamma)\rangle=c\,\|j(s)(\gamma)\|^2$. If this number is $\ge0$, then $c>0$.
	\end{proof}
	
	\subsubsection*{The criterion}
	
	\begin{proposition}\label{prop:twist}
		Let $(A_n,B_n)$ be $C^*$-diagonals with $A_n$ separable (not necessarily unital), with conditional expectations $\EE_n$, and let
		$\varphi_n\colon A_n\to A_{n+1}$ be injective $*$-homomorphisms satisfying \textup{(L1)--(L3)}. Suppose there are subsets
		$S_n\subset\Nor(B_n)$ such that
		\begin{enumerate}[label=\textup{(S\arabic*)}]
			\item $S_n$ is closed under products and adjoints;
			\item the linear span of $S_n$ is dense in $A_n$;
			\item $\EE_n(s^*t)\ge0$ for all $s,t\in S_n$;
			\item $\varphi_n(S_n)\subset S_{n+1}$.
		\end{enumerate}
		Then $B=\lims B_n$ is a $C^*$-diagonal in $A=\lims A_n$ whose Weyl twist is trivial. Consequently there is a second countable,
		locally compact Hausdorff, principal étale groupoid $G$ such that $(A,B)\cong(C^*_r(G),C_0(G^{(0)}))$. The unit space $G^{(0)}$ is
		compact if and only if $A$ is unital.
	\end{proposition}
	
	The conditions should be read as follows. (S1) and (S3) say that $S_n$ is a multiplicative family of normalizers any two of which have ``nonnegative overlap'': by
	Lemma~\ref{lem:overlap}, wherever the open supports of $s$ and $t$ meet, their coefficients point in the same direction. (S2) guarantees that the open supports of the
	elements of $S_n$ cover the whole groupoid. (S4) is the compatibility with the connecting maps discussed above. Before the proof we give two remarks and an example on the
	hypotheses.
	
	\begin{remark}[The criterion is sharp]\label{rem:sharp}
		If $(A,B)=(C^*_r(G),C_0(G^{(0)}))$ with $G$ principal and étale, let $S_G$ be the set of nonnegative functions in $C_c(G)$ supported in open bisections. Then $S_G$
		satisfies (S1)--(S3).
		\begin{itemize}
			\item \emph{(S1).} If $f,g\in S_G$ are supported in the open bisections $U,V$, then $f*g$ is supported in the set $UV=\{\alpha\beta:\alpha\in U,\beta\in V,
			s(\alpha)=r(\beta)\}$, which is an open bisection, and $(f*g)(\alpha\beta)=f(\alpha)g(\beta)\ge0$ because an arrow of $UV$ has exactly one factorization with
			factors in $U$ and $V$. Similarly $f^*(\gamma)=\overline{f(\gamma^{-1})}=f(\gamma^{-1})\ge0$ is supported in the bisection $U^{-1}$.
			\item \emph{(S2).} By a partition of unity, every $h\in C_c(G)$ is a finite sum of functions supported in open bisections (cover the compact support by finitely
			many bisections), and every function supported in a bisection $U$ is a linear combination of four nonnegative functions supported in $U$ (positive and negative
			parts of its real and imaginary parts). Since $C_c(G)$ is dense, (S2) holds.
			\item \emph{(S3).} By \eqref{eq:Ests}, $\EE(f^*g)(x)=\sum_{\gamma\in G_x}f(\gamma)g(\gamma)$, which has at most one nonzero term (in $G_x\cap U$) and is $\ge0$.
		\end{itemize}
		Together with Proposition~\ref{prop:twist} for the constant sequence $(A_n,B_n,\varphi_n)=(A,B,\id)$, this shows: \emph{a separable $C^*$-diagonal has trivial
			Weyl twist if and only if it admits a set of normalizers satisfying (S1)--(S3).} The content of Proposition~\ref{prop:twist} for inductive limits is that (S4)
		carries such sets to the limit.
	\end{remark}
	
	\begin{remark}[Model dependence]\label{rem:modeldep}
		The sets $S_n$ are \emph{part of the data}; ``normalizers with nonnegative coefficients'' depend on the chosen untwisted model,
		not only on the pair $(A_n,B_n)$. Two untwisted models of the same Cartan pair differ by an automorphism of the form
		$f\mapsto c\cdot(f\circ\alpha^{-1})$, with $c\colon G\to\T$ a continuous homomorphism, and this preserves nonnegativity only if
		$c\equiv1$. For example, in $M_2$ the automorphism $\Ad(\diag(1,-1))$ preserves the diagonal, the normalizers and the expectation,
		but sends $e_{12}$ to $-e_{12}$. For the double covering $\pi\colon S^2\to\R P^2$, $R_\pi=\Delta\sqcup\{(y,-y)\}$ with both pieces
		clopen; write $\Gamma=\{(y,-y):y\in S^2\}$ for the second piece. The function $c$ equal to $1$ on $\Delta$ and $-1$ on $\Gamma$ is a continuous
		homomorphism, and $f\mapsto cf$ is an automorphism of $C^*(R_\pi)$ fixing $C(S^2)$ pointwise (in particular the unit $1_\Delta$). The indicator function
		$v=1_\Gamma$ is a nonnegative normalizer, with $v^*=v$ and $v^2=1_\Delta$, and the automorphism sends it to $-v$; and $c$ is not of the
		form $c(y,y')=u(y)\overline{u(y')}$ (that would give an odd map $S^2\to S^1$, contradicting the Borsuk--Ulam theorem). In our
		construction one untwisted model is fixed at every stage (sheet bases and standard matrix units), and the connecting maps respect
		these models; this is what (S4) expresses. If, say, a connecting map $M_2\to M_2$ were the automorphism $\Ad(\diag(1,-1))$ and $S$ the nonnegative monomial
		matrices, (S4) would fail; the limit of the constant sequence with this connecting map is of course still untwisted, but the criterion does not see it with
		these sets $S_n$, and it would be necessary to change the sets.
	\end{remark}
	
	\begin{example}[A twisted pair has no such family]\label{ex:notwistfamily}
		Let $(F,D)$ be the Cartan pair of Remark~\ref{ex:hopf}, with $X=S^2$, $L_1$ trivial and $L_2$ the tautological line bundle. Its Weyl twist over the component
		$X\times\{(1,2)\}$ is the Hopf fibration, which is nontrivial. By Remark~\ref{rem:sharp}, $(F,D)$ admits no set of normalizers with (S1)--(S3). Concretely,
		a normalizer supported on the off-diagonal component is a section of $\Hom(L_2,L_1)$, and (S2) and (S3) would force the normalized nonvanishing parts of such
		sections to agree on overlaps, hence to glue to a nowhere vanishing section of the unit circle bundle of $\Hom(L_2,L_1)\cong\overline{L_2}$ over $S^2$; no such
		section exists.
	\end{example}
	
	\begin{figure}[ht]
		\centering
		\begin{tikzpicture}[>={Stealth[length=2mm]}]
			\draw[thick] (0,0) -- (10,0) node[right]{$G^{(0)}$};
			\fill[gray!8] (0,1.2) .. controls (3,1.9) and (6,0.9) .. (10,1.6) -- (10,2.6) .. controls (6,1.9) and (3,2.9) .. (0,2.2) -- cycle;
			\node[gray] at (10.4,2.1) {$G$};
			\draw[blue!70!black,line width=2.2pt] (1,1.62) .. controls (2.5,1.98) and (4.2,1.6) .. (6,1.55);
			\draw[red!70!black,line width=2.2pt] (4,1.95) .. controls (5.5,1.75) and (7.2,1.7) .. (9,1.98);
			\node[blue!70!black] at (1.5,1.2) {$U_s$};
			\node[red!70!black] at (8.6,2.35) {$U_t$};
			\foreach \x in {4.4,5.0,5.6}{
				\draw (\x,3.4) circle (0.25);
				\draw[->,blue!70!black,thick] (\x,3.4) -- ++(55:0.25);
				\draw[->,red!70!black,thick,dashed] (\x,3.4) -- ++(55:0.25);
				\draw[densely dotted] (\x,3.15) -- (\x,1.85);}
			\node[align=center,font=\small] at (5,4.25) {on $U_s\cap U_t$: $\langle j(s)(\gamma),j(t)(\gamma)\rangle=\EE(s^*t)(s(\gamma))>0$,\\ so the phases agree};
			\foreach \x in {2,3,7,8}{\draw[->,gray] (\x,1.45) -- (\x,0.1);}
			\node[gray,font=\small] at (3.6,0.55) {source map};
		\end{tikzpicture}
		\caption{Two bisections $U_s,U_t$ of the Weyl groupoid. Nonnegativity of $\EE(s^*t)$ forces the unit vectors
			$j(s)/|j(s)|$ and $j(t)/|j(t)|$ to coincide on the overlap; they glue to a section of the twist.}\label{fig:bisections}
	\end{figure}
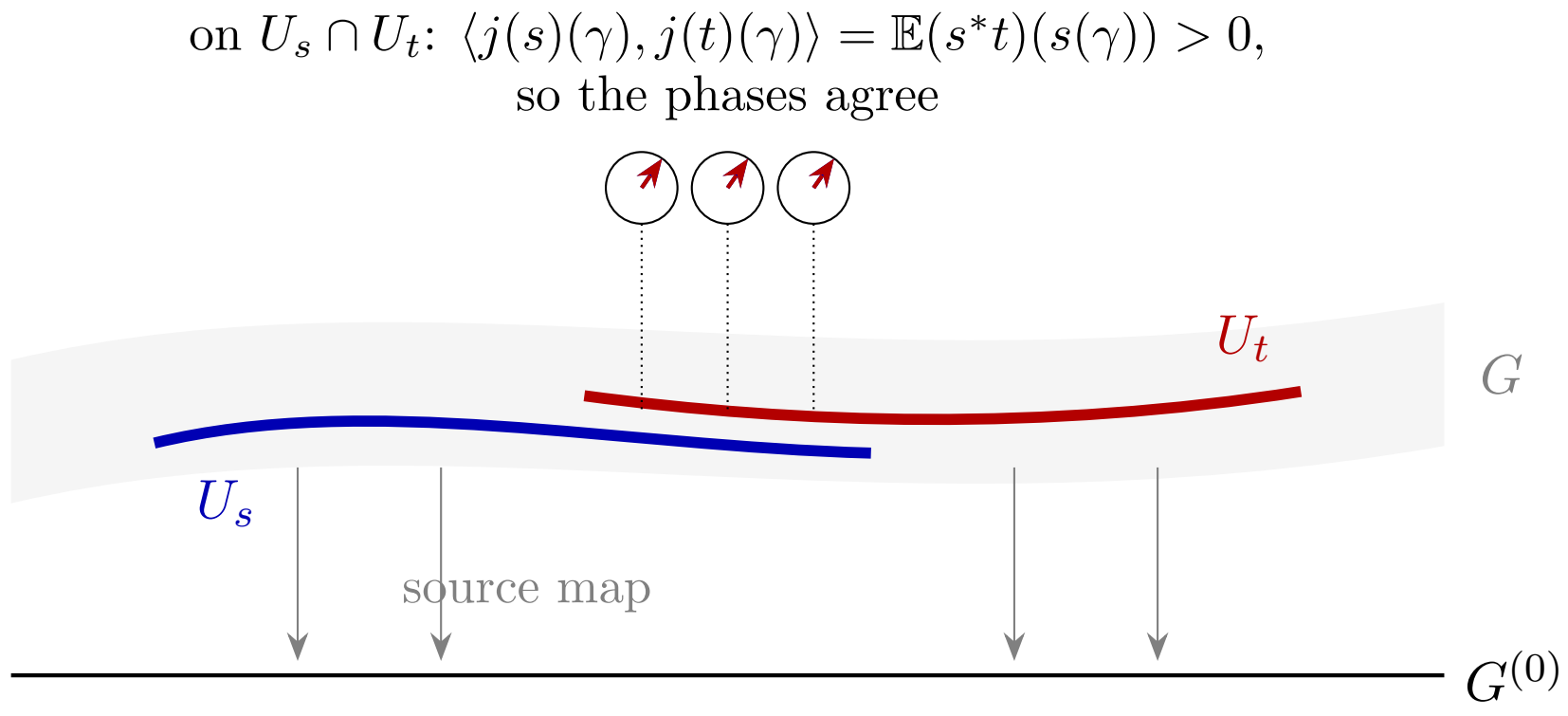
	
	\begin{proof}[Proof of Proposition~\ref{prop:twist}]
		\emph{Step 1: the limit is a $C^*$-diagonal.} By Theorem~\ref{thm:Li}, $B$ is a $C^*$-diagonal in $A$. As an inductive limit of separable $C^*$-algebras, $A$ is separable,
		so by Theorem~\ref{thm:renault}(a),(e) there are a second countable, locally compact Hausdorff, principal étale groupoid $G$ and a twist $\Sigma$ with
		$(A,B)\cong(C^*_r(G,\Sigma),C_0(G^{(0)}))$. By Remark~\ref{rem:pitts}, $B$ contains an approximate unit of $A$; hence $A$ is unital if and only if $B=C_0(G^{(0)})$ is
		unital, i.e.\ if and only if $G^{(0)}$ is compact. We fix this isomorphism and the maps $j$ and $\EE$ of Section~\ref{ss:cartan}, and write $L$ for the line bundle of
		$\Sigma$. In the rest of the proof we construct a continuous section $\sigma\colon G\to\Sigma$ of $p$ which is a groupoid homomorphism.
		
		\emph{Step 2: the limit expectation is the canonical one.} By (L3) and since each $\EE_n$ is contractive, the maps
		$\varphi_{n,\infty}(a)\mapsto\varphi_{n,\infty}(\EE_n(a))$ are well defined and contractive on the dense subalgebra
		$\bigcup_n\varphi_{n,\infty}(A_n)$: if $\varphi_{n,\infty}(a)=\varphi_{m,\infty}(a')$ with $n\le m$, then $\varphi_{n,m}(a)=a'$ by injectivity, and
		$\varphi_{m,\infty}(\EE_m(a'))=\varphi_{m,\infty}(\EE_m(\varphi_{n,m}(a)))=\varphi_{n,\infty}(\EE_n(a))$ by (L3). They extend to a conditional expectation
		$\EE_\infty\colon A\to B$ (positivity, idempotence and the bimodule property pass to the limit). By Theorem~\ref{thm:renault}(c) there is only one conditional
		expectation onto $B$, so $\EE_\infty$ is the canonical expectation $a\mapsto j(a)|_{G^{(0)}}$. (Directly: let $(u_k)$ be an approximate unit of $A$ contained in $B$
		(Remark~\ref{rem:pitts}). For every pure state $\chi$ of $B$ and every conditional expectation $\EE'$, $\chi\circ\EE'$ is a positive functional extending $\chi$, of
		norm $\lim_k\chi(\EE'(u_k))=\lim_k\chi(u_k)=1$; so it is the unique state extending $\chi$. Hence $\chi\circ\EE'=\chi\circ\EE''$ for any two expectations, and pure
		states separate the points of $B$.)
		
		\emph{Step 3: the family $S$.} Let $S=\bigcup_n\varphi_{n,\infty}(S_n)\subset A$. By (S4), $\varphi_{n,\infty}(S_n)=\varphi_{n+1,\infty}(\varphi_n(S_n))\subset
		\varphi_{n+1,\infty}(S_{n+1})$, so the union is increasing. Hence any two elements of $S$ lie in a common $\varphi_{n,\infty}(S_n)$, and $S$ is closed under products and
		adjoints by (S1). Its span is dense in $A$ by (S2), since $\bigcup_n\varphi_{n,\infty}(A_n)$ is dense and $\varphi_{n,\infty}$ is continuous. Every element of $S$
		normalizes $B$: for $s\in S_n$ and $m\ge n$, $\varphi_{n,m}(s)\in S_m\subset\Nor(B_m)$, so $\varphi_{n,\infty}(s)$ normalizes $\varphi_{m,\infty}(B_m)$ for every $m$, hence
		the dense subalgebra $\bigcup_m\varphi_{m,\infty}(B_m)$ of $B$, hence (by continuity of $b\mapsto sbs^*$ and closedness of $B$) all of $B$. Finally, for $s,t\in S$
		choose $n$ and $s_n,t_n\in S_n$ with $s=\varphi_{n,\infty}(s_n)$, $t=\varphi_{n,\infty}(t_n)$; then $\EE_\infty(s^*t)=\varphi_{n,\infty}(\EE_n(s_n^*t_n))\ge0$ by (S3), and
		$\EE_\infty=\EE$ by Step~2. Summarizing: \emph{$S\subset\Nor(B)$ is closed under products and adjoints, has dense span, and $\EE(s^*t)\ge0$ for all $s,t\in S$.}
		
		\emph{Step 4: supports are bisections.} For $s\in S$ let $U_s=\{\gamma\in G:j(s)(\gamma)\ne0\}$. By Lemma~\ref{lem:overlap} this is an open bisection.
		
		\emph{Step 5: the supports cover $G$.} Let $\gamma\in G$. The functional $a\mapsto j(a)(\gamma)\in L_\gamma$ is linear and bounded ($\|j(a)\|_\infty\le\|a\|$).
		It does not vanish identically: by local triviality of $\Sigma$ there are an open bisection $W\ni\gamma$ and a continuous section $W\to\Sigma$ of $p$, hence a
		nowhere vanishing continuous section of $L$ over $W$; multiplying it by a function in $C_c(W)$ which is nonzero at $\gamma$ gives a compactly supported section $f$
		of $L$ with $j(f)(\gamma)=f(\gamma)\ne0$. A bounded linear functional which vanishes on a set with dense span vanishes identically; since it does not, there is
		$s\in S$ with $j(s)(\gamma)\ne0$, i.e.\ $\gamma\in U_s$. Hence $G=\bigcup_{s\in S}U_s$.
		
		\emph{Step 6: the phases agree.} For $s\in S$ and $\gamma\in U_s$ put $\sigma_s(\gamma)=j(s)(\gamma)/\|j(s)(\gamma)\|$, a unit vector in $L_\gamma$; the unit vectors in
		$L_\gamma=\Sigma_\gamma\times_\T\C$ are exactly the classes $[\sigma',1]$ with $\sigma'\in p^{-1}(\gamma)$, so $\sigma_s(\gamma)$ is a point of $\Sigma$ over $\gamma$. If
		$\gamma\in U_s\cap U_t$, Lemma~\ref{lem:overlap} and $\EE(s^*t)\ge0$ (Step~3) show that $j(t)(\gamma)$ is a positive multiple of $j(s)(\gamma)$, i.e.\
		$\sigma_s(\gamma)=\sigma_t(\gamma)$ (Figure~\ref{fig:bisections}).
		
		\emph{Step 7: a continuous section.} By Steps~5 and~6 the maps $\sigma_s$ glue to a well defined map $\sigma\colon G\to\Sigma$ with $p\circ\sigma=\id$. It is
		continuous because it is continuous on each set of the open cover $\{U_s\}$: over a neighbourhood $W\subset U_s$ of a point on which $L$ is trivialized by a
		continuous section of $p$, $j(s)$ is a nowhere vanishing continuous complex function, and $\sigma_s$ corresponds to $j(s)/|j(s)|$, which is continuous.
		
		\emph{Step 8: multiplicativity.} Let $(\gamma,\eta)\in G^{(2)}$, and choose $s,t\in S$ with $\gamma\in U_s$, $\eta\in U_t$ (Step~5). By \eqref{eq:jprod},
		\[
		j(st)(\gamma\eta)=\sum_{\alpha\in G^{r(\gamma)}}j(s)(\alpha)\,j(t)(\alpha^{-1}\gamma\eta).
		\]
		A term can be nonzero only if $\alpha\in U_s$; since $r$ is injective on the bisection $U_s$ and $r(\alpha)=r(\gamma)$, the only such $\alpha$ is $\gamma$, and then
		$\alpha^{-1}\gamma\eta=\eta$. So $j(st)(\gamma\eta)=j(s)(\gamma)\,j(t)(\eta)$, which is nonzero because the multiplication $L_\gamma\times L_\eta\to L_{\gamma\eta}$ is
		isometric. Thus $\gamma\eta\in U_{st}$, and since $st\in S$ (by (S1) for $S$), normalizing gives
		\[
		\sigma(\gamma\eta)=\sigma_{st}(\gamma\eta)=\frac{j(s)(\gamma)\,j(t)(\eta)}{\|j(s)(\gamma)\|\,\|j(t)(\eta)\|}=\sigma_s(\gamma)\sigma_t(\eta)=\sigma(\gamma)\sigma(\eta).
		\]
		Similarly $\sigma(\gamma^{-1})=\sigma(\gamma)^{-1}$, using $s^*\in S$ and $j(s^*)(\gamma^{-1})=\overline{j(s)(\gamma)}$ (Lemma~\ref{lem:jformulas}(a)); and on units
		$\sigma(x)$ is the unit of $\Sigma$ over $x$, because $\sigma(x)=\sigma(x)\sigma(x)$.
		
		\emph{Step 9: trivialization.} The map $\Psi\colon G\times\T\to\Sigma$, $(\gamma,z)\mapsto z\sigma(\gamma)$, is continuous, and it is a groupoid homomorphism by Step~8
		and centrality of $\T$. It is bijective, because $\T$ acts freely and transitively on each fibre $p^{-1}(\gamma)$. Its inverse
		$\sigma'\mapsto\big(p(\sigma'),\sigma'\sigma(p(\sigma'))^{-1}\big)$ is continuous, since $p$, $\sigma$ and the groupoid operations of $\Sigma$ are continuous and
		$\sigma'\sigma(p(\sigma'))^{-1}$ lies in $G^{(0)}\times\T$. Hence $\Psi$ is an isomorphism of extensions, $\Sigma$ is trivial, and
		$(A,B)\cong(C^*_r(G,\Sigma),C_0(G^{(0)}))\cong(C^*_r(G),C_0(G^{(0)}))$.
	\end{proof}
	
	\begin{remark}[Application of the criterion]\label{rem:applytwist}
		In our construction $A_n$ is a subhomogeneous $C^*$-algebra whose untwisted groupoid model $\mathcal G_n$ is given explicitly (Lemma~\ref{lem:groupoid}), and $S_n$ is the
		set $N_{+,n}$ of \emph{nonnegative monomial elements}: the elements all of whose fibre matrices, in the sheet bases of the covering algebras and the standard matrix
		units, are nonnegative monomial matrices. Conditions (S1)--(S3) are verified in Lemma~\ref{lem:normalizers}. Conditions (S1) and (S3) reduce to the finite-dimensional statements that
		products and adjoints of nonnegative monomial matrices are nonnegative monomial, and that the diagonal of $\mu^*\nu$ is nonnegative for nonnegative matrices
		$\mu,\nu$; condition (S2) is a partition of unity argument which takes the endpoint conditions of the building blocks into account. Condition (S4) holds because every connecting map is built from block embeddings, permutations of sheets and point evaluations, all of which are given by
		nonnegative monomial matrices (Lemmas~\ref{lem:basic} and~\ref{lem:labelledprops}(a)). No groupoid model of a connecting map is ever needed.
	\end{remark}
	%======================================================================
	\section{Topology: finite complexes on which torsion is carried by coverings}\label{sec:top}
	
	\subsection{Contents of this section}\label{ss:topoverview}
	
	The building blocks of Section~\ref{sec:blocks} contain a covering algebra $F_\lambda=\Gamma(\End V_\lambda)$ over a finite
	complex $Z_n$, and by Lemma~\ref{lem:KFpi} its unit has $K_0$-class $\deg\lambda+\tc(\lambda)$. To produce torsion in $K_0$ in
	this way, we need spaces on which prescribed torsion classes in reduced $K$-theory are reduced classes of finite coverings,
	and maps between such spaces which induce prescribed homomorphisms. The goal of this section is the following theorem;
	its precise, inductive form is Theorem~\ref{thm:realize} below.
	
	\begin{theorem}[Realization of torsion by coverings, short form]\label{thm:realize-short}
		Let $E_1\xrightarrow{i_1}E_2\xrightarrow{i_2}E_3\to\cdots$ be a sequence of finite abelian groups and homomorphisms. There are
		connected finite CW complexes $T_n$, based maps $g_n\colon T_{n+1}\to T_n$, and isomorphisms $\phi_n\colon
		E_n\to\Kt^0(T_n)$, such that $K^1(T_n)=0$, $g_n^*\circ\phi_n=\phi_{n+1}\circ i_n$, and every element of $\Kt^0(T_n)$ is the
		reduced covering class $\tc(\nu)$ of a finite covering $\nu$ of $T_n$. Moreover, $T_{n+1}$ and $g_n$ can be constructed from
		the data of stage $n$ and from $i_n$ alone.
	\end{theorem}
	
	The last sentence is essential for the inductive constructions of Sections~\ref{ss:choices} and~\ref{ss:SPorder}: the spaces are built one at a time, and at
	each step only finitely many data are known.
	
	\subsubsection*{The use of the theorem}
	The theorem enters the construction of the building blocks in three ways.
	\begin{enumerate}[label=(\roman*)]
		\item \emph{The torsion of $K_0$.} The torsion subgroup of the invariant is presented as $\lims(E_n,i_n)$ (Section~\ref{ss:defects}). At stage $n$ the
		building block contains the covering algebra $F_{\lambda_n}$ over $Z_n=T_n\vee Y_n\vee[0,1]$, and $K_0(F_{\lambda_n})=K^0(Z_n)=\Z\oplus E_n$
		(Lemmas~\ref{lem:KFpi} and~\ref{lem:Zn}). The unit of $F_{\lambda_n}$ has class $\deg\lambda_n+\tc(\lambda_n)$; so the torsion part $\zeta_n$ of the
		class of the unit, which is dictated by the invariant, must be the reduced class of a covering. This is the reason for the requirement that
		\emph{every} element of $\Kt^0(T_n)$ be a covering class.
		\item \emph{The connecting maps.} The connecting maps of the building blocks are built from pullbacks along maps $Z_{n+1}\to Z_n$, and on the
		torsion part of $K_0$ they must induce $i_n$. This is the condition $g_n^*\phi_n=\phi_{n+1}i_n$. Pullback along $g_n$ maps coverings to coverings
		of the same degree (Section~\ref{ss:Kcov}, (C2)), so the covering algebras are preserved by the connecting maps.
		\item \emph{Separation of $K_0$ and $K_1$.} The condition $K^1(T_n)=0$ means that $T_n$ contributes nothing to $K_1$; the group $K_1$ is carried by
		the separate wedge summand $Y_n$ (Section~\ref{ss:K1}), on which no coverings are needed.
	\end{enumerate}
	The reason to use coverings rather than arbitrary vector bundles is the following. A general vector bundle $V$ has no canonical global diagonal
	in $\Gamma(\End V)$. If $V=\bigoplus_iL_i$ is a direct sum of line bundles, the summand diagonal $\bigoplus_i\Gamma(\End L_i)$ of $\Gamma(\End V)$ can have a
	nontrivial Weyl twist (Remark~\ref{ex:hopf}), while for a permutation bundle $V_\pi$ the sheet diagonal gives the pair $(F_\pi,D_\pi)$, which is the untwisted
	groupoid $C^*$-algebra of the covering groupoid $R_\pi$ (Proposition~\ref{prop:coveringalgebra}). The whole point of the present section is to stay inside
	the class of permutation bundles.
	
	\subsubsection*{Origin of the difficulty}
	By Remark~\ref{rem:dim3}, on a complex of dimension at most $3$ every reduced covering class has order at most $2$, and by
	Example~\ref{ex:P2p} the Moore space $P^2(p)$ has torsion $\Z/p$ in $\Kt^0$ but no covering represents a nonzero class.
	More generally, a covering class is always in the image of the map $[T,B\Sigma_N]\to\Kt^0(T)$ induced by the permutation
	representation (Section~\ref{ss:cov3}), and this is a strong restriction: if the monodromy of the covering factors through a finite group $H$, its class is fixed by
	the Adams operation $\psi^k$ for every $k$ prime to $|H|$ (the character of $\psi^kV$ at $h$ is the character of $V$ at $h^k$, and for such $k$ the
	permutations $h$ and $h^k$ have the same fixed points); and the classes are in the image of the stable cohomotopy of $T$ (Section~\ref{ss:BPQ}).
	They need not be fixed by every $\psi^k$: for the double covering $S^2\to\R P^2$ the reduced class is $x=[L]-1\neq0$, where $L$ is the
	complexified tautological line bundle (Example~\ref{ex:RP2}), and $\psi^2(x)=[L^{\otimes2}]-1=0$. On the other hand a simply connected space has no nontrivial
	connected coverings at all, so the classes must be carried by the fundamental group.
	
	So two requirements pull in opposite directions. To control $K$-theory and maps between the spaces we would like simple, highly connected spaces such as
	Moore spaces, whose $K$-theory is computed by their cells and on which maps are constructed by homotopy theory. But coverings need a large fundamental
	group, and the $K$-theory classes they carry must satisfy the restrictions above. The construction reconciles the two by replacing a Moore space $X$ by a
	finite complex $T$ with the same homology but with a large \emph{perfect} fundamental group, connected to $X$ by an \emph{acyclic} map $T\to X$. Such a map
	does not change homology or $K$-theory, but $T$ has many more coverings than $X$.
	
	\subsubsection*{The method}
\begin{figure}[ht]
		\centering
		\begin{tikzcd}[column sep=3.2em,row sep=3.0em]
			& & B\Sigma_\infty \arrow[d,"b"] & B\Sigma_N\arrow[l,hook']\arrow[r,"\text{perm}"] & BU(N)\\
			T_n \arrow[r,"\gamma"] \arrow[rd,"a_n"'] \arrow[rru,bend left=18,"c"] & P \arrow[d,"q"] \arrow[r] & Q_0S^0 \arrow[r,"\Omega^\infty\eta"] & BU &\\
			& X_n=\bigvee_jP^{d_n}(p_j^{r_j}) \arrow[ru,"\alpha"'] & & &
		\end{tikzcd}
		\caption{The construction of $T_n$ (one class shown). The maps $b$, $q$ and $a_n$ are acyclic. $\alpha$ is a stable cohomotopy lift of a class $x\in\Kt^0(X_n)$, $P$ is
			the homotopy pullback of $b$ along $\alpha$, and $T_n$ is a finite complex with an acyclic map to $X_n$ through $P$. The map $c$
			factors through some $B\Sigma_N$ and classifies a covering whose reduced class is $a_n^*x$.}\label{fig:strategy}
	\end{figure}
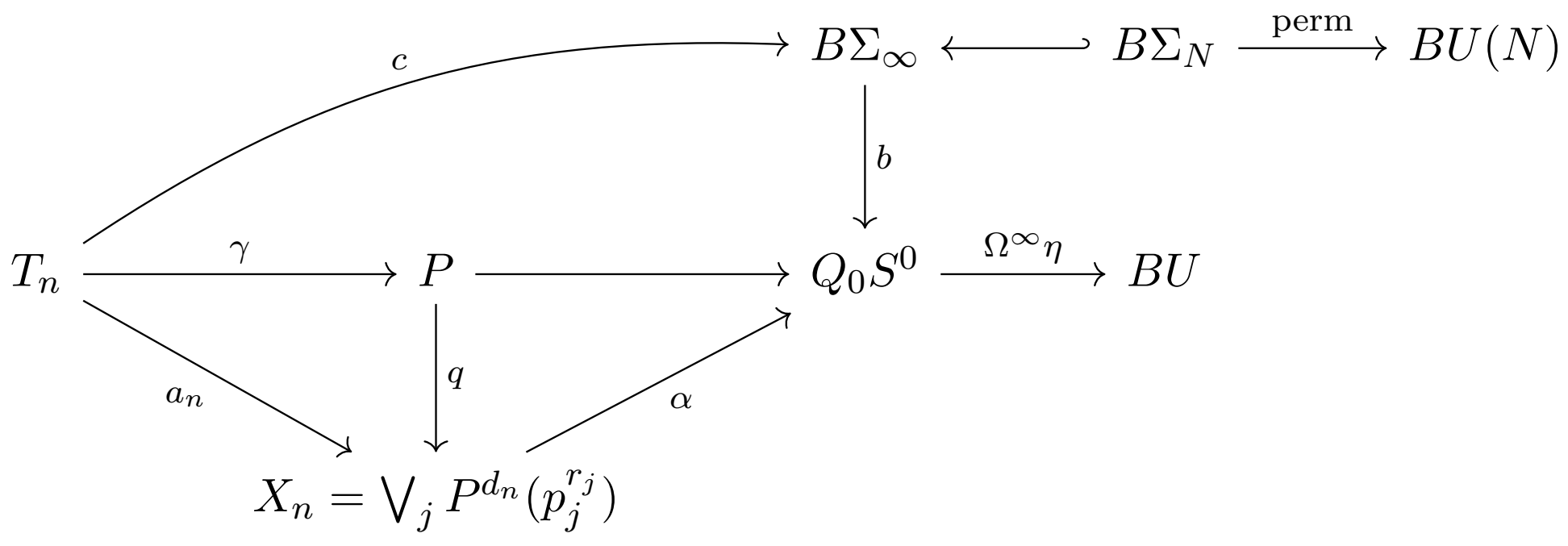

	\begin{enumerate}[label=(\arabic*)]
		\item \emph{Moore spaces whose torsion lifts to stable cohomotopy.} We realize $E_n$ as $\Kt^0(X_n)$ for a wedge $X_n$ of Moore
		spaces $P^{d}(p^r)$ (Section~\ref{ss:moore}). If $d$ is a multiple of Adams' period $q(p,r)$, every class in
		$\Kt^0(P^d(p^r))$ comes from the stable cohomotopy of $P^d(p^r)$ (Lemma~\ref{lem:lift}); this is where Adams' periodic maps
		(Section~\ref{ss:adams}) enter. Maps $X_{n+1}\to X_n$ realizing $i_n$ are built from degree maps and desuspended Adams maps
		(Lemma~\ref{lem:wedgemaps}). The Moore spaces are chosen of even dimension $d\ge4$, so that $K^1(X_n)=0$ and $X_n$ is $2$-connected.
		\item \emph{From stable cohomotopy to coverings.} By the Barratt--Priddy--Quillen theorem (Section~\ref{ss:BPQ}), the
		infinite loop space $Q_0S^0$ of stable cohomotopy is the plus construction of $B\Sigma_\infty$, and the unit map to $K$-theory
		corresponds to the permutation representation. A class in stable cohomotopy is therefore a map $\alpha\colon X_n\to Q_0S^0$, and pulling back
		along the acyclic map $b\colon B\Sigma_\infty\to Q_0S^0$ gives a space $P$ with an acyclic map $q\colon P\to X_n$ and a map $P\to B\Sigma_\infty$ which lifts
		$\alpha$ up to homotopy. The space $P$ is not compact, and a map $P\to B\Sigma_\infty$ need not factor through any finite $B\Sigma_N$; so $P$ does not yet
		carry finite coverings.
		\item \emph{Finiteness.} The Hausmann--Vogel theorem (Section~\ref{ss:HV}) gives a finite complex $T_n$ with a map $\gamma\colon T_n\to P$ such that
		$a_n=q\circ\gamma\colon T_n\to X_n$ is still acyclic. Since $T_n$ is compact, the composite $T_n\to P\to B\Sigma_\infty$ lands in some $B\Sigma_N$
		and classifies an $N$-sheeted covering of $T_n$ whose reduced class is $a_n^*x$. Acyclic maps induce isomorphisms in $K$-theory, so
		$\Kt^0(T_n)\cong\Kt^0(X_n)\cong E_n$.
	\end{enumerate}
	Figure~\ref{fig:strategy} summarizes the construction. To build $T_{n+1}$ together with a map $g_n\colon T_{n+1}\to T_n$ compatible with $f_n\colon
	X_{n+1}\to X_n$, steps (2) and (3) are carried out simultaneously for the classes on $X_{n+1}$ and for the acyclic map $a_n\colon T_n\to X_n$
	(Proposition~\ref{prop:lift}). A worked example, for the groups $\Z/3$ and $\Z/9$, is Example~\ref{ex:smallcase}.
	
	\subsubsection*{Guide to the section}
	Sections~\ref{ss:CW}--\ref{ss:acyclic} collect background from homotopy theory: CW complexes and homotopy groups (\ref{ss:CW}), topological $K$-theory
	(\ref{ss:Ktheory}), classifying spaces and the classification of coverings (\ref{ss:cov3}), Moore spaces (\ref{ss:moore}), spectra and duality (\ref{ss:spectra}),
	Adams' periodic maps (\ref{ss:adams}) and acyclic maps (\ref{ss:acyclic}). This material is standard. The theorems are stated with references and
	explanations of their content; the lemmas, which are formulated for our needs, are proved. The construction itself
	occupies Sections~\ref{ss:BPQ}--\ref{ss:Zn}: the Barratt--Priddy--Quillen theorem and the fundamental group of the homotopy fibre of $b$
	(\ref{ss:BPQ}), the Hausmann--Vogel theorem (\ref{ss:HV}), the simultaneous lifting proposition (\ref{ss:lift}), torsion stages and the realization
	theorem (\ref{ss:realize}), the complexes for $K_1$ (\ref{ss:K1}) and the spaces $Z_n$ used in Section~\ref{sec:blocks} (\ref{ss:Zn}).
	
	\emph{Conventions.} CW complexes are based at a $0$-cell. Maps between CW complexes are continuous; ``based'' means basepoint preserving. For a connected
	finite CW complex $X$ we identify $K^0(X)=\Z\oplus\Kt^0(X)$ via the rank. $\Sigma$ denotes reduced suspension (of spaces or spectra) and, depending on the
	context, the symmetric group $\Sigma_N$; no confusion will arise.
	
	\subsection{CW complexes and homotopy}\label{ss:CW}
	
	We recall the notions and theorems of homotopy theory that are used; standard references are \cite{Hat} and \cite{Wh78}.
	
	\begin{definition}[CW complexes]
		A \emph{CW complex} is a Hausdorff space $X$ with a filtration $X^0\subset X^1\subset\cdots$ by closed subspaces (the
		\emph{skeleta}) with $X=\bigcup_kX^k$, such that $X^0$ is discrete, each $X^k$ is obtained from $X^{k-1}$ by attaching
		$k$-dimensional discs along their boundaries,
		\[
		X^k=X^{k-1}\cup_{\bigsqcup\varphi_\alpha}\textstyle\bigsqcup_\alpha D^k_\alpha,\qquad\varphi_\alpha\colon\partial D^k_\alpha=S^{k-1}\to X^{k-1},
		\]
		and $X$ has the weak topology (a set is closed iff it meets every closed cell in a closed set). The images $e^k_\alpha$ of the
		open discs are the \emph{cells}, the maps $\varphi_\alpha$ the \emph{attaching maps}, and the maps $\chi_\alpha\colon D^k_\alpha\to
		X$ the \emph{characteristic maps}. $X$ is \emph{finite} if it has finitely many cells, and its \emph{dimension} is the largest
		dimension of a cell. A \emph{subcomplex} is a closed union of cells which is itself a CW complex with the induced cells. A map
		$f\colon X\to Y$ of CW complexes is \emph{cellular} if $f(X^k)\subset Y^k$ for all $k$.
	\end{definition}
	
	A finite CW complex is compact and metrizable; CW complexes are locally contractible, in particular locally path-connected
	\cite[Appendix]{Hat}. A connected CW complex is path-connected. A compact subset of a CW complex meets only finitely many cells, and hence lies in a
	finite subcomplex \cite[Proposition~A.1]{Hat}; this will be used to pass from the infinite complex $B\Sigma_\infty$ to a finite $B\Sigma_N$. The product
	of two CW complexes, with the product topology, is a CW complex whose cells are the products of cells if one of the factors is locally compact or if both
	have countably many cells \cite[Theorem~A.6]{Hat}.
	
	\begin{definition}[Based spaces and constructions]
		A \emph{based space} is a space with a chosen point, the \emph{basepoint}; for CW complexes we always choose a $0$-cell.
		\emph{Based maps} preserve basepoints, and \emph{based homotopies} are homotopies through based maps. We write $[A,B]$ for free
		homotopy classes of maps and $[A,B]_*$ for based homotopy classes. For based spaces $(A,a_0)$, $(B,b_0)$:
		\begin{itemize}
			\item the \emph{wedge} $A\vee B=A\sqcup B/(a_0\sim b_0)$;
			\item the \emph{smash product} $A\wedge B=A\times B/A\vee B$;
			\item the \emph{reduced suspension} $\Sigma A=S^1\wedge A$, i.e.\ $A\times[0,1]$ with $A\times\{0,1\}\cup\{a_0\}\times[0,1]$
			collapsed to a point; the \emph{unreduced suspension} $SA$ is $A\times[0,1]$ with $A\times\{0\}$ and $A\times\{1\}$ collapsed to two points;
			\item the \emph{loop space} $\Omega B$ of based loops $[0,1]\to B$, with the compact-open topology; there is a natural
			bijection $[\Sigma A,B]_*\cong[A,\Omega B]_*$.
		\end{itemize}
		If $A,B$ are CW complexes based at $0$-cells, then $A\vee B$, $A\wedge B$ and $\Sigma A$ are CW complexes in a natural way
		(products are taken in the category of compactly generated spaces, which changes nothing for finite complexes). For a CW complex $A$ based at a $0$-cell,
		the quotient map $SA\to\Sigma A$ collapses the contractible subcomplex $\{a_0\}\times[0,1]$ and is therefore a homotopy equivalence
		\cite[Proposition~0.17]{Hat}.
	\end{definition}
	
	$[\Sigma A,B]_*$ is a group (concatenation of the suspension coordinate), abelian if $A$ is itself a suspension. The
	\emph{pinch map} $\mu\colon\Sigma A\to\Sigma A\vee\Sigma A$ collapses $A\times\{\frac12\}$; the sum of two based maps
	$u_1,u_2\colon\Sigma A\to W$ is $u_1+u_2=\nabla\circ(u_1\vee u_2)\circ\mu$, where $\nabla$ is the fold map (Figure~\ref{fig:pinch}).
	
	\begin{definition}[Connectivity]
		A space $B$ is \emph{$n$-connected} if it is nonempty, path-connected and $\pi_i(B)=0$ for $1\le i\le n$. A map $f\colon A\to B$
		is \emph{$n$-connected} if $f_*\colon\pi_i(A)\to\pi_i(B)$ is an isomorphism for $i<n$ and onto for $i=n$, for all basepoints.
	\end{definition}
	
	\begin{definition}[Cofibrations]\label{def:cofibration}
		A map $i\colon A'\to A$ is a \emph{cofibration} if it has the \emph{homotopy extension property}: for every space $Y$, every map $f\colon A\to Y$ and every
		homotopy $h\colon A'\times[0,1]\to Y$ with $h_0=f\circ i$, there is a homotopy $H\colon A\times[0,1]\to Y$ with $H_0=f$ and $H_t\circ i=h_t$ for all $t$. For
		the inclusion of a closed subspace $A'\subset A$ this means that $A\times\{0\}\cup A'\times[0,1]$ is a retract of $A\times[0,1]$ \cite[Chapter~0]{Hat}.
	\end{definition}
	
	\begin{theorem}[Standard results {\cite[\S\S4.1--4.2]{Hat}}]\label{thm:standardhtpy}
		Let $A,B$ be CW complexes.
		\begin{enumerate}[label=\textup{(\alph*)}]
			\item \emph{Cellular approximation.} Every map $A\to B$ is homotopic, relative to any subcomplex on which it is already
			cellular, to a cellular map. In particular, if $B$ has no cells of dimension $1,\dots,n$ apart from a single $0$-cell (e.g.
			$B^n=\pt$), then $B$ is $n$-connected; and if $\dim A<m$ and $B^{m-1}$ is a point, every based map $A\to B$ is null-homotopic.
			\item \emph{Homotopy extension.} A homotopy of a map $A\to B$ defined on a subcomplex $A'\subset A$ extends to a homotopy of the
			whole map. In particular, every map from a connected CW complex to a path-connected space is homotopic to a based map.
			\item \emph{Whitehead's theorem.} A map between connected CW complexes inducing isomorphisms on all homotopy groups is a homotopy
			equivalence. A map between simply connected CW complexes inducing isomorphisms on all integral homology groups is a homotopy
			equivalence.
			\item \emph{Maps into a highly connected map.} If $f\colon B\to B'$ is $n$-connected and $A$ is a CW complex, then
			$f_*\colon[A,B]_*\to[A,B']_*$ is bijective if $\dim A<n$ and onto if $\dim A\le n$.
			\item \emph{CW approximation.} For every space $P$ there are a CW complex $\widetilde P$ and a weak homotopy equivalence
			$\widetilde P\to P$ (a map inducing isomorphisms on all homotopy groups at all basepoints).
		\end{enumerate}
	\end{theorem}
	
	We comment on the content of these statements. (a) says that, up to homotopy, maps between CW complexes respect the skeleta. Its two consequences are
	the ones we use: a complex built with no cells in dimensions $1,\dots,n$ has no homotopy in these dimensions, and a low-dimensional complex cannot map
	essentially into a complex whose cells are all of high dimension. This is how we see that Moore spaces $P^d(m)$ with $d\ge4$ are $2$-connected, and that
	there are no essential maps from lower to higher dimensional Moore spaces (Remark~\ref{rem:monotone}). (b) says that inclusions of subcomplexes are
	cofibrations (Definition~\ref{def:cofibration}); we use it to replace maps and homotopies by based ones. (c) is the criterion for homotopy equivalence; the homological version for simply
	connected complexes is the one used in the proof of Proposition~\ref{prop:lift}, applied to suspensions. (d) says that an $n$-connected map is invisible to
	complexes of dimension less than $n$, and that all homotopy classes from $n$-dimensional complexes lift through it. Applied to the inclusion of a space
	into $\Omega\Sigma$ of itself, it gives the Freudenthal theorem in the form below, and it is the reason why stable maps between Moore spaces of large
	dimension can be desuspended. (e) allows us to apply results stated for CW complexes to spaces such as homotopy pullbacks.
	
	\begin{theorem}[Seifert--van Kampen theorem {\cite[Theorem~1.20, Proposition~1.26]{Hat}}]\label{thm:vankampen}
		\begin{enumerate}[label=\textup{(\alph*)}]
			\item Let $X=U\cup V$ with $U,V\subset X$ open, and suppose that $U$, $V$ and $U\cap V$ are path-connected; let $x_0\in U\cap V$. Then the homomorphism
			$\pi_1(U,x_0)*\pi_1(V,x_0)\to\pi_1(X,x_0)$ induced by the inclusions is surjective, and its kernel is the normal subgroup generated by the elements
			$i_U(\omega)\,i_V(\omega)^{-1}$, $\omega\in\pi_1(U\cap V,x_0)$, where $i_U,i_V$ are induced by the inclusions of $U\cap V$. In other words,
			$\pi_1(X)$ is the amalgamated free product $\pi_1(U)*_{\pi_1(U\cap V)}\pi_1(V)$.
			\item Let $Y$ be obtained from a path-connected space $X$ by attaching $2$-cells along loops $\varphi_\alpha\colon S^1\to X$. Then $\pi_1(Y)\cong\pi_1(X)/N$,
			where $N$ is the normal subgroup generated by the classes of the loops $\varphi_\alpha$ (transported to the basepoint along paths). Attaching cells of
			dimension at least $3$ does not change $\pi_1$.
			\item If $A$ is path-connected, its unreduced suspension $SA$ is simply connected; if $A$ is a connected CW complex based at a $0$-cell, its reduced
			suspension $\Sigma A$ is simply connected.
		\end{enumerate}
	\end{theorem}
	
	Part (a) is the theorem of Seifert and van Kampen; part (b) is its standard application to cell attachments \cite[Proposition~1.26]{Hat}. Part (c) follows
	from (a): $SA$ is the union of the two open cones $U=A\times[0,1)/A\times\{0\}$ and $V=A\times(0,1]/A\times\{1\}$, which are contractible, and $U\cap V=A\times(0,1)$
	is path-connected; so $\pi_1(SA)$ is the amalgamated product of two trivial groups, hence trivial. For a connected CW complex $A$ based at a $0$-cell, $\Sigma A$ is
	homotopy equivalent to $SA$ (see above), hence simply connected. We use (b) in the computation of $\pi_1(P^2(m))$ (Lemma~\ref{lem:flatP2}) and (c) for Moore spaces
	of dimension at least $3$ and in the proof of Proposition~\ref{prop:lift}.
	
	\begin{theorem}[Freudenthal suspension theorem {\cite{Fr}}, {\cite[\S4.2]{Hat}}, {\cite[Chapter~VII]{Wh78}}]\label{thm:freudenthal}
		Let $B$ be an $(n-1)$-connected CW complex, $n\ge2$. Then the natural map $B\to\Omega\Sigma B$ is $(2n-1)$-connected.
		Consequently, for every CW complex $A$ of dimension at most $2n-1$, the suspension
		\[
		\Sigma\colon[A,B]_*\to[\Sigma A,\Sigma B]_*
		\]
		is onto, and it is bijective if $\dim A\le2n-2$.
	\end{theorem}
	
	The first statement is the classical statement that $\pi_i(B)\to\pi_{i+1}(\Sigma B)$ is an isomorphism for $i<2n-1$ and onto for $i=2n-1$. The second
	is obtained from it by Theorem~\ref{thm:standardhtpy}(d) and the adjunction $[\Sigma A,\Sigma B]_*=[A,\Omega\Sigma B]_*$. In words: in the
	\emph{stable range} $\dim A\le2\cdot\mathrm{conn}(B)$, homotopy classes of maps $A\to B$ do not change under suspension. This is what makes the
	passage between maps of spaces and maps of spectra (Section~\ref{ss:spectra}) possible. A consequence used repeatedly: if $B$ is $(n-1)$-connected, then
	$\Sigma B$ is $n$-connected, because $\pi_i(\Sigma B)\cong\pi_{i-1}(\Omega\Sigma B)\cong\pi_{i-1}(B)=0$ for $1\le i-1\le n-1$ and $\pi_1(\Sigma B)=0$ by
	Theorem~\ref{thm:vankampen}(c).
	
	\begin{definition}[Fibrations and homotopy fibres]
		A map $p\colon E\to B$ is a \emph{(Hurewicz) fibration} if it has the homotopy lifting property for all spaces. The
		\emph{homotopy fibre} of an arbitrary map $f\colon A\to B$ over $b\in B$ is
		\[
		\hofib_f(b)=\{(a,\omega)\in A\times B^{[0,1]}:\ \omega(0)=f(a),\ \omega(1)=b\}.
		\]
		The map $A\to B^{[0,1]}\times_BA$, $a\mapsto(a,\text{const}_{f(a)})$, is a homotopy equivalence, and
		$\{(a,\omega)\}\to B$, $(a,\omega)\mapsto\omega(1)$, is a fibration whose fibre over $b$ is $\hofib_f(b)$; so every map is,
		up to homotopy equivalence, a fibration \cite[\S4.3]{Hat}.
	\end{definition}
	
	If $B$ is path-connected, the homotopy fibres over different points are homotopy equivalent (transport along a path), so we may speak of \emph{the}
	homotopy fibre of $f$. If $f$ is itself a fibration, its fibre $f^{-1}(b)$ is homotopy equivalent to $\hofib_f(b)$.
	
	\begin{theorem}[Exact sequence of a fibration {\cite[Theorem~4.41]{Hat}}, {\cite[\S4.1]{Hat}}]\label{thm:fibLES}
		Let $F\to E\xrightarrow{p}B$ be a fibration with $B$ path-connected, $F=p^{-1}(b_0)$ and $e_0\in F$.
		\begin{enumerate}[label=\textup{(\alph*)}]
			\item There is a long exact sequence
			\begin{equation}\label{eq:LES}
				\cdots\to\pi_{i+1}(B)\xrightarrow{\partial}\pi_i(F)\to\pi_i(E)\to\pi_i(B)\xrightarrow{\partial}\pi_{i-1}(F)\to\cdots\to\pi_0(E)\to\pi_0(B),
			\end{equation}
			with basepoints $e_0$ and $b_0$.
			\item The image of $\partial\colon\pi_2(B,b_0)\to\pi_1(F,e_0)$ is central in $\pi_1(F,e_0)$.
		\end{enumerate}
	\end{theorem}
	
	Part (a) is the basic computational tool for homotopy groups. The boundary map $\partial$ is defined by lifting: a based map $(D^{i+1},S^i)\to(B,b_0)$
	representing an element of $\pi_{i+1}(B)$ (with $S^i$ collapsed) lifts to a map $D^{i+1}\to E$ sending the boundary sphere into $F$, and $\partial$ of the
	element is the class of this boundary map in $\pi_i(F)$. The end of the sequence consists of pointed sets: exactness at $\pi_0(E)$ means that two path
	components of $E$ have the same image in $\pi_0(B)$ if and only if they both meet $F$. Part (b) is a refinement in the lowest degree, where $\pi_1(F)$ need not
	be abelian: the boundary map takes values in the centre. This reflects the fact that $\partial$ is equivariant for the action of $\pi_1(F)$ on the relative
	group $\pi_2(E,F)\cong\pi_2(B)$ and its action on itself by conjugation, and that a loop in $F$ acts on $\pi_2(B)$ through its image in $\pi_1(B)$, which is trivial
	\cite[\S4.1]{Hat}. We use (b) to see that the fundamental group of the homotopy fibre of the Barratt--Priddy--Quillen map is a \emph{central} extension of
	$A_\infty$ (Lemma~\ref{lem:Fb}).
	
	\subsection{Topological \texorpdfstring{$K$}{K}-theory}\label{ss:Ktheory}
	
	Standard references are \cite{At67}, \cite{Ka78} and \cite{Hus}.
	
	\begin{definition}
		For a compact Hausdorff space $X$, $K^0(X)$ is the Grothendieck group of the abelian monoid of isomorphism classes of complex
		vector bundles over $X$ under direct sum. By the Serre--Swan theorem \cite{Swan}, $K^0(X)=K_0(C(X))$. We put
		$K^1(X)=K_1(C(X))$; equivalently $K^1(X)=\Kt^0(\Sigma X_+)$, where $X_+$ is $X$ with a disjoint basepoint. For a based compact
		space $(X,x_0)$, the \emph{reduced group} $\Kt^0(X)$ is the kernel of restriction $K^0(X)\to K^0(x_0)=\Z$, and
		$\Kt^1(X)=K^1(X)$ (because $K^1(\pt)=0$).
	\end{definition}
	
	If $X$ is connected, the rank of a vector bundle is constant, the rank homomorphism $K^0(X)\to\Z$ does not depend on the
	point at which it is computed, and $K^0(X)=\Z\oplus\Kt^0(X)$. Every element of $\Kt^0(X)$ is of the form $[V]-\rank V$ for a vector bundle $V$, since
	every vector bundle over a compact space is a direct summand of a trivial bundle. We use the following standard facts \cite{At67}.
	\begin{enumerate}[label=(K\arabic*)]
		\item \emph{Homotopy invariance:} homotopic maps (free homotopy) induce the same homomorphism on $K^*$.
		\item \emph{Exactness:} for a based map $f\colon A\to X$ of finite CW complexes with mapping cone $C_f=X\cup_f CA$, there is an exact sequence
		\[
		\Kt^0(\Sigma X)\xrightarrow{(\Sigma f)^*}\Kt^0(\Sigma A)\xrightarrow{c^*}\Kt^0(C_f)\xrightarrow{\ i^*\ }\Kt^0(X)\xrightarrow{f^*}\Kt^0(A),
		\]
		in which $i\colon X\to C_f$ is the inclusion and $c\colon C_f\to C_f/X=\Sigma A$ the collapse; and similarly with $\Kt^1$. Together these form a long exact sequence.
		\item \emph{Suspension and Bott periodicity:} $\Kt^{j+1}(\Sigma X)\cong\Kt^j(X)$ naturally, and $\Kt^{j+2}\cong\Kt^j$. Hence
		$\Kt^0(S^d)=\Z$, $K^1(S^d)=0$ for even $d\ge2$, and $\Kt^0(S^d)=0$, $K^1(S^d)=\Z$ for odd $d$.
		\item \emph{Wedges:} if $A=A_1\vee\dots\vee A_k$ is a wedge of finite CW complexes, the inclusions $\iota_j$ induce an
		isomorphism $\Kt^*(A)\to\bigoplus_j\Kt^*(A_j)$.
		\item \emph{Representability:} for a finite CW complex $A$, $K^0(A)\cong[A,\Z\times BU]$ and, if $A$ is connected and based,
		$\Kt^0(A)\cong[A,BU]_*$. Here $BU(N)$ is the classifying space of the unitary group $U(N)$ (Section~\ref{ss:cov3}), which classifies rank $N$
		bundles: a rank $N$ bundle $V$ over $A$ is $c^*\gamma_N$ for the universal bundle $\gamma_N$ over $BU(N)$ and a map $c\colon A\to BU(N)$ unique up to homotopy;
		$BU=\operatorname{colim}_NBU(N)$ along the maps induced by the inclusions $U(N)\subset U(N+1)$; and $[V]-N\in\Kt^0(A)$ corresponds to the composite
		$A\to BU(N)\to BU$ \cite[Chapter~1]{At67}.
	\end{enumerate}
	Exactness (K2) and the suspension isomorphism (K3) say that $\Kt^*$ is a cohomology theory on finite complexes; this is what allows us to compute it cell by
	cell, as for ordinary cohomology. (K5) says that a reduced $K$-theory class on a connected complex is the same thing as a homotopy class of maps to $BU$; in
	this form $K$-theory classes can be compared with stable cohomotopy classes, which are maps to $Q_0S^0$ (Definition~\ref{def:stablelift}).
	
	\begin{lemma}\label{lem:K0basic}
		Let $u,u'\colon A\to B$ be maps of connected finite CW complexes.
		\begin{enumerate}[label=\textup{(\roman*)}]
			\item $u^*$ maps $\Kt^0(B)$ into $\Kt^0(A)$, whether or not $u$ is based.
			\item If $u\simeq u'$ by a free homotopy, then $u^*=u'^*$ on $K^*(B)$.
			\item $\Kt^0(A)=[A,BU]_*=[A,BU]$: reduced classes are free homotopy classes of maps $A\to BU$.
			\item If $A=A_1\vee\dots\vee A_k$ is a wedge at $0$-cells, with inclusions $\iota_j$ and collapses $c_j\colon A\to A_j$ (the identity
			on $A_j$ and constant on the other summands), then $(\iota_j^*)_j\colon\Kt^*(A)\to\bigoplus_j\Kt^*(A_j)$ is an isomorphism with
			inverse $(y_j)_j\mapsto\sum_jc_j^*y_j$.
		\end{enumerate}
	\end{lemma}
	
	\begin{proof}
		(i) Pullback of vector bundles preserves ranks pointwise, so $u^*$ maps a class of virtual rank zero to a class of virtual rank zero; this does not require
		$u$ to preserve basepoints, because on connected complexes the rank can be computed at any point. (ii) is (K1). (iii) Since $\pi_1(BU)=\pi_0(U)=0$, the
		space $BU$ is simply connected, and for a simply connected target the action of $\pi_1$ on based homotopy classes, whose orbits are the free homotopy classes,
		is trivial \cite[\S4.A]{Hat}; so $[A,BU]_*=[A,BU]$, and (K5) gives $\Kt^0(A)=[A,BU]$. (iv) The map $(\iota_j^*)_j$ is an isomorphism by (K4). Since
		$c_j\circ\iota_l$ is the identity of $A_j$ for $l=j$ and the constant map for $l\ne j$, and constant maps induce zero on reduced $K$-theory,
		$\iota_l^*\sum_jc_j^*y_j=y_l$; hence $(y_j)_j\mapsto\sum_jc_j^*y_j$ is a right inverse of the isomorphism $(\iota_j^*)_j$, and therefore its inverse.
	\end{proof}
	
	These are reformulations of (K1)--(K5) which free us from keeping track of basepoints. By (iii), in comparisons of $K$-theory classes we may use free
	homotopies, and the maps constructed in Proposition~\ref{prop:lift}, which commute only up to free homotopy, still give equalities of $K$-theory classes. (iv)
	says that a class on a wedge is determined by its restrictions to the summands, and that a family of classes on the summands is realized by pulling them back
	along the collapse maps and adding.
	
	\begin{lemma}[Pinch sums]\label{lem:pinch}
		Let $A=\Sigma A'$ with $A'$ a finite CW complex, $W$ a based finite CW complex, and $u_1,\dots,u_k\colon A\to W$ based maps. Then
		$(u_1+\dots+u_k)^*=u_1^*+\dots+u_k^*$ on $\Kt^*(W)$. For based maps $f_l\colon A_l\to W$ ($l=1,\dots,k$), the wedge map
		$\bigvee_lf_l\colon\bigvee_lA_l\to W$ satisfies $\iota_l^*(\bigvee_lf_l)^*=f_l^*$.
	\end{lemma}
	
	\begin{proof}
		By definition $u_1+u_2=\nabla\circ(u_1\vee u_2)\circ\mu$, so $(u_1+u_2)^*=\mu^*\circ(u_1\vee u_2)^*\circ\nabla^*$. For $y\in\Kt^*(W)$, the class $\nabla^*y\in\Kt^*(W\vee W)$
		restricts to $y$ on both summands, so $\nabla^*y=c_1^*y+c_2^*y$ by Lemma~\ref{lem:K0basic}(iv), where $c_1,c_2\colon W\vee W\to W$ are the collapses. Then
		$(u_1\vee u_2)^*c_k^*y=(c_k\circ(u_1\vee u_2))^*y=\bar c_k^*u_k^*y$, where $\bar c_k\colon\Sigma A'\vee\Sigma A'\to\Sigma A'$ is the collapse onto the $k$th summand, and
		$\bar c_k\circ\mu$ is homotopic to the identity of $\Sigma A'$ (it collapses one half of the suspension coordinate, which is homotopic to the identity by a
		reparametrization). Hence $(u_1+u_2)^*y=\mu^*\bar c_1^*u_1^*y+\mu^*\bar c_2^*u_2^*y=u_1^*y+u_2^*y$, and the case of $k$ summands follows by induction. For the
		second statement, $(\bigvee_lf_l)\circ\iota_l=f_l$, so $\iota_l^*(\bigvee_lf_l)^*=f_l^*$.
	\end{proof}
	
\begin{figure}[ht]
		\centering
		\begin{tikzpicture}[scale=0.9]
			\draw[thick] (0,0) ellipse (1.3 and 0.5);
			\draw[thick] (-1.3,0) -- (0,2) -- (1.3,0); \draw[thick] (-1.3,0) -- (0,-2) -- (1.3,0);
			\draw[red!70!black,thick,dashed] (0,0.5) arc (90:270:0.3 and 0.5);
			\node at (0,-2.5) {$\Sigma A'$};
			\draw[->,thick] (1.8,0) -- node[above]{$\mu$} (3.0,0);
			\begin{scope}[xshift=5cm]
				\draw[thick] (-1.0,1.0) ellipse (0.9 and 0.3); \draw[thick] (-1.9,1.0) -- (-1.0,2.2) -- (-0.1,1.0); \draw[thick] (-1.9,1.0)--(-1.0,0)--(-0.1,1.0);
				\draw[thick] (-1.0,-1.0) ellipse (0.9 and 0.3); \draw[thick] (-1.9,-1.0) -- (-1.0,0) -- (-0.1,-1.0); \draw[thick] (-1.9,-1.0)--(-1.0,-2.2)--(-0.1,-1.0);
				\fill (-1,0) circle (1.5pt);
				\node at (-1.0,-2.6) {$\Sigma A'\vee\Sigma A'$};
			\end{scope}
			\draw[->,thick] (5.4,0) -- node[above]{$u_1\vee u_2$} (6.9,0);
			\node at (7.7,0) {$W\vee W$};
			\draw[->,thick] (8.5,0) -- node[above]{$\nabla$} (9.4,0);
			\node at (9.8,0) {$W$};
		\end{tikzpicture}
		\caption{The pinch sum $u_1+u_2$: collapse the equator of the suspension, then apply $u_1$ to the upper and $u_2$ to the lower
			half.}\label{fig:pinch}
	\end{figure}
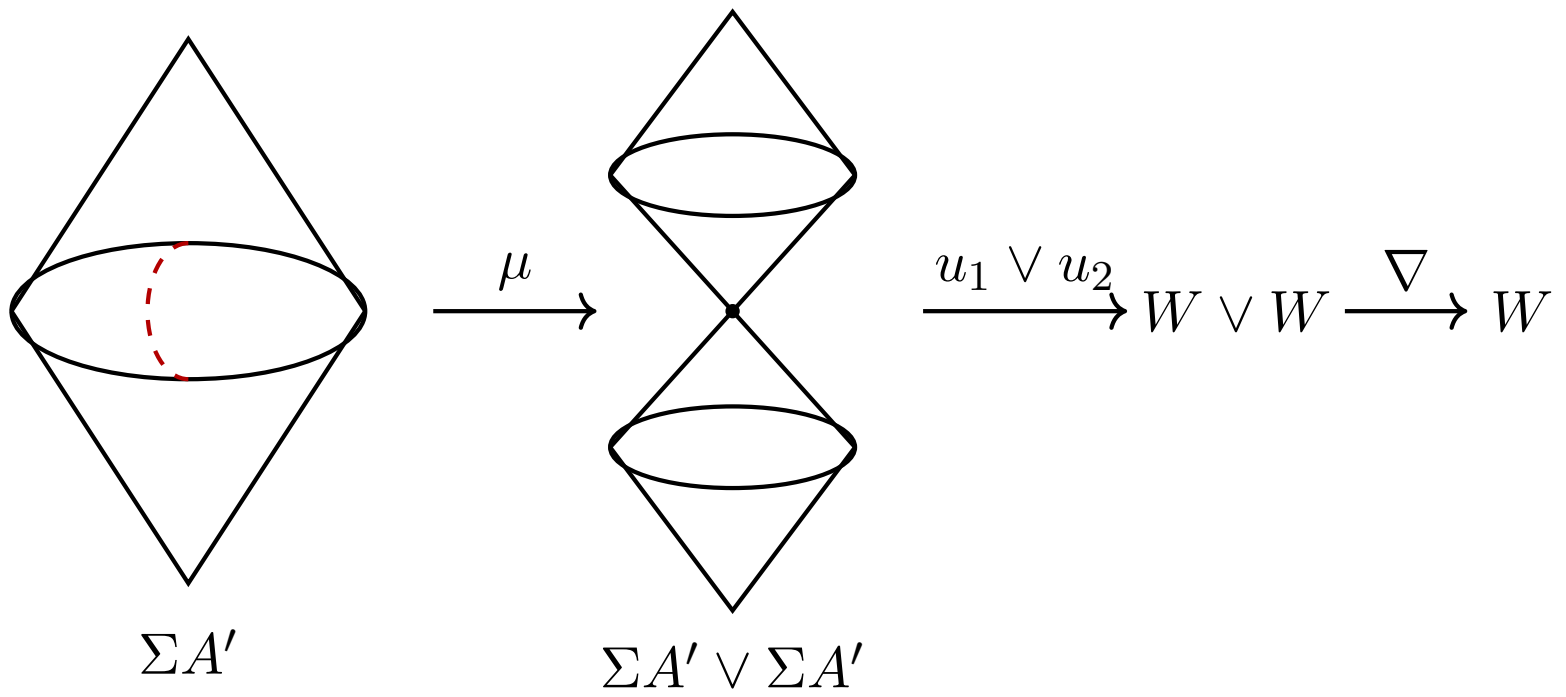

	The lemma says that the group structure of $[\Sigma A',W]_*$, defined geometrically by pinching (Figure~\ref{fig:pinch}), is turned into addition of
	homomorphisms by $K$-theory. It allows us to realize a sum of homomorphisms $\Kt^*(W)\to\Kt^*(A)$ by a single map, which is how maps between wedges of
	Moore spaces with prescribed effect on $K$-theory are assembled (Lemmas~\ref{lem:wedgemaps} and~\ref{lem:K1}). The second statement says that a map out of a
	wedge is determined, on $K$-theory, by its restrictions to the summands. For example, a based map $S^d\to S^d$ of degree $m$ is homotopic to the $m$-fold
	pinch sum of the identity, which also has degree $m$ \cite[Corollary~4.25]{Hat}, so it induces multiplication by $m$ on $\Kt^*(S^d)$.
	
	\subsection{Classifying spaces and finite coverings}\label{ss:cov3}
	
	\subsubsection*{Classifying spaces of discrete groups}
	\begin{definition}[Nerve]\label{def:nerve}
		The \emph{nerve} of a small category $\mathcal C$ is the simplicial set $N\mathcal C$ whose $k$-simplices are the strings
		$c_0\xrightarrow{f_1}c_1\xrightarrow{f_2}\cdots\xrightarrow{f_k}c_k$ of $k$ composable morphisms; the face maps compose two adjacent morphisms (or delete the
		first or the last one), and the degeneracy maps insert an identity. Its geometric realization $|N\mathcal C|$ is a CW complex with one $k$-cell for every
		nondegenerate $k$-simplex, that is, every string containing no identity morphism. A group or monoid $G$ is a category with one object whose morphisms are the
		elements of $G$, composed by multiplication; a $k$-simplex of its nerve is a $k$-tuple $(g_1,\dots,g_k)$ of elements of $G$. For a topological category the
		sets of simplices carry topologies, and the realization is formed accordingly.
	\end{definition}
	
	Let $G$ be a discrete group. We use the following model of its classifying space, the realization of the nerve of $G$ (Definition~\ref{def:nerve}) \cite[Example~1B.7]{Hat}.
	Let $EG$ be the geometric realization of the simplicial set whose $k$-simplices are the $(k+1)$-tuples $(g_0,\dots,g_k)$ of elements of $G$, with face
	maps given by omitting an entry and degeneracy maps given by repeating one. Thus $EG$ is a CW complex with one $k$-cell for every
	nondegenerate $(k+1)$-tuple (one with $g_{i}\ne g_{i+1}$ for all $i$); it is contractible, and $G$ acts on it freely and cellularly by
	$g\cdot(g_0,\dots,g_k)=(gg_0,\dots,gg_k)$. The quotient
	\[
	BG=EG/G
	\]
	is a CW complex, the \emph{classifying space} of $G$, and the quotient map $EG\to BG$ is a covering with group $G$ (a universal principal $G$-bundle). The
	$k$-cells of $BG$ correspond to the tuples $(1,g_1,g_1g_2,\dots,g_1\cdots g_k)$, i.e.\ to the $k$-simplices $(g_1,\dots,g_k)$ of the nerve of $G$ regarded
	as a category with one object. We shall use the following properties.
	\begin{enumerate}[label=(B\arabic*)]
		\item $\pi_1(BG)\cong G$ and $\pi_k(BG)=0$ for $k\ge2$: by the exact sequence \eqref{eq:LES} of the covering $EG\to BG$ with contractible total space and
		discrete fibre $G$. Thus $BG$ is an Eilenberg--MacLane space $K(G,1)$, and its homotopy type is determined by $G$ \cite[Theorem~1B.8]{Hat}.
		\item For a connected CW complex $A$ with basepoint $a_0$, the assignment $c\mapsto\pi_1(c)$ is a bijection from the set of based homotopy classes of based
		maps $c\colon(A,a_0)\to(BG,*)$ onto the set of homomorphisms $\pi_1(A,a_0)\to G$; free homotopy classes of maps $A\to BG$ correspond to conjugacy classes of
		homomorphisms \cite[Proposition~1B.9]{Hat}.
		\item A homomorphism $\theta\colon G\to H$ induces a cellular based map $B\theta\colon BG\to BH$, $(g_0,\dots,g_k)\mapsto(\theta(g_0),\dots,\theta(g_k))$, with
		$\pi_1(B\theta)=\theta$ and $B(\theta'\theta)=B\theta'\circ B\theta$. If $\theta$ is injective, $B\theta$ is a homeomorphism onto a subcomplex; we then regard $BG$
		as a subcomplex of $BH$. If $\theta$ is the inner automorphism $x\mapsto gxg^{-1}$ of $G$, then $B\theta$ is freely homotopic to the identity, by (B2).
		\item If $G$ and $H$ are countable, the natural map $B(G\times H)\to BG\times BH$ is a homeomorphism, the product being a CW complex by \cite[Theorem~A.6]{Hat}
		(the geometric realization commutes with finite products of countable simplicial sets).
		\item If $G=\bigcup_NG_N$ is the union of an increasing sequence of subgroups, then $BG=\bigcup_NBG_N$, with the $BG_N$ subcomplexes by (B3).
	\end{enumerate}
	For the symmetric groups $\Sigma_N$ of the sets $[N]=\{1,\dots,N\}$ we regard $\Sigma_N\subset\Sigma_{N+1}$ via $\sigma\mapsto\sigma\oplus1$ (extension by the
	identity on $N+1$), so that $B\Sigma_N\subset B\Sigma_{N+1}$ is a subcomplex, and we put $\Sigma_\infty=\bigcup_N\Sigma_N$, the group of permutations of $\N$ moving
	finitely many points, with $B\Sigma_\infty=\bigcup_NB\Sigma_N$ (B5). Its alternating subgroup is $A_\infty=\bigcup_NA_N$. Another model of $B\Sigma_N$ is the space of
	$N$-element subsets of $\R^\infty$: the space of ordered $N$-tuples of distinct points of $\R^\infty$ is contractible with a free action of $\Sigma_N$, and the quotient
	is a $K(\Sigma_N,1)$.
	
	\subsubsection*{Classifying spaces of topological groups}
	The same construction applies to a topological group $G$ such as $U(N)$, with the realization of the nerve now formed in the category of compactly generated
	spaces: $EG\to BG$ is a numerable principal $G$-bundle with contractible total space, and for a CW complex $A$ (or any paracompact space) pullback of this bundle
	gives a bijection between $[A,BG]$ and isomorphism classes of principal $G$-bundles over $A$ \cite[Chapter~4]{Hus}. For $G=U(N)$ this is the classification
	of rank $N$ Hermitian vector bundles by maps to $BU(N)$ used in (K5): the universal bundle $\gamma_N$ is $EU(N)\times_{U(N)}\C^N$. From the fibration
	$G\to EG\to BG$ and \eqref{eq:LES}, $\pi_k(BG)\cong\pi_{k-1}(G)$ for all $k\ge1$; in particular $\pi_k(BU)=\pi_{k-1}(U)$, which is $\Z$ for even $k\ge2$ and $0$
	for odd $k$ by Bott periodicity, and $BU(1)$ is a $K(\Z,2)$. These facts were used in the proof of Lemma~\ref{lem:lowdim}.
	
	\subsubsection*{Classification of finite coverings}
	We continue the list (C1), (C2) of Section~\ref{ss:Kcov}.
	\begin{enumerate}[label=(C\arabic*),start=3]
		\item \emph{Homotopy invariance.} If $A$ is a CW complex and $u\simeq u'\colon A\to X$, then $u^*\pi\cong u'^*\pi$ for every finite
		covering $\pi$ of $X$.
		\item \emph{Classification.} Let $U_N=E\Sigma_N\times_{\Sigma_N}[N]\to B\Sigma_N$ be the \emph{universal $N$-sheeted covering}: the quotient of $E\Sigma_N\times[N]$
		($[N]$ discrete) by the diagonal action $\sigma\cdot(e,i)=(\sigma e,\sigma(i))$, mapping to $B\Sigma_N$ by $[e,i]\mapsto[e]$. For a connected CW complex $A$,
		every $N$-sheeted covering of $A$ is isomorphic to $c^*U_N$ for a map $c\colon A\to B\Sigma_N$, unique up to free homotopy. Moreover
		$V_{c^*U_N}=c^*\C^N_\perm$, where $\C^N_\perm=E\Sigma_N\times_{\Sigma_N}\C^N$ is the vector bundle of the permutation representation of $\Sigma_N$ on $\C^N$;
		consequently the reduced class of $c^*U_N$ is the class of the composite
		\[
		A\xrightarrow{\ c\ }B\Sigma_N\xrightarrow{\ \perm\ }BU(N)\longrightarrow BU ,
		\]
		where $\perm=B(\Sigma_N\to U(N))$ is induced by the permutation representation.
	\end{enumerate}
	
	We explain these statements. Item (C3) is a consequence of the covering homotopy property of coverings \cite[Proposition~1.30]{Hat}: a covering of
	$A\times[0,1]$ is isomorphic to the pullback of its restriction to $A\times\{0\}$, so a homotopy $A\times[0,1]\to X$ pulls $\pi$ back to a covering which restricts
	to $u^*\pi$ and to $u'^*\pi$ at the two ends.
	
	For (C4), first note that $U_N\to B\Sigma_N$ is an $N$-sheeted covering: locally over $B\Sigma_N$ the covering $E\Sigma_N\to B\Sigma_N$ is $W\times\Sigma_N\to W$, and
	then $U_N$ is $W\times(\Sigma_N\times[N])/\Sigma_N\cong W\times[N]$. Next recall the classification of coverings by monodromy \cite[\S1.3]{Hat}: for a connected CW
	complex $A$ with universal covering $\widetilde A\to A$, every $N$-sheeted covering of $A$ is isomorphic to one of the form
	\[
	\widetilde A\times_{\rho}[N]=\big(\widetilde A\times[N]\big)/\pi_1(A,a_0),
	\]
	the quotient by the diagonal action of $\pi_1(A,a_0)$ through deck transformations on $\widetilde A$ and through a homomorphism $\rho\colon\pi_1(A,a_0)\to\Sigma_N$ (the
	\emph{monodromy}) on $[N]$; two such coverings $\widetilde A\times_\rho[N]$ and $\widetilde A\times_{\rho'}[N]$ are isomorphic if and only if $\rho$ and $\rho'$ are conjugate
	in $\Sigma_N$. Now $E\Sigma_N\to B\Sigma_N$ is the universal covering of $B\Sigma_N$ (a covering with simply connected total space), with $\Sigma_N=\pi_1(B\Sigma_N)$ acting
	by deck transformations (B1); so $U_N=E\Sigma_N\times_{\Sigma_N}[N]$ is the covering of $B\Sigma_N$ with monodromy the identity of $\Sigma_N$. For a based map
	$c\colon A\to B\Sigma_N$ with $\theta=\pi_1(c)\colon\pi_1(A,a_0)\to\Sigma_N$, the pullback $c^*E\Sigma_N$ is the principal $\Sigma_N$-bundle $\widetilde A\times_\theta\Sigma_N$
	(the map $\widetilde A\times\Sigma_N\to c^*E\Sigma_N$, $(\tilde a,\sigma)\mapsto(p(\tilde a),\sigma\cdot\tilde c(\tilde a))$, where $\tilde c\colon\widetilde A\to E\Sigma_N$ is a lift of
	$c\circ p$, induces an isomorphism), and therefore
	\[
	c^*U_N=c^*E\Sigma_N\times_{\Sigma_N}[N]\cong\widetilde A\times_\theta\Sigma_N\times_{\Sigma_N}[N]=\widetilde A\times_\theta[N]
	\]
	is the covering with monodromy $\theta=\pi_1(c)$. Since by (B2) free homotopy classes of maps $A\to B\Sigma_N$ correspond to conjugacy classes of homomorphisms
	$\pi_1(A,a_0)\to\Sigma_N$, the assignment $c\mapsto c^*U_N$ is a bijection from free homotopy classes of maps $A\to B\Sigma_N$ onto isomorphism classes of $N$-sheeted
	coverings of $A$: this is the classification statement. Finally, $\C^N_\perm$ is obtained from
	$U_N$ by applying the functor $S\mapsto\C[S]$ (the vector space with basis $S$) fibrewise, so $\C^N_\perm=V_{U_N}$ in the notation of Section~\ref{ss:Kcov}, and this
	construction commutes with pullback; thus $V_{c^*U_N}=c^*V_{U_N}=c^*\C^N_\perm$. The composite $B\Sigma_N\to BU(N)\to BU$ classifies $[\C^N_\perm]-N$ by (K5),
	because $\C^N_\perm=B(\Sigma_N\to U(N))^*\gamma_N$ (the associated bundle construction is natural in the group), and pulling back along $c$ gives the last claim.
	This is the form in which coverings enter the Barratt--Priddy--Quillen theorem (Section~\ref{ss:BPQ}).
	
	\begin{remark}\label{rem:simplyconnected}
		A finite covering of a simply connected, locally path-connected space is trivial: its monodromy is trivial, so each sheet
		through a point of the fibre is a section \cite[Proposition~1.33]{Hat}. Hence on a simply connected complex every reduced
		covering class is $0$. The Moore spaces $P^d(m)$ with $d\ge3$ are simply connected, so the coverings we need cannot live on
		them. They will live on complexes $T_n$ with large perfect fundamental groups.
	\end{remark}
	
	\subsection{Moore spaces}\label{ss:moore}
	
	\begin{definition}[Moore spaces]
		For $m\ge2$ let $P^2(m)=S^1\cup_{z^m}e^2$ be the space obtained from the circle by attaching a $2$-cell along the map
		$z\mapsto z^m$, and for $d\ge3$ let $P^d(m)=\Sigma^{d-2}P^2(m)$. Thus $P^d(m)=S^{d-1}\cup_{\rho_m}e^d$ has three cells
		$e^0,e^{d-1},e^d$, its basepoint is $e^0$, and the attaching map $\rho_m\colon\partial e^d=S^{d-1}\to S^{d-1}$ of the top cell
		has degree $m$. We call $S^{d-1}$ the \emph{bottom sphere} and $\partial e^d$ the \emph{attaching sphere}. Let $c\colon
		P^d(m)\to P^d(m)/S^{d-1}=S^d$ be the collapse of the bottom sphere, and let $b_d$ be a generator of $\Kt^0(S^d)$ ($d$ even)
		or of $K^1(S^d)$ ($d$ odd).
	\end{definition}
	
	The cellular chain complex of $P^d(m)$ is $\Z\xrightarrow{m}\Z$ in degrees $d$ and $d-1$, so $\widetilde H_{d-1}(P^d(m))=\Z/m$ and all other reduced homology
	groups vanish; this is the defining property of a Moore space of type $(\Z/m,d-1)$. Examples: $P^2(2)=\R P^2$, and $P^2(m)$ is the space of
	Example~\ref{ex:P2p}. The spaces $P^d(m)$ with $d\ge3$ are suspensions, hence simply connected (Theorem~\ref{thm:vankampen}(c)).
	
	\begin{figure}[ht]
		\centering
		\begin{tikzpicture}[scale=1]
			\draw[thick,fill=blue!6] (0,0) circle (1.4);
			\draw[thick,red!70!black] (0,0) circle (1.4);
			\foreach \a in {0,60,...,300}{\draw[->,red!70!black,thick] ({1.4*cos(\a)},{1.4*sin(\a)}) arc (\a:\a+25:1.4);}
			\node at (0,0) {$e^d$};
			\node[red!70!black] at (2.1,1.1) {$\partial e^d$};
			\draw[->,thick] (2.0,0) -- node[above]{$\rho_m$, degree $m$} (4.4,0);
			\draw[thick,blue!70!black] (5.5,0) circle (0.9);
			\node[blue!70!black] at (5.5,-1.3) {bottom sphere $S^{d-1}$};
			\fill (6.4,0) circle (1.5pt) node[right]{$e^0$};
			\node[align=left,font=\small] at (10.1,0) {$P^d(m)=S^{d-1}\cup_{\rho_m}e^d$\\[2pt]
				cells $e^0,e^{d-1},e^d$\\[2pt] $\widetilde H_{d-1}=\Z/m$, $\widetilde H_d=0$\\[2pt]
				$\Kt^0=\Z/m$ ($d$ even)};
		\end{tikzpicture}
		\caption{The Moore space $P^d(m)$: the boundary of the top cell wraps $m$ times around the bottom sphere.}\label{fig:moore}
	\end{figure}
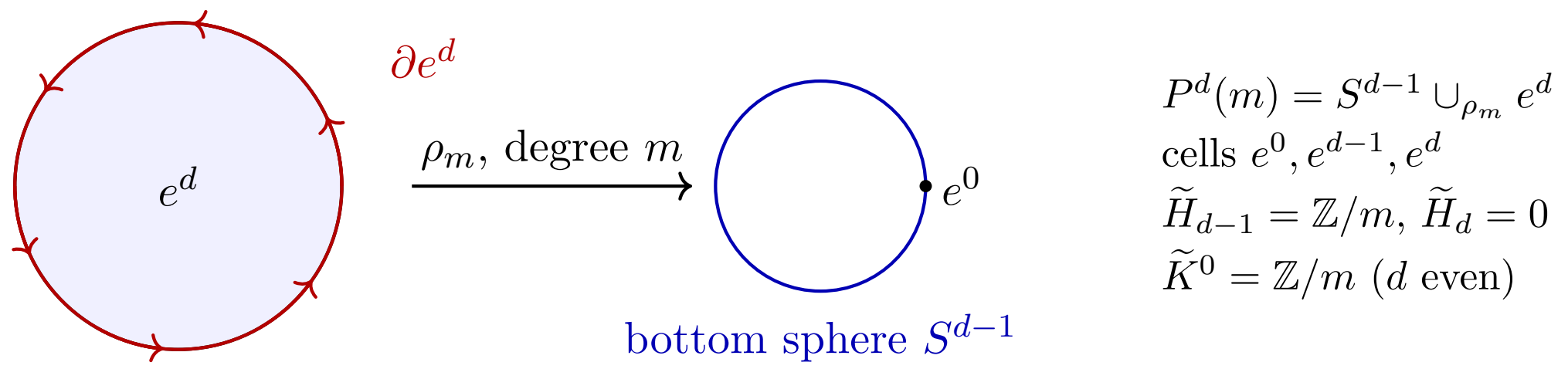
	
	\begin{lemma}\label{lem:mooreK}
		Let $m\ge2$ and $d\ge2$.
		\begin{enumerate}[label=\textup{(\alph*)}]
			\item $P^d(m)$ is $(d-2)$-connected; more generally, a finite wedge (at the basepoints) of Moore spaces $P^{d_j}(m_j)$ with $d_j\ge d$ for all $j$ is
			$(d-2)$-connected. In particular these spaces are $2$-connected if $d\ge4$. For $d\ge3$, $P^d(m)$ is a reduced suspension and simply connected.
			\item For even $d$, $\Kt^0(P^d(m))=\Z/m$, generated by $\gamma_d:=c^*b_d$, and $K^1(P^d(m))=0$. For odd $d$, $\Kt^0(P^d(m))=0$
			and $K^1(P^d(m))=\Z/m$, generated by $c^*b_d$.
		\end{enumerate}
	\end{lemma}
	
	\begin{proof}
		(a) The cells of $P^d(m)$ have dimensions $0$, $d-1$ and $d$. A wedge of CW complexes at $0$-cells is a CW complex whose cells are the common $0$-cell and the
		cells of positive dimension of the summands; so a wedge of Moore spaces $P^{d_j}(m_j)$ with $d_j\ge d$ has a single $0$-cell and no cells of dimension
		$1,\dots,d-2$. By Theorem~\ref{thm:standardhtpy}(a) it is $(d-2)$-connected. For $d\ge3$, $P^d(m)=\Sigma^{d-2}P^2(m)$ is a reduced suspension by definition,
		and it is simply connected by Theorem~\ref{thm:vankampen}(c).
		
		(b) By definition, $P^d(m)$ is the mapping cone $C_{\rho_m}$ of the based map $\rho_m\colon S^{d-1}\to S^{d-1}$ of degree $m$, the inclusion
		$i\colon S^{d-1}\to P^d(m)$ is the inclusion of the bottom sphere, and the collapse $C_{\rho_m}\to C_{\rho_m}/S^{d-1}=\Sigma S^{d-1}=S^d$ is the map $c$. The
		exact sequence (K2), for $\Kt^0$ and for $\Kt^1$, reads
		\[
		\Kt^j(S^d)\xrightarrow{(\Sigma\rho_m)^*}\Kt^j(S^d)\xrightarrow{\ c^*\ }\Kt^j(P^d(m))\xrightarrow{\ i^*\ }\Kt^j(S^{d-1})\xrightarrow{\ \rho_m^*\ }\Kt^j(S^{d-1})\qquad(j=0,1).
		\]
		The maps $\rho_m$ and $\Sigma\rho_m$ have degree $m$, and a based self-map of a sphere of degree $m$ is homotopic to the $m$-fold pinch sum of the identity,
		so it induces multiplication by $m$ on $\Kt^*$ (Lemma~\ref{lem:pinch}; the spheres $S^{d-1}$ and $S^d$ are suspensions since $d\ge2$). Let $d$ be even. By (K3),
		$\Kt^0(S^d)=\Z$, $\Kt^1(S^d)=0$, $\Kt^0(S^{d-1})=0$ and $\Kt^1(S^{d-1})=\Z$. For $j=0$ the sequence becomes $\Z\xrightarrow{m}\Z\xrightarrow{c^*}\Kt^0(P^d(m))\to0$, so
		$c^*$ induces an isomorphism $\Z/m\cong\Kt^0(P^d(m))$ and $\gamma_d=c^*b_d$ is a generator. For $j=1$ it becomes $0\to\Kt^1(P^d(m))\to\Z\xrightarrow{m}\Z$, and
		multiplication by $m$ is injective, so $K^1(P^d(m))=\Kt^1(P^d(m))=0$. For odd $d$ the roles of $j=0$ and $j=1$ are interchanged, and the same argument gives
		$\Kt^0(P^d(m))=0$ and $K^1(P^d(m))=\Z/m$ generated by $c^*b_d$.
	\end{proof}
	
	The important feature of (b) is that the generator is pulled back from the top cell: \emph{the torsion of $P^d(m)$ is carried by the top cell}. Consequently a map
	between Moore spaces acts on $K$-theory through its effect on the top cells, which is made explicit in the next lemma.
	
	\begin{lemma}[Degree squares]\label{lem:degreesquare}
		Let $d\ge2$, $m,m'\ge2$, and $s,t\in\Z$ with $tm'=ms$. There is a based cellular map $\varphi_{t,s}\colon P^d(m')\to P^d(m)$ of
		degree $t$ on the bottom spheres such that $c\circ\varphi_{t,s}=\sigma\circ c'$ for a map $\sigma\colon S^d\to S^d$ of degree $s$
		($c,c'$ the collapses). Hence $\varphi_{t,s}^*(c^*b_d)=s\,c'^*b_d$: $\varphi_{t,s}$ induces $\Z/m\to\Z/m'$, $1\mapsto s$, on
		$\Kt^0$ ($d$ even) or on $K^1$ ($d$ odd). Every homomorphism $\Z/m\to\Z/m'$ arises in this way: $1\mapsto s$ is well defined iff
		$m'\mid ms$, and then one takes $t=ms/m'$.
	\end{lemma}
	
	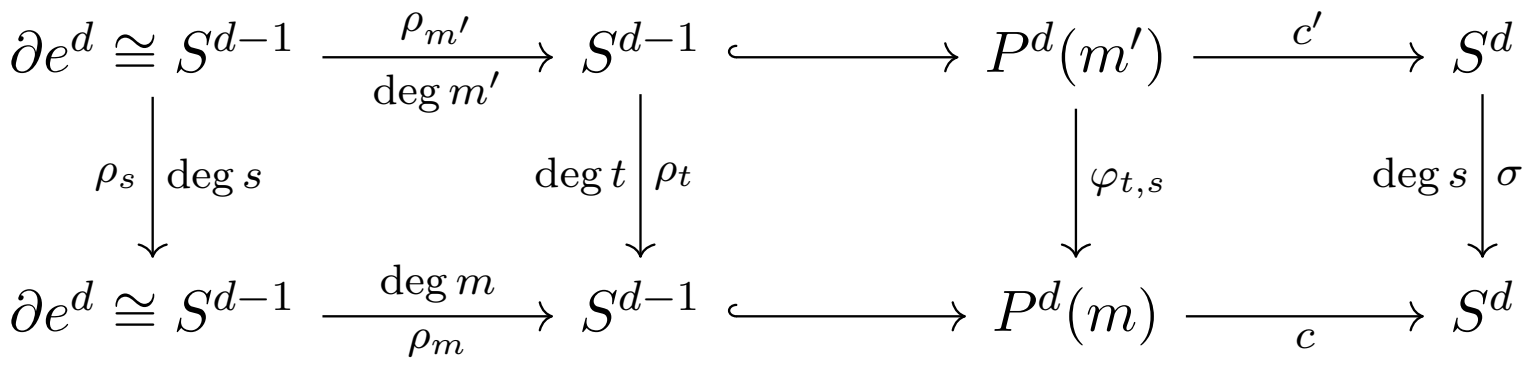
\begin{figure}[ht]
		\centering
		\begin{tikzcd}[column sep=4em,row sep=2.6em]
			\partial e^d\cong S^{d-1}\arrow[r,"\rho_{m'}","\deg m'"']\arrow[d,"\rho_s"',"\deg s"] & S^{d-1}\arrow[d,"\rho_t","\deg t"']\arrow[r,hook] & P^d(m')\arrow[d,"\varphi_{t,s}"]\arrow[r,"c'"] & S^d\arrow[d,"\sigma","\deg s"']\\
			\partial e^d\cong S^{d-1}\arrow[r,"\rho_m"',"\deg m"] & S^{d-1}\arrow[r,hook] & P^d(m)\arrow[r,"c"'] & S^d
		\end{tikzcd}
		\caption{A degree square. The left square commutes up to homotopy exactly when $tm'=ms$, which is the condition for the
			cellular chain map $(t,s)$ between the chain complexes $\Z\xrightarrow{m'}\Z$ and $\Z\xrightarrow{m}\Z$.}\label{fig:degsq}
	\end{figure}
	
	\begin{proof}
		Write $P^d(m)=S^{d-1}\cup_{\rho_m}e^d$ with characteristic map $\chi\colon D^d\to P^d(m)$ of the top cell, so that $\chi|_{\partial D^d}=\rho_m$ followed by the
		inclusion of the bottom sphere, and similarly $P^d(m')=S^{d-1}\cup_{\rho_{m'}}e^d$ with characteristic map $\chi'$. Choose based maps
		$\rho_s,\rho_t\colon S^{d-1}\to S^{d-1}$ of degrees $s$ and $t$ (for instance the $(d-2)$-fold suspensions of $z\mapsto z^s$ and $z\mapsto z^t$ on $S^1$). The
		composites $\rho_t\circ\rho_{m'}$ and $\rho_m\circ\rho_s$ both have degree $tm'=ms$, and based self-maps of $S^{d-1}$ ($d-1\ge1$) of the same degree are based
		homotopic \cite[Corollary~4.25]{Hat}; let $H\colon S^{d-1}\times[0,1]\to S^{d-1}$ be a based homotopy with $H_0=\rho_t\circ\rho_{m'}$ and $H_1=\rho_m\circ\rho_s$.
		
		Define $\varphi=\varphi_{t,s}$ on the bottom sphere as $\rho_t$, and on the top cell by
		\[
		\varphi\big(\chi'(rx)\big)=\begin{cases}H(x,2-2r)&\text{for }x\in S^{d-1},\ \tfrac12\le r\le1,\\[2pt]
			\chi\big(2r\,\rho_s(x)\big)&\text{for }x\in S^{d-1},\ 0\le r\le\tfrac12,\end{cases}
		\]
		where $rx\in D^d$ denotes the point at distance $r$ from the centre in the direction $x$. The two formulas agree at $r=\tfrac12$, since
		$H(x,1)=\rho_m(\rho_s(x))=\chi(\rho_s(x))$, and at $r=1$ the first formula gives $\rho_t(\rho_{m'}(x))=\rho_t(\chi'(x))$, which is the value prescribed on the bottom
		sphere. The second formula is well defined at $r=0$ (constant value $\chi(0)$). Hence $\varphi$ is a well defined continuous map $P^d(m')\to P^d(m)$, by the
		gluing lemma applied to the closed cover of $P^d(m')$ by the bottom sphere and the closed top cell. It is based, because $\rho_t$ is based, and cellular,
		because it maps the bottom sphere to the bottom sphere.
		
		Now apply the collapses. The composite $c\circ\varphi$ sends the bottom sphere and the collar $\{\chi'(rx):\tfrac12\le r\le1\}$ to the basepoint of
		$S^d=P^d(m)/S^{d-1}$, because $H$ takes values in the bottom sphere of $P^d(m)$; on the inner disc it is $\chi'(rx)\mapsto\bar\chi(2r\rho_s(x))$, where
		$\bar\chi\colon D^d\to D^d/\partial D^d=S^d$ is the quotient map. Hence $c\circ\varphi=\sigma\circ c'$, where $\sigma\colon S^d=D^d/\partial D^d\to S^d$ is the map
		which collapses the collar $\{\tfrac12\le r\le1\}$ to the basepoint and is the cone $rx\mapsto 2r\rho_s(x)$ on the inner disc. The map $\sigma$ is homotopic to the
		suspension of $\rho_s$ (compose the cone on $\rho_s$ with the quotient map $D^d\to D^d/\{\frac12\le r\le1\}\cong D^d/\partial D^d$, which is homotopic to the
		identity of $S^d$ through maps of degree one), so $\deg\sigma=\deg\rho_s=s$.
		
		By Lemma~\ref{lem:pinch}, $\sigma^*$ is multiplication by $s$ on $\Kt^*(S^d)$, so $\varphi^*(c^*b_d)=(c\varphi)^*b_d=(\sigma c')^*b_d=c'^*(\sigma^*b_d)=s\,c'^*b_d$. Since
		$c^*b_d$ and $c'^*b_d$ generate $\Kt^0$ (for even $d$) or $K^1$ (for odd $d$) of $P^d(m)$ and $P^d(m')$ (Lemma~\ref{lem:mooreK}(b)), this says that $\varphi^*$ is the
		homomorphism $\Z/m\to\Z/m'$, $1\mapsto s$. Finally, $1\mapsto s$ defines a homomorphism $\Z/m\to\Z/m'$ if and only if $ms\equiv0\bmod m'$, i.e.\ $m'\mid ms$, and then
		$t=ms/m'$ is an integer with $tm'=ms$.
	\end{proof}
	
	\begin{remark}\label{rem:bottomdegree}
		The two spheres play different roles. On $K$-theory only $s$, the degree on the attaching sphere, matters, because the torsion is carried by the top cell;
		on the homology group $\widetilde H_{d-1}$, the map $\varphi_{t,s}$ is multiplication by $t$ (the degree on the bottom sphere). For example, with $m=3$, $m'=9$ and
		$s=3$ one gets $t=1$: the map $\varphi_{1,3}\colon P^d(9)\to P^d(3)$ induces $\Z/9\to\Z/3$, $1\mapsto1$, on $\widetilde H_{d-1}$ and $\Z/3\to\Z/9$, $1\mapsto3$, on $K$-theory.
	\end{remark}
	
	\subsection{Spectra, stable maps and Spanier--Whitehead duality}\label{ss:spectra}
	
	Adams' periodic maps are maps of spectra, and the passage from $K$-theory to stable cohomotopy uses duality. We recall what is
	needed; references are \cite[Part~III]{Ad74}, \cite{Swi}, \cite{SW}.
	
	\begin{definition}[Spectra]\label{def:spectra}
		A \emph{(sequential CW) spectrum} $E$ is a sequence of based CW complexes $E_k$ ($k\ge0$) with cellular inclusions $\Sigma
		E_k\to E_{k+1}$ onto subcomplexes. Its \emph{homotopy groups} are $\pi_j(E)=\operatorname{colim}_k\pi_{j+k}(E_k)$ for $j\in\Z$, the colimit being taken along
		the suspension maps $\pi_{j+k}(E_k)\to\pi_{j+k+1}(\Sigma E_k)\to\pi_{j+k+1}(E_{k+1})$.
		Maps of spectra are defined so that the homotopy category of spectra is a triangulated category with suspension $\Sigma$
		invertible \cite[Part~III, \S\S2--3]{Ad74}; we write $[E,F]$ for the group of homotopy classes of maps $E\to F$ and
		$[E,F]_j=[\Sigma^jE,F]$. Every spectrum is equivalent to an \emph{$\Omega$-spectrum}, one in which the adjoints $E_k\to\Omega E_{k+1}$ of the structure maps are
		weak homotopy equivalences; for an $\Omega$-spectrum $E$ and a based CW complex $X$, $[\Sigma^\infty X,\Sigma^jE]=[X,E_j]_*$. Examples:
		\begin{itemize}
			\item the \emph{suspension spectrum} $\Sigma^\infty X$ of a based CW complex, with $(\Sigma^\infty X)_k=\Sigma^kX$. For finite based CW complexes $A,B$,
			\begin{equation}\label{eq:stablemaps}
				[\Sigma^\infty A,\Sigma^\infty B]=\operatorname{colim}_k[\Sigma^kA,\Sigma^kB]_*,
			\end{equation}
			the colimit over suspension \cite[Part~III, Proposition~2.8]{Ad74}: a map of spectra $\Sigma^\infty A\to\Sigma^\infty B$ is represented by a based map $\Sigma^kA\to\Sigma^kB$ for
			some $k$, and two such maps represent the same map of spectra if and only if they become homotopic after further suspension. Thus maps between suspension
			spectra of finite complexes are \emph{stable} homotopy classes of maps of spaces, and $\Sigma^\infty$ sends mapping cones to cofibres;
			\item the \emph{sphere spectrum} $\Sph=\Sigma^\infty S^0$, whose homotopy groups $\pi_j(\Sph)=\pi_j^s$ are the stable homotopy
			groups of spheres; $\pi_0^s=\Z$ (a map $S^k\to S^k$ corresponds to its degree), $\pi_1^s=\Z/2$ (generated by the class of the Hopf map $S^3\to S^2$), $\pi_2^s=\Z/2$ (generated by the square of that class);
			\item the \emph{Moore spectrum} $\Sph/m$, the cofibre of $m\colon\Sph\to\Sph$ (the $m$-fold sum of the identity), with cells in dimensions $0$ and $1$; there is a
			cofibre sequence $\Sph\xrightarrow{m}\Sph\xrightarrow{i}\Sph/m\xrightarrow{j}\Sigma\Sph$;
			\item the \emph{complex $K$-theory spectrum} $\KU$. Its spaces are $\KU_{2k}=\Z\times BU$ and $\KU_{2k+1}=U$ (replaced by CW approximations), and the
			adjoint structure maps $\Z\times BU\to\Omega U$ and $U\to\Omega(\Z\times BU)$ are the homotopy equivalences given by Bott periodicity and by the fibration
			$U\to EU\to BU$ (the structure maps are converted into cellular inclusions by mapping cylinders, as in \cite[Part~III, \S2]{Ad74}); so $\KU$ is an
			$\Omega$-spectrum whose spaces repeat with period $2$. For a finite CW complex $X$ it represents $K$-theory,
			\begin{gather*}
				\Kt^0(X)=[X,\Z\times BU]_*=[\Sigma^\infty X,\KU],\qquad\Kt^1(X)=[X,U]_*=[\Sigma^\infty X,\Sigma\KU],\\
				\Kt^j(X)=[\Sigma^\infty X,\Sigma^j\KU]\quad(j\in\Z),
			\end{gather*}
			and (K3) is built into it: $[\Sigma^\infty\Sigma X,\Sigma^{j+1}\KU]=[\Sigma^\infty X,\Sigma^j\KU]$. Its homotopy groups are $\pi_{2k}(\KU)=\pi_0(\Z\times BU)=\Z$ and
			$\pi_{2k+1}(\KU)=\pi_0(U)=0$; as a graded ring (with the multiplication induced by the tensor product of vector bundles) $\pi_*(\KU)=\Z[u^{\pm1}]$, where
			$u\in\pi_2(\KU)$ is the Bott element. The spectrum $\KU$ collects the $K$-groups of all degrees into a single object, on which the constructions of
			stable homotopy theory (smash products, duality) can be performed.
		\end{itemize}
	\end{definition}
	
	\begin{lemma}[Stable maps between finite complexes]\label{lem:stablemaps}
		Let $A,B$ be finite based CW complexes, and suppose that $B$ is $(n-1)$-connected with $n\ge2$. If $\dim A\le2n-1$, then the map $[A,B]_*\to[\Sigma^\infty A,\Sigma^\infty B]$,
		sending a map to the stable class it represents, is onto; if $\dim A\le2n-2$, it is bijective.
	\end{lemma}
	
	\begin{proof}
		For $k\ge0$ the complex $\Sigma^kB$ is $(n+k-1)$-connected (apply the last sentence after Theorem~\ref{thm:freudenthal} $k$ times) and
		$\dim\Sigma^kA=\dim A+k$. If $\dim A\le2n-1$, then $\dim\Sigma^kA\le2(n+k)-1$ for every $k$, and Theorem~\ref{thm:freudenthal}, applied to $\Sigma^kA$ and $\Sigma^kB$,
		shows that every map $\Sigma\colon[\Sigma^kA,\Sigma^kB]_*\to[\Sigma^{k+1}A,\Sigma^{k+1}B]_*$ in the colimit \eqref{eq:stablemaps} is onto. Every element of the colimit
		is represented at some stage $k$, and the composite $[A,B]_*\to[\Sigma^kA,\Sigma^kB]_*$ of $k$ surjections is onto, so the map from $[A,B]_*$ to the colimit is onto.
		If $\dim A\le2n-2$, then $\dim\Sigma^kA\le2(n+k)-2$ for every $k$, all maps in the colimit are bijective, and so is the map from $[A,B]_*$ to the colimit.
	\end{proof}
	
	In words: a stable map from a complex of dimension at most $2n-1$ into an $(n-1)$-connected complex is represented by an actual map of spaces. This is how
	Adams' maps are desuspended in Lemma~\ref{lem:unstable}.
	
	\begin{definition}[Smash products, homology and cohomology theories, the Hurewicz map]\label{def:smash}
		The homotopy category of spectra has a \emph{smash product} $\wedge$, which is associative and commutative up to coherent natural equivalence, has the
		sphere spectrum as unit, $\Sph\wedge E\simeq E$, and satisfies $\Sigma^\infty(A\wedge B)\simeq\Sigma^\infty A\wedge\Sigma^\infty B$ for based CW complexes; it also has
		\emph{function spectra} $F(E,F)$, characterized by the natural isomorphism $[W\wedge E,F]\cong[W,F(E,F)]$ \cite[Part~III, \S\S3--5]{Ad74}. For a spectrum $E$ and a
		based CW complex $X$ one writes $E\wedge X=E\wedge\Sigma^\infty X$.
		
		A spectrum $E$ defines a \emph{homology theory} and a \emph{cohomology theory} on spectra, and on based CW complexes via $\Sigma^\infty$:
		\[
		E_j(M)=\pi_j(E\wedge M)=[\Sigma^j\Sph,E\wedge M],\qquad E^j(M)=[M,\Sigma^jE]\qquad(j\in\Z).
		\]
		Both are functorial in $M$ (covariantly, resp.\ contravariantly), they turn cofibre sequences $M\to N\to C\to\Sigma M$ into long exact sequences, and they
		have suspension isomorphisms $E_{j+1}(\Sigma M)\cong E_j(M)$, $E^{j+1}(\Sigma M)\cong E^j(M)$ \cite[Part~III, \S6]{Ad74}. For $E=\Sph$ these are stable homotopy and stable
		cohomotopy; for $E=\KU$ they are \emph{$K$-homology} $\KU_*$ and \emph{$K$-cohomology} $\KU^*$, and for a finite CW complex $X$,
		$\KU^j(\Sigma^\infty X)=[\Sigma^\infty X,\Sigma^j\KU]=\Kt^j(X)$ by Definition~\ref{def:spectra}. The coefficient groups are $E_j(\Sph)=E^{-j}(\Sph)=\pi_j(E)$; for example
		$\KU_j(\Sph/m)$ is $\Z/m$ for even $j$ and $0$ for odd $j$, by the long exact sequence of the cofibre sequence $\Sph\xrightarrow{m}\Sph\to\Sph/m$.
		
		The \emph{unit} of $\KU$ is the map of spectra $\eta\colon\Sph\to\KU$ corresponding to $1\in\pi_0(\KU)=[\Sph,\KU]$; it is the unit of the ring spectrum structure
		of $\KU$ given by the tensor product of bundles. For a spectrum $M$ it induces the \emph{$\KU$-Hurewicz map}
		\[
		h\colon\pi_j(M)=\pi_j(\Sph\wedge M)\longrightarrow\pi_j(\KU\wedge M)=\KU_j(M),\qquad\beta\mapsto(\eta\wedge1_M)\circ\beta ,
		\]
		a homomorphism natural in $M$. It is the analogue of the classical Hurewicz homomorphism $\pi_j(M)\to H_j(M;\Z)$, which is obtained in the same way from the
		unit $\Sph\to H\Z$ of the Eilenberg--MacLane spectrum; it records what $K$-homology sees of a homotopy class. We shall use it to transfer information about
		Adams' maps, which live in the homotopy of Moore spectra, into $K$-theory (Lemma~\ref{lem:lift}).
	\end{definition}
	
	\begin{definition}[Spanier--Whitehead duality {\cite{SW}, \cite[Part~III, \S5]{Ad74}}]\label{def:SWduality}
		A spectrum is \emph{finite} if it is equivalent to $\Sigma^{-k}\Sigma^\infty A$ for some finite CW complex $A$ and $k\ge0$. The \emph{(Spanier--Whitehead) dual}
		of a finite spectrum $M$ is the function spectrum $DM=F(M,\Sph)$. It is again finite, and the \emph{evaluation} $DM\wedge M\to\Sph$ (adjoint to the identity of $DM$)
		has a \emph{coevaluation} $\Sph\to M\wedge DM$ such that the two composites
		\begin{gather*}
			M\simeq\Sph\wedge M\xrightarrow{\ \text{coev}\wedge1\ }M\wedge DM\wedge M\xrightarrow{\ 1\wedge\text{ev}\ }M\wedge\Sph\simeq M,\\
			DM\simeq DM\wedge\Sph\xrightarrow{\ 1\wedge\text{coev}\ }DM\wedge M\wedge DM\xrightarrow{\ \text{ev}\wedge1\ }\Sph\wedge DM\simeq DM
		\end{gather*}
		are the identities. Consequently $D$ is an exact contravariant equivalence of the homotopy category of finite spectra with itself (it takes cofibre sequences to
		cofibre sequences and commutes with $\Sigma$ up to sign in the degree: $D(\Sigma M)\simeq\Sigma^{-1}DM$), with $DDM\simeq M$ and $D\Sph=\Sph$, and for all spectra $W$
		and $E$ there are isomorphisms
		\begin{equation}\label{eq:duality}
			[W\wedge DM,E]\cong[W,M\wedge E],\qquad\text{in particular}\qquad[\Sigma^dDM,E]\cong\pi_d(M\wedge E),
		\end{equation}
		natural in $W$, $E$ and $M$; the first is $[W\wedge DM,E]\cong[W,F(DM,E)]$ together with the equivalence $M\wedge E\simeq F(DM,E)$ which expresses finiteness,
		and the second is the case $W=\Sigma^d\Sph$.
	\end{definition}
	
	Duality exchanges homology and cohomology: by \eqref{eq:duality}, $E^j(DM)=[DM,\Sigma^jE]\cong[\Sph,M\wedge\Sigma^jE]=E_{-j}(M)$ for every spectrum $E$. We need two consequences of duality.
	
	\begin{lemma}[Duality for Moore spectra]\label{lem:mooredual}
		$D(\Sph/m)\simeq\Sigma^{-1}\Sph/m$, and $\Sigma^\infty P^d(m)\simeq\Sigma^{d-1}\Sph/m$ for $d\ge2$. Hence
		$\Sigma^\infty P^d(m)\simeq\Sigma^dD(\Sph/m)$.
	\end{lemma}
	
	\begin{proof}
		Apply $D$ to the cofibre sequence $\Sph\xrightarrow{m}\Sph\xrightarrow{i}\Sph/m\xrightarrow{j}\Sigma\Sph$. Since $D$ is exact and contravariant, and $D\Sph=\Sph$,
		$D(\Sigma\Sph)=\Sigma^{-1}\Sph$, we obtain a cofibre sequence
		\[
		\Sigma^{-1}\Sph\xrightarrow{Dj}D(\Sph/m)\xrightarrow{Di}\Sph\xrightarrow{Dm}\Sph .
		\]
		Here $Dm=m$, because $D$ is additive on morphisms and $D(\id_\Sph)=\id_\Sph$. Rotating the triangle, $\Sph\xrightarrow{m}\Sph\to\Sigma D(\Sph/m)\to\Sigma\Sph$ is a cofibre
		sequence beginning with $m\colon\Sph\to\Sph$, so $\Sigma D(\Sph/m)$ is equivalent to the cofibre $\Sph/m$ of $m$ (the third term of a cofibre sequence is determined up
		to equivalence by the first map). Hence $D(\Sph/m)\simeq\Sigma^{-1}\Sph/m$.
		
		For the second statement, $P^2(m)$ is the mapping cone of the degree $m$ map $\rho_m\colon S^1\to S^1$, and $\Sigma^\infty$ takes mapping cones to cofibres
		(Definition~\ref{def:spectra}). The stable class of $\rho_m$ in $[\Sigma^\infty S^1,\Sigma^\infty S^1]=[\Sigma\Sph,\Sigma\Sph]=\pi_0^s=\Z$ is its degree $m$, i.e.\
		$\Sigma^\infty\rho_m=\Sigma(m\cdot\id_\Sph)$. Hence $\Sigma^\infty P^2(m)\simeq\operatorname{cofibre}(\Sigma m)=\Sigma(\Sph/m)$, since $\Sigma$ is exact. Then
		$\Sigma^\infty P^d(m)=\Sigma^\infty\Sigma^{d-2}P^2(m)=\Sigma^{d-2}\Sigma^\infty P^2(m)\simeq\Sigma^{d-1}\Sph/m$. Combining the two statements,
		$\Sigma^dD(\Sph/m)\simeq\Sigma^{d-1}\Sph/m\simeq\Sigma^\infty P^d(m)$.
	\end{proof}
	
	\begin{lemma}\label{lem:KUhomtocohom}
		A map of finite spectra which induces an isomorphism on $\KU_*$ induces an isomorphism on $\KU^*$.
	\end{lemma}
	
	\begin{proof}
		Let $f\colon M\to N$ be a map of finite spectra inducing an isomorphism on $\KU_*$, and let $C$ be its cofibre, a finite spectrum. The long exact sequence of
		$\KU_*$ for the cofibre sequence $M\xrightarrow{f}N\to C\to\Sigma M$ shows that $\KU_*(C)=\pi_*(\KU\wedge C)=0$; a spectrum with vanishing homotopy groups is
		contractible, so $\KU\wedge C\simeq*$. We claim that also $\KU\wedge DC\simeq*$. Indeed, by the second triangle identity of Definition~\ref{def:SWduality},
		$DC$ is a retract of $DC\wedge C\wedge DC$; smashing with $\KU$, the spectrum $\KU\wedge DC$ is a retract of
		$\KU\wedge DC\wedge C\wedge DC\simeq(\KU\wedge C)\wedge DC\wedge DC\simeq*$, hence contractible. Now by \eqref{eq:duality} with $W=\Sph$, $M$ replaced by $DC$ (so that $DM=DDC\simeq C$) and
		$E=\Sigma^j\KU$,
		\[
		\KU^j(C)=[C,\Sigma^j\KU]=[\Sph\wedge DDC,\Sigma^j\KU]\cong[\Sph,DC\wedge\Sigma^j\KU]=\pi_{-j}(\KU\wedge DC)=0
		\]
		for every $j$. The long exact sequence of $\KU^*$ for the cofibre sequence $M\to N\to C\to\Sigma M$ then shows that
		$f^*\colon\KU^j(N)\to\KU^j(M)$ is an isomorphism for every $j$.
	\end{proof}
	
	In our application, Adams' maps are constructed as $\KU_*$-isomorphisms (Theorem~\ref{thm:adams}), while we need isomorphisms on the $K$-theory of spaces,
	i.e.\ on $\KU^*$; the lemma bridges the two.
	
	\begin{definition}[Stable cohomotopy and stable cohomotopy lifts]\label{def:stablelift}
		Let $QS^0=\operatorname{colim}_k\Omega^kS^k$, the colimit along the maps $\Omega^kS^k\to\Omega^{k+1}S^{k+1}$ which send a map $S^k\to S^k$ to its suspension. This is
		the zeroth space of an $\Omega$-spectrum equivalent to $\Sph$, the \emph{infinite loop space} $\Omega^\infty\Sph$ of the sphere spectrum, and for a finite based CW
		complex $X$ (which is compact, so that maps from $X$ into the colimit factor through a finite stage)
		\[
		[X,QS^0]_*=\operatorname{colim}_k[X,\Omega^kS^k]_*=\operatorname{colim}_k[\Sigma^kX,S^k]_*=[\Sigma^\infty X,\Sph]=:\pi^0_s(X),
		\]
		the \emph{stable cohomotopy} group of $X$. The components of $QS^0$ are indexed by the degree, $\pi_0(QS^0)=\operatorname{colim}_k\pi_k(S^k)=\Z$, and $Q_0S^0$
		denotes the component of degree $0$; a based map from a connected space into $QS^0$ takes values in a single component. Applying the functor $\Omega^\infty$
		(zeroth space of the associated $\Omega$-spectrum) to the unit $\eta\colon\Sph\to\KU$ gives an infinite loop map
		\[
		\Omega^\infty\eta\colon QS^0\to\Omega^\infty\KU=\Z\times BU ,
		\]
		which on $\pi_0$ is $\eta_*\colon\pi_0(\Sph)=\Z\to\pi_0(\KU)=\Z$, the identity, because $\eta$ represents $1$; so it maps $Q_0S^0$ to $\{0\}\times BU=BU$. Under the
		natural identifications $[X,QS^0]_*=[\Sigma^\infty X,\Sph]$ and $[X,\Z\times BU]_*=[\Sigma^\infty X,\KU]$ (Definition~\ref{def:spectra}), composition with
		$\Omega^\infty\eta$ corresponds to composition with $\eta$, that is, to the homomorphism $\eta_*\colon\pi^0_s(X)\to K^0(X)$.
		
		Let $X$ be a connected finite based CW complex and $x\in\Kt^0(X)=[X,BU]_*$. A \emph{stable cohomotopy lift} of $x$ is a based map $\alpha\colon X\to Q_0S^0$ such that
		$(\Omega^\infty\eta)\circ\alpha$ represents $x$; equivalently, $x=\eta_*(\alpha^\#)$ for the adjoint stable map $\alpha^\#\colon\Sigma^\infty X\to\Sph$. Thus $x$ has a
		stable cohomotopy lift if and only if it lies in the image of $\eta_*\colon\pi^0_s(X)\to\Kt^0(X)$.
	\end{definition}
	
	Since $\eta_*$ is a homomorphism, the classes with stable cohomotopy lifts form a subgroup of $\Kt^0(X)$: if $\alpha$ and $\beta$ lift $x$ and $y$, then $\alpha+\beta$
	and $-\alpha$, formed in the group $[X,Q_0S^0]_*=\pi^0_s(X)$, lift $x+y$ and $-x$. For a based map $u\colon X'\to X$, if $\alpha$ lifts $x$, then $\alpha\circ u$ lifts $u^*x$, by
	naturality of $\eta_*$. The reason for introducing stable cohomotopy lifts is the Barratt--Priddy--Quillen theorem (Section~\ref{ss:BPQ}): a class is a covering class
	if and only if it is represented by a map $X\to Q_0S^0$ which can be lifted, up to homotopy, along $b\colon B\Sigma_\infty\to Q_0S^0$ and then compressed into a finite
	$B\Sigma_N$. The existence of a stable cohomotopy lift is the homotopy-theoretic half of this condition; the other half is achieved by changing the space
	(Proposition~\ref{prop:lift}). Figure~\ref{fig:strategy} shows the two halves.
	
	\subsection{Adams operations and Adams' periodic maps}\label{ss:adams}
	
	\begin{definition}[Adams operations {\cite{Ad62}, \cite[Chapter~III]{At67}}]\label{def:adamsops}
		For an integer $k\ge1$ the \emph{Adams operation} $\psi^k$ is the natural ring endomorphism of $K^0(X)$, $X$ compact, characterized by
		$\psi^k[L]=[L]^k=[L^{\otimes k}]$ for line bundles $L$; it is defined by $\psi^k(x)=Q_k(\lambda^1x,\dots,\lambda^kx)$, where $\lambda^i$ are the exterior power
		operations and $Q_k$ is the Newton polynomial expressing the power sum $\sum t_i^k$ in terms of the elementary symmetric functions. (One also puts
		$\psi^0(x)=\rank(x)$, $\psi^{-1}(x)=\bar x$, the class of the complex conjugate bundle, and $\psi^{-k}=\psi^{-1}\psi^k$; only $k\ge1$ is used below.) The
		operations satisfy $\psi^k\psi^l=\psi^{kl}$, they commute with pullbacks and preserve $\Kt^0$, and on $\Kt^0(S^{2t})=\Z$ the operation $\psi^k$ is
		multiplication by $k^t$ \cite[Chapter~III]{At67}; $\psi^1$ is the identity. For $k\ge2$ the operation $\psi^k$ is not stable: on $\Kt^0(S^2)$ it is
		multiplication by $k$, while on $\Kt^0(S^0)=\Z$ it is the identity, so it does not commute with the Bott periodicity isomorphism
		$\Kt^0(S^0)\cong\Kt^0(S^2)$, and it does not define a self-map of the spectrum $\KU$, whose structure maps are given by Bott periodicity. The obstruction is
		that the image $k\beta$ of the Bott element $\beta\in\pi_2(\KU)$ is not a unit. After inverting $k$, the operation extends to a map of ring spectra
		$\psi^k\colon\KU[1/k]\to\KU[1/k]$, where $\KU[1/k]$ is the localization of $\KU$ with $\pi_*(\KU[1/k])=\pi_*(\KU)\otimes\Z[1/k]$ \cite{Ad74}. It sends
		$\beta$ to $k\beta$, so it acts on $\pi_{2t}(\KU[1/k])=\Z[1/k]$ by multiplication by $k^t$; in particular it is still not linear over multiplication by
		$\beta$ (it commutes with the suspension isomorphisms, as every map of spectra does, but not with the Bott periodicity isomorphism). It induces the
		operations $\psi^k$ on $K^0(X)\otimes\Z[1/k]$. In particular, if $M$ is a spectrum
		whose $\KU$-homology is $p$-torsion for a prime $p\nmid k$, then multiplication by $k$ is invertible on $\KU_*(M)$, so $\KU[1/k]_*(M)=\KU_*(M)$ and $\psi^k$ acts on
		$\KU_{2t}(M)$, naturally in $M$.
	\end{definition}
	
	The \emph{period} that will be used throughout is
	\begin{equation}\label{eq:period}
		q(p,r)=\begin{cases}2(p-1)p^{r-1},&p\text{ odd},\\ \max\{8,2^{r-1}\},&p=2,\end{cases}\qquad(p\text{ prime},\ r\ge1).
	\end{equation}
	For $p=2$ the values are $8,8,8,8,16,32,\dots$; for odd $p$, $q(p,r+1)=p\,q(p,r)$. In particular $q(p,r)$ is even, and the periods
	of a fixed prime form a divisibility chain: $q(p,c)\mid q(p,r)$ for $1\le c\le r$.
	
	\begin{theorem}[{\cite[\S12]{Ad66}}; see also {\cite[\S6]{ACL}}, {\cite[Chapter~1]{Rav}}]\label{thm:adams}
		For every prime $p$ and $r\ge1$ there is a map of spectra
		\[
		v=v_{p,r}\colon\Sigma^{q(p,r)}\Sph/p^r\longrightarrow\Sph/p^r
		\]
		inducing an isomorphism on $\KU_*$.
	\end{theorem}
	
	\begin{lemma}[Desuspension in the Freudenthal range]\label{lem:unstable}
		Let $p$ be a prime, $r\ge1$ and $q=q(p,r)$.
		\begin{enumerate}[label=\textup{(\roman*)}]
			\item If $d\ge q+3$, there is a based map $w\colon P^{d+q}(p^r)\to P^d(p^r)$ with $\Sigma^\infty w=\Sigma^{d-1}v_{p,r}$ under
			Lemma~\ref{lem:mooredual}; $w^*$ is an isomorphism on $\Kt^*$.
			\item If moreover $d$ is even and $k\ge0$, there is a based map $w_{d,k}\colon P^{d+kq}(p^r)\to P^d(p^r)$ with
			$w_{d,k}^*\gamma_d=\gamma_{d+kq}$.
		\end{enumerate}
	\end{lemma}
	
	\begin{proof}
		(i) Put $A=P^{d+q}(p^r)$ and $B=P^d(p^r)$. $B$ is $(d-2)$-connected (Lemma~\ref{lem:mooreK}(a)), i.e.\ $(n-1)$-connected with
		$n=d-1\ge2$, and $\dim A=d+q\le2n-1=2d-3$ iff $d\ge q+3$. By Lemma~\ref{lem:stablemaps}, $[A,B]_*\to[\Sigma^\infty A,\Sigma^\infty
		B]$ is onto, and we choose $w$ mapping to $\Sigma^{d-1}v$. Since $\Kt^*(Y)=\KU^*(\Sigma^\infty Y)$ naturally for finite $Y$,
		$w^*=(\Sigma^{d-1}v)^*$ on $\Kt^*$, which is an isomorphism by Theorem~\ref{thm:adams} and Lemma~\ref{lem:KUhomtocohom}.
		
		(ii) For $k=0$ take the identity. For $k\ge1$ compose the maps of (i) at the dimensions $d+jq$, $j=k-1,\dots,0$,
		\[
		P^{d+kq}(p^r)\xrightarrow{w}P^{d+(k-1)q}(p^r)\xrightarrow{w}\cdots\xrightarrow{w}P^{d+q}(p^r)\xrightarrow{w}P^d(p^r);
		\]
		each of these dimensions is at least $d\ge q+3$, so each step exists. The composite $w'$ induces an isomorphism
		$\Z/p^r=\Kt^0(P^d)\to\Kt^0(P^{d+kq})=\Z/p^r$, so $w'^*\gamma_d=u\gamma_{d+kq}$ for a unit $u$ modulo $p^r$. With $s\equiv
		u^{-1}\bmod p^r$ put $w_{d,k}=\varphi_{s,s}\circ w'$, where $\varphi_{s,s}\colon P^d(p^r)\to P^d(p^r)$ is the degree square of
		Lemma~\ref{lem:degreesquare} with $m=m'=p^r$. Then $w_{d,k}^*\gamma_d=w'^*(s\gamma_d)=su\,\gamma_{d+kq}=\gamma_{d+kq}$.
	\end{proof}
	
	\begin{remark}
		The inequality $d\ge q+3$ is exactly what the Freudenthal theorem needs: $d+q\le2(d-2)+1$ iff $d\ge q+3$. The iterate in (ii) must
		be built one step at a time: desuspending the $k$-fold stable iterate $\Sigma^{kq}\Sph/p^r\to\Sph/p^r$ in one step would require
		$d\ge kq+3$.
	\end{remark}
	
	\begin{lemma}[Stable cohomotopy lifts on Moore spaces]\label{lem:lift}
		If $d$ is a positive multiple of $q(p,r)$, every element of $\Kt^0(P^d(p^r))$ has a stable cohomotopy lift.
	\end{lemma}
	
	\begin{proof}
		Let $M=\Sph/p^r$, $i\colon\Sph\to M$ the inclusion of the bottom cell, $q=q(p,r)$ and $d=kq$. The composite
		\[
		v^{(k)}=v\circ\Sigma^qv\circ\dots\circ\Sigma^{(k-1)q}v\colon\Sigma^dM\to M
		\]
		is a $\KU_*$-isomorphism, being a composite of $\KU_*$-isomorphisms (Theorem~\ref{thm:adams}), and $i_*\colon\KU_0(\Sph)=\Z\to\KU_0(M)=\Z/p^r$ is onto by the long exact
		sequence of the cofibre sequence $\Sph\xrightarrow{p^r}\Sph\xrightarrow{i}M$. Hence $\beta=v^{(k)}\circ\Sigma^di\in\pi_d(M)$
		has a $\KU$-Hurewicz image which generates $\KU_d(M)=\Z/p^r$: by naturality of the Hurewicz map, $h(\beta)=v^{(k)}_*\Sigma^d(i_*h(\iota))$, where $h(\iota)\in\KU_0(\Sph)$
		is the Hurewicz image of the identity, which is $1$. Let $\alpha'\in[\Sigma^dDM,\Sph]$ correspond to $\beta$ under
		\eqref{eq:duality} (with $E=\Sph$). By naturality of \eqref{eq:duality} in $E$, applied to $\eta\colon\Sph\to\KU$, the class
		$\eta\circ\alpha'\in[\Sigma^dDM,\KU]$ corresponds to the Hurewicz image of $\beta$ in $\pi_d(M\wedge\KU)=\KU_d(M)$; so it generates
		$[\Sigma^dDM,\KU]\cong\Z/p^r$. By Lemma~\ref{lem:mooredual}, $\Sigma^dDM\simeq\Sigma^\infty P^d(p^r)$, so we have found
		$\alpha\colon\Sigma^\infty P^d(p^r)\to\Sph$ such that $\eta\alpha$ generates $\KU^0(\Sigma^\infty P^d(p^r))=\Kt^0(P^d(p^r))$. Its
		adjoint $P^d(p^r)\to Q_0S^0$ (it lands in the zero component because $P^d(p^r)$ is connected) is a stable cohomotopy lift of a
		generator. Multiples of the generator lift because lifts form a subgroup.
	\end{proof}
	
	\begin{remark}[The hypothesis is almost necessary]\label{rem:liftnecessary}
		The Adams operation $\psi^k$ is a map of ring spectra $\KU[1/k]\to\KU[1/k]$, so $\psi^k\eta=\eta$ after inverting $k$; hence a class with a stable
		cohomotopy lift is fixed by $\psi^k$ in $\KU^0(-)[1/k]$. On $p$-primary torsion this gives invariance under $\psi^k$ for every $k$ prime to $p$
		(and no information for $p\mid k$; indeed, on $\Kt^0(P^{2t}(p^r))$ the operation $\psi^p$ is multiplication by $p^t$, which fixes no generator). On
		$\Kt^0(P^d(p^r))=\Z/p^r$, $\psi^k$ acts by $k^{d/2}$, so for odd $p$ a generator can lift only if $k^{d/2}\equiv1\bmod p^r$ for all $k$ prime to $p$,
		i.e.\ only if $q(p,r)\mid d$. For $p=2$ a lift must also factor through
		real $K$-theory. For instance, the generator of $\Kt^0(P^4(2))=\Z/2$ has no lift, whereas on $P^2(2)=\R P^2$ the generator does
		lift: it is the reduced class of the covering $S^2\to\R P^2$ (Example~\ref{ex:RP2}). This is why we take the dimensions in
		Theorem~\ref{thm:realize} to be multiples of the periods: this choice is sufficient, uniform in the primes involved, and, for odd primes, forced
		by the obstruction just described; it is not necessary in every case, as $P^2(2)$ shows.
	\end{remark}
	
	\begin{corollary}[Wedges of Moore spaces]\label{cor:wedgelift}
		Let $X=\bigvee_{j\in J}P^d(p_j^{r_j})$ with $J$ finite ($X$ a point if $J=\emptyset$), $d$ even, and $q(p_j,r_j)\mid d$ for all $j$.
		Let $\gamma^{(j)}$ be the generator $\gamma_d$ of the $j$th summand and $x_j=c_j^*\gamma^{(j)}$ (with $c_j$ the collapse onto the $j$th
		summand). Then $\Kt^0(X)=\bigoplus_j\Z/p_j^{r_j}$ with $x_j$ the $j$th standard generator, $K^1(X)=0$, and every element of
		$\Kt^0(X)$ has a stable cohomotopy lift.
	\end{corollary}
	
	\begin{proof}
		Lemma~\ref{lem:K0basic}(iv) and Lemma~\ref{lem:mooreK}(b) give the groups, and $\iota_l^*x_j=(c_j\iota_l)^*\gamma^{(j)}$ is
		$\gamma^{(j)}$ for $l=j$ and $0$ otherwise. If $\alpha_j$ lifts $\gamma^{(j)}$ (Lemma~\ref{lem:lift}), then $\alpha_j\circ c_j$
		lifts $x_j$. Sums lift because lifts form a subgroup.
	\end{proof}
	
	\subsection{Acyclic maps and the plus construction}\label{ss:acyclic}
	
	\begin{definition}[Perfect and locally perfect groups]
		A group $\Pi$ is \emph{perfect} if $\Pi=[\Pi,\Pi]$, i.e.\ every element is a product of commutators; equivalently,
		$H_1(\Pi;\Z)=\Pi^{ab}=0$. It is \emph{locally perfect} if every finite subset of $\Pi$ is contained in a finitely generated perfect
		subgroup.
	\end{definition}
	
	\begin{example}\label{ex:Ainfty}
		The alternating groups $A_N$, $N\ge5$, are perfect (they are simple and nonabelian). The infinite alternating group
		$A_\infty=\bigcup_NA_N\subset\Sigma_\infty=\bigcup_N\Sigma_N$ is locally perfect: a finite set of even permutations of $\N$ moves
		only finitely many points, so it lies in some $A_N$ with $N\ge5$, which is finite and perfect. Another example is the group
		$SL(2,5)=SL(2,\mathbb F_5)$ of $2\times2$ matrices of determinant $1$ over the field with five elements. It has order $5(5^2-1)=120$, its centre is
		$\{\pm1\}$, and $SL(2,5)/\{\pm1\}=PSL(2,5)\cong A_5$; it is the \emph{binary icosahedral group}, the preimage of the icosahedral rotation group $A_5\subset SO(3)$
		under the double covering $S^3=SU(2)\to SO(3)$. It is perfect (as is $SL(2,\mathbb F_q)$ for every $q\ge4$), and, being a finite subgroup of the group
		$S^3$, it acts freely on $S^3$ by left multiplication; the quotient $S^3/SL(2,5)$ is the Poincaré homology sphere, a closed $3$-manifold with the integral
		homology of $S^3$ and fundamental group $SL(2,5)$.
	\end{example}
	
	\begin{lemma}\label{lem:locperf}
		\begin{enumerate}[label=\textup{(\alph*)}]
			\item Locally perfect groups are perfect; finitely generated perfect groups, and finite products of locally perfect groups, are
			locally perfect.
			\item \emph{(Central extensions.)} If $\Pi$ is perfect and $\Pi\twoheadrightarrow\Lambda$ is a central extension of a locally perfect
			group $\Lambda$ (i.e.\ the kernel is central in $\Pi$), then $\Pi$ is locally perfect.
		\end{enumerate}
	\end{lemma}
	
	\begin{proof}
		(a) If $g$ lies in a perfect subgroup $H$, then $g\in[H,H]\subset[\Pi,\Pi]$. A finitely generated perfect group $\Pi$ is its own
		finitely generated perfect subgroup containing any finite subset. If $S\subset\Pi_1\times\Pi_2$ is finite, its projections lie
		in finitely generated perfect subgroups $H_1,H_2$, and $S\subset H_1\times H_2$, which is finitely generated and perfect.
		
		(b) Let $S\subset\Pi$ be finite. Each $s\in S$ is a finite product of commutators of elements of $\Pi$; let $Y$ be the finite set
		of all elements so used. The image of $Y$ in $\Lambda$ lies in a finitely generated perfect subgroup $\Lambda_0$. Let
		$\mathcal E\le\Pi$ be the preimage of $\Lambda_0$, $\mathcal D=[\mathcal E,\mathcal E]$, and $\mathcal A=\ker(\Pi\to\Lambda)$,
		which is central. Then $S\subset\mathcal D$ (all commutators used are commutators of elements of $Y\subset\mathcal E$), and
		$\mathcal D$ maps onto $[\Lambda_0,\Lambda_0]=\Lambda_0$, so $\mathcal E=\mathcal D\mathcal A$. Since $\mathcal A$ is central,
		commutators in $\mathcal E=\mathcal D\mathcal A$ are commutators in $\mathcal D$, so $\mathcal D=[\mathcal E,\mathcal E]=[\mathcal
		D,\mathcal D]$. Choose $h_1,\dots,h_m\in\mathcal D$ mapping onto generators of $\Lambda_0$ and let $H=\langle h_1,\dots,h_m\rangle\le
		\mathcal D$. For $x\in\mathcal D$ choose $h\in H$ with the same image; then $h^{-1}x\in\mathcal D\cap\mathcal A$. So $\mathcal
		D=H(\mathcal D\cap\mathcal A)$ with $\mathcal D\cap\mathcal A$ central, and as before $\mathcal D=[\mathcal D,\mathcal D]=[H,H]\le H$.
		Thus $\mathcal D=H$ is finitely generated and perfect and contains $S$.
	\end{proof}
	
	Both hypotheses of (b) are needed: $A_5\times\Z/2\to A_5$ is a central extension of a perfect group by a group which is not
	perfect, and the kernel of an extension of a locally perfect group need not be central.
	
	\begin{definition}[Acyclic spaces and maps]
		A space $F$ is \emph{acyclic} if it is nonempty and $\widetilde H_*(F;\Z)=0$; in particular it is path-connected. A map $q\colon
		P\to X$ is \emph{acyclic} if all its homotopy fibres $\hofib_q(x)$, $x\in X$, are acyclic.
	\end{definition}
	
	An acyclic map is a map which is invisible to homology, not only globally but fibrewise: it may change the fundamental group drastically, but only by killing
	perfect subgroups, and it does not change homology with any coefficients coming from the target.
	
	\begin{definition}[The plus construction {\cite{HH}, \cite[\S4.2]{Hat}}]\label{def:plus}
		Let $P$ be a connected CW complex and $\Pi\trianglelefteq\pi_1(P)$ a perfect normal subgroup. A \emph{plus construction} of $P$
		with respect to $\Pi$ is a map $P\to P^+_\Pi$ to a CW complex which induces the quotient map $\pi_1(P)\to\pi_1(P)/\Pi$ and is
		acyclic. It exists and is unique up to homotopy equivalence under $P$. Construction: attach $2$-cells to $P$ along loops
		representing generators of $\Pi$; this kills $\Pi$ (Theorem~\ref{thm:vankampen}(b)) but creates new classes in $H_2$; since $\Pi$ is perfect, each new $2$-cell's
		attaching loop is a product of commutators, and one can show that the new homology classes are spherical in the universal cover
		of the intermediate complex, so they can be killed by attaching $3$-cells, which does not change $\pi_1$; the result is $P^+_\Pi$ \cite[\S4.2]{Hat},
		\cite{HH}.
	\end{definition}
	
	\begin{theorem}[Properties of acyclic maps {\cite{HH}}, {\cite[\S4.1]{Ra}}]\label{thm:acyclic}
		\begin{enumerate}[label=\textup{(A\arabic*)}]
			\item Acyclic maps induce isomorphisms on integral homology, and indeed on homology with any local coefficient system pulled back
			from the target. Conversely, a homology isomorphism to a simply connected space is acyclic.
			\item Products and homotopy pullbacks of acyclic maps are acyclic; maps homotopic to acyclic maps, and composites of acyclic maps
			with weak homotopy equivalences on either side, are acyclic.
			\item If $P$ is a connected CW complex, the acyclic maps $q\colon P\to X$ to spaces $X$ of the homotopy type of a CW complex are, up to homotopy
			equivalence under $P$, exactly the plus constructions with respect to perfect normal subgroups of $\pi_1(P)$; the subgroup is $\ker\pi_1(q)$, and
			$\pi_1(q)$ is onto.
			\item If $q\colon P\to X$ is acyclic and $X$ is path-connected, then $P$ is nonempty and path-connected.
		\end{enumerate}
	\end{theorem}
	
	We explain what these statements say. (A1) is the defining property of acyclic maps seen from homology: by the Serre spectral sequence, a map whose homotopy
	fibres have the homology of a point is a homology isomorphism, even with twisted coefficients. Its converse for simply connected targets says that for such
	targets there is no difference between homology isomorphisms and acyclic maps. (A2) lists the operations under which acyclic maps are stable. The most
	important for us is the stability under homotopy pullback: a homotopy pullback of $q$ along any map has the same homotopy fibres as $q$; this is what makes the
	pullback $P$ of Figure~\ref{fig:strategy} acyclic over $X_n$. The stability under products uses the Künneth theorem (a product of acyclic spaces is acyclic). (A3)
	identifies acyclic maps with plus constructions; in particular, an acyclic map is surjective on $\pi_1$ with perfect kernel, and conversely the kernel of $\pi_1$
	determines the acyclic map up to homotopy equivalence. We use (A3) to translate between Raptis' formulation of the Hausmann--Vogel theorem and ours. (A4) is the
	statement in degree zero: an acyclic map to a path-connected space has path-connected source.
	
	\begin{example}
		(1) The Poincaré homology sphere $\Sigma^3=S^3/SL(2,5)$ of Example~\ref{ex:Ainfty} has the integral homology of $S^3$, and a degree one map
		$\Sigma^3\to S^3$ is a homology isomorphism to a simply connected space, hence acyclic by (A1). It is the plus construction
		with respect to the perfect group $\pi_1(\Sigma^3)=SL(2,5)$. Here the source is a finite complex with a large perfect
		fundamental group which maps acyclically to a simply connected finite complex. This is the shape of the maps $a_n\colon T_n\to
		X_n$ we shall construct: $X_n$ is a wedge of Moore spaces, $T_n$ has the same homology, and the perfect group $\pi_1(T_n)$
		carries the coverings.
		
		(2) The prototype of an acyclic map with infinite complexes is the Barratt--Priddy--Quillen map of the next subsection.
	\end{example}
	
	\subsection{The Barratt--Priddy--Quillen theorem}\label{ss:BPQ}
	
	This subsection provides the link between stable cohomotopy and finite coverings. We first define the map $b\colon B\Sigma_\infty\to Q_0S^0$, then state the
	Barratt--Priddy--Quillen theorem in the two forms we need and explain how they follow from the literature, and finally compute the fundamental
	group of the homotopy fibre of $b$ (Lemma~\ref{lem:Fb}), which is the input needed for the Hausmann--Vogel theorem.
	
	\subsubsection*{The monoid of symmetric groups}
	For $N,N'\ge0$ the \emph{block sum} $\Sigma_N\times\Sigma_{N'}\to\Sigma_{N+N'}$, $(\sigma,\tau)\mapsto\sigma\oplus\tau$, lets $\sigma$ act on $\{1,\dots,N\}$ and $\tau$ on
	$\{N+1,\dots,N+N'\}$ (shifted by $N$). It is an injective homomorphism, and it is strictly associative. By (B3) and (B4) of Section~\ref{ss:cov3} it induces maps
	$B\Sigma_N\times B\Sigma_{N'}=B(\Sigma_N\times\Sigma_{N'})\to B\Sigma_{N+N'}$, and
	\[
	M=\bigsqcup_{N\ge0}B\Sigma_N
	\]
	becomes a topological monoid, with unit the point $B\Sigma_0$ and with $\pi_0(M)=\N$. It is homotopy commutative: $\sigma\oplus\tau$ and $\tau\oplus\sigma$ are conjugate by a
	fixed permutation, and conjugation by a fixed element induces a self-map of a classifying space which is homotopic to the identity (B3). The inclusion
	$B\Sigma_N\subset B\Sigma_{N+1}$ of Section~\ref{ss:cov3} is the map induced by $\sigma\mapsto\sigma\oplus1$, i.e.\ right multiplication in $M$ with the point $B\Sigma_1$.
	
	\emph{The Pontryagin ring and its localization.} The multiplication $\mu\colon M\times M\to M$ makes $H_*(M;\Z)$ a graded ring, the \emph{Pontryagin ring}, with product
	$x\cdot y=\mu_*(x\times y)$, where $\times$ is the homology cross product. For $N\in\pi_0(M)=\N$ let $[N]\in H_0(M)$ be the class of a point of $B\Sigma_N$; then
	$[N]\cdot[N']=[N+N']$, so $\pi_0(M)$ is a multiplicative subset of $H_*(M)$, generated by $t=[1]$. Since $M$ is homotopy commutative, these elements are central, and
	\[
	H_*(M;\Z)[\pi_0(M)^{-1}]
	\]
	denotes the ring of fractions $x/[N]$ ($x\in H_*(M)$, $N\in\N$), where $x/[N]=y/[N']$ if and only if $x\cdot[N'+K]=y\cdot[N+K]$ for some $K\in\N$. Multiplication by $t$ maps
	$H_*(B\Sigma_N)$ to $H_*(B\Sigma_{N+1})$ by the homomorphism induced by the inclusion $B\Sigma_N\subset B\Sigma_{N+1}$. Hence, as a graded abelian group, the localization is
	the colimit of the sequence $H_*(M)\xrightarrow{t}H_*(M)\xrightarrow{t}\cdots$, that is,
	\begin{equation}\label{eq:localization}
		H_*(M;\Z)[\pi_0(M)^{-1}]\cong\bigoplus_{k\in\Z}\operatorname{colim}_NH_*(B\Sigma_N;\Z)=H_*(\Z\times B\Sigma_\infty;\Z),
	\end{equation}
	where the summand indexed by $k$ consists of the fractions $x/[N]$ with $x\in H_*(B\Sigma_{N+k})$, and where we have used that homology commutes with the
	colimit of the subcomplexes $B\Sigma_N$ of $B\Sigma_\infty$. For example, in degree $0$ the localization is $\Z[t,t^{-1}]$, the localization of $H_0(M)=\Z[\N]=\Z[t]$.
	
	The monoid $M$ has a classifying space $BM$ (the realization of the nerve of $M$ regarded as a topological category with one object) and a natural
	\emph{group completion map} $M\to\Omega BM$, which is a map of $H$-spaces, is the inclusion $\N\to\Z$ on $\pi_0$, and has as target a group up to homotopy (loop
	multiplication has homotopy inverses). Its degree $N$ component $\Omega_NBM$ receives $B\Sigma_N$.
	
	\begin{theorem}[{\cite{BP}}, {\cite{Pr71}}, {\cite[Proposition~3.5]{Se74}}]\label{thm:BPQbasic}
		There is a homotopy equivalence $\Omega BM\simeq QS^0$ which is the identity on $\pi_0=\Z$. The group completion map induces an isomorphism of rings
		\[
		H_*(M;\Z)[\pi_0(M)^{-1}]\xrightarrow{\ \cong\ }H_*(\Omega BM;\Z),\qquad x/[N]\mapsto(\text{group completion})_*(x)\cdot[\ell]^N ,
		\]
		where $[\ell]\in H_0(\Omega_{-1}BM)$ is the class of a point of degree $-1$ and $\cdot$ is the Pontryagin product of $H_*(\Omega BM)$.
	\end{theorem}
	
	The first statement says that $QS^0$ is the group completion of the monoid of finite sets; in Segal's language, the $K$-theory spectrum of the category of finite
	sets and bijections is the sphere spectrum \cite[Proposition~3.5]{Se74}. The second statement is the \emph{group completion theorem} for $M$ \cite{BP}: the homology
	of the group completion is obtained from the Pontryagin ring of $M$ by inverting the components. Combined with \eqref{eq:localization} it says that the group completion
	map, followed by translation into the degree-zero component, induces an isomorphism $H_*(B\Sigma_\infty)=\operatorname{colim}_NH_*(B\Sigma_N)\to H_*(\Omega_0BM)$; this is made
	precise in Definition~\ref{def:bmap} and Theorem~\ref{thm:BPQ}(i).
	
	\begin{definition}[The standard group-completion map]\label{def:bmap}
		The group completion map sends $B\Sigma_N$ into the component $\Omega_NBM$ of degree $N$, so it does not define a map from $B\Sigma_\infty$ into a single
		component. To obtain one, we translate each component to degree $0$. Fix a point $\ell\in\Omega_{-1}BM$ (a loop of degree $-1$, i.e.\ a homotopy inverse of the
		point $B\Sigma_1$ in the group $\Omega BM$), and for $N\ge0$ let $\tau_N\colon\Omega_NBM\to\Omega_0BM$ be loop multiplication by $\ell^N$; it is a homotopy equivalence,
		with homotopy inverse the multiplication by $(B\Sigma_1)^N$. The maps
		\[
		\tau_N\circ(\text{group completion})|_{B\Sigma_N}\colon B\Sigma_N\to\Omega_0BM\qquad(N\ge0)
		\]
		are compatible up to homotopy: the restriction of the $(N+1)$st map to $B\Sigma_N$ is homotopic to the $N$th map, because the inclusion $B\Sigma_N\subset B\Sigma_{N+1}$
		is right multiplication by the point $B\Sigma_1$, the group completion map is a homotopy homomorphism, and $\ell\cdot B\Sigma_1$ is homotopic to the unit. Since
		$B\Sigma_N$ is a subcomplex of $B\Sigma_{N+1}$, homotopy extension (Theorem~\ref{thm:standardhtpy}(b)) allows us to replace these maps inductively by maps
		$\beta_N\colon B\Sigma_N\to\Omega_0BM$ with $\beta_N\simeq\tau_N\circ(\text{group completion})|_{B\Sigma_N}$ and $\beta_{N+1}|_{B\Sigma_N}=\beta_N$: having chosen $\beta_N$, the
		homotopy from the restriction of the $(N+1)$st map to $\beta_N$ on $B\Sigma_N$ extends to a homotopy on $B\Sigma_{N+1}$, whose end is $\beta_{N+1}$. The $\beta_N$ define a
		map
		\[
		b\colon B\Sigma_\infty=\bigcup_NB\Sigma_N\to\Omega_0BM\simeq Q_0S^0 ,
		\]
		the \emph{standard group-completion map} \cite{BP}, \cite{Se74}. We fix such a map $b$; it depends on the choices of $\ell$ and of the homotopies, but by
		construction $b|_{B\Sigma_N}\simeq\tau_N\circ(\text{group completion})|_{B\Sigma_N}$ for every $N$, and only this property is used.
	\end{definition}
	
	\begin{theorem}[{\cite{BP}}, {\cite{Pr71}}, {\cite{Se74}}]\label{thm:BPQ}
		Let $\Sigma_\infty=\bigcup_N\Sigma_N$, with alternating subgroup $A_\infty$.
		\begin{enumerate}[label=\textup{(\roman*)}]
			\item The map $b$ induces an isomorphism on integral homology. It is acyclic, and $\pi_1(b)\colon\Sigma_\infty\to\pi_1(Q_0S^0)=\pi^s_1
			=\Z/2$ is the sign homomorphism; so $b$ is the plus construction of $B\Sigma_\infty$ with respect to $A_\infty$:
			$Q_0S^0\simeq B\Sigma_\infty^+$.
			\item For every $N\ge1$, every finite CW complex $A$ and every map $c\colon A\to B\Sigma_N$, the composite
			\[
			A\xrightarrow{c}B\Sigma_N\xrightarrow{\iota_N}B\Sigma_\infty\xrightarrow{b}Q_0S^0\xrightarrow{\Omega^\infty\eta}BU
			\]
			represents $[c^*\C^N_\perm]-N\in\Kt^0(A)$, where $\iota_N\colon B\Sigma_N\subset B\Sigma_\infty$. By (C4), this is the reduced covering class $\tc(c^*U_N)$ of the
			covering classified by $c$.
		\end{enumerate}
	\end{theorem}
	
	\subsubsection*{Explanation of part (i)}
	\emph{Homology.} The isomorphism of Theorem~\ref{thm:BPQbasic} sends the fraction $x/[N]$, $x\in H_*(B\Sigma_N)$, to $(\tau_N\circ\text{group completion})_*(x)$, and by
	\eqref{eq:localization} these fractions form the summand $H_*(B\Sigma_\infty)$ of index $k=0$, which is mapped onto $H_*(\Omega_0BM)$. Since $b$ restricted to $B\Sigma_N$
	is homotopic to $\tau_N\circ(\text{group completion})$ (Definition~\ref{def:bmap}), $H_*(b)$ is this isomorphism $H_*(B\Sigma_\infty)\to H_*(\Omega_0BM)=H_*(Q_0S^0)$.
	
	\emph{Fundamental group.} $Q_0S^0$ is a loop space, so $\pi_1(Q_0S^0)$ is abelian and the Hurewicz map $\pi_1(Q_0S^0)\to H_1(Q_0S^0)$ is an isomorphism; and
	$\pi_1(Q_0S^0)=\pi_1^s=\Z/2$. The Hurewicz map of $B\Sigma_\infty$ is the abelianization $\Sigma_\infty\to\Sigma_\infty^{ab}$. Since $[\Sigma_N,\Sigma_N]=A_N$ for all $N$, we
	have $[\Sigma_\infty,\Sigma_\infty]=A_\infty$ and $\Sigma_\infty^{ab}=\Z/2$, the abelianization being the sign. By naturality of the Hurewicz map and since $H_1(b)$ is an
	isomorphism, $\pi_1(b)$ is surjective with kernel $A_\infty$. There is only one surjection $\Sigma_\infty\to\Z/2$, so $\pi_1(b)$ is the sign.
	
	\emph{Acyclicity.} The target $Q_0S^0$ is not simply connected, so the converse part of (A1) does not apply directly; the acyclicity of $b$ is Priddy's theorem
	$Q_0S^0\simeq B\Sigma_\infty^+$ \cite{Pr71}. In outline: the kernel $A_\infty$ of $\pi_1(b)$ is perfect (Example~\ref{ex:Ainfty}); if $j\colon B\Sigma_\infty\to B\Sigma_\infty^+$ is a plus
	construction with respect to $A_\infty$ (Definition~\ref{def:plus}), the universal property of the plus construction \cite{HH} gives a map $h\colon B\Sigma_\infty^+\to Q_0S^0$
	with $h\circ j\simeq b$; since $j$ and $b$ are homology isomorphisms, so is $h$, and a homology isomorphism between $H$-spaces of the homotopy type of CW complexes
	(here $B\Sigma_\infty^+$ carries the $H$-space structure induced by block sum) is a homotopy equivalence. Hence $b\simeq h\circ j$ is acyclic by (A2).
	
	\subsubsection*{Explanation of part (ii)}
	Part (ii) says that, under the equivalence $Q_0S^0\simeq B\Sigma_\infty^+$, the unit $\eta$ of $K$-theory is induced by the permutation representations. This follows
	from Segal's theory of categories with sum \cite{Se74}, as follows. Segal associates to every (topological) category $\mathcal C$ with a sum, i.e.\ a symmetric monoidal
	structure, a connective spectrum $\mathbf K(\mathcal C)$ whose zeroth space is the group completion $\Omega B|\mathcal C|$ of the classifying space $|\mathcal C|$; the
	construction is functorial for functors preserving sums, and for such a functor $\Phi$ the infinite loop map $\Omega^\infty\mathbf K(\Phi)$ restricts, along the group
	completion maps, to $|\Phi|$ up to homotopy \cite[\S\S1--2]{Se74}. For the category $\mathcal F$ of finite sets and bijections (with disjoint union) one has
	$|\mathcal F|\simeq M$ and $\mathbf K(\mathcal F)\simeq\Sph$ (Theorem~\ref{thm:BPQbasic}); for the category $\mathcal V$ of finite-dimensional Hermitian vector spaces and
	unitary isomorphisms (with orthogonal sum), $|\mathcal V|\simeq\bigsqcup_NBU(N)$, $\mathbf K(\mathcal V)$ is the connective $K$-theory spectrum with zeroth space
	$\Z\times BU$, and the group completion map $BU(N)\to\{N\}\times BU$ classifies the universal bundle $\gamma_N$ \cite[\S2]{Se74}. The functor $S\mapsto\C[S]$
	preserves sums and induces a map of spectra $\Sph=\mathbf K(\mathcal F)\to\mathbf K(\mathcal V)$, which is determined up to homotopy by the element $[\C[\pt]]=1$ of
	$\pi_0$ that it represents; composed with the map from connective to periodic $K$-theory it is therefore the unit $\eta$. Restricting along the group completion
	maps, we obtain for every $N$ a square
	\[
	\begin{tikzcd}[column sep=large]
		B\Sigma_N\arrow[r]\arrow[d,"\text{perm}"'] & Q_NS^0\arrow[d,"\Omega^\infty\eta"]\\
		BU(N)\arrow[r] & \{N\}\times BU
	\end{tikzcd}
	\]
	which commutes up to homotopy, with horizontal arrows the group completion maps and left arrow $B(\Sigma_N\to U(N))$, which classifies $\C^N_\perm$ (C4). Finally,
	on $B\Sigma_N$ the map $b$ is the group completion map followed by loop multiplication with $\ell^N$ (Definition~\ref{def:bmap}); as $\Omega^\infty\eta$ is an $H$-map, it
	carries this translation to the translation of $\Z\times BU$ by the image of $\ell^N$, which is the point $-N\in\pi_0(\Z\times BU)$. Hence $\Omega^\infty\eta\circ b\circ\iota_N$
	classifies $[\C^N_\perm]-N\in\Kt^0$ on every finite subcomplex of $B\Sigma_N$, and composing with $c\colon A\to B\Sigma_N$ gives (ii).
	
	\subsubsection*{The homotopy fibre of $b$}
	The following lemma is the group-theoretic input of Proposition~\ref{prop:lift}.
	
	\begin{lemma}[The homotopy fibre of $b$]\label{lem:Fb}
		Let $F_b$ be the homotopy fibre of $b$ over a point of $Q_0S^0$.
		\begin{enumerate}[label=\textup{(\alph*)}]
			\item $F_b$ is acyclic, in particular path-connected.
			\item There is a central extension
			\[
			1\to\Z/2\to\pi_1(F_b)\to A_\infty\to1 .
			\]
			\item $\pi_1(F_b)$ is perfect and locally perfect, but not finitely generated.
		\end{enumerate}
	\end{lemma}
	
	\begin{proof}
		(a) is the acyclicity of $b$ (Theorem~\ref{thm:BPQ}(i)); acyclic spaces are path-connected by definition. The homotopy type of $F_b$ does not depend on the point, as
		$Q_0S^0$ is path-connected.
		
		(b) Replace $b$ by the fibration $E_b=\{(y,\omega)\in B\Sigma_\infty\times(Q_0S^0)^{[0,1]}:\omega(0)=b(y)\}\to Q_0S^0$, $(y,\omega)\mapsto\omega(1)$; its fibre is $F_b$ and
		$E_b\simeq B\Sigma_\infty$ by $(y,\omega)\mapsto y$. The exact sequence \eqref{eq:LES} gives
		\[
		\pi_2(B\Sigma_\infty)\to\pi_2(Q_0S^0)\xrightarrow{\partial}\pi_1(F_b)\to\pi_1(B\Sigma_\infty)\xrightarrow{\pi_1(b)}\pi_1(Q_0S^0).
		\]
		Here $\pi_2(B\Sigma_\infty)=0$: every map $S^k\to B\Sigma_\infty$ and every homotopy between such maps has compact image, which lies in some $B\Sigma_N$ (Section~\ref{ss:CW}), so
		$\pi_k(B\Sigma_\infty)=\operatorname{colim}_N\pi_k(B\Sigma_N)$, which is $0$ for $k\ge2$ and $\Sigma_\infty$ for $k=1$ by (B1). Further $\pi_2(Q_0S^0)=\pi_2^s=\Z/2$, and by part (i) of
		Theorem~\ref{thm:BPQ} the kernel of $\pi_1(b)$ is $A_\infty$. So the sequence becomes $1\to\Z/2\to\pi_1(F_b)\to A_\infty\to1$. The image of $\partial$ is central by
		Theorem~\ref{thm:fibLES}(b).
		
		(c) Since $F_b$ is path-connected, $\pi_1(F_b)^{ab}\cong H_1(F_b)=0$ by the Hurewicz theorem and (a); so $\pi_1(F_b)$ is perfect. By (b) it is a central extension of the locally
		perfect group $A_\infty$ (Example~\ref{ex:Ainfty}), hence locally perfect by Lemma~\ref{lem:locperf}(b). It is not finitely generated, because its quotient $A_\infty$ is not: a
		finitely generated subgroup of $A_\infty$ lies in some $A_N$.
	\end{proof}
	
	The lemma explains why the Hausmann--Vogel theorem has to be applied with an infinitely generated kernel: the fundamental group of the homotopy pullback $P$ of
	Figure~\ref{fig:strategy} contains copies of $\pi_1(F_b)$.
	
	\subsection{The Hausmann--Vogel theorem}\label{ss:HV}
	
	The homotopy pullback in Figure~\ref{fig:strategy} has the right homology and carries a map to $B\Sigma_\infty$, but it is far from finite. A
	finite complex is needed: the spaces $Z_n$ of Section~\ref{sec:blocks} must be compact, carry chosen coverings, and admit
	surjections from $[0,1]$. The naive attempt, replacing $P$ by a finite subcomplex of a CW approximation, does not work: a finite subcomplex has the wrong
	homology in general, and the restriction of $q$ to it is not acyclic. The Hausmann--Vogel theorem says that a finite complex with the right homology exists
	after all, provided the target is finite and the kernel of $\pi_1(q)$ is locally perfect.
	
	\begin{definition}
		A space is \emph{homotopy finite} if it is homotopy equivalent to a finite CW complex, and \emph{finitely dominated} if it is a
		homotopy retract of a finite CW complex. Wall showed that a finitely dominated space $X$ is homotopy finite iff an obstruction in
		$\widetilde K_0(\Z[\pi_1X])$ vanishes \cite{Wall}.
	\end{definition}
	
	\begin{theorem}[{\cite[Theorem~3.1]{HV}}, {\cite[Theorem~4.1]{Ra}}]\label{thm:HV}
		Let $P$ be a connected CW complex and $q\colon P\to X$ an acyclic map such that $\ker\pi_1(q)$ is locally perfect, and suppose
		that $X$ is homotopy finite. Then there are a finite CW complex $L$ and a map $\gamma\colon L\to P$ such that $q\circ\gamma$ is
		acyclic.
	\end{theorem}
	
	In words: an acyclic map onto a homotopy finite space whose $\pi_1$-kernel is locally perfect can be \emph{refined} by a finite complex, in the sense that a finite complex
	$L$ maps to the source, and the composite to the target is still acyclic. Note that $L$ is not required to be a subcomplex of $P$, and $\gamma$ is not required to be injective
	in any sense; what matters is that $L$ is finite, that $L$ has the same homology as $X$ (by (A1)), and that $L$ maps to $P$, so that everything $P$ carries (in our case,
	maps to $B\Sigma_\infty$ and homotopies to the given stable cohomotopy classes) is pulled back to $L$.
	
	The theorem is due to Hausmann and Vogel; we use it in the form stated by Raptis: \emph{Let $X$ be a path-connected space and let $f\colon X\to X^+_P$ be the plus construction with
		respect to a locally perfect normal subgroup $P$ of $\pi_1(X,x)$. Suppose that $X^+_P$ is homotopy finite. Then there is a finite CW
		complex $L$ and a map $\gamma\colon L\to X$ such that the composition $L\to X\to X^+_P$ is again acyclic} \cite[Theorem~4.1]{Ra}, where it is recorded as a special case of
	\cite[Theorem~3.1]{HV}. (In this quotation $X$ and $P$ denote Raptis' space and subgroup, not ours.)
	
	\begin{proof}[Derivation of Theorem~\ref{thm:HV} from Raptis' formulation]
		Let $P$, $q$ and $X$ be as in Theorem~\ref{thm:HV}, and put $\Pi=\ker\pi_1(q)$. As $X$ is homotopy finite, it has the homotopy type of a CW complex, so by (A3) there are a plus
		construction $j\colon P\to P^+_\Pi$ with respect to the perfect normal subgroup $\Pi$ and a homotopy equivalence $h\colon P^+_\Pi\to X$ with $h\circ j\simeq q$. The space $P^+_\Pi$ is
		homotopy finite, because it is homotopy equivalent to $X$, and $\Pi$ is locally perfect by hypothesis. Raptis' formulation, applied to the path-connected space $P$ and to $j$,
		gives a finite CW complex $L$ and a map $\gamma\colon L\to P$ such that $j\circ\gamma$ is acyclic. Then $q\circ\gamma$ is homotopic to $h\circ(j\circ\gamma)$, the composite of an
		acyclic map with a homotopy equivalence, and it is acyclic by (A2).
	\end{proof}
	
	\begin{remark}[On the hypotheses and the conclusion]
		\begin{enumerate}[label=(\arabic*)]
			\item Raptis' formulation allows an arbitrary path-connected source. We use the theorem only for CW complexes, after a CW approximation (Step~4 of
			the proof of Proposition~\ref{prop:lift}). The kernel $\ker\pi_1(q)$ need not be finitely generated; only local perfectness is used. This occurs in our
			application, where $\pi_1(P)$ contains the infinitely generated group $\pi_1(F_b)$ (Lemma~\ref{lem:Fb}).
			\item The point of the theorem is that $L$ is \emph{finite}, not merely finitely dominated. The usual trick for removing a Wall
			obstruction, taking a product with $S^1$, is not available: it destroys acyclicity over $X$.
			\item A finite complex $L$ has a finitely presented fundamental group, whereas $\ker\pi_1(q)$ may be infinitely generated;
			local perfectness says that the kernel is the union of its finitely generated perfect subgroups, which is what makes finite
			approximations possible. Conversely, the conclusion forces perfectness of the relevant subgroups: by (A3) applied to $q\circ\gamma$, the group $\pi_1(L)$ maps onto $\pi_1(X)$
			with perfect kernel. In our application $\pi_1(X)=0$, so $\pi_1(L)$ itself is a finitely presented perfect group.
			\item The hypothesis that $X$ be homotopy finite cannot be dropped: if $L\to X$ is acyclic with $L$ finite, then $H_*(X;\Z)\cong H_*(L;\Z)$ is finitely generated and
			vanishes in large degrees.
		\end{enumerate}
		We do not reproduce the proof of the theorem and refer to \cite[\S3]{HV} for the original argument and to \cite[\S4]{Ra} for the form used here.
	\end{remark}
	
	\subsection{Simultaneous finite lifting}\label{ss:lift}
	
	We now combine Sections~\ref{ss:adams}--\ref{ss:HV}. The following proposition produces, from an acyclic map $a\colon Z\to X$ and
	a map $f\colon X'\to X$, a finite complex $Z'$ over $X'$ which maps to $Z$ compatibly with $f$ and carries coverings
	representing prescribed classes. In the inductive construction, $X\to X'$ will be the passage from the wedge of Moore spaces of stage $n$ to that of stage $n+1$, $Z=T_n$ the
	complex already constructed, and $Z'=T_{n+1}$ the new complex. The word \emph{simultaneous} refers to the fact that the two requirements, coverings for the classes $x_j$ on
	$X'$ and a map to the old complex $Z$, are met by one and the same finite complex. Neither requirement alone would suffice: coverings without a map $Z'\to Z$ would not give
	connecting maps, and a map to $Z$ without coverings would not give the right unit classes.
	
	We first record the elementary properties of homotopy pullbacks which are used.
	
	\begin{lemma}[Homotopy pullbacks]\label{lem:hpb}
		Let $F\colon A\to B$ and $G\colon X'\to B$ be maps, and
		\[
		P=\{(y,\omega,x')\in A\times B^{[0,1]}\times X':\ \omega(0)=F(y),\ \omega(1)=G(x')\},
		\]
		with the maps $q\colon P\to X'$, $(y,\omega,x')\mapsto x'$, and $\mathrm{pr}_A\colon P\to A$, $(y,\omega,x')\mapsto y$.
		\begin{enumerate}[label=\textup{(\alph*)}]
			\item $q$ is a fibration, and for every $x'\in X'$ the fibre $q^{-1}(x')$ is homeomorphic to $\hofib_F(G(x'))$.
			\item The map $P\times[0,1]\to B$, $((y,\omega,x'),t)\mapsto\omega(t)$, is a homotopy from $F\circ\mathrm{pr}_A$ to $G\circ q$.
			\item If $A=A_1\times A_2$, $B=B_1\times B_2$ and $F=F_1\times F_2$, then $\hofib_F((b_1,b_2))\cong\hofib_{F_1}(b_1)\times\hofib_{F_2}(b_2)$.
			\item A product of finitely many acyclic spaces is acyclic. Consequently, if $F$ is a product of acyclic maps, then $q$ is acyclic.
		\end{enumerate}
	\end{lemma}
	
	\begin{proof}
		Throughout, path spaces carry the compact-open topology, so that the evaluation map $B^{[0,1]}\times[0,1]\to B$ is continuous and a map $W\to B^{[0,1]}$ is continuous if and
		only if the corresponding map $W\times[0,1]\to B$ is continuous (the interval being locally compact Hausdorff).
		
		(a) Let $E_F=\{(y,\omega)\in A\times B^{[0,1]}:\omega(0)=F(y)\}$ and $e\colon E_F\to B$, $e(y,\omega)=\omega(1)$. We first show that $e$ is a fibration. Let $W$ be a space,
		$h\colon W\times[0,1]\to B$ a homotopy, and $\tilde h_0\colon W\to E_F$, $w\mapsto(y(w),\omega_w)$, a lift of $h(\cdot,0)$, so that $\omega_w(1)=h(w,0)$. For $t\in[0,1]$ let $\omega_w^t$ be the path
		which first traverses $\omega_w$ and then the path $h(w,\cdot)|_{[0,t]}$:
		\[
		\omega_w^t(u)=\begin{cases}\omega_w\big(u(1+t)\big)&0\le u\le\tfrac1{1+t},\\[2pt] h\big(w,u(1+t)-1\big)&\tfrac1{1+t}\le u\le1.\end{cases}
		\]
		The two formulas agree at $u=\frac1{1+t}$, since $\omega_w(1)=h(w,0)$; moreover $\omega_w^0=\omega_w$, $\omega_w^t(0)=F(y(w))$ and $\omega_w^t(1)=h(w,t)$, and
		$(w,t)\mapsto(y(w),\omega_w^t)$ is continuous because the corresponding map $W\times[0,1]\times[0,1]\to B$ is. Thus $\tilde h(w,t)=(y(w),\omega^t_w)$ is a lift
		$W\times[0,1]\to E_F$ of $h$ extending $\tilde h_0$, and $e$ has the homotopy lifting property \cite[\S4.3]{Hat}.
		
		The space $P$ is the pullback $\{(x',(y,\omega)):G(x')=e(y,\omega)\}$ of $e$ along $G$, and $q$ is the projection to $X'$. The pullback of a fibration is a fibration: a
		homotopy $h\colon W\times[0,1]\to X'$ with a lift $\tilde h_0\colon W\to P$ of $h(\cdot,0)$ gives the homotopy $G\circ h$ in $B$ with the lift $\mathrm{pr}_{E_F}\circ\tilde h_0$ to
		$E_F$, which extends to a lift $H\colon W\times[0,1]\to E_F$ of $G\circ h$; then $(h,H)\colon W\times[0,1]\to P$ is a lift of $h$ extending $\tilde h_0$. The fibre of $q$ over $x'$
		is $\{(y,\omega,x'):\omega(0)=F(y),\ \omega(1)=G(x')\}$, which is $\hofib_F(G(x'))\times\{x'\}$ by the definition of the homotopy fibre.
		
		(b) The map is continuous by the remark on evaluation at the beginning of the proof, and at $t=0$ and $t=1$ it takes the values $\omega(0)=F(y)=F(\mathrm{pr}_A(y,\omega,x'))$
		and $\omega(1)=G(x')=G(q(y,\omega,x'))$.
		
		(c) A path in $B_1\times B_2$ is a pair of paths, and $(B_1\times B_2)^{[0,1]}=B_1^{[0,1]}\times B_2^{[0,1]}$ as topological spaces (the compact-open topology of a product is the
		product of the compact-open topologies). The conditions $\omega(0)=F(y)$ and $\omega(1)=(b_1,b_2)$ are the pairs of conditions $\omega_i(0)=F_i(y_i)$ and $\omega_i(1)=b_i$. So
		$((y_1,y_2),(\omega_1,\omega_2))\mapsto((y_1,\omega_1),(y_2,\omega_2))$ is a homeomorphism $\hofib_F((b_1,b_2))\to\hofib_{F_1}(b_1)\times\hofib_{F_2}(b_2)$.
		
		(d) Let $F_1,F_2$ be acyclic spaces. By the Künneth theorem,
		\[
		H_n(F_1\times F_2)\cong\bigoplus_{i+j=n}H_i(F_1)\otimes H_j(F_2)\ \oplus\bigoplus_{i+j=n-1}\Tor\big(H_i(F_1),H_j(F_2)\big).
		\]
		Since $H_0(F_1)=H_0(F_2)=\Z$ and $H_n(F_1)=H_n(F_2)=0$ for $n>0$, the right side is $\Z$ for $n=0$ and $0$ for $n>0$; so $F_1\times F_2$ is acyclic, and finite products are acyclic by
		induction. Now let $F=F_1\times\dots\times F_r$ be a product of acyclic maps $F_i\colon A_i\to B_i$. By (c), applied $r-1$ times, every homotopy fibre of $F$ is a product of
		homotopy fibres of the $F_i$, hence a product of acyclic spaces, hence acyclic. By (a), the fibres of the fibration $q$ are homeomorphic to homotopy fibres of $F$, hence
		acyclic. Finally, for a fibration the fibre over a point is homotopy equivalent to the homotopy fibre over that point \cite[\S4.3]{Hat} (the inclusion
		$q^{-1}(x')\to\hofib_q(x')$, $z\mapsto(z,\mathrm{const}_{x'})$, is a homotopy equivalence); so all homotopy fibres of $q$ are acyclic, and $q$ is acyclic.
	\end{proof}
	
	\begin{proposition}[Simultaneous finite lifting]\label{prop:lift}
		Let $X,X'$ be $2$-connected finite CW complexes, $Z$ a connected finite CW complex, $a\colon Z\to X$ an acyclic map, $f\colon
		X'\to X$ a map, and $x_1,\dots,x_s\in\Kt^0(X')$ ($s\ge0$) classes with stable cohomotopy lifts. Then there are a connected finite
		CW complex $Z'$, a based acyclic map $a'\colon Z'\to X'$, a map $g\colon Z'\to Z$, and finite coverings $\nu_1,\dots,\nu_s$ of $Z'$
		of degree at least $1$, such that
		\begin{equation}\label{eq:liftprop}
			a\circ g\simeq f\circ a'\ \text{(free homotopy)},\qquad\tc(\nu_j)=a'^*x_j\quad(j=1,\dots,s),
		\end{equation}
		and $a'^*\colon K^*(X')\to K^*(Z')$ is an isomorphism. The cases $X=Z=\pt$, $a=\id$, and $X'=\pt$, $s=0$, are allowed.
	\end{proposition}
	
	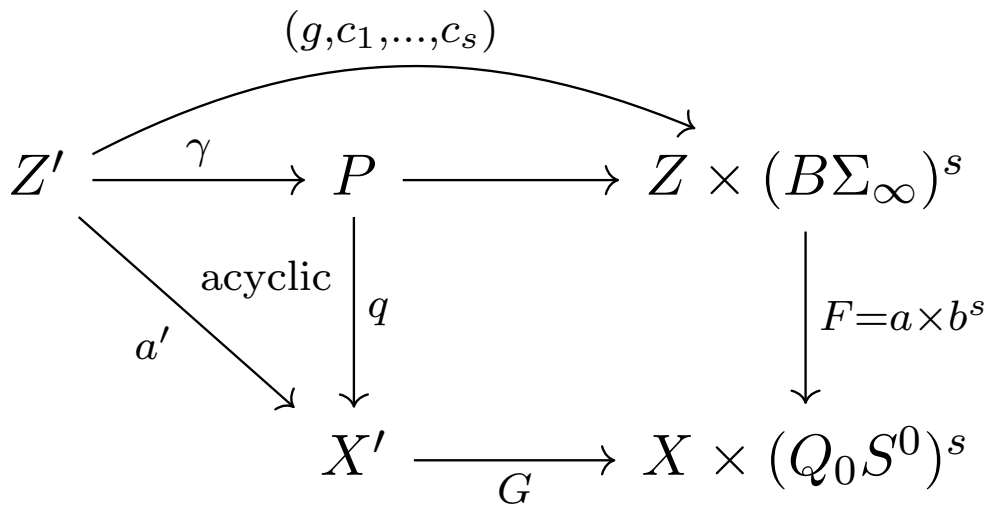
\begin{figure}[ht]
		\centering
		\begin{tikzcd}[column sep=3.5em,row sep=3em]
			Z'\arrow[r,"\gamma"]\arrow[rd,"a'"',"\text{acyclic}"]\arrow[rr,bend left=25,"{(g,c_1,\dots,c_s)}"] & P\arrow[r]\arrow[d,"q"] & Z\times(B\Sigma_\infty)^s\arrow[d,"F=a\times b^s"]\\
			& X'\arrow[r,"G"'] & X\times(Q_0S^0)^s
		\end{tikzcd}
		\caption{The homotopy pullback square in the proof of Proposition~\ref{prop:lift}, with $G=(f,\alpha_1,\dots,\alpha_s)$. The square on the right commutes up to the
			canonical homotopy which is part of the definition of $P$.}\label{fig:pullback}
	\end{figure}
	
	\begin{proof}
		\emph{Step 1: the homotopy pullback.} By hypothesis there are based maps $\alpha_j\colon X'\to Q_0S^0$ lifting the $x_j$
		(Definition~\ref{def:stablelift}). Put
		\begin{gather*}
			A=Z\times(B\Sigma_\infty)^s,\qquad B=X\times(Q_0S^0)^s,\\
			F=a\times b^s\colon A\to B,\qquad G=(f,\alpha_1,\dots,\alpha_s)\colon X'\to B,
		\end{gather*}
		with $b$ the map of Theorem~\ref{thm:BPQ}, and let $P$ and $q\colon P\to X'$ be as in Lemma~\ref{lem:hpb} (Figure~\ref{fig:pullback}). By Lemma~\ref{lem:hpb}(a),(c), $q$ is a
		fibration whose fibre over $x'$ is
		\[
		\hofib_F(G(x'))\cong\hofib_a(f(x'))\times\prod_{j=1}^s\hofib_b(\alpha_j(x')).
		\]
		The factors are acyclic: $a$ is acyclic by hypothesis and $b$ by Theorem~\ref{thm:BPQ}(i). By Lemma~\ref{lem:hpb}(d), $q$ is acyclic.
		
		\emph{Step 2: $P$ is path-connected.} By (A4), since $X'$ is path-connected and $q$ is acyclic.
		
		\emph{Step 3: $\pi_1(P)$ is locally perfect.} Write $F_a$ and $F_b$ for the homotopy fibres of $a$ and $b$; their homotopy types do not depend on the point, since $X$ and
		$Q_0S^0$ are path-connected, and they are path-connected, being acyclic. Fix a basepoint $x'_0\in X'$ and let $\Phi=q^{-1}(x'_0)\cong F_a\times F_b^s$.
		
		\emph{The group $\pi_1(F_a)$.} The exact sequence \eqref{eq:LES} for $a$ (replaced by a fibration with fibre $F_a$ and total space homotopy equivalent to $Z$) reads
		$\pi_2(X)\to\pi_1(F_a)\to\pi_1(Z)\to\pi_1(X)$. Since $X$ is $2$-connected, $\pi_1(F_a)\cong\pi_1(Z)$. This group is finitely generated, since $Z$ is a finite complex, and it is
		perfect, since $\pi_1(Z)^{ab}=H_1(Z)\cong H_1(X)=0$ by (A1) and the Hurewicz theorem. So $\pi_1(F_a)$ is locally perfect (Lemma~\ref{lem:locperf}(a)).
		
		\emph{The group $\pi_1(F_b)$} is locally perfect by Lemma~\ref{lem:Fb}(c).
		
		\emph{The group $\pi_1(P)$.} The exact sequence of the fibration $q$ gives $\pi_2(X')\to\pi_1(\Phi)\to\pi_1(P)\to\pi_1(X')$, and both ends vanish since $X'$ is $2$-connected. Hence
		\[
		\pi_1(P)\cong\pi_1(\Phi)\cong\pi_1(F_a)\times\pi_1(F_b)^s,
		\]
		which is locally perfect by Lemma~\ref{lem:locperf}(a). Since $\pi_1(X')=0$, $\ker\pi_1(q)=\pi_1(P)$.
		
		\emph{Step 4: finite refinement.} Let $\kappa\colon\widetilde P\to P$ be a CW approximation (Theorem~\ref{thm:standardhtpy}(e)). Since $\kappa$ is a weak homotopy equivalence, $\widetilde P$
		is path-connected and $\pi_1(\kappa)$ is an isomorphism; $\tilde q=q\circ\kappa$ is acyclic by (A2); and $\ker\pi_1(\tilde q)=\pi_1(\widetilde P)\cong\pi_1(P)$ is locally perfect by
		Step~3. The target $X'$ is a finite CW complex. Theorem~\ref{thm:HV}, applied to $\tilde q$, gives a finite CW complex $Z'$ and a map $\gamma\colon Z'\to\widetilde P$ such that
		$a':=\tilde q\circ\gamma\colon Z'\to X'$ is acyclic. By (A4), $Z'$ is nonempty and connected.
		
		\emph{Step 5: the map $g$ and the homotopies.} Write $\kappa\gamma(z')=(y(z'),\omega_{z'},a'(z'))\in P$, where
		\[
		y(z')=(g(z'),c_1(z'),\dots,c_s(z'))\in Z\times(B\Sigma_\infty)^s .
		\]
		The maps $g\colon Z'\to Z$ and $c_j\colon Z'\to B\Sigma_\infty$ are continuous, being composites of $\kappa\gamma$ with projections. By Lemma~\ref{lem:hpb}(b), $(z',t)\mapsto\omega_{z'}(t)$
		is a homotopy from $F\circ y$ to $G\circ a'$. Its projection to $X$ is a free homotopy from $a\circ g$ to $f\circ a'$, and its projection to the $j$th factor $Q_0S^0$ is a free
		homotopy from $b\circ c_j$ to $\alpha_j\circ a'$.
		
		\emph{Step 6: coverings.} $Z'$ is compact, so $c_j(Z')$ is a compact subset of the CW complex $B\Sigma_\infty$ and lies in a finite subcomplex \cite[Proposition~A.1]{Hat}. Each
		cell of $B\Sigma_\infty=\bigcup_NB\Sigma_N$ lies in some $B\Sigma_N$, and the $B\Sigma_N$ increase; so the finite subcomplex, and hence $c_j(Z')$, lies in some $B\Sigma_{N_j}$, and
		we may take $N_j\ge1$. Thus $c_j=\iota_{N_j}\circ c_j'$ for a map $c_j'\colon Z'\to B\Sigma_{N_j}$, which is continuous because $B\Sigma_{N_j}$ is a closed subspace. Put
		$\nu_j=c_j'^*U_{N_j}$, an $N_j$-sheeted covering with $V_{\nu_j}=c_j'^*\C^{N_j}_\perm$ (C4). By Theorem~\ref{thm:BPQ}(ii), $\tc(\nu_j)=[V_{\nu_j}]-N_j$ is represented by
		$(\Omega^\infty\eta)\circ b\circ\iota_{N_j}\circ c'_j=(\Omega^\infty\eta)\circ b\circ c_j$. By Step~5 this map is freely homotopic to $(\Omega^\infty\eta)\circ\alpha_j\circ a'$, which
		represents $a'^*x_j$ because $\alpha_j$ lifts $x_j$. By Lemma~\ref{lem:K0basic}(iii), free homotopy classes of maps to $BU$ are reduced $K$-theory classes, so
		$\tc(\nu_j)=a'^*x_j$.
		
		\emph{Step 7: $a'$ is a $K$-isomorphism.} By homotopy extension (Theorem~\ref{thm:standardhtpy}(b)) and path-connectedness of $X'$,
		$a'$ is homotopic to a based map; this changes neither acyclicity (A2) nor \eqref{eq:liftprop} nor $a'^*$
		(Lemma~\ref{lem:K0basic}(ii)). So let $a'$ be based. It induces isomorphisms on $H_*(-;\Z)$ by (A1), and therefore so does the
		reduced suspension $\Sigma a'\colon\Sigma Z'\to\Sigma X'$, by the natural suspension isomorphism $\widetilde H_{k+1}(\Sigma Y)\cong\widetilde H_k(Y)$. Reduced suspensions of
		connected CW complexes are simply connected (Theorem~\ref{thm:vankampen}(c)), so $\Sigma a'$ is a homotopy equivalence by Whitehead's theorem
		(Theorem~\ref{thm:standardhtpy}(c)), and $(\Sigma a')^*$ is an isomorphism on $\Kt^*$. By the natural suspension isomorphism (K3), $a'^*\colon\Kt^*(X')\to\Kt^*(Z')$ is an
		isomorphism; it preserves ranks, so it is an isomorphism on $K^0=\Z\oplus\Kt^0$, and $K^1=\Kt^1$.
		
		\emph{The degenerate cases.} If $X=Z=\pt$ and $a=\id$, then $F_a$ is a point and $f$ is constant; the proof is the same. If $X'=\pt$ and $s=0$, then $P=\hofib_a(f(\pt))=F_a$ is
		acyclic with $\pi_1(P)\cong\pi_1(Z)$, and $Z'$ is a finite acyclic complex with a map $g\colon Z'\to Z$.
	\end{proof}
	
	\begin{remark}
		\begin{enumerate}[label=(\arabic*)]
			\item The proof only uses that $X$ is simply connected and $X'$ is connected (for $X'$ one then uses $\ker\pi_1(q)$ instead of
			$\pi_1(P)$). We keep the stronger hypotheses, which hold in the application, because they make the computation transparent.
			\item The suspension in Step~7 can also be taken unreduced; then no basepoint is needed.
			\item The nontrivial perfect group $\pi_1(Z')$ is essential: by Remark~\ref{rem:simplyconnected} the coverings $\nu_j$ with
			nonzero class cannot exist on a simply connected space. The permutation actions of the $\nu_j$ are the homomorphisms
			$\pi_1(Z')\to\pi_1(B\Sigma_{N_j})=\Sigma_{N_j}$ induced by the maps $c_j'$ of Step~6, which classify the coverings; no further information on
			$\pi_1(Z')$ is needed.
			\item The output is not unique, and no choice is canonical: the lifts $\alpha_j$, the CW approximation and the complex given by Theorem~\ref{thm:HV} all involve choices. This
			does not matter, because only the properties listed in the proposition are used later.
		\end{enumerate}
	\end{remark}
	
	\subsection{Torsion stages and the realization theorem}\label{ss:realize}
	
	We now organize the construction as an induction. At each stage we keep the algebraic data (the group and its decomposition), the topological model (a wedge of Moore
	spaces), the finite complex constructed so far with its acyclic map to the model, and coverings representing generators. The induction of Theorem~\ref{thm:realize}
	carries the following data from one stage to the next.
	
	\begin{definition}[Torsion stage]\label{def:stage}
		A \emph{torsion stage} is a tuple $\mathsf S=(E,\Theta,d,X,\iota,T,a,\phi,(\nu_j)_{j\in J})$ such that:
		\begin{enumerate}[label=(T\arabic*)]
			\item $E$ is a finite abelian group and $\Theta\colon\bigoplus_{j\in J}\Z/p_j^{r_j}\to E$ an isomorphism ($J$ finite, possibly
			empty; $p_j$ prime, $r_j\ge1$); $e_j=\Theta(1_j)$ has order $o_j=p_j^{r_j}$;
			\item $d$ is even, $d\ge4$, and $q(p_j,r_j)\mid d$ and $d\ge q(p_j,r_j)+3$ for every $j\in J$;
			\item $X=\bigvee_{j\in J}P^d(p_j^{r_j})$, and $\iota\colon E\to\Kt^0(X)$ is the isomorphism $e_j\mapsto x_j$ of
			Corollary~\ref{cor:wedgelift};
			\item $T$ is a connected finite CW complex, based at a $0$-cell, and $a\colon T\to X$ is acyclic with $a^*\colon K^*(X)\to K^*(T)$
			an isomorphism;
			\item $\phi=a^*\circ\iota\colon E\to\Kt^0(T)$, and the $\nu_j$ ($j\in J$) are finite coverings of $T$ with $\tc(\nu_j)=\phi(e_j)$ and
			$N_j=\deg\nu_j\ge1$; we call them the \emph{generator coverings}.
		\end{enumerate}
		The number $d$ is the \emph{dimension} of the stage. We write $\mathsf S_n$, $E_n$, $d_n$, $X_n$, $T_n$, \dots\ for the data of the
		$n$th stage.
	\end{definition}
	
	The conditions have the following roles. (T1): every finite abelian group is a direct sum of cyclic groups of prime power order, but the decomposition is not unique, and
	since maps between wedges of Moore spaces are built summand by summand (Lemma~\ref{lem:wedgemaps}), the decomposition $\Theta$ is part of the data. (T2): evenness gives
	$K^1(X)=0$ and $\Kt^0(X)\cong E$ (Lemma~\ref{lem:mooreK}); $d\ge4$ makes $X$ $2$-connected, as required in Proposition~\ref{prop:lift}; divisibility by the periods gives stable
	cohomotopy lifts (Corollary~\ref{cor:wedgelift}); and $d\ge q+3$ is the Freudenthal range in which Adams' maps \emph{into} $X$ can be desuspended (Lemma~\ref{lem:unstable}).
	(T3): $X$ is the homotopy-theoretic model of $E$, with a fixed identification $\iota$. (T4): $T$ is the finite complex on which coverings live; the acyclic map $a$ identifies its
	$K$-theory with that of the model. (T5): the generator coverings represent the generators $e_j$; together with Corollary~\ref{cor:bound} they give covering representatives of all
	elements with explicitly bounded degrees.
	
	By (T4) and Corollary~\ref{cor:wedgelift}, $\phi$ is an isomorphism and $K^1(T)\cong K^1(X)=0$. The finite groups which are
	realized are arbitrary: every finite abelian group is a direct sum of cyclic groups of prime power order, and the choice of the
	decomposition $\Theta$ is part of the data.
	
	\begin{lemma}[Maps between wedges]\label{lem:wedgemaps}
		Let $(E,\Theta,d)$ and $(E',\Theta',d')$ satisfy \textup{(T1)--(T2)}, with $d'\ge d$, and let $X,\iota$ and $X',\iota'$ be as in
		\textup{(T3)}. Every homomorphism $i\colon E\to E'$ is induced by a based map $f\colon X'\to X$, in the sense that $f^*\circ\iota=
		\iota'\circ i$.
	\end{lemma}
	
	\begin{proof}
		Write $E=\bigoplus_j\Z/p_j^{a_j}$ and $E'=\bigoplus_k\Z/p_k'^{\,b_k}$ via $\Theta,\Theta'$, and let $h_{kj}\colon\Z/p_j^{a_j}\to
		\Z/p_k'^{\,b_k}$ be the components of $\Theta'^{-1}i\Theta$, so that $i(e_j)=\sum_kh_{kj}(1)e'_k$. We construct, for each pair $(k,j)$, a based map $f_{kj}\colon
		P^{d'}(p_k'^{\,b_k})\to P^d(p_j^{a_j})$ with $f_{kj}^*\gamma^{(j)}=h_{kj}(1)\gamma'^{(k)}$, where $\gamma^{(j)}$ and $\gamma'^{(k)}$ are the generators of the summands
		(Lemma~\ref{lem:mooreK}(b)).
		
		\emph{Different primes.} If $p_j\ne p'_k$ then $h_{kj}=0$, since a homomorphism between cyclic groups of coprime orders is zero, and $f_{kj}$ is the constant map.
		
		\emph{Same prime.} Write $p=p_j=p'_k$, $a=a_j$, $b=b_k$, $c=\min\{a,b\}$, $h=h_{kj}$, and choose an integer $s$ representing $h(1)$;
		then $p^b\mid p^as$, since $p^a\cdot1=0$ in $\Z/p^a$. Factor $h=h_2\circ h_1$ with
		\[
		h_1\colon\Z/p^a\to\Z/p^c,\ 1\mapsto1,\qquad h_2\colon\Z/p^c\to\Z/p^b,\ 1\mapsto s .
		\]
		These are well defined: $h_1$ because $c\le a$; $h_2$ trivially if $c=b$, and because $p^b\mid p^as$ if $c=a$. Realize $h_1$ by the
		degree square $\varphi_1=\varphi_{p^{a-c},1}\colon P^d(p^c)\to P^d(p^a)$ and $h_2$ by $\varphi_2=\varphi_{p^cs/p^b,\,s}\colon
		P^{d'}(p^b)\to P^{d'}(p^c)$ (Lemma~\ref{lem:degreesquare}; the bottom degrees $p^{a-c}$ and $p^cs/p^b$ are integers, and the conditions $p^{a-c}\cdot p^c=p^a\cdot1$ and
		$(p^cs/p^b)\cdot p^b=p^c\cdot s$ of the lemma hold). The two squares live in different dimensions; to connect them we use an Adams map. Let $q=q(p,c)$. Since $c\le a$ and
		$c\le b$, the divisibility chain of periods gives $q\mid q(p,a)\mid d$ and $q\mid q(p,b)\mid d'$, so $q\mid d'-d$; and $d\ge q(p,a)+3\ge q+3$.
		Lemma~\ref{lem:unstable}(ii) gives $w=w_{d,(d'-d)/q}\colon P^{d'}(p^c)\to P^d(p^c)$ with $w^*\gamma_d=\gamma_{d'}$ ($w=\id$ if
		$d'=d$). Put
		\[
		f_{kj}=\varphi_1\circ w\circ\varphi_2\colon P^{d'}(p^b)\xrightarrow{\varphi_2}P^{d'}(p^c)\xrightarrow{w}P^d(p^c)\xrightarrow{\varphi_1}P^d(p^a).
		\]
		Then $f_{kj}^*\gamma_d=\varphi_2^*w^*\varphi_1^*\gamma_d=\varphi_2^*w^*\gamma_d=\varphi_2^*\gamma_{d'}=s\gamma_{d'}$, so $f_{kj}$
		induces $h$ on $\Kt^0$.
		
		\emph{Assembly.} Since $d'\ge4$, every summand $P^{d'}(p_k'^{\,b_k})=\Sigma P^{d'-1}(p_k'^{\,b_k})$ is a suspension, and we may form pinch sums of maps out of it. Let
		$f_k\colon P^{d'}(p_k'^{\,b_k})\to X$ be the pinch sum $\sum_j\iota_j\circ f_{kj}$ ($f_k$ is constant if $J=\emptyset$), and $f=\bigvee_kf_k\colon X'\to X$. By Lemmas~\ref{lem:pinch}
		and~\ref{lem:K0basic}(iv), and since $\iota_l^*x_j$ is $\gamma^{(j)}$ for $l=j$ and $0$ otherwise,
		\[
		\iota_k'^*f^*x_j=f_k^*x_j=\sum_lf_{kl}^*\iota_l^*x_j=f_{kj}^*\gamma^{(j)}=h_{kj}(1)\gamma'^{(k)}.
		\]
		By Lemma~\ref{lem:K0basic}(iv) a class on $X'$ is determined by its restrictions to the summands, so $f^*x_j=\sum_kh_{kj}(1)x'_k$, that is, $f^*\iota(e_j)=\iota'(i(e_j))$ for all
		$j$. All maps are based.
	\end{proof}
	
	\begin{figure}[ht]
		\centering
		\begin{tikzpicture}[>={Stealth[length=2mm]},node distance=1.2cm]
			\node (A) at (0,0) {$P^{d'}(p^b)$};
			\node (B) at (3.6,0) {$P^{d'}(p^c)$};
			\node (C) at (7.2,0) {$P^{d}(p^c)$};
			\node (D) at (10.6,0) {$P^{d}(p^a)$};
			\draw[->] (A) -- node[above,font=\small]{$\varphi_2$} node[below,font=\scriptsize]{$1\mapsto s$} (B);
			\draw[->] (B) -- node[above,font=\small]{$w$ (Adams)} node[below,font=\scriptsize]{$\gamma_d\mapsto\gamma_{d'}$} (C);
			\draw[->] (C) -- node[above,font=\small]{$\varphi_1$} node[below,font=\scriptsize]{$1\mapsto1$} (D);
			\node[font=\small] at (5.3,-1.2) {on $\Kt^0$: \quad $\Z/p^a\xrightarrow{\ 1\mapsto1\ }\Z/p^c\xrightarrow{\ \cong\ }\Z/p^c\xrightarrow{\ 1\mapsto s\ }\Z/p^b$};
		\end{tikzpicture}
		\caption{The component $f_{kj}$ of Lemma~\ref{lem:wedgemaps} for $c=\min\{a,b\}$. The Adams map changes the dimension; the degree
			squares change the order.}\label{fig:fkj}
	\end{figure}
	
	\begin{theorem}[Realization of finite torsion systems]\label{thm:realize}
		\begin{enumerate}[label=\textup{(\alph*)}]
			\item \emph{(Base step.)} For every finite abelian group $E_1$ with a decomposition $\Theta_1$, and every $d_1$ satisfying
			\textup{(T2)}, there is a torsion stage $\mathsf S_1$ with data $(E_1,\Theta_1,d_1)$.
			\item \emph{(Inductive step.)} Let $\mathsf S_n$ be a torsion stage, $E_{n+1}$ a finite abelian group with a decomposition $\Theta_{n+1}$,
			$i_n\colon E_n\to E_{n+1}$ a homomorphism, and $d_{n+1}\ge d_n$ an integer satisfying \textup{(T2)} for $(E_{n+1},\Theta_{n+1})$.
			Then there are a torsion stage $\mathsf S_{n+1}$ with data $(E_{n+1},\Theta_{n+1},d_{n+1})$ and a based map $g_n\colon T_{n+1}\to
			T_n$ with $g_n^*\circ\phi_n=\phi_{n+1}\circ i_n$.
			\item \emph{(Dimensions.)} Such a $d_{n+1}$ exists: if $E_{n+1}=0$ take $d_{n+1}=d_n$; otherwise take the least multiple of
			$L=\lcm\{q(p,r):\Z/p^r\text{ a summand of }E_{n+1}\}$ which is at least $d_n$ and at least $q(p,r)+3$ for every summand.
			\item In every torsion stage $K^1(T)=0$, and every element of $E$ is $\phi^{-1}\tc(\nu)$ for a finite covering $\nu$ of $T$
			(Corollary~\ref{cor:bound}).
		\end{enumerate}
		Consequently, for every sequence $(E_n,i_n)_{n\ge1}$ of finite abelian groups and homomorphisms $i_n\colon E_n\to E_{n+1}$ there are connected finite CW complexes
		$T_n$, based maps $g_n\colon T_{n+1}\to T_n$ and isomorphisms $\phi_n\colon E_n\to\Kt^0(T_n)$ with $K^1(T_n)=0$ and
		$g_n^*\phi_n=\phi_{n+1}i_n$, such that every element of $E_n$ is the reduced class of a finite covering of $T_n$. The stage
		$\mathsf S_{n+1}$ and the map $g_n$ are constructed from $\mathsf S_n$, $(E_{n+1},\Theta_{n+1})$ and $i_n$ alone.
	\end{theorem}
	
	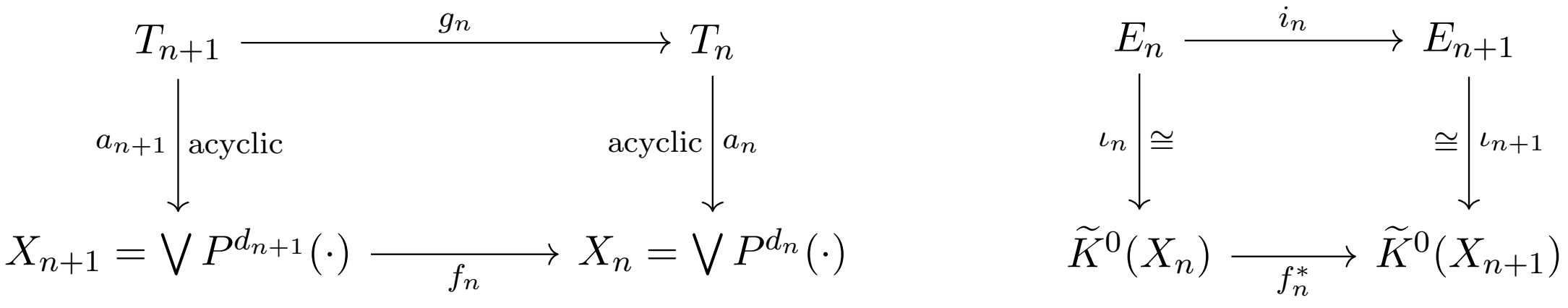
\begin{figure}[ht]
		\centering
		\begin{tikzcd}[column sep=4.5em,row sep=3.2em]
			T_{n+1}\arrow[r,"g_n"]\arrow[d,"a_{n+1}"',"\text{acyclic}"] & T_n\arrow[d,"a_n","\text{acyclic}"']\\
			X_{n+1}=\bigvee P^{d_{n+1}}(\cdot)\arrow[r,"f_n"'] & X_n=\bigvee P^{d_n}(\cdot)
		\end{tikzcd}
		\qquad\qquad
		\begin{tikzcd}[column sep=3em,row sep=3.2em]
			E_n\arrow[r,"i_n"]\arrow[d,"\iota_n"',"\cong"] & E_{n+1}\arrow[d,"\iota_{n+1}","\cong"']\\
			\Kt^0(X_n)\arrow[r,"f_n^*"'] & \Kt^0(X_{n+1})
		\end{tikzcd}
		\caption{One step of Theorem~\ref{thm:realize}. The left square commutes up to free homotopy; the right square commutes by
			Lemma~\ref{lem:wedgemaps}.}\label{fig:step}
	\end{figure}
	
	The theorem is the inductive form of Theorem~\ref{thm:realize-short}. The essential point is (b): a new stage is obtained from the old one by solving the lifting problem
	of Proposition~\ref{prop:lift}, in which the old complex $T_n$ plays the role of $Z$ and the new wedge $X_{n+1}$ the role of $X'$. Part (c) shows that the hypotheses of (b)
	can always be met, and that the choice of $d_{n+1}$ only depends on $d_n$ and $E_{n+1}$.
	
	\begin{proof}
		(c) Each period $q(p,r)$ is even, so $L$ is even, and every multiple of $L$ is divisible by every period $q(p,r)$ of a summand of $E_{n+1}$. The least multiple of $L$ which is at
		least $\max\{d_n,q(p,r)+3\}$ (maximum over the summands) exists, is even, is at least $d_n\ge4$, and satisfies all conditions of (T2) for $(E_{n+1},\Theta_{n+1})$. If
		$E_{n+1}=0$, then (T2) is the condition that $d_{n+1}$ be even and at least $4$, which holds for $d_{n+1}=d_n$.
		
		(b) \emph{Step 1: the map of models.} Lemma~\ref{lem:wedgemaps}, applied to $(E_n,\Theta_n,d_n)$, $(E_{n+1},\Theta_{n+1},d_{n+1})$ and $i_n$, gives a
		based map $f_n\colon X_{n+1}\to X_n$ with $f_n^*\iota_n=\iota_{n+1}i_n$. Here $d_{n+1}\ge d_n$ is used.
		
		\emph{Step 2: the lifting problem.} Apply Proposition~\ref{prop:lift} with $X=X_n$, $Z=T_n$, $a=a_n$, $X'=X_{n+1}$, $f=f_n$ and the classes $x'_k$ ($k\in
		J_{n+1}$) of Corollary~\ref{cor:wedgelift}. The hypotheses hold: $X_n$ and $X_{n+1}$ are $2$-connected finite CW complexes by
		Lemma~\ref{lem:mooreK}(a) ($d_n,d_{n+1}\ge4$; points if the index set is empty); $T_n$ is a connected finite CW complex and $a_n$ is
		acyclic by (T4) for $\mathsf S_n$; the $x'_k$ have stable cohomotopy lifts by Corollary~\ref{cor:wedgelift}, because $d_{n+1}$ is even and divisible by the periods of
		$E_{n+1}$. The output is a connected finite CW complex $T_{n+1}=Z'$, an acyclic map $a_{n+1}=a'$ with $a_{n+1}^*$ an isomorphism on $K^*$, a map $g'\colon T_{n+1}\to T_n$ with
		$a_ng'\simeq f_na_{n+1}$, and coverings $\nu'_k$ of degree at least $1$ with $\tc(\nu'_k)=a_{n+1}^*x'_k$.
		
		\emph{Step 3: basing.} Choose a $0$-cell of $T_{n+1}$ as basepoint, and arrange that $a_{n+1}$ is based for it (Step~7 of the proof of Proposition~\ref{prop:lift}). As $T_n$ is
		path-connected, homotopy extension (Theorem~\ref{thm:standardhtpy}(b)) gives a homotopy from $g'$ to a based map $g_n$.
		
		\emph{Step 4: the identity $g_n^*\phi_n=\phi_{n+1}i_n$.} Put $\phi_{n+1}=a_{n+1}^*\iota_{n+1}$. Since $a_ng_n\simeq a_ng'\simeq f_na_{n+1}$ (free homotopies),
		Lemma~\ref{lem:K0basic}(ii) gives
		\[
		g_n^*\phi_n=g_n^*a_n^*\iota_n=(a_ng_n)^*\iota_n=(f_na_{n+1})^*\iota_n=a_{n+1}^*f_n^*\iota_n=a_{n+1}^*\iota_{n+1}i_n=\phi_{n+1}i_n .
		\]
		
		\emph{Step 5: the new stage.} Let $\mathsf S_{n+1}=(E_{n+1},\Theta_{n+1},d_{n+1},X_{n+1},\iota_{n+1},T_{n+1},a_{n+1},\phi_{n+1},(\nu'_k)_{k\in J_{n+1}})$. Conditions (T1) and (T2) hold by
		hypothesis, and (T3) by definition of $X_{n+1}$ and $\iota_{n+1}$. (T4) holds by Steps~2 and~3. For (T5), $\phi_{n+1}=a_{n+1}^*\iota_{n+1}$ by definition, and
		$\tc(\nu'_k)=a_{n+1}^*x'_k=a_{n+1}^*\iota_{n+1}(e'_k)=\phi_{n+1}(e'_k)$ with $\deg\nu'_k\ge1$. Only $\mathsf S_n$, $(E_{n+1},\Theta_{n+1})$, $d_{n+1}$ and $i_n$ were used, and
		$d_{n+1}$ was chosen by (c) from $d_n$ and $E_{n+1}$.
		
		(a) The \emph{point stage} $\mathsf S_0=(0,\emptyset,4,\pt,0,\pt,\id,0,())$ is a torsion stage: (T1)--(T3) and (T5) are vacuous for $J=\emptyset$, and $\id_\pt$ is acyclic and a
		$K$-isomorphism. Apply (b) to $\mathsf S_0$, $(E_1,\Theta_1)$, $i_0=0$ and $d_1\ge4=d_0$; this is Proposition~\ref{prop:lift} with $X=Z=\pt$ and $f_0$ constant. Discard $g_0$.
		
		(d) $K^1(T)\cong K^1(X)=0$ by (T4) and Corollary~\ref{cor:wedgelift} ($d$ even); representability is Corollary~\ref{cor:bound}.
		
		\emph{The final statement.} Let $(E_n,i_n)_{n\ge1}$ be given, and choose decompositions $\Theta_n$ of the groups $E_n$ (they exist by the structure theorem for finite abelian
		groups). We construct torsion stages $\mathsf S_n$ with data $(E_n,\Theta_n,d_n)$ and based maps $g_n\colon T_{n+1}\to T_n$ by induction on $n$. For $n=1$, choose $d_1$ by (c)
		applied to $d_0=4$ and $E_1$, and let $\mathsf S_1$ be the stage given by (a). If $\mathsf S_1,\dots,\mathsf S_n$ and $g_1,\dots,g_{n-1}$ have been constructed, choose
		$d_{n+1}$ by (c) applied to $d_n$ and $E_{n+1}$; then $d_{n+1}\ge d_n$ and $d_{n+1}$ satisfies (T2) for $(E_{n+1},\Theta_{n+1})$, so (b) provides $\mathsf S_{n+1}$ and
		$g_n$ with $g_n^*\phi_n=\phi_{n+1}i_n$. The complexes $T_n$ are connected and finite by (T4), the maps $g_n$ are based by (b), and $\phi_n\colon E_n\to\Kt^0(T_n)$ is an isomorphism
		because $\iota_n$ is an isomorphism (Corollary~\ref{cor:wedgelift}) and $a_n^*$ is an isomorphism (T4). By (d), $K^1(T_n)=0$ and every element of $E_n$ is $\phi_n^{-1}\tc(\nu)$ for
		a finite covering $\nu$ of $T_n$, i.e.\ every element of $\Kt^0(T_n)=\phi_n(E_n)$ is a reduced covering class. Finally, the construction of $\mathsf S_{n+1}$ and $g_n$ in (b) and
		(c) uses only $\mathsf S_n$, $(E_{n+1},\Theta_{n+1})$ and $i_n$; in particular it does not depend on the groups $E_m$ and maps $i_m$ for $m>n+1$, which need not be known when
		stage $n+1$ is built.
	\end{proof}
	
	\begin{corollary}[Covering representatives and their degrees]\label{cor:bound}
		In a torsion stage, write $x\in E$ uniquely as $x=\sum_jc_je_j$ with $0\le c_j<o_j$. For $x\ne0$, the disjoint union $\nu_x$ of $c_j$
		copies of $\nu_j$ ($j\in J$) has $\tc(\nu_x)=\phi(x)$ by (C1), and degree $\sum_jc_jN_j$; for $x=0$ take $\nu_0=\id_T$, the trivial
		covering of degree one. Hence every $x\in E$ has a covering representative of degree at most
		\begin{equation}\label{eq:ccbound}
			\cc(\mathsf S)=\max\Big\{1,\ \sum_j(o_j-1)N_j\Big\},
		\end{equation}
		a number determined by $T$ and its generator coverings, hence known as soon as $T$ has been constructed. For every map $u\colon A\to
		T$, $u^*\nu_x$ has the same degree (C2). The bound is attained by $x=\sum_j(o_j-1)e_j$ within this family.
	\end{corollary}
	
	\begin{proof}
		The representation $x=\sum_jc_je_j$ with $0\le c_j<o_j$ exists and is unique by (T1). By (C1), degrees and reduced classes add under disjoint union, so
		$\deg\nu_x=\sum_jc_jN_j$ and $\tc(\nu_x)=\sum_jc_j\tc(\nu_j)=\sum_jc_j\phi(e_j)=\phi(x)$ by (T5) and additivity of $\phi$. If $x\ne0$, some $c_j\ge1$, and $N_j\ge1$ by (T5), so
		$1\le\deg\nu_x\le\sum_j(o_j-1)N_j$. The trivial covering $\id_T$ has degree $1$ and class $0=\phi(0)$. Pullback preserves degrees by (C2).
	\end{proof}
	
	The bound \eqref{eq:ccbound} is used in Sections~\ref{sec:blocks} and~\ref{ss:choices}: the coverings of $T_n$ which enter the homogeneous connecting map at stage $n$ have degree
	at most $\cc(\mathsf S_n)$, and the multiplicities of the connecting maps are chosen at least $\cc(\mathsf S_n)$. This is possible because the bound only involves the stage
	$\mathsf S_n$, which is known when these multiplicities are chosen, and because pullback preserves degrees.
	
	\begin{remark}[The reason for $d_{n+1}\ge d_n$]\label{rem:monotone}
		The hypothesis $d_{n+1}\ge d_n$ in Theorem~\ref{thm:realize}(b) cannot be dropped. If $d_{n+1}<d_n$, then $d_{n+1}\le d_n-2$, both numbers being even; since $X_n$ has no
		cells in dimensions $1,\dots,d_n-2$, every based map $X_{n+1}\to X_n$ is null-homotopic by cellular approximation (Theorem~\ref{thm:standardhtpy}(a)), so only $i_n=0$ can
		be realized by a map of models. The hypothesis couples consecutive stages only through this inequality: (T2) involves one stage at a time, and the conditions needed
		for the map $f_n$ in Lemma~\ref{lem:wedgemaps} follow from (T2) at the stages $n$ and $n+1$ and the divisibility chain of periods. In the degenerate case $E=0$, the
		model $X$ is a point, $T$ is a finite acyclic complex and $\cc(\mathsf S)=1$; if $i_n=0$, the map $f_n$ may be taken to be constant.
	\end{remark}
	
	\begin{example}[The groups $\Z/3$ and $\Z/9$]\label{ex:smallcase}
		Let $E_1=\Z/3$ and $E_2=\Z/9$, with $i_1(1)=3$ and the obvious decompositions. The periods are $q(3,1)=4$ and $q(3,2)=12$.
		
		\emph{Dimensions.} For $\Z/3$ alone, dimension $4$ would suffice: $q(3,1)=4$ divides $4$, so the generator of $\Kt^0(P^4(3))=\Z/3$ has a stable cohomotopy lift
		(Lemma~\ref{lem:lift}), and Proposition~\ref{prop:lift} applies to $X'=P^4(3)$. The condition $d\ge q+3$ of (T2) is needed only for maps between stages: realizing $i_1$
		requires an Adams map which changes the dimension, and it desuspends into $P^{d_1}(3)$ only if $d_1\ge q(3,1)+3=7$ (Lemma~\ref{lem:unstable}). So (T2) and
		Theorem~\ref{thm:realize}(c) give $d_1=8$ and $d_2=24$ (a multiple of $12$, at least $15$ and at least $d_1$).
		
		\emph{Stage $1$.} $X_1=P^8(3)$, and Proposition~\ref{prop:lift} with $X=Z=\pt$, $X'=P^8(3)$ and the generator $x_1=\gamma_8$ gives a finite complex $T_1$ with an acyclic map
		$a_1\colon T_1\to P^8(3)$ and a generator covering $\nu_1$ of some degree $N_1$ with $\tc(\nu_1)=a_1^*\gamma_8$; here $\cc(\mathsf S_1)=2N_1$, realized by $\nu_1\sqcup\nu_1$. Since
		$P^8(3)$ is $2$-connected, the homotopy pullback $P$ in the proof has $\pi_1(P)\cong\pi_1(F_b)$, the infinitely generated central extension of $A_\infty$ by $\Z/2$
		(Lemma~\ref{lem:Fb}), and $T_1$ has the integral homology of $P^8(3)$ but a finitely presented perfect fundamental group which carries $\nu_1$. No Moore space could carry
		such a covering: on $P^2(3)$ no covering represents a nonzero class (Example~\ref{ex:P2p}), and $P^d(3)$ is simply connected for $d\ge3$ (Remark~\ref{rem:simplyconnected}).
		
		\emph{The map of models.} $f_1\colon P^{24}(9)\to P^8(3)$ is constructed as in Lemma~\ref{lem:wedgemaps} with $a=1$, $b=2$, $c=1$, $s=3$: the degree square
		$\varphi_2=\varphi_{1,3}\colon P^{24}(9)\to P^{24}(3)$ (bottom degree $3\cdot3/9=1$, top degree $3$), followed by the iterated Adams map $w_{8,4}\colon P^{24}(3)\to P^8(3)$, the composite of
		four desuspended Adams maps $P^{24}(3)\to P^{20}(3)\to P^{16}(3)\to P^{12}(3)\to P^8(3)$ of period $4$ (each target dimension is at least $7=q(3,1)+3$), followed by $\varphi_1=\id$.
		It induces $\Z/3\to\Z/9$, $1\mapsto3$.
		
		\emph{Stage $2$.} Proposition~\ref{prop:lift} with $X=P^8(3)$, $Z=T_1$, $a=a_1$, $X'=P^{24}(9)$, $f=f_1$ and the generator $x'_1=\gamma_{24}$ gives $T_2$, an acyclic map
		$a_2\colon T_2\to P^{24}(9)$, a based map $g_1\colon T_2\to T_1$ and a generator covering $\nu'_1$ with $\tc(\nu'_1)=a_2^*\gamma_{24}$. By Theorem~\ref{thm:realize},
		$g_1^*\phi_1(1)=\phi_2(3)$: the pulled back covering $g_1^*\nu_1$ of $T_2$ has the same degree $N_1$ as $\nu_1$, and its reduced class is three times the generator.
	\end{example}

	\subsection{The spaces for \texorpdfstring{$K_1$}{K1}}\label{ss:K1}
	
	The group $K_1$ will be carried by Moore spaces of odd dimension. They are kept out of Proposition~\ref{prop:lift}, which they do
	not need (no coverings are required for $K_1$, because $K_1$ is not part of the class of the unit), and to which they would not apply ($\pi_2(P^3(m))=\Z/m$, so they are not
	$2$-connected). Since all spaces for $K_1$ have dimension $3$, no Adams maps are needed either: maps between them are built from degree maps and degree squares alone.
	
	\begin{lemma}\label{lem:K1}
		Let $K$ be a finitely generated abelian group with a decomposition $K=\Z^a\oplus\bigoplus_l\Z/m_l$ ($m_l\ge2$), and
		$Y(K)=(\bigvee^aS^3)\vee\bigvee_lP^3(m_l)$ (a point if $K=0$). Then $Y(K)$ is a simply connected finite CW complex with
		$\Kt^0(Y(K))=0$ and $K^1(Y(K))=K$, and every homomorphism $\chi\colon K\to K'$ is induced by a single based map $Y(K')\to Y(K)$.
	\end{lemma}
	
	\begin{proof}
		\emph{The $K$-theory.} By (K3), $\Kt^0(S^3)=0$ and $K^1(S^3)=\Z$, generated by $b_3$; by Lemma~\ref{lem:mooreK}(b) with $d=3$, $\Kt^0(P^3(m))=0$ and $K^1(P^3(m))=\Z/m$, generated
		by $c^*b_3$. By Lemma~\ref{lem:K0basic}(iv), $\Kt^0(Y(K))=0$ and $K^1(Y(K))=\Z^a\oplus\bigoplus_l\Z/m_l=K$, where the $j$th standard generator corresponds to the generator of the
		$j$th summand pulled back along the collapse onto that summand. We use this identification.
		
		\emph{Simple connectivity.} $Y(K)$ is a finite CW complex whose only $0$-cell is the wedge point and which has no $1$-cells, so it is simply connected by cellular approximation
		(Theorem~\ref{thm:standardhtpy}(a)).
		
		\emph{The components.} Write $K=\bigoplus_jC_j$ and $K'=\bigoplus_kC'_k$ for the given decompositions into cyclic groups, $Y_j$ and $Y'_k$ for the corresponding summands of $Y(K)$
		and $Y(K')$, and $\chi_{kj}\colon C_j\to C'_k$ for the components of $\chi$. Each $\chi_{kj}$ is realized by a based map $f_{kj}\colon Y'_k\to Y_j$ with $f_{kj}^*=\chi_{kj}$ on $K^1$:
		\begin{itemize}
			\item $\Z\to\Z$, $1\mapsto n$: a based map $S^3\to S^3$ of degree $n$, which induces multiplication by $n$ on $K^1(S^3)$ (Lemma~\ref{lem:pinch});
			\item $\Z\to\Z/m'$, $1\mapsto n$: the composite $P^3(m')\xrightarrow{c}S^3\xrightarrow{n}S^3$, which maps $b_3$ to $n\,c^*b_3$;
			\item $\Z/m\to\Z$: the only homomorphism is $0$, realized by the constant map;
			\item $\Z/m\to\Z/m'$, $1\mapsto s$ with $m'\mid ms$: the degree square $\varphi_{ms/m',s}\colon P^3(m')\to P^3(m)$
			(Lemma~\ref{lem:degreesquare} with $d=3$, acting on $K^1$).
		\end{itemize}
		
		\emph{Assembly.} As $S^3=\Sigma S^2$ and $P^3(m')=\Sigma P^2(m')$, each summand $Y'_k$ is a suspension. Let $f_k\colon Y'_k\to Y(K)$ be the pinch sum $\sum_j\iota_j\circ f_{kj}$, and
		$f=\bigvee_kf_k\colon Y(K')\to Y(K)$. By Lemmas~\ref{lem:pinch} and~\ref{lem:K0basic}(iv), exactly as in the proof of Lemma~\ref{lem:wedgemaps}, the restriction of $f^*y_j$ to $Y'_k$ is
		$f_{kj}^*$ of the generator of $Y_j$, i.e.\ $\chi_{kj}$ of the $j$th generator. So $f^*=\chi$ under the identifications above. The map $f$ is based.
	\end{proof}
	
	These are the Moore spaces used for $K_1$ in Li's construction \cite[\S3]{Li}. Since $Y(K)$ is simply connected, every finite covering of it is trivial
	(Remark~\ref{rem:simplyconnected}); this is used in Lemma~\ref{lem:Zn}(b).
	
	\subsection{The spaces \texorpdfstring{$Z_n$}{Zn} and coverings over them}\label{ss:Zn}
	
	We now assemble the spaces used in Section~\ref{sec:blocks}. The space $Z_n$ has three wedge summands with three different tasks: $T_n$ carries the torsion of $K_0$ and the
	coverings, $Y_n$ carries $K_1$, and an interval carries a surjection $Z_{n+1}\to Z_n$ which is invisible in $K$-theory and is used to make the connecting maps injective
	(Section~\ref{ss:homog}). Let $\mathsf S_n,\mathsf S_{n+1}$ be torsion stages and $g_n\colon
	T_{n+1}\to T_n$ a based map with $g_n^*\phi_n=\phi_{n+1}i_n$ (Theorem~\ref{thm:realize}); let $Y_n=Y(K_n)$ and $k_n\colon Y_{n+1}\to Y_n$ a
	based map inducing $\chi_n$ (Lemma~\ref{lem:K1}). Put
	\begin{equation}\label{eq:Zn}
		Z_n=T_n\vee Y_n\vee[0,1],
	\end{equation}
	the interval attached at $0$, with wedge point $\theta_n$. Let $\ret_n\colon Z_n\to T_n$ be the collapse onto $T_n$ (the identity on $T_n$
	and constant $\theta_n$ elsewhere), and $\rho_n\colon Z_n\to[0,1]$ the collapse onto the interval.
	
	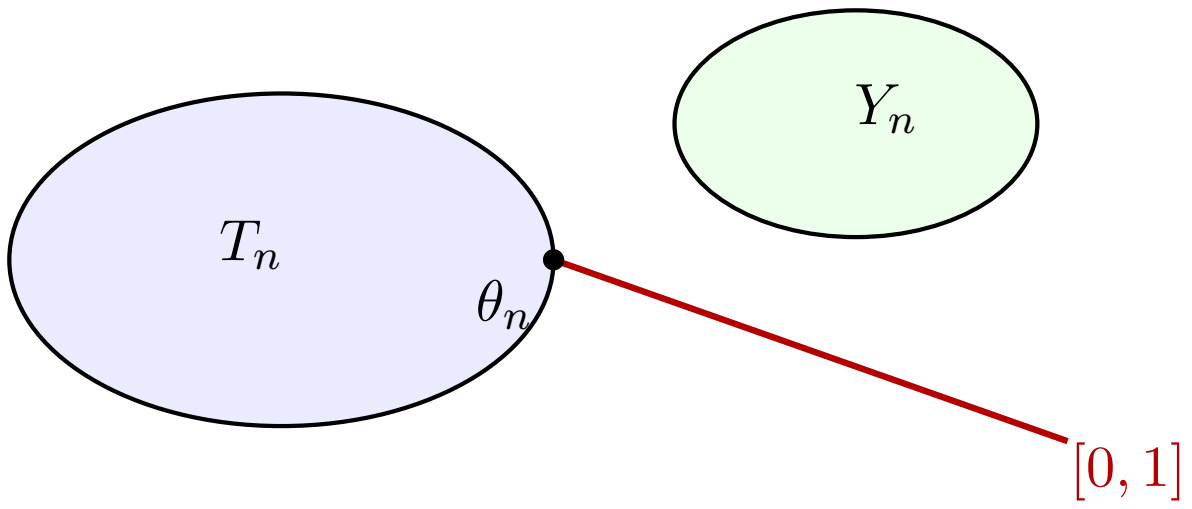
\begin{figure}[ht]
		\centering
		\begin{tikzpicture}[scale=1]
			\draw[thick,fill=blue!8] (-2.2,0) ellipse (1.8 and 1.1); \node at (-2.4,0.1) {$T_n$};
			\draw[thick,fill=green!8] (1.6,0.9) ellipse (1.2 and 0.75); \node at (1.8,1.0) {$Y_n$};
			\draw[very thick,red!70!black] (-0.4,0) -- (3.0,-1.2); \node[red!70!black] at (3.4,-1.4) {$[0,1]$};
			\fill (-0.4,0) circle (2pt); \node[below left] at (-0.4,0) {$\theta_n$};
			\node[align=left,font=\small] at (8.0,0.2) {$\Kt^0(Z_n)\cong\Kt^0(T_n)\cong E_n$ (torsion, by coverings)\\[2pt]
				$K^1(Z_n)\cong K^1(Y_n)\cong K_n$\\[2pt]
				the interval carries the surjection $h_n$};
		\end{tikzpicture}
		\caption{The space $Z_n$.}\label{fig:Zn}
	\end{figure}
	
	\begin{lemma}\label{lem:Zn}
		\begin{enumerate}[label=\textup{(\alph*)}]
			\item $Z_n$ is a connected finite CW complex, hence compact, metrizable and locally connected. Restriction
			gives isomorphisms $\Kt^0(Z_n)\to\Kt^0(T_n)\cong E_n$ (via $\phi_n^{-1}$), with inverse $\ret_n^*$, and $K^1(Z_n)\to K^1(Y_n)=K_n$; the
			latter uses $K^1(T_n)=0$.
			\item For a finite covering $\nu$ of $T_n$, $\ret_n^*\nu$ is a finite covering of $Z_n$ which is $\nu$ over $T_n$ and trivial with fibre
			$\nu^{-1}(\theta_n)$ over $Y_n\vee[0,1]$; it has the same degree, its fibre over $\theta_n$ is $\nu^{-1}(\theta_n)$, and
			$\tc(\ret_n^*\nu)=\ret_n^*\tc(\nu)$. Every finite covering $\mu$ of $Z_n$ is isomorphic to $\ret_n^*(\mu|_{T_n})$ by an isomorphism which is
			the identity over $T_n$. So the extension of coverings across $Y_n\vee[0,1]$ by trivial sheets is well defined and unique up to
			canonical isomorphism; it needs no labelling.
			\item $\kappa_n=g_n\vee k_n\vee\mathrm{const}_{\theta_n}\colon Z_{n+1}\to Z_n$ is well defined, continuous and based (because $g_n$ and
			$k_n$ are based); under (a), $\kappa_n^*=i_n$ on $\Kt^0$ and $\kappa_n^*=\chi_n$ on $K^1$; and $\ret_n\circ\kappa_n=g_n\circ \ret_{n+1}$.
			\item For a finite covering $\nu$ of $T_n$, $\kappa_n^*(\ret_n^*\nu)=\ret_{n+1}^*(g_n^*\nu)$; the covering $g_n^*\nu$ has degree $\deg\nu$ and
			fibre over $\theta_{n+1}$ canonically equal to $\nu^{-1}(\theta_n)$; and $\tc(\ret_{n+1}^*g_n^*\nu)$ corresponds to
			$i_n(\phi_n^{-1}\tc(\nu))$. With Corollary~\ref{cor:bound}: for every $e\in E_n$ there is $\nu$ on $T_n$ with
			$\deg\nu\le\cc(\mathsf S_n)$ and $\tc(\nu)=\phi_n(e)$, and then $\ret_{n+1}^*g_n^*\nu$ has degree at most $\cc(\mathsf S_n)$ and class
			$i_n(e)$.
			\item If $h=p\circ\rho_{n+1}$ for a path $p\colon[0,1]\to Z_n$ with $p(0)=\theta_n$, then $h$ is based, $h^*=0$ on $\Kt^*$, and for
			every finite covering $\mu$ of $Z_n$ the covering $h^*\mu$ is trivial, with fibre over $\theta_{n+1}$ canonically $\mu^{-1}(\theta_n)$.
			\item $Z_n$ is path-connected, and fibres of a finite covering of $Z_n$ can be transported along paths (unique path lifting).
		\end{enumerate}
	\end{lemma}
	
	\begin{proof}
		(a) The wedge is formed at $0$-cells: the basepoints of $T_n$ and $Y_n$ and the endpoint $0$ of $[0,1]$ (with its CW structure of two $0$-cells and one $1$-cell). A wedge of finitely
		many finite CW complexes at $0$-cells is a finite CW complex, and it is connected because the summands are; so $Z_n$ is compact, metrizable and locally connected
		(Section~\ref{ss:CW}). By Lemma~\ref{lem:K0basic}(iv), restriction to the summands is an isomorphism
		\[
		\Kt^*(Z_n)\to\Kt^*(T_n)\oplus\Kt^*(Y_n)\oplus\Kt^*([0,1]),
		\]
		with inverse given by pulling back along the collapses and adding. Here $\Kt^*([0,1])=0$ by (K1), since $[0,1]$ is contractible; $\Kt^0(Y_n)=0$ by Lemma~\ref{lem:K1}; and $K^1(T_n)=0$
		by Theorem~\ref{thm:realize}(d). So on $\Kt^0$ restriction to $T_n$ is an isomorphism whose inverse is pullback along the collapse $\ret_n=c_{T_n}$, and on $K^1$ restriction to $Y_n$ is
		an isomorphism.
		
		(b) Let $\nu\colon W\to T_n$. Then $\ret_n^*\nu=\{(z,w)\in Z_n\times W:\ret_n(z)=\nu(w)\}$. Over $T_n$, where $\ret_n$ is the identity, it is $\nu$; over $Y_n\vee[0,1]$, where $\ret_n$
		is constant with value $\theta_n$, it is $(Y_n\vee[0,1])\times\nu^{-1}(\theta_n)$, the trivial covering. Degree, fibre over $\theta_n$ and reduced class are given by (C2). For the
		uniqueness statement, let $\mu\colon W'\to Z_n$ be a finite covering and $F=\mu^{-1}(\theta_n)$. The space $Y_n\vee[0,1]$ is simply connected and locally path-connected ($Y_n$ is
		simply connected by Lemma~\ref{lem:K1} and $[0,1]$ is contractible), so $\mu$ is trivial over it (Remark~\ref{rem:simplyconnected}): there is a homeomorphism
		$\tau\colon\mu^{-1}(Y_n\vee[0,1])\to(Y_n\vee[0,1])\times F$ over $Y_n\vee[0,1]$ which is the identity on $F$, and it is unique, because a section of a covering over a connected
		space is determined by its value at one point. Define $\Phi\colon W'\to\ret_n^*(\mu|_{T_n})$ as the identity over $T_n$ and as $\tau$ over $Y_n\vee[0,1]$. The two closed sets
		$\mu^{-1}(T_n)$ and $\mu^{-1}(Y_n\vee[0,1])$ cover $W'$ and meet in $F$, where both definitions are the identity; so $\Phi$ is well defined and continuous by the gluing lemma, and
		the same argument applies to $\Phi^{-1}$. Thus $\Phi$ is an isomorphism of coverings which is the identity over $T_n$.
		
		(c) A map out of a wedge which is continuous on each summand and agrees at the wedge point is continuous, and
		$g_n(\theta_{n+1})=\theta_n=k_n(\theta_{n+1})$ because $g_n$ and $k_n$ are based; on the interval, $\kappa_n$ is constant with value $\theta_n$. For $x\in\Kt^0(Z_n)$, the restriction of
		$\kappa_n^*x$ to $T_{n+1}$ is $g_n^*(x|_{T_n})$, because $\kappa_n|_{T_{n+1}}=g_n$; so under the identifications of (a), $\kappa_n^*=\phi_{n+1}^{-1}g_n^*\phi_n=i_n$. Similarly, the
		restriction of $\kappa_n^*y$ to $Y_{n+1}$ is $k_n^*(y|_{Y_n})$, which gives $\chi_n$ on $K^1$. Finally, on $T_{n+1}$ both $\ret_n\kappa_n$ and $g_n\ret_{n+1}$ equal $g_n$, and on
		$Y_{n+1}\vee[0,1]$ both are constant with value $\theta_n$: the first because $\kappa_n$ maps this set to $Y_n\vee\{\theta_n\}$, which $\ret_n$ collapses to $\theta_n$, the second because
		$\ret_{n+1}$ collapses it to $\theta_{n+1}$ and $g_n(\theta_{n+1})=\theta_n$ (here basedness of $g_n$ is used).
		
		(d) By (c) and functoriality of pullback, $\kappa_n^*(\ret_n^*\nu)\cong(\ret_n\kappa_n)^*\nu=(g_n\ret_{n+1})^*\nu\cong\ret_{n+1}^*(g_n^*\nu)$. By (C2), $g_n^*\nu$ has degree $\deg\nu$, and its fibre
		over $\theta_{n+1}$ is $\{\theta_{n+1}\}\times\nu^{-1}(g_n(\theta_{n+1}))=\nu^{-1}(\theta_n)$. Its reduced class is $g_n^*\tc(\nu)$, so by (a) the class of $\ret_{n+1}^*g_n^*\nu$
		corresponds to $\phi_{n+1}^{-1}g_n^*\tc(\nu)=\phi_{n+1}^{-1}g_n^*\phi_n\big(\phi_n^{-1}\tc(\nu)\big)=i_n(\phi_n^{-1}\tc(\nu))$. The last sentence follows from Corollary~\ref{cor:bound}.
		
		(e) $h(\theta_{n+1})=p(\rho_{n+1}(\theta_{n+1}))=p(0)=\theta_n$. On $\Kt^*$, $h^*=\rho_{n+1}^*p^*$ factors through $\Kt^*([0,1])=0$. For a covering $\mu$ of $Z_n$,
		$h^*\mu\cong\rho_{n+1}^*(p^*\mu)$, and $p^*\mu$ is a covering of the simply connected space $[0,1]$, hence trivial, with fibre $\mu^{-1}(p(0))=\mu^{-1}(\theta_n)$ over $0$; a pullback of a
		trivial covering is trivial, and the fibre over $\theta_{n+1}$ is the fibre of $p^*\mu$ over $\rho_{n+1}(\theta_{n+1})=0$, which is $\mu^{-1}(\theta_n)$ by (C2).
		
		(f) Connected CW complexes are path-connected (Section~\ref{ss:CW}); unique path lifting is \cite[Proposition~1.30]{Hat}.
	\end{proof}
	
	\begin{remark}[Labels]\label{rem:labels}
		A \emph{labelling} of a covering $\mu$ of $Z_n$ of degree $r$ is a bijection $[r]\to\mu^{-1}(\theta_n)$. By (C2) and
		Lemma~\ref{lem:Zn}(d),(e), the fibre over $\theta_{n+1}$ of each pulled back covering used in Section~\ref{ss:homog} is canonically a fibre
		over $\theta_n$ (of the distinguished covering $\lambda_n$, of a generator covering, or of a trivial covering), so labels at $\theta_n$
		transport to labels at $\theta_{n+1}$. This is why all eigenvalue maps in Section~\ref{ss:homog} must be \emph{based}: if one of them
		sent $\theta_{n+1}$ to a point $z\ne\theta_n$, the value of the connecting map at $\theta_{n+1}$ would involve values away from
		$\theta_n$, and the evaluation at the basepoint would not be compatible with the connecting maps.
	\end{remark}
	%======================================================================
	\section{An algebraic presentation with delayed torsion corrections}\label{sec:alg}
	
	In this section we construct, from the invariant to be realized, an inductive system of finitely generated abelian groups which
	presents the $K_0$-group, and in which every torsion correction occurs one step later than in the usual presentation. This delay is what
	makes the inductions of Sections~\ref{ss:choices} and~\ref{ss:SPorder} possible. The construction is the same in the unital and in the
	stably projectionless case, except for the choice of an auxiliary order unit $u$; the consequences of this choice are discussed in
	Section~\ref{ss:twocases}. Everything in this section is elementary algebra, apart from the Effros--Handelman--Shen theorem and Edwards'
	separation theorem, which we state.
	
	\subsection{Ordered groups and the invariant}\label{ss:ordered}
	
	\begin{definition}[Ordered groups {\cite{Goo}}, {\cite{Eff}}]
		A \emph{partially ordered abelian group} is a pair $(G,G^+)$ of an abelian group and a subset $G^+$ (the \emph{positive cone}) with
		$G^++G^+\subset G^+$, $G^+\cap(-G^+)=\{0\}$ and $G=G^+-G^+$. We write $x\le y$ if $y-x\in G^+$, and $x<y$ if moreover $x\neq y$. An
		\emph{order unit} is $u\in G^+$ such that for every $x\in G$ there is $N$ with $x\le Nu$. A \emph{state} of $(G,G^+,u)$ is a
		homomorphism $s\colon G\to\R$ with $s(G^+)\subset[0,\infty)$ and $s(u)=1$; the states form a compact convex set $S(G)$ in the
		topology of pointwise convergence. $(G,G^+)$ is
		\begin{itemize}
			\item \emph{unperforated} if $nx\in G^+$ for some $n\ge1$ implies $x\in G^+$;
			\item an \emph{interpolation group} (has the \emph{Riesz interpolation property}) if whenever $x_1,x_2\le y_1,y_2$ there is $z$ with
			$x_1,x_2\le z\le y_1,y_2$;
			\item \emph{simple} if every nonzero positive element is an order unit;
			\item a \emph{dimension group} if it is isomorphic, as an ordered group, to an inductive limit
			$\lims(\Z^{k_n},\Z_{\ge0}^{k_n};A_n)$ of simplicially ordered groups with connecting maps given by nonnegative integer matrices
			$A_n$; the positive cone of the limit is the union of the images of the cones $\Z_{\ge0}^{k_n}$.
		\end{itemize}
	\end{definition}
	
	\begin{theorem}[{\cite{EHS}}]\label{thm:EHS}
		A countable partially ordered abelian group is a dimension group if and only if it is unperforated and has the Riesz
		interpolation property.
	\end{theorem}
	
	\begin{definition}[Choquet simplices and affine functions]
		A metrizable compact convex subset $T$ of a locally convex space is a \emph{Choquet simplex} if every point is the barycentre of a
		unique probability measure supported on the extreme boundary. $\Aff(T)$ is the real Banach space of continuous affine functions
		$T\to\R$ with the supremum norm, ordered by the strict order: $f\gg0$ means $f(\tau)>0$ for all $\tau$.
	\end{definition}
	
	\begin{theorem}[Edwards' separation theorem {\cite{Edw}}, {\cite[\S II.3]{Alf}}]\label{thm:edwards}
		Let $T$ be a metrizable Choquet simplex, $f$ a continuous convex function and $g$ a continuous concave function on $T$ with $f\le g$.
		Then there is $a\in\Aff(T)$ with $f\le a\le g$. In particular, $\Aff(T)$ has the Riesz interpolation property: if $f_1,f_2\le
		g_1,g_2$ in $\Aff(T)$, apply this to $\max\{f_1,f_2\}$ and $\min\{g_1,g_2\}$.
	\end{theorem}
	
	Recall also that the states of the order-unit space $\Aff(T)$, i.e.\ the positive linear functionals with value $1$ at the constant
	function $1$, are exactly the point evaluations at points of $T$ \cite[Chapter~II]{Alf}, \cite[\S7]{Goo}.
	
	\paragraph{The data to be realized.} Throughout Sections~\ref{sec:alg}--\ref{sec:blocks} we fix
	\[
	(G_0,\,T,\,\rho,\,G_1)\quad\text{and an element }u ,
	\]
	where $G_0,G_1$ are countable abelian groups, $T$ is a nonempty metrizable Choquet simplex, $\rho\colon G_0\to\Aff(T)$ is a homomorphism,
	and $u$ is specified in one of the following two cases.
	\begin{itemize}
		\item \emph{Case (U), unital.} $(G_0,G_0^+,u,T,\rho,G_1)$ is the Elliott invariant of a unital classifiable stably finite $C^*$-algebra $A$:
		$G_j=K_j(A)$, $u=[1_A]\in G_0$, $T=T(A)$ and $\rho(x)(\tau)=\tau_*(x)$. In particular $\rho(u)=1$.
		\item \emph{Case (SP), stably projectionless.} $\rho$ satisfies
		\begin{equation}\label{eq:P}
			\rho(G_0)\cap\Aff_{++}(T)=\emptyset,\qquad\Aff_{++}(T)=\{f\in\Aff(T):f(\tau)>0\text{ for all }\tau\in T\},\tag{P}
		\end{equation}
		and $u$ is the auxiliary element $(0,1)$ of the group $H$ defined in \eqref{eq:H} below; it does not lie in $G_0$.
	\end{itemize}
	In Case (U), $T(A)$ is the tracial state space with the weak$^*$ topology. It is a metrizable Choquet simplex (a standard fact; see e.g.
	\cite{Sak}), and it is nonempty: $A$ is exact, simple, unital and stably finite, so it has a quasitrace \cite{BH}, and quasitraces on
	exact $C^*$-algebras are traces \cite{Haa}. The map $\tau_*$ is the homomorphism $K_0(A)\to\R$ induced by $\tau$ (on projections in
	$M_k(A)$, $\tau_*[p]=(\tau\otimes\Tr_k)(p)$, with the unnormalized trace $\Tr_k$ on $M_k$; see Section~\ref{ss:conventions}). The groups $G_0,G_1$ are countable
	because $A$ is separable. Case (SP) arises from classifiable stably projectionless $C^*$-algebras with continuous scale (Section~\ref{sec:SPprelim}).
	The letter $T$ always denotes the trace simplex; the complexes of Section~\ref{sec:top} carry an index, $T_n$.
	
	\subsection{The dense affine enlargement}\label{ss:enlarge}
	
	Choose a countable subgroup $G'\subset\Aff(T)$ which is dense in the norm topology and contains the constant functions $1$ and $\sqrt2$
	(two $\Q$-linearly independent constants). Such a subgroup exists because $\Aff(T)$ is separable ($T$ is metrizable). Put
	\begin{equation}\label{eq:H}
		H=G_0\oplus G',\qquad H^+=\{0\}\cup\{(g,f):\ \rho(g)(\tau)+f(\tau)>0\text{ for every }\tau\in T\},
	\end{equation}
	and define $\theta\colon H\to\Aff(T)$ by $\theta(g,f)=\rho(g)+f$. In Case (U) we regard $u$ as $(u,0)\in H$; in Case (SP) we put $u=(0,1)$,
	which lies in $H$ because $1\in G'$. In both cases
	\begin{equation}\label{eq:thetau}
		\theta(u)=1,
	\end{equation}
	the constant function; this is all that is used in Sections~\ref{ss:enlarge}--\ref{ss:delayed}. Let $R=\Tor(H)=\Tor(G_0)\oplus0$ be the
	torsion subgroup, $F=H/R$, with the image $F^+$ of $H^+$ as positive cone and the image $\bar u$ of $u$. Since $\rho$ is a
	homomorphism to a torsion-free group, $\theta(R)=0$, and $\theta$ induces a map $\theta\colon F\to\Aff(T)$. The condition
	$\theta(h)\gg0$ defining $H^+\setminus\{0\}$ depends only on the coset $h+R$, so $F^+=\{0\}\cup\{x\in F:\theta(x)\gg0\}$.
	
	This is the \emph{dense affine enlargement} of the unital range construction \cite[\S2.1]{Li}, \cite[\S13]{GLN}. Its purpose: the
	pairing $\rho$ may be far from injective and $\rho(G_0)$ may be far from dense; adding $G'$ makes the order on $F$ ``as rich as
	$\Aff(T)$'', so that $F$ is a simple dimension group with state space exactly $T$, and the AF algebra with this $K_0$ has trace
	simplex $T$. The summand $G'$ will later be removed again by the interval algebras, which cut $K_0$ down to $G_0$
	(Section~\ref{ss:unit}).
	
	\begin{theorem}\label{thm:Fdimgroup}
		$(F,F^+,\bar u)$ is a countable, simple, noncyclic dimension group with order unit, of rank at least $2$. Its state space is
		canonically $T$: the map $\tau\mapsto s_\tau$, $s_\tau(x)=\theta(x)(\tau)$, is an affine homeomorphism $T\to S(F)$.
	\end{theorem}
	
	\begin{proof}
		\emph{Ordered group with order unit.} $F^++F^+\subset F^+$ is clear. If $x,-x\in F^+$ with $x\ne0$, then $\theta(x)\gg0$ and
		$\theta(-x)\gg0$, which is impossible; so $F^+\cap(-F^+)=0$. Let $x\in F$. Since $T$ is compact, $\theta(x)$ is bounded, so
		$\theta(N\bar u-x)=N-\theta(x)\gg0$ for large $N$, i.e.\ $x\le N\bar u$; in particular $x=N\bar u-(N\bar u-x)\in F^+-F^+$, and $\bar u$ is
		an order unit (note $\theta(\bar u)=1$ by \eqref{eq:thetau}).
		
		\emph{Countable, rank, noncyclic.} $G_0$ and $G'$ are countable. The composite $G'\to H\to F$ is injective (as $G'\cap R=0$), so $F$
		contains $G'$, which contains $1$ and $\sqrt2$, two elements independent over $\Q$; hence $F$ has rank at least $2$ and is not cyclic.
		
		\emph{Unperforated.} If $nx\in F^+$, $n\ge1$, $x\ne0$, then $nx\ne0$ ($F$ is torsion-free), $\theta(x)=\theta(nx)/n\gg0$, so $x\in F^+$.
		
		\emph{Riesz interpolation.} Let $x_1,x_2\le y_1,y_2$. If $x_i=y_j$ for some $i,j$, then $z=x_i$ works. Otherwise $\theta(y_j-x_i)\gg0$
		for all $i,j$, and by compactness of $T$ there is $\varepsilon>0$ with $\theta(x_i)+3\varepsilon\le\theta(y_j)$ for all $i,j$.
		Theorem~\ref{thm:edwards}, applied to $f_i=\theta(x_i)+\varepsilon$ and $g_j=\theta(y_j)-\varepsilon$, gives $a\in\Aff(T)$ with
		$\theta(x_i)+\varepsilon\le a\le\theta(y_j)-\varepsilon$. By density choose $f\in G'$ with $\|f-a\|<\varepsilon$, and let $z$ be the image of
		$(0,f)$. Then $\theta(z-x_i)\gg0$ and $\theta(y_j-z)\gg0$, so $x_i<z<y_j$. Here the hypothesis that $T$ is a Choquet simplex is used;
		for a general compact convex set, $\Aff(T)$ need not have interpolation (e.g.\ for a square).
		
		\emph{Dimension group.} By Theorem~\ref{thm:EHS}.
		
		\emph{Simple.} If $x\in F^+\setminus\{0\}$ then $\theta(x)\ge\delta$ for some $\delta>0$ by compactness, and for $y\in F$,
		$\theta(Nx-y)\gg0$ for large $N$; so $y\le Nx$.
		
		\emph{State space.} Each $s_\tau$ is a state. Conversely let $s$ be a state. If $\theta(x)=\theta(y)$, then $\theta(\bar u\pm n(x-y))=1\gg0$,
		so $\bar u\pm n(x-y)\ge0$ and $|s(x-y)|\le1/n$ for all $n$; thus $s=s'\circ\theta$ for a function $s'$ on $\theta(F)$. If $p,m$ are
		integers with $m\ge1$ and $p/m>\|\theta(x)\|$, then $\theta(p\bar u\mp mx)\gg0$, so $m|s(x)|\le p$; hence $|s'(\theta(x))|\le\|\theta(x)\|$.
		So $s'$ is additive and contractive on the subgroup $\theta(F)$, which contains the dense subgroup $G'$; it extends uniquely to a
		continuous additive, hence $\R$-linear, functional on $\Aff(T)$, which is positive and takes the value $1$ at the function $1$. Such a
		functional is evaluation at a unique $\tau\in T$, and then $s=s_\tau$. The map $\tau\mapsto s_\tau$ is affine, injective (as $G'$ is
		dense and $\Aff(T)$ separates points of $T$) and continuous for the topology of pointwise convergence; as $T$ is compact and $S(F)$
		Hausdorff, it is a homeomorphism.
	\end{proof}
	
	\subsection{A positive presentation with large multiplicities}\label{ss:positive}
	
	\begin{definition}
		A \emph{positive presentation} of a dimension group with order unit $(F,F^+,\bar u)$ is a sequence of finite sets $I_n$, nonnegative
		integer matrices $C_n\colon\Z^{I_n}\to\Z^{I_{n+1}}$ and vectors $v_n\in\Z_{\ge0}^{I_n}$ with $v_{n+1}=C_nv_n$, together with an
		isomorphism of ordered groups with order unit $(F,F^+,\bar u)\cong\lims(\Z^{I_n},\Z_{\ge0}^{I_n},v_n;C_n)$. The matrix entry
		$(C_n)_{j,i}$ is the coefficient of the basis vector $h_j$ in $C_nh_i$. We write $\mu_n\colon\Z^{I_n}\to F$ for the canonical maps and
		$C_{n,m}=C_{m-1}\cdots C_n$.
	\end{definition}
	
	\begin{theorem}\label{thm:presentation}
		There is a positive presentation of $(F,F^+,\bar u)$ with the following properties.
		\begin{enumerate}[label=\textup{(\alph*)}]
			\item Every $v_n$ is strictly positive, every $\mu_n(h_i)$ is nonzero, and $|I_n|\ge2$.
			\item For every $n$ there is $m>n$ such that $C_{n,m}$ is strictly positive; and for every $n$ and every bound $b$ there is $m>n$ such that
			every entry of $C_{n,m}$ is at least $b$.
		\end{enumerate}
		Consequently, the constructions of Sections~\ref{ss:choices} and~\ref{ss:SPorder} may use a \emph{telescoped} presentation, chosen successively: having fixed
		stages up to $n$ (that is, an index $m_n$ of the presentation above), and having determined a bound $b_n$ from the data of the first $n$
		stages, one chooses $m_{n+1}>m_n$ with all entries of $C_{m_n,m_{n+1}}$ at least $b_n$, and takes $C_{m_n,m_{n+1}}$ as the next
		connecting matrix.
	\end{theorem}
	This is the form of the presentation used in the unital range construction \cite[13.24]{GLN}, where it is derived from \cite{EHS} using three $\Q$-linearly
	independent elements of the dense subgroup; we give a complete proof.

	\begin{proof}
		\emph{Step 0: some presentation.} By Theorem~\ref{thm:Fdimgroup} and Theorem~\ref{thm:EHS}, $(F,F^+)\cong\lims(\Z^{k_n},\Z_{\ge0}^{k_n};A_n)$
		for nonnegative integer matrices $A_n$. Since $\bar u\in F^+$ is the image of a positive vector at some stage, after deleting the
		first stages we may assume $\bar u=\mu_1(v_1)$ with $v_1\ge0$, and put $v_n=A_{1,n}v_1$. Recall two facts about inductive limits:
		an element $\mu_n(z)$ is $0$ iff $A_{n,m}z=0$ for some $m$; and $\mu_n(z)\in F^+$ iff $A_{n,m}z\ge0$ for some $m$ (because $F^+$ is the
		union of the images of the cones, and two representatives of the same element agree at a late stage).
		
		\emph{Step 1: strictly positive units.} Let $S_n=\{e:\ v_n(e)>0\}$ be the support of $v_n$. If $e\in S_n$, then $h_e\le v_n$, so
		$A_nh_e\le A_nv_n=v_{n+1}$ and $\supp(A_nh_e)\subset S_{n+1}$. So the coordinate subgroups $\Z^{S_n}$ form a subsystem with
		connecting maps the restrictions of the $A_n$, and with strictly positive units $v_n$. Its limit maps into $F$ injectively (the
		restrictions are honest restrictions, so the criterion for zero is the same). It maps onto $F^+$, hence onto $F$: if $y\in F^+$, write
		$y=\mu_n(z)$ with $z\ge0$; since $\bar u$ is an order unit, $N\bar u-y\ge0$ for some $N$, so $A_{n,m}(Nv_n-z)\ge0$ for some $m$, i.e.\
		$0\le A_{n,m}z\le Nv_m$, and $A_{n,m}z\in\Z_{\ge0}^{S_m}$ represents $y$. The same argument shows that the positive cone of the
		subsystem limit is $F^+$. So we may assume all $v_n$ strictly positive.
		
		\emph{Step 2: nonzero images.} Let $N_n=\{e:\ \mu_n(h_e)=0\}$. If $e\in N_n$, then $0=\mu_{n+1}(A_nh_e)=\sum_f(A_n)_{fe}\mu_{n+1}(h_f)$ is
		a sum of positive elements, so each term vanishes; as $F$ is torsion-free, $\mu_{n+1}(h_f)=0$ whenever $(A_n)_{fe}>0$. Thus
		$A_n(\Z^{N_n})\subset\Z^{N_{n+1}}$, and the $A_n$ induce maps of the quotients $\Z^{k_n}/\Z^{N_n}=\Z^{[k_n]\setminus N_n}$, given by the
		submatrices on the coordinates outside $N$. The maps $\mu_n$ factor through these quotients, the induced map of limits is onto and
		injective (an element of the quotient with zero image in $F$ is the image of some $z$ with $A_{n,m}z=0$), and the positive cones
		correspond. The images of the $v_n$ remain strictly positive. So we may assume $\mu_n(h_e)\neq0$ for every coordinate.
		
		\emph{Step 3: at least two coordinates.} If $k_n=1$ for infinitely many $n$, then $F$ would be an increasing union of cyclic groups,
		hence of rank at most $1$, contradicting Theorem~\ref{thm:Fdimgroup}. So $k_n\ge2$ for all large $n$, and we delete the first stages.
		This proves (a).
		
		\emph{Step 4: strict positivity after telescoping.} Fix $n$ and a coordinate $e$. By simplicity, $\mu_n(h_e)\in F^+\setminus\{0\}$ is an
		order unit, so $N\mu_n(h_e)-\bar u\ge0$ for some $N$, i.e.\ $A_{n,m}(Nh_e-v_n)\ge0$ for some $m$, and then
		$NA_{n,m}h_e\ge A_{n,m}v_n=v_m\gg0$, so $A_{n,m}h_e$ is strictly positive; the same holds for all larger $m$, because $A_{m,m'}$ maps strictly
		positive vectors to strictly positive vectors (every row of $A_{m,m'}$ is nonzero, because $A_{m,m'}v_m=v_{m'}$ is strictly
		positive). Taking $m$ large enough for all the finitely many $e$ gives $A_{n,m}$ strictly positive.
		
		\emph{Step 5: large entries.} If $X$ and $Y$ are strictly positive integer matrices and the middle dimension is at least $2$, then every
		entry of $YX$ is a sum of at least two terms, each at least $1$, so it is at least $2$; inductively, a product of $t$ strictly
		positive integer matrices with all intermediate dimensions at least $2$ has all entries at least $2^{t-1}$. By Step~4 we can write
		$A_{n,m}$, for $m$ large, as a product of $t$ strictly positive matrices, with $t$ as large as we like. This proves (b), and the last
		assertion is a reformulation.
	\end{proof}
	
	\begin{example}
		Both conditions in (a) are needed for (b). For $A=\left(\begin{smallmatrix}2&1\\0&0\end{smallmatrix}\right)$ the limit $\lims(\Z^2,A)\cong\Z[\tfrac12]$
		is simple, and no basis vector has zero image ($A^ke_2=2^{k-1}e_1\neq0$); but the second row of every power of $A$ vanishes, so no unit vector at
		a later stage is strictly positive, and no product of these matrices is strictly positive. This illustrates the first condition in (a). For
		$B=\left(\begin{smallmatrix}2&0\\1&0\end{smallmatrix}\right)$ with unit $(1,1)$ the unit vectors $B^k(1,1)=2^{k-1}(2,1)$ are strictly positive, but
		$Be_2=0$, so $e_2$ has zero image, and again no power of $B$ is strictly positive. This illustrates the second condition.
	\end{example}
	
	\subsection{Lifts, defects and the torsion groups}\label{ss:defects}
	
	Fix the telescoped presentation $(\Z^{I_n},C_n)$ of Theorem~\ref{thm:presentation} (its connecting matrices will be chosen successively in
	Sections~\ref{ss:choices} and~\ref{ss:SPorder}). Choose lifts $\ell_n\colon\Z^{I_n}\to H$ of the maps $\mu_n\colon\Z^{I_n}\to F=H/R$, i.e.\ choose for each
	basis vector $h_i$ an element $\ell_n(h_i)\in H$ with image $\mu_n(h_i)$. The \emph{defect}
	\begin{equation}\label{eq:defect}
		t_n=\ell_n-\ell_{n+1}C_n\colon\Z^{I_n}\to H
	\end{equation}
	takes values in $R$, because $\mu_{n+1}C_n=\mu_n$. Choose an enumeration $(\mathfrak r_k)_{k\ge1}$ of the countable torsion group $R$ and define
	finite subgroups $E_n\subset R$ inductively by $E_0=0$ and
	\[
	E_{n}=\big\langle E_{n-1},\ t_{n-1}(\Z^{I_{n-1}}),\ \mathfrak r_1,\dots,\mathfrak r_n\big\rangle .
	\]
	A finitely generated torsion abelian group is finite, so each $E_n$ is finite; the $E_n$ increase and exhaust $R$, and
	$t_{n-1}(\Z^{I_{n-1}})\subset E_n$. Write $i_n\colon E_n\hookrightarrow E_{n+1}$ for the inclusions.
	
\begin{figure}[ht]
		\centering
		\begin{tikzpicture}[>={Stealth[length=2mm]},x=1cm,y=1cm]
			\foreach \x/\l in {0/{stage $n$},4.5/{stage $n+1$},9/{stage $n+2$}}{\draw[thick] (\x,0) -- (\x,-4.4); \node[above] at (\x,0) {\l};}
			\node[draw,rounded corners,fill=blue!6,align=center,font=\small] (a) at (2.25,-0.8) {choose $C_n$ (bound uses $\cc(\mathsf S_n)$)};
			\node[draw,rounded corners,fill=blue!6,align=center,font=\small] (b) at (2.25,-1.9) {lifts $\ell_{n+1}$, defect $t_n$, group $E_{n+1}$};
			\node[draw,rounded corners,fill=green!8,align=center,font=\small] (c) at (4.5,-3.0) {$T_{n+1}$, coverings, $\cc(\mathsf S_{n+1})$};
			\node[draw,rounded corners,fill=blue!6,align=center,font=\small] (d) at (6.75,-4.0) {$C_{n+1}$: $m_{n+1}\ge\cc(\mathsf S_{n+1})$;\\ correction $t_n$ paid out of $x_i$};
			\draw[->] (a) -- (b); \draw[->] (b) -- (c); \draw[->] (c) -- (d);
		\end{tikzpicture}
		\caption{The order of choices. The correction $t_n$ becomes known after $C_n$ is chosen, and the coverings realizing it live on
			$T_{n+1}$, whose degrees are known only after $T_{n+1}$ is built. The delayed presentation pays for $t_n$ in the map from stage $n+1$
			to stage $n+2$, whose multiplicity is chosen after $T_{n+1}$.}\label{fig:timeline}
	\end{figure}
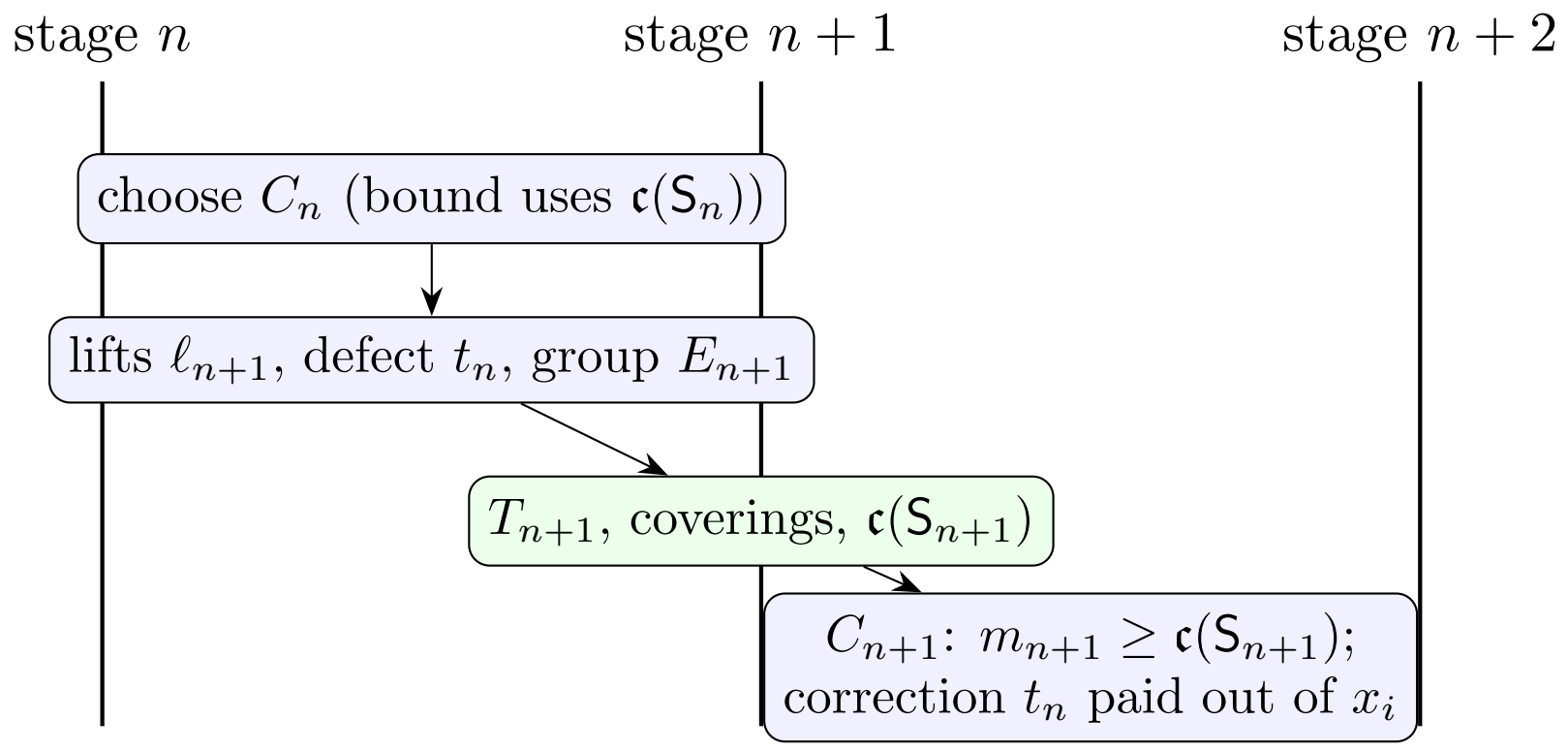

	\paragraph{The failure of the usual presentation.} The group $H$ is the limit of the system
	$(\Z^{I_n}\oplus E_n,(h,e)\mapsto(C_nh,\ i_n(e)+t_n(h)))$, via $(h,e)\mapsto\ell_n(h)+e$. In a $C^*$-algebraic model of this map,
	the torsion correction $t_n(h_i)\in E_{n+1}$ has to be realized by a covering over the target complex $T_{n+1}$, and the degree of
	this covering has to be absorbed by the matrix multiplicity $(C_n)_{\ib_{n+1},i}$ (see Section~\ref{ss:homog}). But $t_n$ depends on
	$\ell_{n+1}$ and on $C_n$ itself, and $T_{n+1}$ and its coverings are constructed only after $E_{n+1}$, hence after $t_n$, is known; the
	covering degrees of $T_{n+1}$ cannot be bounded in advance. So the multiplicity would have to be chosen before the number it must
	dominate is known. We remove this circularity by realizing the correction $t_n$ one step later, out of an extra coordinate, by a
	covering over the \emph{source} complex $T_{n+1}$ of the next map, whose degree is known at that time (Figure~\ref{fig:timeline}).
	
	\subsection{The delayed presentation}\label{ss:delayed}
	
	Designate an index $\ib_n\in I_n$ for each $n$ (the \emph{distinguished} coordinate; it will correspond to the covering algebra). For
	$n\ge1$ introduce one extra coordinate $x_i$ for each $i\in I_{n-1}$, and write
	\[
	\mathcal X_n=\{x_i:i\in I_{n-1}\},\qquad J_n=I_n\sqcup\mathcal X_n .
	\]
	(The letters $x_i$ denote basis vectors of $\Z^{J_n}$; the index $n$ is suppressed.) Let $\pi_n\colon\Z^{J_n}\to\Z^{I_n}$ be the identity on
	the original coordinates and send every extra coordinate to $h_{\ib_n}$. For $l\in J_n$ write $\pi_n(l)\in I_n$ for the corresponding
	index ($\pi_n(l)=l$ for $l\in I_n$, $\pi_n(x_i)=\ib_n$). Define $\delta_n\colon\Z^{J_n}\to E_n$ by
	\[
	\delta_n(h_i)=0\ (i\in I_n),\qquad\delta_n(x_i)=t_{n-1}(h_i)\in E_n\quad(i\in I_{n-1}).
	\]
	The connecting maps are
	\begin{align}
		M_nh_l&=\sum_{j\in I_{n+1}}(C_n)_{j,\pi_n(l)}h_j-h_{\ib_{n+1}}+x_{\pi_n(l)}\qquad(l\in J_n),\label{eq:Mn}\\
		\gt_n(h,e)&=\big(M_nh,\ i_n(e+\delta_n(h))\big)\qquad(h\in\Z^{J_n},\ e\in E_n).\label{eq:gt}
	\end{align}
	Here $x_{\pi_n(l)}$ is the extra coordinate of stage $n+1$ indexed by $\pi_n(l)\in I_n$. In words: the column of $l$ is the column of
	$\pi_n(l)$ in $C_n$, with one unit moved from the distinguished target coordinate to the extra coordinate $x_{\pi_n(l)}$. The entries
	of $C_n$ will always be at least $2$, so $M_n$ is a nonnegative integer matrix. It has sparse extra rows: the row of $x_i\in\mathcal X_{n+1}$ has
	entry $1$ in the columns $l$ with $\pi_n(l)=i$ and $0$ elsewhere. Define also
	\begin{equation}\label{eq:Mprime}
		M_n'\colon\Z^{I_n}\to\Z^{J_{n+1}},\qquad M_n'h_i=\sum_{j\in I_{n+1}}(C_n)_{j,i}h_j-h_{\ib_{n+1}}+x_i .
	\end{equation}
	It is nonnegative (indeed, its original rows are strictly positive, because $(C_n)_{\ib_{n+1},i}\ge2$).
	
	\begin{example}
		Let $I_{n-1}=I_n=I_{n+1}=\{1,2\}$ with $\ib=1$, so $J_n=J_{n+1}=\{1,2,x_1,x_2\}$, and $C_n=\left(\begin{smallmatrix}c_{11}&c_{12}\\
			c_{21}&c_{22}\end{smallmatrix}\right)$. Then $\pi_n(1)=\pi_n(x_1)=\pi_n(x_2)=1$, $\pi_n(2)=2$, and in the ordered basis $(1,2,x_1,x_2)$
		\[
		M_n=\begin{pmatrix}c_{11}-1&c_{12}-1&c_{11}-1&c_{11}-1\\ c_{21}&c_{22}&c_{21}&c_{21}\\ 1&0&1&1\\ 0&1&0&0\end{pmatrix},\qquad
		\delta_n=\big(0,\ 0,\ t_{n-1}(h_1),\ t_{n-1}(h_2)\big).
		\]
		The extra coordinates $x_1,x_2$ of stage $n$ behave like copies of the distinguished coordinate $1$ for the free part, but carry the
		torsion corrections $t_{n-1}(h_1),t_{n-1}(h_2)$; the extra rows of stage $n+1$ receive exactly one unit from each column.
	\end{example}
	
	\begin{lemma}\label{lem:delayed}
		The inductive limit of $(\Z^{J_n}\oplus E_n,\gt_n)$ is $H$, via the maps
		\begin{equation}\label{eq:jn}
			j_n(h,e)=\ell_n(\pi_nh)+e+\delta_n(h).
		\end{equation}
		Its free quotient $\lims(\Z^{J_n},M_n)$, with the coordinatewise positive cones, has ordered limit $(F,F^+)$, via the maps induced by
		the $\pi_n$.
	\end{lemma}
	
	\begin{proof}
		\emph{Intertwining relations.} The cancellations in \eqref{eq:Mn} give
		\begin{equation}\label{eq:intertwine}
			M_n=M_n'\pi_n,\qquad\pi_{n+1}M_n'=C_n,\qquad\pi_{n+1}M_n=C_n\pi_n,\qquad\delta_{n+1}M_n=t_n\pi_n .
		\end{equation}
		Indeed, $M_nh_l$ depends only on $\pi_n(l)$, which gives the first relation; $\pi_{n+1}$ sends $x_i$ to $h_{\ib_{n+1}}$, which cancels the
		term $-h_{\ib_{n+1}}$, giving the second; the third follows from the first two; and $\delta_{n+1}$ vanishes on original coordinates and sends
		$x_{\pi_n(l)}$ to $t_n(h_{\pi_n(l)})$, giving the fourth.
		
		\emph{Compatibility of the $j_n$.} Using \eqref{eq:intertwine} and \eqref{eq:defect},
		\[
		j_{n+1}\gt_n(h,e)=\ell_{n+1}C_n\pi_nh+e+\delta_n(h)+t_n\pi_nh=\ell_n\pi_nh+e+\delta_n(h)=j_n(h,e).
		\]
		So the $j_n$ induce a homomorphism $j\colon\lims(\Z^{J_n}\oplus E_n,\gt_n)\to H$.
		
		\emph{Free part.} By \eqref{eq:intertwine}, the nonnegative maps
		\[
		\Z^{J_n}\xrightarrow{\pi_n}\Z^{I_n}\xrightarrow{M_n'}\Z^{J_{n+1}}\xrightarrow{\pi_{n+1}}\Z^{I_{n+1}}\xrightarrow{M_{n+1}'}\cdots
		\]
		form one inductive system of ordered groups, whose subsystems at the $J$-places and at the $I$-places are $(\Z^{J_n},M_n)$ and
		$(\Z^{I_n},C_n)$. A cofinal subsystem of an inductive system of ordered groups has the same ordered limit, so both have the same ordered
		limit, and the identification is induced by the $\pi_n$. This identifies the ordered free limit with $(F,F^+)$.
		
		\emph{Torsion part.} The short exact sequences $0\to E_n\to\Z^{J_n}\oplus E_n\to\Z^{J_n}\to0$ are compatible with the connecting maps,
		and inductive limits of abelian groups are exact. The map $j$ is the identity on the limiting torsion group $\bigcup_nE_n=R$, and on the
		quotient it is the isomorphism just established, because $j_n(h,e)$ modulo $R$ is $\mu_n(\pi_nh)$. Comparing with $0\to R\to H\to F\to0$,
		the five lemma shows that $j$ is an isomorphism.
	\end{proof}
	
	\subsection{The unit vectors and the quotient equations}\label{ss:unit}
	
	Prepend an algebraic stage $I_0=\{0\}$, $E_0=0$, $\ell_0(1)=u$, with $C_0\colon\Z\to\Z^{I_1}$, $C_0(1)=v_1$. Using
	Theorem~\ref{thm:presentation}, choose $v_1$ strictly positive with distinguished coordinate at least $2$, and include
	$t_0(1)=u-\ell_1(v_1)$ in $E_1$ (it lies in $R$ because $\mu_1(v_1)=\bar u$). The $C^*$-algebra construction starts at stage $1$; then
	$\mathcal X_1=\{x_0\}$, $J_1=I_1\sqcup\{x_0\}$, and the \emph{unit vector} is
	\begin{equation}\label{eq:unit1}
		(r_1,\zeta_1)=(v_1-h_{\ib_1}+x_0,\ 0)\in\Z^{J_1}\oplus E_1 .
	\end{equation}
	Every coordinate of $r_1$ is positive, and $j_1(r_1,0)=\ell_1(v_1-h_{\ib_1}+h_{\ib_1})+\delta_1(x_0)=\ell_1(v_1)+t_0(1)=u$. In particular, the
	torsion component $\zeta_1$ is $0$, so the first covering for the distinguished algebra is trivial. The potentially large covering
	degree needed for $t_0$ is not required at this first stage: it is paid for by the next map, out of the size-one extra summand $x_0$.
	Transport the unit vector by \eqref{eq:gt}:
	\begin{equation}\label{eq:unitn}
		r_{n+1}=M_nr_n,\qquad\zeta_{n+1}=i_n\big(\zeta_n+\delta_n(r_n)\big).
	\end{equation}
	Every coordinate of every $r_n$ is positive: the original rows of $M_n$ are strictly positive (as $(C_n)_{j,i}\ge2$), and each extra row
	of $x_i$ has entry $1$ in the column of the original coordinate $i$, where $r_n$ is positive. By Lemma~\ref{lem:delayed},
	$j_n(r_n,\zeta_n)=u$ for every $n$. The pair $(r_n,\zeta_n)$ will be the $K_0$-class of the unit of the homogeneous algebra $F_n$ of
	Section~\ref{sec:blocks}; in Case (U) it will also be the class of the unit of the building block, in Case (SP) it will not.
	
	\paragraph{The subgroups $G_n$.} Let $G_n=j_n^{-1}(G_0)\subset\Z^{J_n}\oplus E_n$, where $G_0$ is embedded in $H$ as the first summand. It
	contains $E_n$ (as $E_n\subset R\subset G_0$). Since $j_{n+1}\gt_n=j_n$, $\gt_n$ maps $G_n$ into $G_{n+1}$, and by exactness of inductive limits
	$\lims(G_n,\gt_n)=G_0$. The quotient $(\Z^{J_n}\oplus E_n)/G_n$ embeds into $H/G_0\cong G'$, a torsion-free group, and it is finitely generated,
	so it is free abelian, of rank $l_n$ say. Choose an identification with $\Z^{l_n}$ and denote the quotient map, restricted to $\Z^{J_n}$, by
	$\Delta_n\colon\Z^{J_n}\to\Z^{l_n}$. It is onto (the quotient is generated by the image of $\Z^{J_n}$, since $E_n\subset G_n$). Then
	\begin{equation}\label{eq:quot}
		G_n=\ker(\Delta_n\oplus0),\qquad\Delta_{n+1}M_n=L_n\Delta_n
	\end{equation}
	for an integer matrix $L_n\colon\Z^{l_n}\to\Z^{l_{n+1}}$, namely the map induced by $\gt_n$ on the quotients. Neither $\Delta_n$ nor $L_n$ is
	required to be positive. Since the images of the second components of the $\ell_n$ increase and exhaust $G'\neq0$, $l_n$ is nondecreasing and
	eventually positive. Put
	\begin{equation}\label{eq:dn}
		\mathbf d_n:=\Delta_nr_n\in\Z^{l_n},
	\end{equation}
	the coordinate vector of the image of $u$ in $H/G_0\cong G'$, written in the quotient basis at stage $n$.
	
	\subsection{The two cases}\label{ss:twocases}
	
	\paragraph{Case (U).} Here $u\in G_0$, so $(r_n,\zeta_n)\in G_n$ and
	\begin{equation}\label{eq:dnzero}
		\mathbf d_n=\Delta_nr_n=0\qquad\text{for all }n .
	\end{equation}
	We allow $l_n=0$ (this happens at early stages if $\ell_n(\Z^{I_n})$ lies in $G_0\oplus0$ modulo torsion); then there is no interval summand at
	stage $n$, and every statement about interval summands at that stage is vacuous.
	
	\paragraph{Case (SP).} Here $u=(0,1)$ and its image in $H/G_0=G'$ is the function $1\neq0$. Hence
	\begin{equation}\label{eq:dnnonzero}
		\mathbf d_n=\Delta_nr_n\neq0\qquad\text{for all }n ;
	\end{equation}
	in particular $l_n\ge1$ for every $n$. This is the algebraic origin of the non-unitality of the building blocks in this case. The order
	structure has a further consequence, which is the algebraic shadow of stable projectionlessness.
	
	\begin{lemma}\label{lem:nopositive}
		In Case (SP), $G_0\cap H^+=\{0\}$; and if $(h,e)\in G_n$ with $h\ge0$ coordinatewise, then $h=0$.
	\end{lemma}
	
	\begin{proof}
		If $(g,0)\in H^+\setminus\{0\}$, then $\rho(g)=\theta(g,0)\gg0$ on $T$, contradicting \eqref{eq:P}. Now suppose $(h,e)\in G_n$ with $h\ge0$ and
		$h\ne0$. The original rows of $M_n$ are strictly positive, so $M_nh$ is nonzero and nonnegative, and inductively the image of $h$ at every later
		stage is nonzero and nonnegative. Hence the class of $h$ in $F=\lims(\Z^{J_m},M_m)$ is a nonzero element of $F^+$ (Lemma~\ref{lem:delayed}). Since the
		class of $j_n(h,e)$ in $F$ is the class of $\pi_nh$, the element $j_n(h,e)$ is a nonzero element of $H^+$ (the condition defining $H^+$ only depends
		on the class modulo $R$). But $j_n(h,e)\in G_0$, a contradiction.
	\end{proof}
	
	If $p$ is a nonzero projection over a homogeneous algebra whose $K_0$-group is $\Z^{J_n}\oplus E_n$ with ranks in the first summand, its class
	$(h,e)$ has $h\ge0$ and $h\neq0$. Lemma~\ref{lem:nopositive} therefore says that no such class lies in $G_n$ (compare \cite[\S2.2]{Li}).
	
	\begin{remark}
		In Section~\ref{sec:blocks}, $\Z^{J_n}\oplus E_n$ will be $K_0$ of the homogeneous part $F_n$ of the building block (with the class of its unit
		$(r_n,\zeta_n)$), the map $\partial_n\colon\Z^{J_n}\oplus E_n\to\Z^{l_n}$, $\partial_n(h,e)=\Delta_nh$ (that is, $\Delta_n\oplus0$), will be the boundary map of the
		extension defining the building block $A_n$, and $G_n=\ker\partial_n=(\ker\Delta_n)\oplus E_n$ will be $K_0(A_n)$.
		So the interval algebras remove the summand $G'$ which was added in \eqref{eq:H}. The vector $\mathbf d_n$ measures the difference of the ranks of the two
		endpoint maps on the unit; it vanishes in Case (U) and not in Case (SP).
	\end{remark}
	%======================================================================
	\section{The building blocks}\label{sec:blocks}
	
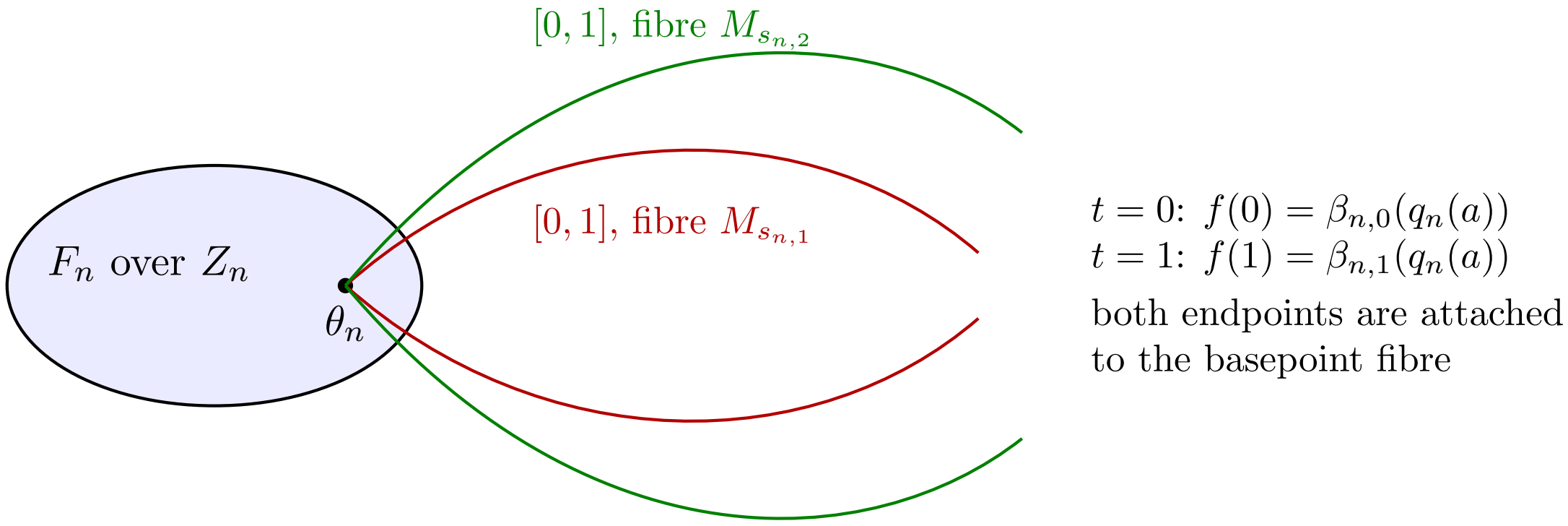
\begin{figure}[ht]
		\centering
		\begin{tikzpicture}[scale=1]
			\draw[thick,fill=blue!8] (0,0) ellipse (1.9 and 1.1);
			\node at (-0.6,0.2) {$F_n$ over $Z_n$};
			\fill (1.2,0) circle (2pt); \node[below] at (1.2,-0.05) {$\theta_n$};
			\draw[thick,red!70!black] (1.2,0) .. controls (3,1.6) and (5.5,1.6) .. (7,0.3);
			\draw[thick,red!70!black] (1.2,0) .. controls (3,-1.6) and (5.5,-1.6) .. (7,-0.3);
			\draw[thick,green!50!black] (1.2,0) .. controls (3.4,2.6) and (6,2.5) .. (7.4,1.4);
			\draw[thick,green!50!black] (1.2,0) .. controls (3.4,-2.6) and (6,-2.5) .. (7.4,-1.4);
			\node[red!70!black,font=\small] at (4.2,0.55) {$[0,1]$, fibre $M_{s_{n,1}}$};
			\node[green!50!black,font=\small] at (4.2,2.35) {$[0,1]$, fibre $M_{s_{n,2}}$};
			\node[font=\small,align=left] at (10.2,0) {$t=0$: $f(0)=\beta_{n,0}(q_n(a))$\\ $t=1$: $f(1)=\beta_{n,1}(q_n(a))$\\[3pt]
				both endpoints are attached\\ to the basepoint fibre};
		\end{tikzpicture}
		\caption{A building block $A_n$: a homogeneous algebra $F_n$ over $Z_n$ (together with matrix summands, not shown) and interval algebras
			whose two endpoint values are determined by the value of $F_n$ at $\theta_n$, through two different block embeddings.}\label{fig:block}
	\end{figure}

	In this section we construct, at every stage $n$, a building block $A_n$ with a $C^*$-diagonal $B_n$, an explicit untwisted groupoid
	model, and $K_0(A_n)\cong G_n$, $K_1(A_n)\cong K_n$. The building block is a one-dimensional analogue of Li's blocks \cite[\S\S2--3]{Li}: a
	homogeneous part $F_n$ over a finite complex, glued to finitely many interval algebras whose endpoints are attached to the fibre of $F_n$ over
	a basepoint (Figure~\ref{fig:block}). The homogeneous part and its connecting maps $\psi_n$ are the same in both cases. The endpoint maps of
	the interval algebras are unital in Case (U) and have \emph{zero blocks} in Case (SP); we treat them uniformly, allowing zero blocks
	throughout. The connecting maps of the full building blocks are constructed separately in the two cases (Sections~\ref{sec:unitalmaps}
	and~\ref{sec:SPmaps}).
	
	\subsection{The homogeneous part}\label{ss:homog}
	
	\subsubsection*{The spaces}
	Choose a presentation $G_1=\lims(K_n,\chi_n)$ of $K_1(A)$ by finitely generated abelian groups $K_n$ with chosen decompositions. Apply
	Theorem~\ref{thm:realize} to the system $(E_n,i_n)$ of Section~\ref{ss:defects} and Lemma~\ref{lem:K1} to the $K_n$. This gives torsion
	stages $\mathsf S_n$ with complexes $T_n$, based maps $g_n\colon T_{n+1}\to T_n$, and $Y_n=Y(K_n)$ with based maps $k_n\colon Y_{n+1}\to Y_n$
	inducing $\chi_n$. Let $Z_n=T_n\vee Y_n\vee[0,1]$ with basepoint $\theta_n$ and $\kappa_n\colon Z_{n+1}\to Z_n$ as in Section~\ref{ss:Zn}. By
	Lemma~\ref{lem:Zn}, $\Kt^0(Z_n)\cong E_n$ and $K^1(Z_n)\cong K_n$, and $\kappa_n^*$ is $i_n$ on $\Kt^0$ and $\chi_n$ on $K^1$. We
	abbreviate $\cc_n=\cc(\mathsf S_n)$ (Corollary~\ref{cor:bound}).
	
	\subsubsection*{A surjection from the interval}
	We need one more map $Z_{n+1}\to Z_n$, which is surjective and invisible in $K$-theory; it will make the connecting maps injective.
	
	\begin{theorem}[Hahn--Mazurkiewicz theorem {\cite[\S31]{Wil}}, {\cite[Chapter~VIII]{Na}}]\label{thm:HM}
		A Hausdorff space is a continuous image of $[0,1]$ if and only if it is compact, connected, locally connected and metrizable.
	\end{theorem}
	
	A connected finite CW complex has all four properties, so there is a continuous surjection $[0,1]\to Z_n$. For CW complexes a surjection
	can be written down directly, which shows what the theorem gives here. For every $k\ge1$ there is a continuous surjection $[0,1]\to
	[0,1]^k\cong D^k$ (a Peano or Hilbert curve). For each closed cell of $Z_n$, compose such a curve with the characteristic map; then
	concatenate these finitely many paths, joined by paths in the connected space $Z_n$ through the basepoint. The result is a continuous
	surjection $p_n\colon[0,1]\to Z_n$, and by prefixing a path from $\theta_n$ we may require $p_n(0)=\theta_n$. Let
	\[
	h_n=p_n\circ\rho_{n+1}\colon Z_{n+1}\to Z_n ,
	\]
	where $\rho_{n+1}\colon Z_{n+1}\to[0,1]$ is the collapse onto the attached interval. So $h_n$ is $p_n$ on the interval and constant with value
	$\theta_n=p_n(0)$ on $T_{n+1}\vee Y_{n+1}$.
	
	\begin{lemma}\label{lem:kappah}
		The maps $\kappa_n,h_n\colon Z_{n+1}\to Z_n$ are based and continuous, $h_n$ is surjective, and
		\begin{equation}\label{eq:kappah}
			\kappa_n^*|_{\Kt^0}=i_n,\qquad\kappa_n^*|_{K^1}=\chi_n,\qquad h_n^*|_{\Kt^*}=0 .
		\end{equation}
		For every finite covering $\mu$ of $Z_n$, $h_n^*\mu$ is trivial, and $\kappa_n^*\mu$ is computed by Lemma~\ref{lem:Zn}(d).
	\end{lemma}
	
	\begin{proof}
		Lemma~\ref{lem:Zn}(c),(e). Surjectivity of $h_n$ is surjectivity of $p_n$.
	\end{proof}
	
	The surjection is used for injectivity of the connecting maps, not for the invariant; no surjectivity of $g_n$ is needed.
	
	\subsubsection*{The homogeneous algebra}
	Recall from Section~\ref{sec:alg} the index set $J_n=I_n\sqcup\mathcal X_n$, the distinguished index $\ib_n\in I_n$, the unit
	$(r_n,\zeta_n)\in\Z^{J_n}\oplus E_n$ with all $r_{n,l}\ge1$, and the matrix $M_n$. Suppose a finite covering
	\[
	\lambda_n\colon\widetilde Y_n\to Z_n,\qquad\deg\lambda_n=r_{n,\ib_n},\qquad\tc(\lambda_n)=\zeta_n,
	\]
	has been chosen (inductively below; for $n=1$ it is the trivial covering, since $\zeta_1=0$), together with a \emph{labelling}, i.e.\ a
	bijection $[r_{n,\ib_n}]\to\lambda_n^{-1}(\theta_n)$. Put
	\begin{equation}\label{eq:Fn}
		F_n=F_{\lambda_n}\oplus\bigoplus_{l\in J_n\setminus\{\ib_n\}}M_{r_{n,l}},\qquad\widehat F_n=\bigoplus_{l\in J_n}M_{r_{n,l}} .
	\end{equation}
	Here $F_{\lambda_n}=\Gamma(\End V_{\lambda_n})\cong C^*(R_{\lambda_n})$ is the covering algebra of Section~\ref{ss:coverings}. Evaluation of the
	distinguished summand at $\theta_n$, written as a matrix in the labelled sheet basis, with the other summands unchanged, is a surjective
	$*$-homomorphism
	\[
	q_n\colon F_n\twoheadrightarrow\widehat F_n,\qquad q_n(a)_{\ib_n}=a_{\ib_n}(\theta_n),\quad q_n(a)_l=a_l\ (l\ne\ib_n).
	\]
	
	\begin{lemma}\label{lem:KFn}
		In \emph{Morita coordinates}, i.e.\ identifying $K_0(M_{r})=\Z$ via the class of a minimal projection and $K_0(F_{\lambda_n})=K^0(Z_n)=
		\Z\oplus\Kt^0(Z_n)=\Z\oplus E_n$ via Lemma~\ref{lem:KFpi} and Lemma~\ref{lem:Zn}(a),
		\[
		K_0(F_n)=\Z^{J_n}\oplus E_n,\qquad[1_{F_n}]=(r_n,\zeta_n),\qquad K_1(F_n)=K^1(Z_n)=K_n,
		\]
		and $(q_n)_*\colon K_0(F_n)\to K_0(\widehat F_n)=\Z^{J_n}$ is the projection onto the first summand.
	\end{lemma}
	
	\begin{proof}
		$K_*$ is additive. By Lemma~\ref{lem:KFpi}, $[1_{F_{\lambda_n}}]=[V_{\lambda_n}]=\deg\lambda_n+\tc(\lambda_n)=r_{n,\ib_n}+\zeta_n$, and
		$[1_{M_r}]=r$. Evaluation at a point maps the class of a vector bundle to its rank, so $(q_n)_*$ kills reduced classes.
	\end{proof}
	
	The matrix summands carry their standard bases, and the covering algebra its sheet bases. Every element $a\in F_n$ therefore has
	canonical \emph{fibre matrices}: $a_{\ib_n}(z)$ in the sheet basis of $\lambda_n^{-1}(z)$, for $z\in Z_n$, and $a_l$ for $l\ne\ib_n$. We call
	$a$ \emph{diagonal}, \emph{monomial}, or \emph{nonnegative monomial} if all its fibre matrices are (Definition~\ref{def:monomial}).
	
	\subsubsection*{The homogeneous connecting map $\psi_n$}
	We now define $\psi_n\colon F_n\to F_{n+1}$; the covering $\lambda_{n+1}$ will be part of the output. Write
	\[
	m_n=(C_n)_{\ib_{n+1},\ib_n}-1 .
	\]
	Recall that $(M_n)_{\ib_{n+1},l}=(C_n)_{\ib_{n+1},\pi_n(l)}-1$, so $(M_n)_{\ib_{n+1},l}=m_n$ for $l=\ib_n$ and for every extra coordinate $l=x_i$,
	and $(M_n)_{\ib_{n+1},l}=(C_n)_{\ib_{n+1},l}-1\ge1$ for the other original coordinates.
	
	\emph{Into the matrix targets.} For a target index $l'\in J_{n+1}\setminus\{\ib_{n+1}\}$, put
	\[
	\psi_n(a)_{l'}=\diag\Big(\underbrace{q_n(a)_{l_1},\dots,q_n(a)_{l_1}}_{(M_n)_{l',l_1}},\ \underbrace{q_n(a)_{l_2},\dots}_{(M_n)_{l',l_2}},\ \dots\Big)
	\in M_{r_{n+1,l'}},
	\]
	the block-diagonal matrix containing $(M_n)_{l',l}$ copies of the basepoint value $q_n(a)_l$, for $l\in J_n$ in a fixed order (the
	\emph{standard block order}). The sizes match because $r_{n+1}=M_nr_n$.
	
	\emph{Into the distinguished target.} Here $\psi_n(a)_{\ib_{n+1}}$ is the direct sum, over the source indices $l\in J_n$, of basic
	homomorphisms \eqref{eq:basic}, with coverings $\nu_l$ of $Z_{n+1}$ and maps $g_l$, as follows.
	\begin{enumerate}[label=(\roman*)]
		\item \emph{From an original matrix source} $l\in I_n\setminus\{\ib_n\}$: the trivial covering of $Z_{n+1}$ of degree
		$(C_n)_{\ib_{n+1},l}-1$, with the matrix-source version of \eqref{eq:basic}. The image is
		$a_l\otimes1\in M_{r_{n,l}}\otimes C(Z_{n+1},M_{(C_n)_{\ib_{n+1},l}-1})$: constant copies of $a_l$.
		\item \emph{From an extra matrix source} $x_i$, $i\in I_{n-1}$: at the source stage $n$ choose a covering $\nu_{n,i}$ of $T_n$ with
		\[
		\tc(\nu_{n,i})=\phi_n(\delta_n(x_i)),\qquad\deg\nu_{n,i}\le\cc_n,
		\]
		by Corollary~\ref{cor:bound}; if $\delta_n(x_i)=0$, take $\nu_{n,i}$ to be the trivial covering of degree one. Pull it back along $g_n$
		and extend it over $Y_{n+1}\vee[0,1]$ by trivial sheets, i.e.\ form $\ret_{n+1}^*g_n^*\nu_{n,i}$ (Lemma~\ref{lem:Zn}(b),(d)); then add trivial
		sheets until the degree is $m_n$. Call the result $\nu'_{n,i}$; it is a covering of $Z_{n+1}$ of degree $m_n$ with
		$\tc(\nu'_{n,i})=i_n(\delta_n(x_i))$ (identifying $\Kt^0(Z_{n+1})$ with $E_{n+1}$). Use the matrix-source version of \eqref{eq:basic}: the
		image of $a_{x_i}$ is $a_{x_i}\otimes1\in M_{r_{n,x_i}}\otimes F_{\nu'_{n,i}}$, a covering algebra of degree $r_{n,x_i}m_n$, whose fibre at $z$
		consists of copies of the constant matrix $a_{x_i}$, permuted by the monodromy of $\nu'_{n,i}$.
		\item \emph{From the distinguished source}: the trivial $m_n$-sheeted covering $Z_{n+1}\times[m_n]$, with the maps $\kappa_n$ on the first
		sheet, $h_n$ on the second sheet, and the constant map $\theta_n$ on the remaining $m_n-2$ sheets. By the first special case after
		Definition~\ref{def:basic}, the image is
		\[
		z\mapsto\diag\big(a_{\ib_n}(\kappa_n(z)),\ a_{\ib_n}(h_n(z)),\ a_{\ib_n}(\theta_n),\dots,a_{\ib_n}(\theta_n)\big),
		\]
		written in the sheet bases of the pullback coverings $\kappa_n^*\lambda_n$, $h_n^*\lambda_n$, $\theta_n^*\lambda_n$.
	\end{enumerate}
	The conditions $m_n\ge\cc_n$ (needed in (ii)) and $m_n\ge3$ (so that both $\kappa_n$ and $h_n$ occur) are imposed when $C_n$ is chosen
	(Sections~\ref{ss:choices} and~\ref{ss:SPorder}); they only involve data of the completed stage $n$.
	
	\begin{definition}[The covering $\lambda_{n+1}$]\label{def:lambda}
		$\lambda_{n+1}$ is the disjoint union of the image coverings of the entries (i)--(iii):
		\[
		\lambda_{n+1}=\kappa_n^*\lambda_n\ \sqcup\ h_n^*\lambda_n\ \sqcup\ (m_n-2)\cdot\theta_n^*\lambda_n\ \sqcup\ \bigsqcup_{l\in I_n\setminus\{\ib_n\}}r_{n,l}\cdot
		\big(Z_{n+1}\times[(C_n)_{\ib_{n+1},l}-1]\big)\ \sqcup\ \bigsqcup_{i\in I_{n-1}}r_{n,x_i}\cdot\nu'_{n,i},
		\]
		where $k\cdot\nu$ denotes $k$ disjoint copies. Its fibre over $\theta_{n+1}$ is the disjoint union of the corresponding fibres, each of
		which is canonically a fibre over $\theta_n$ (of $\lambda_n$, of a trivial covering, or of a generator covering; Remark~\ref{rem:labels}).
		We label it by concatenating the labels of these fibres in a fixed order of the entries (the \emph{block order} of
		$\widehat\psi_n$ below).
	\end{definition}
	
	\begin{lemma}\label{lem:lambdaclass}
		$\deg\lambda_{n+1}=r_{n+1,\ib_{n+1}}$ and $\tc(\lambda_{n+1})=\zeta_{n+1}$. Thus $\psi_n$ is a unital $*$-homomorphism $F_n\to F_{n+1}$.
	\end{lemma}
	
	\begin{proof}
		The degree is $\sum_l(M_n)_{\ib_{n+1},l}\,r_{n,l}=(M_nr_n)_{\ib_{n+1}}=r_{n+1,\ib_{n+1}}$. By (C1), (C2) and \eqref{eq:kappah}, the reduced class is
		\[
		\kappa_n^*\zeta_n+h_n^*\zeta_n+0+0+\sum_ir_{n,x_i}\,i_n\delta_n(x_i)=i_n(\zeta_n)+i_n\Big(\sum_ir_{n,x_i}\delta_n(x_i)\Big)=i_n\big(\zeta_n+\delta_n(r_n)\big)
		=\zeta_{n+1},
		\]
		using that $\delta_n$ vanishes on original coordinates and \eqref{eq:unitn}; trivial coverings and pullbacks along constant maps have
		reduced class $0$. Each entry is a unital $*$-homomorphism into the corresponding corner of $F_{\lambda_{n+1}}$ (Lemma~\ref{lem:basicK}), and
		the matrix targets are clearly unital.
	\end{proof}
	
	The degree requirement on $\nu_{n,i}$ concerns a covering of the \emph{source} complex $T_n$, known at stage $n$; this is the point of the
	delay (Figure~\ref{fig:timeline}).
	
	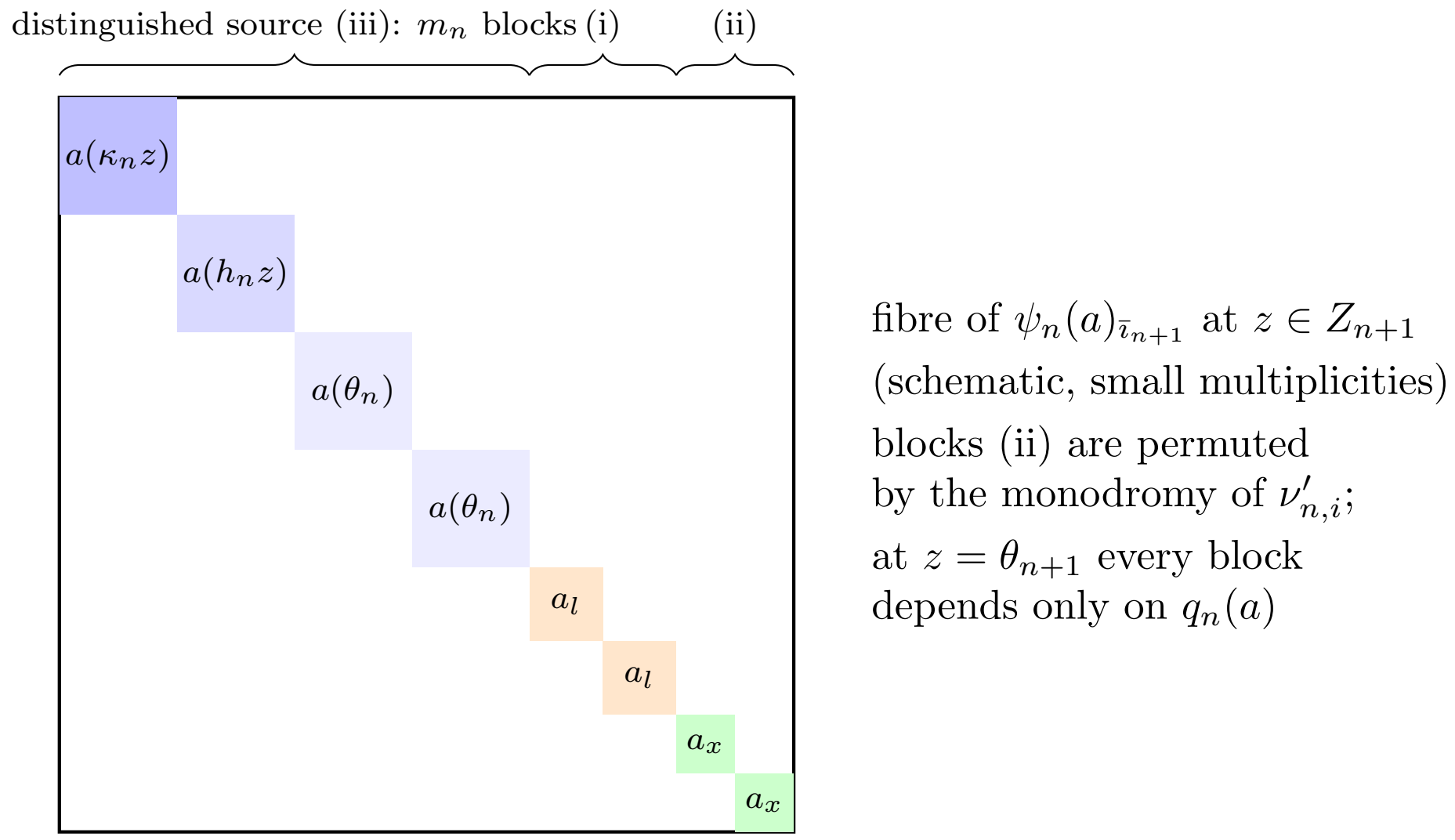
\begin{figure}[ht]
		\centering
		\begin{tikzpicture}[scale=0.62]
			\draw[thick] (0,0) rectangle (10,-10);
			\fill[blue!25] (0,0) rectangle (1.6,-1.6); \node[font=\scriptsize] at (0.8,-0.8) {$a(\kappa_n z)$};
			\fill[blue!15] (1.6,-1.6) rectangle (3.2,-3.2); \node[font=\scriptsize] at (2.4,-2.4) {$a(h_nz)$};
			\fill[blue!8] (3.2,-3.2) rectangle (4.8,-4.8); \node[font=\scriptsize] at (4.0,-4.0) {$a(\theta_n)$};
			\fill[blue!8] (4.8,-4.8) rectangle (6.4,-6.4); \node[font=\scriptsize] at (5.6,-5.6) {$a(\theta_n)$};
			\fill[orange!20] (6.4,-6.4) rectangle (7.4,-7.4); \node[font=\scriptsize] at (6.9,-6.9) {$a_l$};
			\fill[orange!20] (7.4,-7.4) rectangle (8.4,-8.4); \node[font=\scriptsize] at (7.9,-7.9) {$a_l$};
			\fill[green!20] (8.4,-8.4) rectangle (9.2,-9.2); \node[font=\scriptsize] at (8.8,-8.8) {$a_{x}$};
			\fill[green!20] (9.2,-9.2) rectangle (10,-10); \node[font=\scriptsize] at (9.6,-9.6) {$a_{x}$};
			\draw[decorate,decoration={brace,amplitude=5pt}] (0,0.3) -- node[above=4pt,font=\scriptsize]{distinguished source (iii): $m_n$ blocks} (6.4,0.3);
			\draw[decorate,decoration={brace,amplitude=5pt}] (6.4,0.3) -- node[above=4pt,font=\scriptsize,align=center]{(i)} (8.4,0.3);
			\draw[decorate,decoration={brace,amplitude=5pt}] (8.4,0.3) -- node[above=4pt,font=\scriptsize]{(ii)} (10,0.3);
			\node[font=\small,align=left] at (15.0,-5) {fibre of $\psi_n(a)_{\ib_{n+1}}$ at $z\in Z_{n+1}$\\[3pt]
				(schematic, small multiplicities)\\[3pt]
				blocks (ii) are permuted\\ by the monodromy of $\nu'_{n,i}$;\\[3pt]
				at $z=\theta_{n+1}$ every block\\ depends only on $q_n(a)$};
		\end{tikzpicture}
		\caption{The distinguished component of $\psi_n(a)$ at a point $z$: block-diagonal in the sheet basis of $\lambda_{n+1}$.}\label{fig:psi}
	\end{figure}
	
	\begin{lemma}\label{lem:psi}
		The map $\psi_n$ is injective and induces $\gt_n$ on $K_0$ and $\chi_n$ on $K_1$ (in the coordinates of Lemma~\ref{lem:KFn}). It commutes
		with basepoint evaluation:
		\begin{equation}\label{eq:qpsi}
			q_{n+1}\circ\psi_n=\widehat\psi_n\circ q_n,
		\end{equation}
		where $\widehat\psi_n\colon\widehat F_n\to\widehat F_{n+1}$ is, on each target summand $l'$, a labelled block homomorphism
		(Definition~\ref{def:labelled}) of type $((M_n)_{l',l})_{l\in J_n}$; for $l'\ne\ib_{n+1}$ its labelling is the standard one. Moreover $\psi_n$
		maps diagonal, monomial and nonnegative monomial elements to elements of the same kind, and it commutes with fibrewise diagonal compression.
		Finally, the distinguished target alone is faithful: $\psi_n(a)_{\ib_{n+1}}=0$ implies $a=0$.
	\end{lemma}
	
	\begin{proof}
		\emph{$K_0$.} Use Lemma~\ref{lem:basicK} entry by entry. Into a matrix target $l'$, a minimal projection of the source summand $l$ goes to
		$(M_n)_{l',l}$ minimal projections, and reduced classes of the distinguished source go to $0$ (they are killed by $q_n$). Into the
		distinguished target: a minimal projection of an original matrix source $l$ goes to $[V_{\nu_l}]=(M_n)_{\ib_{n+1},l}+0$; a minimal projection
		of an extra source $x_i$ goes to $[V_{\nu'_{n,i}}]=m_n+i_n\delta_n(x_i)$; and a class $(k,e)\in\Z\oplus E_n=K_0(F_{\lambda_n})$ goes to
		$\sum_{\text{sheets}}g_{\text{sheet}}^*(k,e)=(m_nk,\ \kappa_n^*e+h_n^*e+0)=(m_nk,\ i_n(e))$. Altogether, $(h,e)\mapsto\big(M_nh,\
		i_n(e)+\sum_ih_{x_i}i_n\delta_n(x_i)\big)=\gt_n(h,e)$.
		
		\emph{$K_1$.} The matrix summands have $K_1=0$, and on $K_1(F_{\lambda_n})=K^1(Z_n)$ the distinguished entry induces
		$\kappa_n^*+h_n^*+(m_n-2)\cdot\mathrm{const}^*=\chi_n+0+0$ by \eqref{eq:kappah}.
		
		\emph{Basepoint evaluation.} All maps $\kappa_n$, $h_n$, $\theta_n$ send $\theta_{n+1}$ to $\theta_n$. Hence at $z=\theta_{n+1}$ every block of
		$\psi_n(a)_{\ib_{n+1}}$ is either a copy of $a_{\ib_n}(\theta_n)=q_n(a)_{\ib_n}$ or a copy of a constant matrix $a_l=q_n(a)_l$, written in the sheet
		basis of the corresponding fibre over $\theta_{n+1}$; and these fibres are canonically identified with labelled fibres over $\theta_n$. With the
		labelling of Definition~\ref{def:lambda}, the matrix of $\psi_n(a)_{\ib_{n+1}}(\theta_{n+1})$ is therefore obtained from $q_n(a)$ by a labelled
		block homomorphism of type $((M_n)_{\ib_{n+1},l})_l$. Into the matrix targets $\psi_n$ is by definition a standard block homomorphism of $q_n(a)$.
		This proves \eqref{eq:qpsi}. (If an eigenvalue map were not based, the value at $\theta_{n+1}$ would involve values of $a$ away from $\theta_n$,
		and \eqref{eq:qpsi} would fail.)
		
		\emph{Faithfulness and injectivity.} Suppose $\psi_n(a)_{\ib_{n+1}}=0$. The block $a_{\ib_n}\circ h_n$ vanishes, and $h_n$ is surjective, so
		$a_{\ib_n}=0$. Every matrix source $l\ne\ib_n$ occurs in the distinguished target at least once, because $(M_n)_{\ib_{n+1},l}\ge1$; so $a_l=0$. Hence
		$a=0$, and in particular $\psi_n$ is injective.
		
		\emph{Diagonals and monomials.} Each entry is a basic homomorphism, so Lemma~\ref{lem:basic} applies, and block-diagonal sums preserve the three
		classes and commute with diagonal compression (Lemma~\ref{lem:monomialprops}(d)).
	\end{proof}
	
	\subsubsection*{The AF algebra of basepoint values}
	The $C^*$-algebra $\widehat F=\lims(\widehat F_n,\widehat\psi_n)$ is a unital AF algebra. By Lemma~\ref{lem:psi} the connecting maps $\widehat\psi_n$ have
	multiplicity matrices $M_n$, so by Lemma~\ref{lem:delayed} its ordered $K_0$-group is
	\[
	\lims\big(\Z^{J_n},\Z_{\ge0}^{J_n},M_n\big)=(F,F^+),
	\]
	and the class of its unit is the image of $r_n$, which is $\bar u$ by \eqref{eq:unitn}. Hence $\widehat F$ is simple (Theorem~\ref{thm:Fdimgroup}). Its
	traces will be identified in Lemma~\ref{lem:AFtraces}: $T(\widehat F)\cong S(F,\bar u)\cong T$. In both cases the traces of the limit $C^*$-algebra will be
	compared with those of $\widehat F$.
	
	\subsection{Labelled block homomorphisms}\label{ss:interval}
	
	The endpoints of the interval algebras are attached through block embeddings of $\widehat F_n$ in which the order of the blocks and of the basis
	vectors inside them matters, and in which some coordinates may be left empty; we encode this by a \emph{labelling}.
	
	\begin{definition}\label{def:labelled}
		Let $\widehat F=\bigoplus_{l\in J}M_{r_l}$, let $(m_l)_{l\in J}$ be nonnegative integers and $s\ge\sum_lm_lr_l$. A \emph{labelling of type $(m_l)$
			into $[s]$} is an injection
		\[
		\iota\colon\bigsqcup_{l\in J}[m_l]\times[r_l]\to[s],\qquad(l,c,k)\mapsto\iota(l,c,k)
		\]
		(the $k$th basis vector of the $c$th copy of the summand $l$ is placed at position $\iota(l,c,k)$). The associated \emph{labelled block
			homomorphism} $\beta=\beta_\iota\colon\widehat F\to M_s$ is
		\[
		\beta(x)_{\iota(l,c,k),\,\iota(l,c,k')}=(x_l)_{k,k'},\qquad\text{all other entries being }0 .
		\]
		Its \emph{support} is $S(\beta)=\iota\big(\bigsqcup_l[m_l]\times[r_l]\big)$, of cardinality $\sum_lm_lr_l$, and its \emph{zero block} is
		$Z(\beta)=[s]\setminus S(\beta)$. Its \emph{label map} is $S(\beta)\to\bigsqcup_l[r_l]$, $\iota(l,c,k)\mapsto(l,k)$. The \emph{standard labelling}
		(when $s=\sum_lm_lr_l$) places the copies in the standard block order: first the $m_{l_1}$ copies of $M_{r_{l_1}}$, each with its basis in the
		standard order, then those of $M_{r_{l_2}}$, and so on.
	\end{definition}
	
	\begin{example}
		Let $\widehat F=M_2\oplus M_1$ with $m=(1,2)$. For $s=4$, the standard labelling gives $\beta(x\oplus y)=\diag(x,y,y)$, and the labelling
		$\iota(1,1,1)=1$, $\iota(2,1,1)=2$, $\iota(1,1,2)=3$, $\iota(2,2,1)=4$ gives
		\[
		\beta(x\oplus y)=\begin{pmatrix}x_{11}&0&x_{12}&0\\0&y&0&0\\x_{21}&0&x_{22}&0\\0&0&0&y\end{pmatrix}
		=\sigma\,\diag(x,y,y)\,\sigma^{-1}
		\]
		for the permutation matrix $\sigma$ of the transposition of the second and third basis vectors. For $s=5$ and the standard positions in the first
		four coordinates, $\beta(x\oplus y)=\diag(x,y,y,0)$: the fifth coordinate is the zero block, and $\beta(1)=\diag(1,1,1,1,0)$ is not the unit.
	\end{example}
	
	\begin{lemma}\label{lem:labelledprops}
		\begin{enumerate}[label=\textup{(\alph*)}]
			\item A labelled block homomorphism $\beta$ is a $*$-homomorphism; $\beta(1)$ is the diagonal projection onto the coordinates in $S(\beta)$, so
			$\beta$ is unital if and only if $Z(\beta)=\emptyset$. It maps diagonal, monomial and nonnegative monomial elements to matrices of the same kind,
			and it commutes with diagonal compression.
			\item Two labelled block homomorphisms of the same type into the same $M_s$ differ by conjugation with a permutation matrix: if
			$\iota,\iota'$ are labellings of the same type, any permutation of $[s]$ which extends $\iota'\circ\iota^{-1}$ on $S(\beta_\iota)$ (and maps
			$Z(\beta_\iota)$ bijectively onto $Z(\beta_{\iota'})$) conjugates $\beta_\iota$ into $\beta_{\iota'}$. Conversely, $\Ad(\sigma)\circ\beta_\iota=\beta_{\sigma\circ\iota}$
			is labelled of the same type.
			\item Direct sums and composites of labelled block homomorphisms are labelled block homomorphisms; the type of a composite is the product of the
			types (as multiplicity matrices).
			\item On $K_0$, a labelled block homomorphism of type $(m_l)$ induces $\Z^J\to\Z$, $h\mapsto\sum_lm_lh_l$, whether or not it is unital.
		\end{enumerate}
	\end{lemma}
	
	\begin{proof}
		(a) $\beta(x)$ is the matrix of $\bigoplus_{l,c}x_l\oplus0_{Z(\beta)}$ in the basis reordered by $\iota$, so $\beta$ is a $*$-homomorphism with
		$\beta(1)$ as stated. A diagonal (monomial, nonnegative monomial) $x$ gives a block-diagonal sum of such matrices and zeros, conjugated by a
		permutation (Lemma~\ref{lem:monomialprops}(d)). (b) Both are the same block sum in two orders. (c) Reordering a block sum of block sums (with
		zero blocks) is a block sum. (d) A minimal projection of $M_{r_l}$ goes to $m_l$ orthogonal rank-one projections.
	\end{proof}
	
	\subsection{The endpoint data and the building blocks}\label{ss:An}
	
	Choose strictly positive integer matrices $\Lambda_{n,0},\Lambda_{n,1}\colon\Z^{J_n}\to\Z^{l_n}$ with
	\begin{equation}\label{eq:Bn}
		\Lambda_{n,1}-\Lambda_{n,0}=\Delta_n .
	\end{equation}
	For instance, write $\Delta_n=\Delta_n^+-\Delta_n^-$ with $\Delta_n^\pm\ge0$ entrywise, choose a strictly positive matrix $Q_n$ (the \emph{common
		part}), and put $\Lambda_{n,0}=\Delta_n^-+Q_n$, $\Lambda_{n,1}=\Delta_n^++Q_n$. (This is how the matrices are chosen at stage $1$; at later stages they
	are produced, together with the rest of the endpoint data, in the construction of the connecting maps.) Put
	\[
	\varsigma_{n,i}:=\Lambda_{n,i}r_n\in\Z_{>0}^{l_n}\quad(i=0,1),\qquad\text{so that}\qquad\varsigma_{n,1}-\varsigma_{n,0}=\Delta_nr_n=\mathbf d_n
	\]
	by \eqref{eq:dn}. The \emph{endpoint data} at stage $n$ consist of
	\begin{itemize}
		\item sizes $s_{n,q}\ge\max\{\varsigma_{n,0,q},\varsigma_{n,1,q}\}$ ($1\le q\le l_n$), and the $C^*$-algebra $E_n'=\bigoplus_qM_{s_{n,q}}$ (the letter $q$ indexes
		interval summands from now on);
		\item for $i=0,1$ a homomorphism $\beta_{n,i}\colon\widehat F_n\to E_n'$ whose $q$th component is a labelled block homomorphism of type
		$((\Lambda_{n,i})_{q,l})_{l\in J_n}$ into $M_{s_{n,q}}$, with support $S^q_{n,i}$ (of cardinality $\varsigma_{n,i,q}$) and zero block $Z^q_{n,i}$, subject to
		\begin{equation}\label{eq:ZC}
			\textup{(Z)}\quad Z^q_{n,0}\cap Z^q_{n,1}=\emptyset,\qquad\qquad\textup{(C)}\quad S^q_{n,0}\cap S^q_{n,1}\neq\emptyset\qquad\text{for every }q.
		\end{equation}
	\end{itemize}
	Condition (Z) says that every coordinate of $[s_{n,q}]$ lies in the support of at least one endpoint map, and (C) says that some coordinate lies in
	both; in particular $\max_i\varsigma_{n,i,q}\le s_{n,q}<\varsigma_{n,0,q}+\varsigma_{n,1,q}$.
	
	\begin{itemize}
		\item In \emph{Case (U)}, $\mathbf d_n=0$ by \eqref{eq:dnzero}, so $\varsigma_{n,0}=\varsigma_{n,1}$. We always take $s_n=\varsigma_{n,0}=\varsigma_{n,1}$; then both endpoint
		maps are unital, all zero blocks are empty, and (Z), (C) hold trivially.
		\item In \emph{Case (SP)}, $\mathbf d_n\ne0$ by \eqref{eq:dnnonzero}, so for some $q$ the two supports have different cardinalities and some zero block is
		nonempty. At stage $1$ we take $s_{1,q}=\max_i\varsigma_{1,i,q}$ and nested supports; then one support is all of $[s_{1,q}]$ and the other is
		nonempty, so (Z) and (C) hold. At later stages they are inherited from the connecting maps (Lemma~\ref{lem:intervaldata}).
	\end{itemize}
	The labellings $\iota^q_{n,i}$ are fixed during the construction: for $n=1$ standard (or nested), for $n\ge2$ chosen in the construction of the
	connecting map $A_{n-1}\to A_n$.
	
	\begin{definition}[The building block]
		\begin{equation}\label{eq:An}
			A_n=\{(f,a)\in C([0,1],E_n')\oplus F_n:\ f(0)=\beta_{n,0}(q_n(a)),\ f(1)=\beta_{n,1}(q_n(a))\}.
		\end{equation}
		The \emph{fibre matrices} of $(f,a)\in A_n$ are the matrices $f_q(t)\in M_{s_{n,q}}$ ($t\in[0,1]$, standard basis) and the fibre matrices of $a$. We
		call $(f,a)$ \emph{diagonal}, \emph{monomial} or \emph{nonnegative monomial} if all its fibre matrices are. Let $B_n\subset A_n$ be the set of diagonal
		elements, $\EE_n\colon A_n\to B_n$ fibrewise diagonal compression, and $N_{+,n}$ the set of nonnegative monomial elements. Let
		$\mathcal I_n=C_0((0,1),E_n')$, an ideal of $A_n$, and write $\bar q_n\colon A_n\to\widehat F_n$, $\bar q_n(f,a)=q_n(a)$. Put
		\begin{equation}\label{eq:en}
			e_n:=(f_n,1_{F_n})\in A_n,\qquad f_n(t)=(1-t)\beta_{n,0}(1)+t\,\beta_{n,1}(1).
		\end{equation}
	\end{definition}
	
	Since $\beta_{n,i}$ and $q_n$ commute with diagonal compression (Lemmas~\ref{lem:labelledprops}(a) and~\ref{lem:basic}), $\EE_n$ maps $A_n$ onto $B_n$; it
	is a conditional expectation (fibrewise it is the diagonal compression of matrices, which is positive, contractive, idempotent and $B_n$-bimodular),
	and it is faithful (if $\EE_n(x^*x)=0$ then every fibre matrix of $x^*x$ has zero trace, so $x=0$). Note that the endpoint conditions only involve the
	basepoint value $q_n(a)$. In Case (U), $e_n=(1,1)$ is the unit of $A_n$.
	
	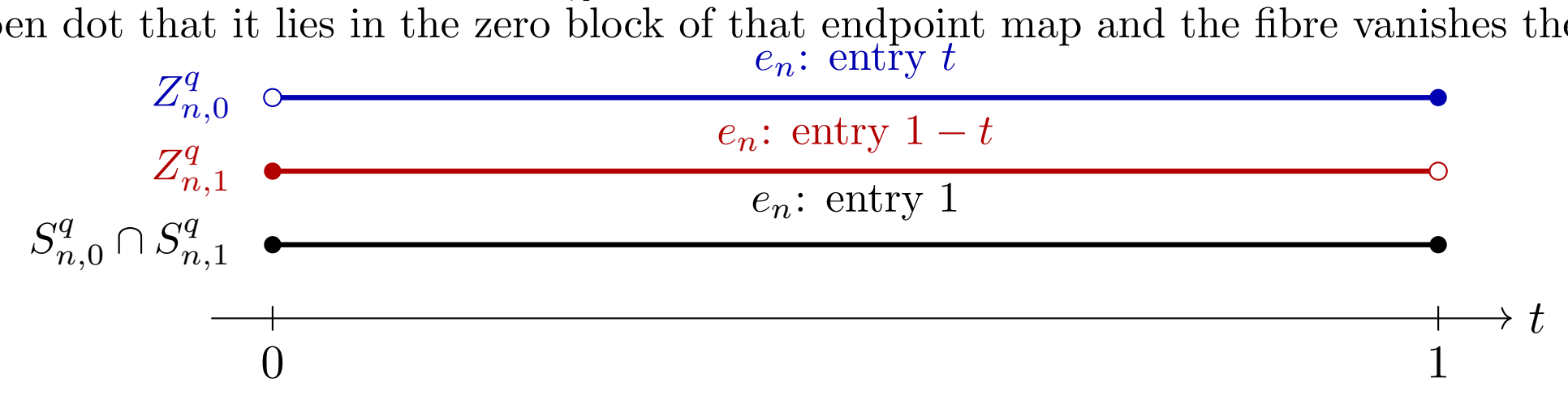
\begin{figure}[ht]
		\centering
		\begin{tikzpicture}[scale=1.0]
			\draw[->] (0,0) -- (10.6,0) node[right]{$t$};
			\draw (0.5,0.1) -- (0.5,-0.1) node[below]{$0$}; \draw (10,0.1) -- (10,-0.1) node[below]{$1$};
			\foreach \y/\lab/\col in {0.6/{S^q_{n,0}\cap S^q_{n,1}}/black,1.2/{Z^q_{n,1}}/red!70!black,1.8/{Z^q_{n,0}}/blue!70!black}{
				\node[left,font=\small,\col] at (0.3,\y) {$\lab$};}
			\draw[very thick] (0.5,0.6) -- (10,0.6);
			\draw[very thick,red!70!black] (0.5,1.2) -- (10,1.2);
			\draw[very thick,blue!70!black] (0.5,1.8) -- (10,1.8);
			\fill[black] (0.5,0.6) circle (2pt); \fill[black] (10,0.6) circle (2pt);
			\fill[red!70!black] (0.5,1.2) circle (2pt); \draw[red!70!black,fill=white] (10,1.2) circle (2pt);
			\draw[blue!70!black,fill=white] (0.5,1.8) circle (2pt); \fill[blue!70!black] (10,1.8) circle (2pt);
			\node[font=\small,align=left] at (5.25,2.6) {coordinates of an interval fibre $M_{s_{n,q}}$; a filled dot means the coordinate carries the value of
				$\beta_{n,i}q_n(a)$,\\ an open dot that it lies in the zero block of that endpoint map and the fibre vanishes there};
			\node[font=\small] at (5.25,0.95) {$e_n$: entry $1$};
			\node[font=\small,red!70!black] at (5.25,1.5) {$e_n$: entry $1-t$};
			\node[font=\small,blue!70!black] at (5.25,2.1) {$e_n$: entry $t$};
		\end{tikzpicture}
		\caption{The three kinds of coordinates of an interval fibre in Case (SP). In Case (U) every coordinate lies in the common support
			$S^q_{n,0}\cap S^q_{n,1}$ (bottom row) and has entry $1$.}\label{fig:zeroblocks}
	\end{figure}
	
	\begin{lemma}\label{lem:An}\label{lem:irreps}
		\begin{enumerate}[label=\textup{(\alph*)}]
			\item $A_n$ is a separable subhomogeneous $C^*$-algebra. It is unital if and only if all zero blocks are empty; this is the case in Case (U), and
			not in Case (SP).
			\item Every irreducible representation of $A_n$ is unitarily equivalent to exactly one of: a \emph{homogeneous evaluation}
			$\ev_z\colon(f,a)\mapsto a_{\ib_n}(z)\in M_{r_{n,\ib_n}}$, $z\in Z_n$ (in a sheet basis); a \emph{matrix summand} $(f,a)\mapsto a_l$, $l\in J_n\setminus\{\ib_n\}$;
			an \emph{interior interval evaluation} $\ev_{q,t}\colon(f,a)\mapsto f_q(t)$, $t\in(0,1)$. The endpoint evaluations $(f,a)\mapsto f_q(i)=\beta_{n,i}(q_n(a))_q$
			are not irreducible: they factor through $\bar q_n$, are direct sums of the basepoint representations $(f,a)\mapsto q_n(a)_l$ (and zeros on $Z^q_{n,i}$).
			\item $e_n$ is a strictly positive element of $A_n$ and lies in $B_n$. Its image in every irreducible representation is invertible, and its fibre matrix at
			$t\in(0,1)$ is diagonal with entries $1$ on $S^q_{n,0}\cap S^q_{n,1}$, $1-t$ on $Z^q_{n,1}$ and $t$ on $Z^q_{n,0}$ (Figure~\ref{fig:zeroblocks}). More generally,
			$(f,a)\in(A_n)_+$ is strictly positive if and only if $a$ is invertible in $F_n$ and $f_q(t)$ is invertible for all $q$ and all $t\in(0,1)$.
			\item The map $K_0(A_n)\to K_0(F_n)=\Z^{J_n}\oplus E_n$ induced by $(f,a)\mapsto a$ is injective with image $\ker(\Delta_n\oplus0)=G_n$, and $K_1(A_n)\to
			K_1(F_n)=K_n$ is an isomorphism:
			\begin{equation}\label{eq:KAn}
				K_0(A_n)\cong G_n,\qquad K_1(A_n)\cong K_n;\qquad\text{in Case (U), }[1_{A_n}]=(r_n,\zeta_n).
			\end{equation}
			\item In Case (SP), $A_n$ is stably projectionless.
		\end{enumerate}
	\end{lemma}
	
	\begin{proof}
		(a) Separability and subhomogeneity are clear (irreducible representations have dimension at most $\max\{s_{n,q},r_{n,l}\}$ by (b)). If all zero blocks are
		empty, $(1,1)$ is a unit. If some $Z^q_{n,i}$ is nonempty, a unit $(f,a)$ would have $f_q(t)=1$ for $t\in(0,1)$, hence $f_q(i)=1\neq\beta_{n,i}(1)_q$.
		In Case (SP) some zero block is nonempty since $\mathbf d_n\ne0$.
		
		(b) Let $\pi$ be an irreducible representation on $H_\pi$. If $\pi(\mathcal I_n)=0$, then $\pi$ factors through $A_n/\mathcal I_n$, which is $F_n$ via
		$(f,a)\mapsto a$ (onto, because $L(a)$ of Lemma~\ref{lem:normalizers} lifts $a$). The irreducible representations of a direct sum are those of the
		summands; those of $M_r$ are the identity, and those of the covering algebra $F_{\lambda_n}=\Gamma(\End V_{\lambda_n})$ are the fibre evaluations, because it is
		locally isomorphic to $C_0(U,M_r)$ and its centre is $C(Z_n)$ (an irreducible representation restricts to a character of the centre, i.e.\ a point $z$, and
		factors through the fibre $\End(V_z)\cong M_r$). Otherwise the closed span of $\pi(\mathcal I_n)H_\pi$ is a nonzero invariant subspace, hence $H_\pi$; an
		operator commuting with $\pi(\mathcal I_n)$ commutes with $\pi(A_n)$, because $\pi(x)\pi(j)=\pi(xj)$ for $j\in\mathcal I_n$; so $\pi|_{\mathcal I_n}$ is an
		irreducible representation of $\mathcal I_n=\bigoplus_qC_0((0,1),M_{s_{n,q}})$, i.e.\ an evaluation at some $(q,t)$, and $\pi$ is its canonical extension
		$(f,a)\mapsto f_q(t)$. Different points give inequivalent representations (different kernels). The description of the endpoint evaluations is
		\eqref{eq:An}; they contain $(\Lambda_{n,i})_{q,l}\ge1$ copies of $q_n(\cdot)_l$ for each of the at least two indices $l$, so they are reducible.
		
		(c) The fibre description follows from \eqref{eq:en} and (Z): a coordinate not in $S^q_{n,0}\cap S^q_{n,1}$ lies in exactly one of $Z^q_{n,0}$, $Z^q_{n,1}$.
		Hence $f_n(t)$ is invertible for $t\in(0,1)$, and the $F_n$-component $1$ is invertible. Necessity in the criterion is clear (a strictly positive element
		is nonzero under every state, in particular invertible in every finite-dimensional irreducible representation), and sufficiency is
		Lemma~\ref{lem:fibrewise}.
		
		(d) Use the six-term exact sequence of $0\to\mathcal I_n\to A_n\to F_n\to0$. Since $K_0(\mathcal I_n)=K_1(E_n')=0$ and $K_1(\mathcal I_n)\cong K_0(E_n')=\Z^{l_n}$,
		it reads
		\[
		0\to K_0(A_n)\to K_0(F_n)\xrightarrow{\ \partial\ }K_1(\mathcal I_n)\to K_1(A_n)\to K_1(F_n)\to K_0(\mathcal I_n)=0 .
		\]
		To compute $\partial$, compare with the extension $0\to\mathcal I_n\to C([0,1],E_n')\to E_n'\oplus E_n'\to0$ via the map $(f,a)\mapsto f$, which is the
		identity on $\mathcal I_n$ and induces $a\mapsto(\beta_{n,0}q_n(a),\beta_{n,1}q_n(a))$ on the quotients. In the second extension
		$K_*(C([0,1],E_n'))=K_*(E_n')$ maps diagonally to $K_*(E_n'\oplus E_n')$, so its boundary map $\partial'\colon K_0(E_n')^2\to K_1(\mathcal I_n)$ is onto with
		kernel the diagonal, i.e.\ $\partial'(y_0,y_1)=\omega(y_1-y_0)$ for an isomorphism $\omega\colon K_0(E_n')\to K_1(\mathcal I_n)$. By naturality of the
		boundary map,
		\[
		\partial=\partial'\circ\big((\beta_{n,0}q_n)_*,(\beta_{n,1}q_n)_*\big)=\omega\circ\big((\beta_{n,1})_*-(\beta_{n,0})_*\big)\circ(q_n)_*=\omega\circ(\Delta_n\oplus0),
		\]
		by Lemma~\ref{lem:KFn}, Lemma~\ref{lem:labelledprops}(d) and \eqref{eq:Bn}; unitality of the $\beta_{n,i}$ plays no role. Since $\Delta_n$ is onto,
		$\partial$ is onto, so $K_1(A_n)\to K_1(F_n)=K_n$ is an isomorphism and $K_0(A_n)\cong\ker\partial=\ker(\Delta_n\oplus0)=G_n$ by \eqref{eq:quot}. In
		Case (U), the unit of $A_n$ maps to $[1_{F_n}]=(r_n,\zeta_n)$.
		
		(e) Let $p\in M_k(A_n)$ be a projection. Its image in $M_k(F_n)$ is a projection whose class $(h,e)$ has $h\ge0$ ($h$ lists fibre ranks, constant along
		the connected space $Z_n$) and lies in the image of $K_0(A_n)$, that is, in $G_n$. By Lemma~\ref{lem:nopositive}, $h=0$, so the image of $p$ in $M_k(F_n)$ is
		$0$. Thus $p\in M_k(\mathcal I_n)$ is a continuous projection-valued function on $(0,1)$ vanishing at the endpoints; its rank is locally constant on the
		connected interval and tends to $0$, so $p=0$.
	\end{proof}
	
	\begin{lemma}\label{lem:fibrewise}
		Let $h=(h_J,h_F)\in(A_n)_+$ be such that $h_F$ is invertible in $F_n$ and $h_J(t)$ is invertible for every $t\in(0,1)$. Then $h$ is strictly positive in $A_n$.
	\end{lemma}
	
	\begin{proof}
		A positive element is strictly positive if and only if it is nonzero under every state \cite[3.10.5]{Ped}. Suppose $\varphi(h)=0$ for a state $\varphi$ of
		$A_n$. By the Cauchy--Schwarz inequality, $|\varphi(h^{1/2}xh^{1/2})|\le\varphi(h)^{1/2}\varphi(h^{1/2}x^*xh^{1/2})^{1/2}=0$ for all $x$, so $\varphi$ vanishes on
		$h^{1/2}\mathcal I_nh^{1/2}$. For $k\in C_c((0,1),E_n')$ the function $t\mapsto h_J(t)^{-1/2}k(t)h_J(t)^{-1/2}$, extended by $0$, lies in $\mathcal I_n$, so
		$k\in h^{1/2}\mathcal I_nh^{1/2}$. Since $C_c((0,1),E_n')$ is dense in $\mathcal I_n$, $\varphi$ vanishes on $\mathcal I_n$ and factors through $A_n/\mathcal
		I_n\cong F_n$ as a state $\bar\varphi$; then $\varphi(h)=\bar\varphi(h_F)\ge\|h_F^{-1}\|^{-1}>0$, a contradiction.
	\end{proof}
	
	\begin{remark}[The role of the intervals]
		The quotient $F_n$ has $K_0=\Z^{J_n}\oplus E_n$, which presents $H$ rather than $G_0$ (Lemma~\ref{lem:delayed}). The interval algebras cut $K_0$ down to the
		kernel of $\Delta_n=\Lambda_{n,1}-\Lambda_{n,0}$, which presents $G_0$: a projection over $F_n$ lifts to $A_n$ only if its two endpoint images have the same rank.
		This is the mechanism of the one-dimensional building blocks of \cite{GLN} and \cite{Li}. In Case (SP) the unit of $F_n$ does not lift, because its two
		endpoint images have ranks $\varsigma_{n,0}\ne\varsigma_{n,1}$; the zero blocks are what remains of it.
	\end{remark}
	
	\subsection{The groupoid model}\label{ss:groupoid}
	
	\begin{lemma}\label{lem:groupoid}
		There is a second countable, locally compact Hausdorff, principal étale groupoid $\mathcal G_n$ with unit space $\Omega_n$ and an isomorphism $A_n\cong
		C^*_r(\mathcal G_n)$ carrying $B_n$ onto $C_0(\Omega_n)$ and $\EE_n$ onto restriction to the unit space. In particular $(A_n,B_n)$ is a $C^*$-diagonal with
		expectation $\EE_n$, and its Weyl twist is trivial. The unit space $\Omega_n$ is compact if and only if $A_n$ is unital. Moreover, the basepoint fibre
		\[
		W_n^\theta:=\lambda_n^{-1}(\theta_n)\sqcup\bigsqcup_{l\neq\ib_n}[r_{n,l}]\;=\;\bigsqcup_{l\in J_n}[r_{n,l}]\quad\text{(via the fixed labels)}
		\]
		is a finite closed subset of $\Omega_n$, the restriction $\widehat R_n=\bigsqcup_l[r_{n,l}]^2$ of $\mathcal G_n$ to it is a finite closed subgroupoid with
		$C^*(\widehat R_n)=\widehat F_n$, and under the isomorphism the map $\bar q_n$ is restriction of functions to $\widehat R_n$; in particular $\bar q_n|_{B_n}$ is
		restriction of functions to $W^\theta_n$.
	\end{lemma}
	
	We write $\mathcal G_n$ for the groupoid, to avoid confusion with the groups $G_n=K_0(A_n)$. Arrows of pair groupoids are written $(x,y)$, with range $x$
	and source $y$.
	
	\begin{proof}
		\emph{Step 1: the groupoid of $F_n$.} Let $\varpi\colon\widetilde Y_n\to Z_n$ be the covering $\lambda_n$, and let
		\[
		W_n=\widetilde Y_n\sqcup\bigsqcup_{l\ne\ib_n}[r_{n,l}],\qquad R^F_n=(\widetilde Y_n\times_{Z_n}\widetilde Y_n)\sqcup\bigsqcup_{l\ne\ib_n}[r_{n,l}]^2 .
		\]
		By Lemma~\ref{lem:Rpietale} and Proposition~\ref{prop:coveringalgebra}, $R^F_n$ is a compact principal étale groupoid with unit space $W_n$, $C^*(R^F_n)=F_n$
		and $C(W_n)$ the diagonal; the value of $a\in F_n$ at an arrow is the corresponding entry of a fibre matrix. The basepoint fibre $W_n^\theta$ is identified,
		via the fixed labels, with the finite set $\bigsqcup_{l\in J_n}[r_{n,l}]$. It is \emph{invariant}: an arrow with source or range in $W_n^\theta$ has both in
		$W_n^\theta$ (it is a full fibre of $\varpi$ together with whole matrix pieces). Hence $\widehat R_n=s^{-1}(W_n^\theta)=r^{-1}(W_n^\theta)$; it is finite, hence
		closed and discrete, and $C^*(\widehat R_n)=\widehat F_n$, with restriction of functions equal to $q_n$.
		
		\emph{Step 2: the endpoint groupoids.} For each $q$ and $i\in\{0,1\}$, the labelling $\iota=\iota^q_{n,i}$ defines the subgroupoid
		\[
		R^q_i=\{(\iota(l,c,k),\iota(l,c,k')):\ l\in J_n,\ c\in[(\Lambda_{n,i})_{q,l}],\ k,k'\in[r_{n,l}]\}\subset(S^q_{n,i})^2\subset[s_{n,q}]^2
		\]
		(pairs of positions in the same copy), whose unit space is the support $S^q_{n,i}$, and the groupoid homomorphism
		\[
		b^q_i\colon R^q_i\to\widehat R_n,\qquad(\iota(l,c,k),\iota(l,c,k'))\mapsto((l,k),(l,k')),
		\]
		which on units is the label map. It is \emph{fibrewise bijective}: for each unit $y=\iota(l,c,k')$, $b^q_i$ maps the arrows of $R^q_i$ with source $y$,
		namely $(\iota(l,c,k),y)$ for $k\in[r_{n,l}]$, bijectively onto the arrows of $\widehat R_n$ with source $(l,k')$; similarly for ranges. Since $\Lambda_{n,i}$ is
		strictly positive, $b^q_i$ is onto on units. In matrix language: $\beta_{n,i}(x)_q$ is the function $x\circ b^q_i$ on $R^q_i$, extended by zero to
		$[s_{n,q}]^2$. The coordinates in the zero block $Z^q_{n,i}$ carry no units and no arrows at the level $t=i$.
		
		\emph{Step 3: the groupoid.} For each $q$ let
		\[
		\check{\mathcal G}^q=\{(t,x,y)\in[0,1]\times[s_{n,q}]^2:\ (x,y)\in R^q_i\text{ if }t=i\in\{0,1\}\}.
		\]
		It is the complement of a finite set in $[0,1]\times[s_{n,q}]^2$, hence open, locally compact and Hausdorff; with $(t,x,y)(t,y,z)=(t,x,z)$ it is an étale
		groupoid (fibrewise a pair groupoid, and a subgroupoid at the endpoints). Let $\check{\mathcal G}_n=R^F_n\sqcup\bigsqcup_q\check{\mathcal G}^q$, and let
		\[
		\mathcal G_n=\check{\mathcal G}_n/{\sim},\qquad(i,x,y)\sim b^q_i(x,y)\quad\text{for }(x,y)\in R^q_i,
		\]
		with the quotient topology. As an abstract groupoid, $\mathcal G_n$ is the disjoint union of $R^F_n$ and the pair groupoids $\{t\}\times[s_{n,q}]^2$, $t\in(0,1)$;
		the identification affects only the topology. (It is a groupoid because the $b^q_i$ are homomorphisms, and arrows at interior points $t\in(0,1)$ are never
		composable with arrows of $R^F_n$.) Every arrow is determined by its range and source, so $\mathcal G_n$ is principal. Its unit space is
		\[
		\Omega_n=\Big(W_n\sqcup\bigsqcup_q\big([0,1]\times[s_{n,q}]\setminus\{(i,x):x\in Z^q_{n,i}\}\big)\Big)\Big/\!\sim ,
		\]
		the endpoint $(i,x)$, $x\in S^q_{n,i}$, being identified with its label in $W_n^\theta$.
		
		\begin{figure}[ht]
			\centering
			\begin{tikzpicture}[scale=0.95]
				\draw[thick,fill=blue!6] (0,0) ellipse (1.6 and 2.0); \node at (-0.5,1.3) {$W_n$};
				\foreach \y/\c in {0.9/red!70!black,0.3/blue!70!black,-0.3/green!50!black,-0.9/orange!80!black}{\fill[\c] (0.8,\y) circle (2pt);}
				\node[font=\scriptsize] at (0.45,-1.45) {$W_n^\theta$};
				\foreach \k/\c/\a/\b in {1/red!70!black/0.9/0.3, 2/blue!70!black/0.3/-0.9, 3/green!50!black/-0.3/0.9}{
					\draw[thick,\c] (0.8,\a) .. controls (2.5,\a+0.5*\k) and (5.5,\b+0.5*\k) .. (7.2,\b);}
				\draw[thick,orange!80!black] (0.8,-0.9) .. controls (2.5,-2.4) and (5.5,-2.2) .. (7.2,-0.3);
				\draw[thick,violet,dashed] (2.0,-2.6) .. controls (3.5,-2.6) and (5.5,-2.4) .. (7.2,-1.5);
				\draw[violet,fill=white] (2.0,-2.6) circle (2pt);
				\draw[thick,dashed] (7.2,1.6) -- (7.2,-1.6);
				\node[font=\small,align=left] at (10.4,0.2) {unit space $\Omega_n$: each interval\\ $[0,1]\times\{x\}$ starts at the label of $x$\\ under $\iota_{n,0}^q$ and ends at its\\
					label under $\iota_{n,1}^q$ (the right ends\\ are glued back to $W_n^\theta$); a coordinate\\ in a zero block (violet) has an\\ open end: $\Omega_n$ is not compact};
			\end{tikzpicture}
			\caption{The unit space of $\mathcal G_n$ (schematic).}\label{fig:unitspace}
		\end{figure}
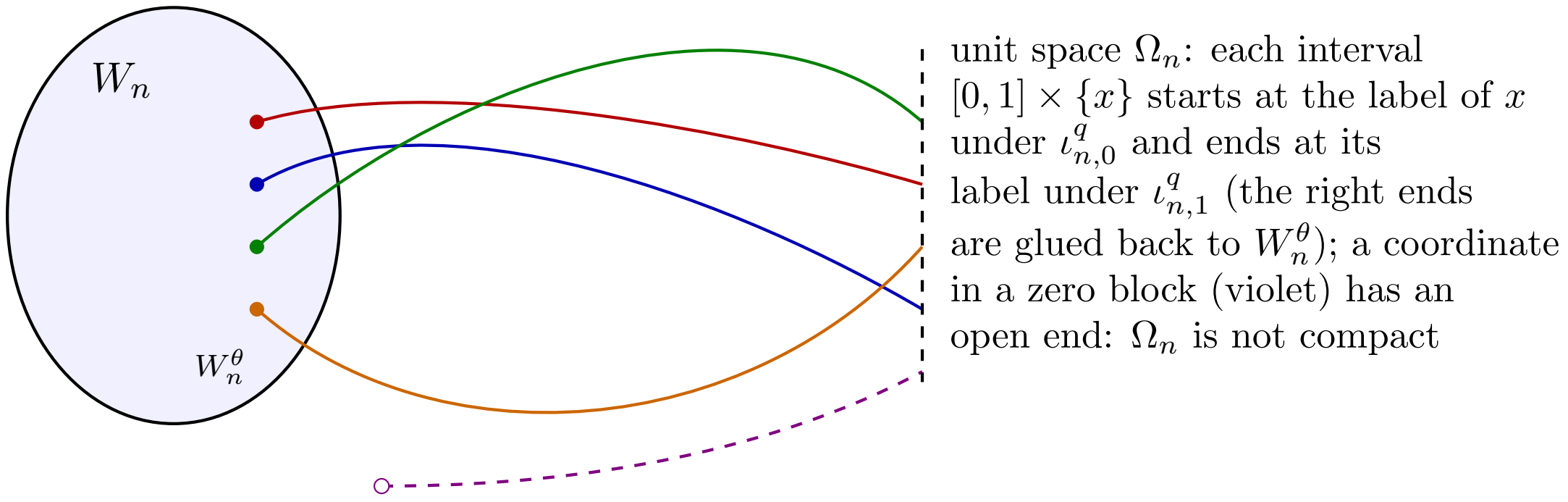
		
		\emph{Step 4: neighbourhoods.} Points not in $\widehat R_n$ keep their neighbourhoods: open subsets of $R^F_n$ disjoint from $\widehat R_n$, and open subsets of
		$(0,1)\times[s_{n,q}]^2$; their preimages in $\check{\mathcal G}_n$ are open, so they are open in $\mathcal G_n$. For $\gamma\in\widehat R_n$, the sets
		\[
		N(V,\varepsilon)=V\cup\bigcup_{q,i}\ O_i(\varepsilon)\times(b^q_i)^{-1}(\gamma),\qquad O_0(\varepsilon)=(0,\varepsilon),\ O_1(\varepsilon)=(1-\varepsilon,1),
		\]
		with $V$ an open bisection of $R^F_n$ such that $V\cap\widehat R_n=\{\gamma\}$, and $0<\varepsilon<\frac12$, form a neighbourhood basis of $\gamma$: the preimage of
		$N(V,\varepsilon)$ in $\check{\mathcal G}_n$ is $V$ together with the half-open intervals $[0,\varepsilon)\times(b^q_0)^{-1}(\gamma)$ and
		$(1-\varepsilon,1]\times(b^q_1)^{-1}(\gamma)$, which is open; and every open set containing $\gamma$ contains some $N(V,\varepsilon)$, because its preimage is open and
		contains the finitely many points $(i,x,y)$ with $b^q_i(x,y)=\gamma$.
		
		\emph{Step 5: Hausdorff, locally compact, second countable.} Two distinct points of $\widehat R_n$ have disjoint neighbourhoods $N(V,\varepsilon)$,
		$N(V',\varepsilon)$, because $R^F_n$ is Hausdorff and the fibres $(b^q_i)^{-1}(\gamma)$, $(b^q_i)^{-1}(\gamma')$ are disjoint. A point of $\widehat R_n$ and an interior
		point $(t,x,y)$ are separated by choosing $\varepsilon<\min\{t,1-t\}$; the remaining cases are clear. The closure of $N(V,\varepsilon/2)$ is contained in the image of
		the compact set $\overline V\sqcup\bigcup_{q,i}[\text{closed intervals of length }\varepsilon/2]\times(b^q_i)^{-1}(\gamma)$, so it is compact; hence $\mathcal G_n$ is
		locally compact. A countable basis is obtained from a countable basis of the compact metrizable space $R^F_n$ and rational intervals.
		
		\emph{Step 6: continuity of the operations.} Inversion maps $N(V,\varepsilon)$ onto $N(V^{-1},\varepsilon)$ for $\gamma^{-1}$. For multiplication, let
		$\gamma=\gamma_1\gamma_2$ with $\gamma_1,\gamma_2\in\widehat R_n$ (the other cases are clear, by invariance of $W_n^\theta$: if a product lies in $\widehat R_n$, so do both
		factors). Given $N(V,\varepsilon)$, choose bisections $V_1\ni\gamma_1$, $V_2\ni\gamma_2$ of $R^F_n$ meeting $\widehat R_n$ only in $\gamma_1,\gamma_2$, with $V_1V_2\subset V$. A
		composable pair in $N(V_1,\varepsilon)\times N(V_2,\varepsilon)$ lies either in $V_1\times V_2$, with product in $V$, or in a single interval level $\{t\}\times[s_{n,q}]^2$
		with $t\in O_i(\varepsilon)$, as $(t,x,z)$, $(t,z,w)$ with $b^q_i(x,z)=\gamma_1$, $b^q_i(z,w)=\gamma_2$; then $(x,w)\in R^q_i$ and $b^q_i(x,w)=\gamma_1\gamma_2=\gamma$, so the
		product $(t,x,w)$ lies in $N(V,\varepsilon)$.
		
		\emph{Step 7: étale.} Let $\gamma\in\widehat R_n$ and $N(V,\varepsilon)$ as above. The range map is injective on $N(V,\varepsilon)$: on $V$ because $V$ is a bisection,
		and on $O_i(\varepsilon)\times(b^q_i)^{-1}(\gamma)$ because, by fibrewise bijectivity, for each unit $x\in S^q_{n,i}$ with label $r(\gamma)$ exactly one arrow of $R^q_i$
		with range $x$ is mapped to $\gamma$. Its image is
		\[
		r(V)\cup\bigcup_{q,i}O_i(\varepsilon)\times\{x\in S^q_{n,i}:\ \text{label of }x\text{ under }\iota^q_{n,i}\text{ is }r(\gamma)\},
		\]
		which is a basic open neighbourhood of $r(\gamma)$ in $\Omega_n$ (note $r(V)\cap W^\theta_n=\{r(\gamma)\}$ by invariance). The neighbourhood bases correspond, so
		$r|_{N(V,\varepsilon)}$ is a homeomorphism onto an open set. The same holds for $s$. At points not in $\widehat R_n$ this is clear. So $\mathcal G_n$ is étale.
		
		\emph{Step 8: compactness of the unit space.} If all zero blocks are empty, $\Omega_n$ is a quotient of the compact space $W_n\sqcup\bigsqcup_q[0,1]\times[s_{n,q}]$,
		hence compact. If $x\in Z^q_{n,i}$, a sequence $(t_k,x)$ with $t_k\to i$ has no convergent subsequence in $\Omega_n$ (the point $(i,x)$ is missing, and $x$ has no
		label under $\iota^q_{n,i}$), so $\Omega_n$ is not compact. By Lemma~\ref{lem:An}(a) the second case is exactly the non-unital one.
		
		\emph{Step 9: the $C^*$-algebra.} Every source fibre of $\mathcal G_n$ has at most $N=\max\{s_{n,q},r_{n,l},\deg\lambda_n\}$ elements, so
		$\|\Upsilon\|_\infty\le\|\Upsilon\|_r\le N\|\Upsilon\|_\infty$ on $C_c(\mathcal G_n)$; hence $C^*_r(\mathcal G_n)$ is $C_0(\mathcal G_n)$ with convolution (remark after
		Definition~\ref{def:Cred}). For $(f,a)\in A_n$ define $\Upsilon(f,a)\colon\mathcal G_n\to\C$ by
		\[
		\Upsilon(f,a)=a\ \text{on }R^F_n,\qquad\Upsilon(f,a)(t,x,y)=f_q(t)_{x,y}\ \ (t\in(0,1)).
		\]
		\begin{itemize}
			\item \emph{Continuity} at $\gamma\in\widehat R_n$: for $(x,y)\in(b^q_i)^{-1}(\gamma)$ and $t\to i$, $f_q(t)_{x,y}\to f_q(i)_{x,y}=\beta_{n,i}(q_na)_{x,y}=a(b^q_i(x,y))=a(\gamma)$.
			\item \emph{Vanishing at infinity}: the arrows $(t,x,y)$ with $(x,y)\notin R^q_i$ and $t\to i$ leave every compact set, and there $f_q(t)_{x,y}\to f_q(i)_{x,y}=0$,
			because $\beta_{n,i}(\cdot)_q$ vanishes outside $R^q_i$; this includes the case that $x$ or $y$ lies in $Z^q_{n,i}$.
			\item \emph{Surjectivity}: for $\Upsilon\in C_0(\mathcal G_n)$, put $a=\Upsilon|_{R^F_n}$ and $f_q(t)_{x,y}=\Upsilon(t,x,y)$ for $t\in(0,1)$, extended to $t=i$ by
			$f_q(i)_{x,y}=a(b^q_i(x,y))$ for $(x,y)\in R^q_i$ and $0$ otherwise; continuity and vanishing at infinity of $\Upsilon$ say exactly that $f$ is continuous on
			$[0,1]$, and then $f(i)=\beta_{n,i}(q_na)$ by Step~2, so $(f,a)\in A_n$.
			\item \emph{Multiplicativity}: by invariance of $W_n^\theta$, convolution at an arrow of $R^F_n$ only involves arrows of $R^F_n$, where it is the product of $F_n$; at
			$(t,x,z)$ it is $\sum_y\Upsilon(t,x,y)\Upsilon'(t,y,z)$, the matrix product. Involution is the adjoint.
			\item \emph{Isometry}: the regular representation at a unit of $W_n$ is the corresponding fibre representation of $F_n$, and at an interior unit $(t,y)$ it is
			$f_q(t)$ acting on $\C^{s_{n,q}}$. So $\|\Upsilon(f,a)\|_r=\max\{\sup_{t,q}\|f_q(t)\|,\|a\|\}=\|(f,a)\|$.
		\end{itemize}
		Thus $\Upsilon$ is an isometric $*$-isomorphism $A_n\cong C^*_r(\mathcal G_n)$. An element is diagonal iff $\Upsilon(f,a)$ is supported on the units, so $\Upsilon$
		carries $B_n$ onto $C_0(\Omega_n)$ and $\EE_n$ onto restriction to $\Omega_n$. By Theorem~\ref{thm:renault}(b),(e) (or Example~\ref{ex:untwisted}), $(A_n,B_n)$ is a
		$C^*$-diagonal whose Weyl groupoid is $\mathcal G_n$ and whose twist is trivial. The complement of $\widehat R_n$ in $\mathcal G_n$ is the union of
		$R^F_n\setminus\widehat R_n$ and the interior levels, both open, so $\widehat R_n$ and $W_n^\theta$ are closed; and $\Upsilon(f,a)$ restricted to $\widehat R_n$ is
		$a|_{\widehat R_n}=q_n(a)=\bar q_n(f,a)$.
	\end{proof}
	
	\subsection{Normalizers}\label{ss:normalizers}
	
	\begin{lemma}\label{lem:normalizers}
		\begin{enumerate}[label=\textup{(\alph*)}]
			\item The normalizers of $B_n$ in $A_n$ are exactly the monomial elements.
			\item The set $N_{+,n}$ is closed under products and adjoints, its linear span is dense in $A_n$, and $\EE_n(s^*t)\ge0$ for $s,t\in N_{+,n}$.
		\end{enumerate}
		Thus $S_n=N_{+,n}$ satisfies \textup{(S1)--(S3)} of Proposition~\ref{prop:twist}.
	\end{lemma}
	
	\begin{proof}
		(a) If $m$ is monomial and $b\in B_n$, every fibre matrix of $mbm^*$ is $\mu\beta\mu^*$ with $\mu$ monomial and $\beta$ diagonal (in the same basis), which is
		diagonal by Lemma~\ref{lem:monomialprops}(b); similarly for $m^*bm$. Conversely, let $m$ be a normalizer. Consider a fibre which is not an endpoint fibre: a
		point $z\in Z_n$ (sheet basis), a matrix summand $l\ne\ib_n$, or an interior interval point $t\in(0,1)$ of a summand $q$. Its coordinates are distinct points of the
		unit space $\Omega_n$ of $\mathcal G_n$. For a coordinate $k$, choose by Urysohn's lemma $b\in C_0(\Omega_n)=B_n$ equal to $1$ at the unit $k$ and $0$ at the other
		units of this fibre; its fibre matrix there is $e_{kk}$. The fibre matrix of $mbm^*$ there is $\mu e_{kk}\mu^*$, which must be diagonal, so the $k$th column of
		$\mu$ has at most one nonzero entry; using $m^*bm$, the same holds for the $k$th row. By Lemma~\ref{lem:monomialprops}(c), $\mu$ is monomial. The endpoint fibre
		matrices are $f_q(i)=\beta_{n,i}(q_n(a))_q$, labelled block images of fibre matrices of $a$ at $\theta_n$, hence monomial too (Lemma~\ref{lem:labelledprops}(a)). (This
		also follows from Theorem~\ref{thm:renault}(d): normalizers are the elements whose open support is a bisection, and a function on $\mathcal G_n$ is supported in a
		bisection iff all its fibre matrices are monomial.)
		
		(b) Products and adjoints are computed fibrewise, and products and adjoints of nonnegative monomial matrices are nonnegative monomial
		(Lemma~\ref{lem:monomialprops}(a)). For $s,t\in N_{+,n}$ all entries of the fibre matrices of $s^*t$ are nonnegative, so its diagonal compression is a nonnegative
		function on $\Omega_n$: $\EE_n(s^*t)\ge0$.
		
		For density, fix $\chi_0,\chi_1\in C([0,1],[0,1])$ with disjoint supports, $\chi_0(0)=1$ and $\chi_1(1)=1$, and define the linear map
		\[
		L\colon F_n\to A_n,\qquad L(a)=\big(\chi_0\,\beta_{n,0}q_n(a)+\chi_1\,\beta_{n,1}q_n(a),\ a\big).
		\]
		$L(a)$ satisfies the endpoint conditions, and $L$ maps nonnegative monomial elements into $N_{+,n}$ (at each $t$ at most one of the two terms is nonzero). The
		nonnegative monomial elements of $F_n=C(R^F_n)$ span $F_n$: cover the compact groupoid $R^F_n$ by finitely many open bisections $U_1,\dots,U_k$ and choose a
		subordinate partition of unity $(\phi_j)$; for $a\in C(R^F_n)$, each $\phi_ja$ is supported in $U_j$, and it is a combination of four nonnegative functions supported in
		$U_j$ (positive and negative parts of real and imaginary parts). A nonnegative function supported in a bisection has nonnegative monomial fibre matrices, because a
		bisection contains at most one arrow with given range and at most one with given source. Next, the elements $\phi\,e_{xy}$ (in summand $q$) with $\phi\in C_c((0,1))$
		nonnegative lie in $N_{+,n}$, and their span is dense in $\mathcal I_n$. Finally, every $(f,a)\in A_n$ is $L(a)+\big((f,a)-L(a)\big)$, and the second term lies in
		$\mathcal I_n$. So the span of $N_{+,n}$ is dense.
	\end{proof}
	
	\begin{remark}[Model dependence]
		By Remark~\ref{rem:modeldep}, the set $N_{+,n}$ depends on the untwisted model, not only on the Cartan pair $(A_n,B_n)$. Here the model is fixed: sheet bases of
		the coverings (which are canonical), standard bases of the matrix summands and of the interval fibres. All connecting maps constructed below are compatible with
		these models, in the sense that they map nonnegative monomial elements to nonnegative monomial elements.
	\end{remark}
	%======================================================================
	\part*{Part II: The unital case}
	
	Throughout Part~II we are in Case (U): $u=[1_A]$, $\mathbf d_n=0$, the endpoint maps $\beta_{n,i}$ are unital, $s_n=\Lambda_{n,0}r_n=\Lambda_{n,1}r_n$, and the building
	blocks $A_n$ are unital with $[1_{A_n}]=(r_n,\zeta_n)$.
	
	%======================================================================
	\section{The unital case: exact connecting maps}\label{sec:unitalmaps}\label{ss:exactint}
	
	We now extend $\psi_n$ to a unital homomorphism $\varphi_n\colon A_n\to A_{n+1}$ which induces the maps $\gt_n$ of Section~\ref{sec:alg} and preserves the
	diagonals, the expectations and the nonnegative monomial elements.
	
	The image of $(f,a)$ must satisfy the endpoint conditions of $A_{n+1}$:
	its interval part at $t=i$ must equal $\beta_{n+1,i}(q_{n+1}(\psi_n(a)))=\beta_{n+1,i}(\widehat\psi_n(q_n(a)))$ by \eqref{eq:qpsi}. As a function of
	$x=q_n(a)$ the right-hand side is a labelled block homomorphism of type $\Lambda_{n+1,i}M_n$ (Lemma~\ref{lem:labelledprops}(c)). In the usual range
	construction the multiplicity matrices are strictly positive, and one may absorb discrepancies by spare point evaluations. Here $M_n$ has zero
	entries (the sparse extra rows), and we need an exact matrix identity instead.
	
	\subsubsection*{The matrix identity}
	Choose integer matrices $P_n,N_n\colon\Z^{l_n}\to\Z^{l_{n+1}}$ with strictly positive entries and $P_n-N_n=L_n$; adding the same strictly positive
	matrix to the positive and negative parts of $L_n$ achieves this. Define
	\[
	D_{n,0}=P_n\Lambda_{n,0}+N_n\Lambda_{n,1},\qquad D_{n,1}=P_n\Lambda_{n,1}+N_n\Lambda_{n,0}.
	\]
	Then, by \eqref{eq:quot} and \eqref{eq:Bn},
	\[
	D_{n,1}-D_{n,0}=(P_n-N_n)(\Lambda_{n,1}-\Lambda_{n,0})=L_n\Delta_n=\Delta_{n+1}M_n=(\Lambda_{n+1,1}-\Lambda_{n+1,0})M_n .
	\]
	Hence
	\begin{equation}\label{eq:Rn}
		R_n:=\Lambda_{n+1,0}M_n-D_{n,0}=\Lambda_{n+1,1}M_n-D_{n,1}.
	\end{equation}
	The matrices $\Lambda_{n+1,i}=\Delta_{n+1}^\mp+Q_{n+1}$ are at our disposal through the common part $Q_{n+1}$. Adding $c$ to every entry of $Q_{n+1}$ adds
	$c$ times the column sums of $M_n$ to the entries of $R_n$; every column of $M_n$ is nonzero (its original entries are positive), so for $c$
	large \emph{every entry of $R_n$ is strictly positive}. We choose $Q_{n+1}$ in this way (and subject to \eqref{eq:wsmall} below).
	
	\subsubsection*{The interval map}
	In the target interval summand $q\in[l_{n+1}]$ set, in standard block order,
	\begin{equation}\label{eq:Phi}
		\Phi_n(f,a)_q(t)=\diag\Big(\big(f_p(t)\big)^{\oplus(P_n)_{qp}},\ \big(f_p(1-t)\big)^{\oplus(N_n)_{qp}},\ \big(q_n(a)_l\big)^{\oplus(R_n)_{ql}}\Big)_{p\in[l_n],\,l\in J_n}.
	\end{equation}
	So the fibre contains $(P_n)_{qp}$ forward copies and $(N_n)_{qp}$ backward copies of each source interval summand $p$, and $(R_n)_{ql}$ constant
	copies of each basepoint matrix $q_n(a)_l$. The size is $\sum_p(P_n+N_n)_{qp}s_{n,p}+\sum_l(R_n)_{ql}r_{n,l}=(D_{n,0}r_n+R_nr_n)_q=(\Lambda_{n+1,0}M_nr_n)_q=s_{n+1,q}$.
	
	At $t=0$ the right side depends only on $x=q_n(a)$, because $f_p(0)=\beta_{n,0}(x)_p$ and $f_p(1)=\beta_{n,1}(x)_p$; as a function of $x$ it is a labelled
	block homomorphism (Lemma~\ref{lem:labelledprops}(c)) of type
	\[
	P_n\Lambda_{n,0}+N_n\Lambda_{n,1}+R_n=D_{n,0}+R_n=\Lambda_{n+1,0}M_n
	\]
	in row $q$. Similarly at $t=1$ the type is $D_{n,1}+R_n=\Lambda_{n+1,1}M_n$. On the other hand, let $\beta^{\mathrm{std}}_{q,i}$ be the labelled block homomorphism
	of $\widehat F_{n+1}$ of type $((\Lambda_{n+1,i})_{q,l'})_{l'}$ with the standard labelling; then $\beta^{\mathrm{std}}_{q,i}\circ\widehat\psi_n$ is a labelled block
	homomorphism of $\widehat F_n$ of the same type $(\Lambda_{n+1,i}M_n)_q$ (Lemma~\ref{lem:psi} and Lemma~\ref{lem:labelledprops}(c)). By
	Lemma~\ref{lem:labelledprops}(b) the two differ by conjugation with a permutation matrix $\sigma_{q,i}$:
	\[
	\Phi_n(f,a)_q(i)=\Ad(\sigma_{q,i})\big(\beta^{\mathrm{std}}_{q,i}(\widehat\psi_n(q_n(a)))\big).
	\]
	\emph{We define} the $q$th component of $\beta_{n+1,i}$ to be $\Ad(\sigma_{q,i})\circ\beta^{\mathrm{std}}_{q,i}$, again a labelled block homomorphism of type
	$((\Lambda_{n+1,i})_{q,l'})_{l'}$. This fixes the labellings $\iota^q_{n+1,i}$ of stage $n+1$. With this choice and \eqref{eq:qpsi},
	\[
	\Phi_n(f,a)(i)=\beta_{n+1,i}\widehat\psi_n q_n(a)=\beta_{n+1,i}q_{n+1}\psi_n(a),\qquad i=0,1 .
	\]
	Therefore
	\begin{equation}\label{eq:phin}
		\varphi_n(f,a)=\big(\Phi_n(f,a),\psi_n(a)\big)
	\end{equation}
	lies in $A_{n+1}$, and $\varphi_n\colon A_n\to A_{n+1}$ is a unital $*$-homomorphism.
	
	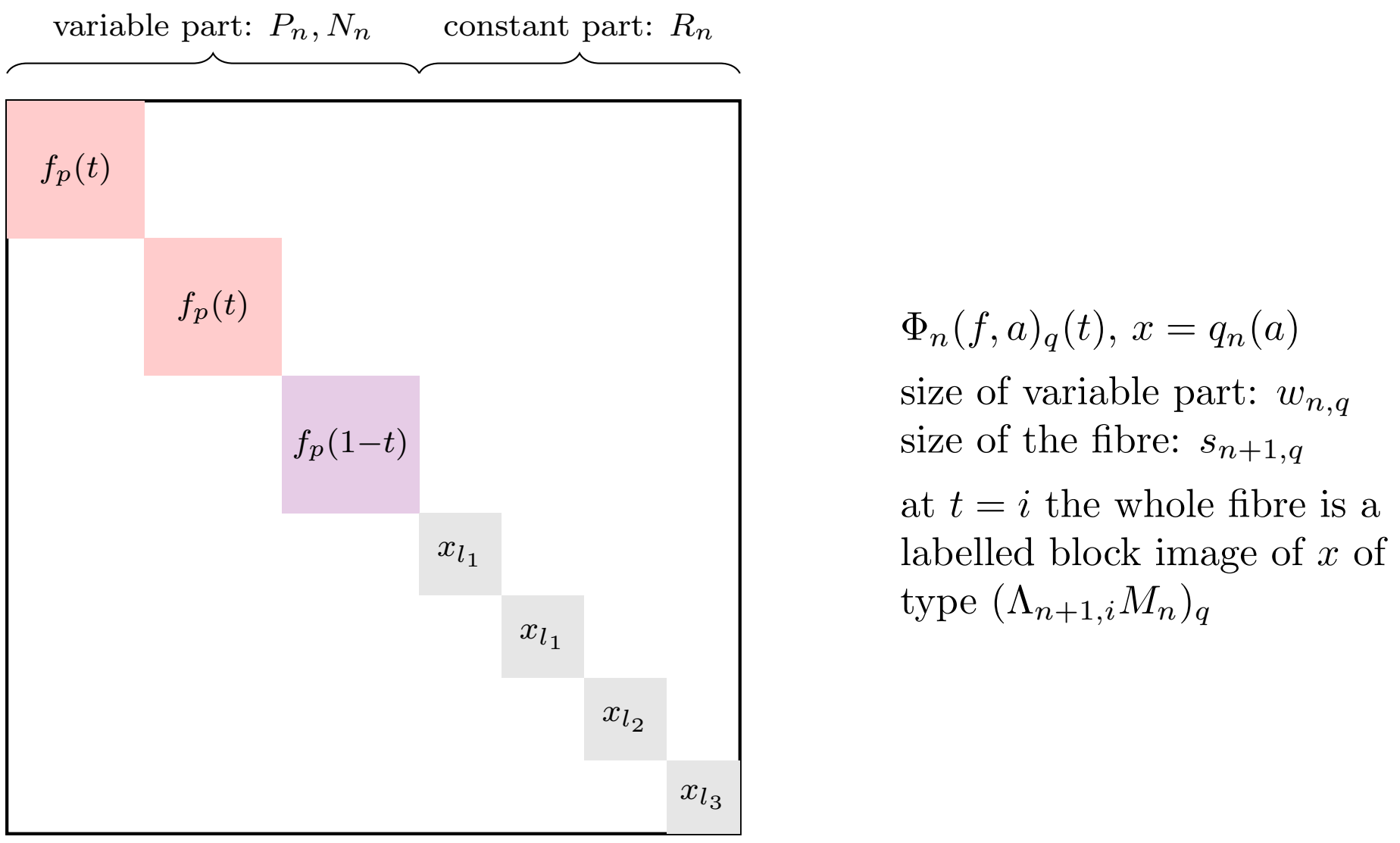
\begin{figure}[ht]
		\centering
		\begin{tikzpicture}[scale=0.8]
			\draw[thick] (0,0) rectangle (8,-8);
			\fill[red!20] (0,0) rectangle (1.5,-1.5); \node[font=\scriptsize] at (0.75,-0.75) {$f_p(t)$};
			\fill[red!20] (1.5,-1.5) rectangle (3,-3); \node[font=\scriptsize] at (2.25,-2.25) {$f_p(t)$};
			\fill[violet!20] (3,-3) rectangle (4.5,-4.5); \node[font=\scriptsize] at (3.75,-3.75) {$f_p(1{-}t)$};
			\fill[gray!20] (4.5,-4.5) rectangle (5.4,-5.4); \node[font=\scriptsize] at (4.95,-4.95) {$x_{l_1}$};
			\fill[gray!20] (5.4,-5.4) rectangle (6.3,-6.3); \node[font=\scriptsize] at (5.85,-5.85) {$x_{l_1}$};
			\fill[gray!20] (6.3,-6.3) rectangle (7.2,-7.2); \node[font=\scriptsize] at (6.75,-6.75) {$x_{l_2}$};
			\fill[gray!20] (7.2,-7.2) rectangle (8,-8); \node[font=\scriptsize] at (7.6,-7.6) {$x_{l_3}$};
			\draw[decorate,decoration={brace,amplitude=5pt}] (0,0.3) -- node[above=4pt,font=\scriptsize]{variable part: $P_n,N_n$} (4.5,0.3);
			\draw[decorate,decoration={brace,amplitude=5pt}] (4.5,0.3) -- node[above=4pt,font=\scriptsize]{constant part: $R_n$} (8,0.3);
			\node[font=\small,align=left] at (12.4,-4) {$\Phi_n(f,a)_q(t)$, $x=q_n(a)$\\[4pt]
				size of variable part: $w_{n,q}$\\ size of the fibre: $s_{n+1,q}$\\[4pt]
				at $t=i$ the whole fibre is a\\ labelled block image of $x$ of\\ type $(\Lambda_{n+1,i}M_n)_q$};
		\end{tikzpicture}
		\caption{The fibre of the interval part of $\varphi_n(f,a)$ in a target summand $q$ (schematic).}\label{fig:Phi}
	\end{figure}
	
	\begin{lemma}\label{lem:phin}
		The map $\varphi_n$ is injective, maps diagonal, monomial and nonnegative monomial elements to elements of the same kind, and satisfies
		$\EE_{n+1}\varphi_n=\varphi_n\EE_n$. In particular it maps $B_n$ into $B_{n+1}$, $\Nor(B_n)$ into $\Nor(B_{n+1})$ and $N_{+,n}$ into $N_{+,n+1}$. On
		$K$-theory it induces $\gt_n|_{G_n}\colon G_n\to G_{n+1}$ and $\chi_n\colon K_n\to K_{n+1}$.
	\end{lemma}
	
	\begin{proof}
		\emph{Injectivity.} If $\varphi_n(f,a)=0$, then $\psi_n(a)=0$, so $a=0$ (Lemma~\ref{lem:psi}); and then the forward copies in \eqref{eq:Phi} give
		$f_p(t)=0$ for all $p$ and $t$, because every entry of $P_n$ is positive.
		
		\emph{Fibre matrices.} The fibre matrices of $\varphi_n(f,a)$ are: those of $\psi_n(a)$, which are block-diagonal sums of fibre matrices of $a$ in the
		sheet bases (Lemma~\ref{lem:psi}); and, in the interval summands, block-diagonal sums of fibre matrices $f_p(t)$, $f_p(1-t)$ and $q_n(a)_l$ of $(f,a)$, in
		standard block order. By Lemma~\ref{lem:monomialprops}(d), diagonal, monomial and nonnegative monomial elements are preserved, and diagonal compression
		commutes with $\varphi_n$. Since $\EE_n$ is fibrewise diagonal compression, $\EE_{n+1}\varphi_n=\varphi_n\EE_n$. Lemma~\ref{lem:normalizers}(a) converts the
		statement on monomial elements into the statement on normalizers.
		
		\emph{$K$-theory.} The quotient maps $A_n\to F_n$ are injective on $K_0$ and isomorphisms on $K_1$ (Lemma~\ref{lem:An}(d)), and $\varphi_n$
		covers $\psi_n$; so Lemma~\ref{lem:psi} gives the claim.
	\end{proof}
	
	\subsubsection*{Smallness of the variable part}
	The variable part of the $q$th interval fibre has size
	\[
	w_{n,q}=\sum_p(P_n+N_n)_{qp}\,s_{n,p}.
	\]
	It is determined by $P_n,N_n$ and $s_n$, which are fixed before the common part $Q_{n+1}$ is enlarged. Enlarging $Q_{n+1}$ by $c$ increases
	$s_{n+1,q}=(\Lambda_{n+1,0}r_{n+1})_q$ by $c\sum_{l'}r_{n+1,l'}$, without bound. Consequently, for any prescribed $\varepsilon_n>0$ we may arrange
	\begin{equation}\label{eq:wsmall}
		w_{n,q}/s_{n+1,q}\le\varepsilon_n\qquad\text{for every }q .
	\end{equation}
	At the homogeneous step we also require
	\begin{equation}\label{eq:msmall}
		2/m_n\le\varepsilon_n .
	\end{equation}
	These two estimates control the traces (Section~\ref{sec:traces}): they say that, in every target fibre, the part of $\varphi_n(f,a)$ which does not
	factor through the basepoint values $q_n(a)$ is small. The strictly positive spare blocks will control simplicity (Section~\ref{sec:simple}).
	
	%======================================================================
	\section{The unital case: making the limit simple}\label{sec:simple}
	
	\subsection{The need for a modification}\label{ss:whymodify}
	The limit of the maps $\varphi_n$ is never simple. Indeed, let
	\[
	\mathfrak J_n=\{(0,a)\in A_n:\ a_l=0\text{ for }l\ne\ib_n,\ a_{\ib_n}(\theta_n)=0\},
	\]
	the elements supported on the distinguished summand and vanishing at the basepoint. It is a nonzero closed ideal of $A_n$ (the endpoint
	conditions hold because $q_n(a)=0$). If $(0,a)\in \mathfrak J_n$, then $q_n(a)=0$, so $\Phi_n(0,a)=0$ and $\psi_n(a)_{l'}=0$ for $l'\ne\ib_{n+1}$, while
	$\psi_n(a)_{\ib_{n+1}}$ consists of the blocks $a_{\ib_n}\circ\kappa_n$, $a_{\ib_n}\circ h_n$ and $a_{\ib_n}(\theta_n)=0$, which all vanish at
	$\theta_{n+1}$ because $\kappa_n$ and $h_n$ are based. Hence $\varphi_n(\mathfrak J_n)\subset \mathfrak J_{n+1}$, and the closure of $\bigcup_n\varphi_{n,\infty}(\mathfrak J_n)$ is a
	nonzero ideal of $\lims(A_n,\varphi_n)$. It is proper, because the quotients $A_n/\mathfrak J_n$ form a compatible system of nonzero unital $C^*$-algebras. The
	reason is that the matrix and interval targets only see the basepoint values $q_n(a)$.
	
	We therefore modify the maps by adding finitely many \emph{point evaluations away from the basepoint} at every step, following the mechanism of
	the unital range construction \cite[13.39--13.43]{GLN}: the connecting maps are changed within their homotopy classes, so that $K$-theory is
	unchanged, and by a small amount on traces. We give the modification explicitly, because two properties are needed later that are not part of
	\cite{GLN}: the modified maps must preserve nonnegative monomial elements (for the twist), and their effect on traces must be controlled in the
	dual norm (for Section~\ref{sec:traces}).
	
	\subsection{Irreducible representations and evaluation paths}\label{ss:paths}
	
	Recall from Lemma~\ref{lem:An}(b) that the irreducible representations of $A_n$ are the homogeneous evaluations $\ev_z$ ($z\in Z_n$), the matrix summands
	$(f,a)\mapsto a_l$ ($l\neq\ib_n$) and the interior interval evaluations $\ev_{q,t}$ ($t\in(0,1)$); in Case (U) the endpoint evaluations are unital block sums of the
	basepoint representations.
	
\begin{figure}[ht]
		\centering
		\begin{tikzpicture}[scale=0.95,>={Stealth[length=2mm]}]
			\draw[thick,fill=blue!6] (0,0) ellipse (2.2 and 1.3); \node at (-1.3,0.7) {$Z_n$};
			\fill (-0.8,-0.3) circle (2pt) node[below]{$\theta_n$};
			\fill (1.3,0.4) circle (2pt) node[above]{$z$};
			\draw[->,thick,blue!70!black] (-0.8,-0.3) .. controls (0,0.9) and (0.6,-0.6) .. (1.3,0.4);
			\node[blue!70!black,font=\small] at (0.3,-0.75) {$z(u)$};
			\node[font=\small] at (0,-1.8) {type (a): transport the labels along $z(u)$};
			\begin{scope}[xshift=7cm]
				\draw[thick] (-2,0) -- (2,0);
				\fill (-2,0) circle (2pt) node[below]{$0$}; \fill (2,0) circle (2pt) node[below]{$1$};
				\fill (0.9,0) circle (2pt) node[above]{$s$};
				\draw[->,thick,red!70!black] (-2,0.25) -- node[above,font=\small]{$us$} (0.9,0.25);
				\node[font=\small,align=center] at (0,-1.0) {type (c): $u\mapsto v_p^*f_p(us)v_p$;\\ at $u=0$: $\beta^{\mathrm{std}}_p(q_n(\cdot))$};
			\end{scope}
		\end{tikzpicture}
		\caption{Evaluation paths: each selected irreducible representation is connected, through representations which preserve the canonical bases up to
			reordering, to a standard block sum of basepoint representations.}\label{fig:paths}
	\end{figure}
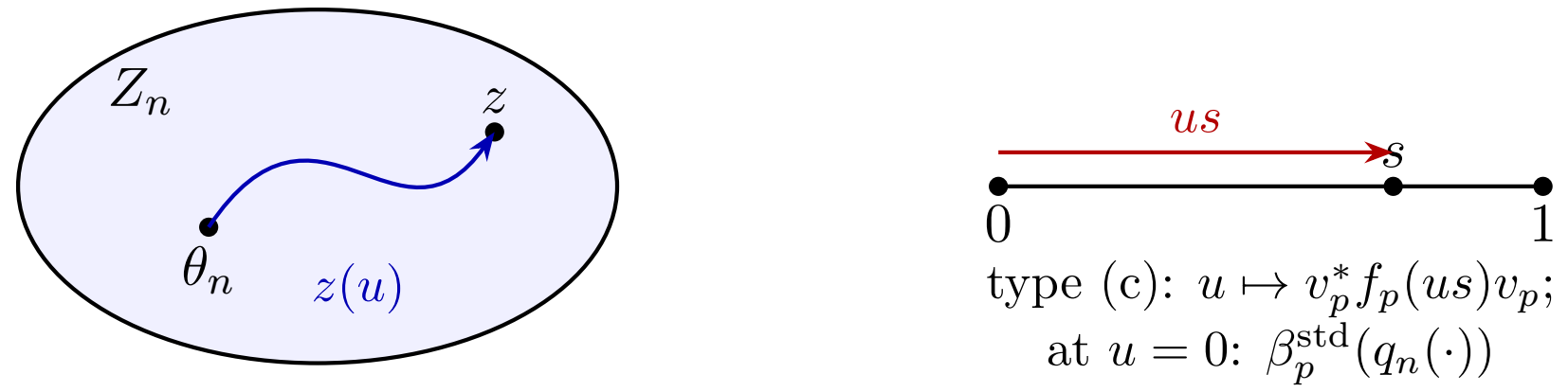

	At a completed stage $n$ we choose finitely many nonzero positive elements of $A_n$ to be ``detected'' by the next map (the \emph{tests}; their
	choice is specified in Proposition~\ref{prop:simple}), and for each test an irreducible representation in which it is nonzero. Call the selected
	representations $\rho_1,\dots,\rho_k$. To each $\rho_j$ we attach a continuous path of unital $*$-homomorphisms
	\[
	\rho_j^u\colon A_n\to M_{\dim\rho_j},\qquad u\in[0,1],
	\]
	with $\rho_j^1\cong\rho_j$ and $\rho_j^0$ a standard block sum of basepoint representations; continuity means that $u\mapsto\rho_j^u(x)$ is
	norm-continuous for each $x$. The paths are as follows (Figure~\ref{fig:paths}).
	\begin{itemize}
		\item \emph{Type (a), $\rho_j=\ev_z$.} Choose a path $z(u)$ in $Z_n$ from $z(0)=\theta_n$ to $z(1)=z$ (Lemma~\ref{lem:Zn}(f)). Each labelled point
		$y_k\in\lambda_n^{-1}(\theta_n)$, $k\in[r_{n,\ib_n}]$, has a unique lift $y_k(u)$ of $z(u)$ with $y_k(0)=y_k$ (unique path lifting). Put
		\[
		\rho_j^u(f,a)=\big(a_{\ib_n}(z(u))_{y_k(u),\,y_{k'}(u)}\big)_{k,k'\in[r_{n,\ib_n}]},
		\]
		the matrix of $a_{\ib_n}(z(u))$ in the transported basis. Its entries are values of the continuous function $k_a$ on $R_{\lambda_n}$
		(Proposition~\ref{prop:coveringalgebra}) along continuous paths, so $\rho_j^u$ is continuous in $u$. Then $\rho_j^0=q_n(\cdot)_{\ib_n}$, and $\rho_j^1$ is
		$\ev_z$ written in the basis $(y_k(1))_k$; we take this ordering when selecting $\rho_j$.
		\item \emph{Type (b), $\rho_j=(f,a)\mapsto a_l$.} The constant path; $\rho_j^0=\rho_j=q_n(\cdot)_l$.
		\item \emph{Type (c), $\rho_j=\ev_{p,s}$.} By Lemma~\ref{lem:labelledprops}(b) there is a permutation matrix $v_p$ with
		$(\beta_{n,0})_p=\Ad(v_p)\circ\beta^{\mathrm{std}}_p$, where $\beta^{\mathrm{std}}_p$ has the same type $((\Lambda_{n,0})_{p,l})_l$ and the standard labelling. Put
		\[
		\rho_j^u(f,a)=v_p^*\,f_p(us)\,v_p .
		\]
		Then $\rho_j^0(f,a)=v_p^*\beta_{n,0}(q_n(a))_pv_p=\beta^{\mathrm{std}}_p(q_n(a))$, a standard block sum of basepoint representations, and $\rho_j^1=\Ad(v_p^*)\circ\ev_{p,s}$
		is unitarily equivalent to $\rho_j$. (Without the conjugation by $v_p$, the path would end at the permuted block sum $\beta_{n,0}(q_n(\cdot))_p$, which
		does not match the standard block order used below; this matching is what makes $H_n(0)=\psi_n$ hold exactly.)
	\end{itemize}
	Each $\rho^u_j$ maps diagonal, monomial and nonnegative monomial elements to matrices of the same kind and commutes with diagonal compression: in
	all three cases $\rho^u_j(x)$ is a fibre matrix of $x$, written in a reordering of its canonical basis (transport of sheets, or conjugation by a
	permutation matrix). Each $\rho_j^0$ is a direct sum, in standard block order, of copies of the basepoint representations $x\mapsto q_n(x)_l$, with
	multiplicities
	\[
	\mu_l(\rho_j)=\begin{cases}\delta_{l,\ib_n}&\text{type (a)},\\ \delta_{l,l_0}&\text{type (b), }\rho_j=(\cdot)_{l_0},\\ (\Lambda_{n,0})_{p,l}&\text{type (c), }\rho_j=\ev_{p,s}.\end{cases}
	\]
	
	\subsubsection*{Reserving copies in one matrix target}
	Choose an original target coordinate $j_*\in I_{n+1}\setminus\{\ib_{n+1}\}$ (it exists because $|I_{n+1}|\ge2$). In this matrix summand,
	$\psi_n(\cdot)_{j_*}$ is the standard block homomorphism of type $((M_n)_{j_*,l})_l=((C_n)_{j_*,\pi_n(l)})_l$ composed with $q_n$. We require
	\begin{equation}\label{eq:reserve}
		(C_n)_{j_*,\pi_n(l)}\ge\sum_{j=1}^k\mu_l(\rho_j)+1\qquad\text{for every }l\in J_n .
	\end{equation}
	Then the copies of $q_n(\cdot)_l$ in $\psi_n(\cdot)_{j_*}$ suffice to accommodate all the $\rho^0_j$, with at least one copy of each $q_n(\cdot)_l$ left
	over. By Lemma~\ref{lem:labelledprops}(b) there are a permutation matrix $w$ and a standard block homomorphism $\rho_{\mathrm{rest}}$ (containing
	at least one copy of every basepoint representation) with
	\[
	\psi_n(f,a)_{j_*}=\Ad(w)\big(\rho^0_1(f,a)\oplus\cdots\oplus\rho^0_k(f,a)\oplus\rho_{\mathrm{rest}}(f,a)\big),
	\]
	because both sides are labelled block homomorphisms of $q_n(a)$ of the same type. Define $H_n(u)\colon A_n\to F_{n+1}$, $u\in[0,1]$, to agree with
	$\psi_n$ (composed with $(f,a)\mapsto a$) in every summand other than $j_*$, and
	\begin{equation}\label{eq:Hn}
		H_n(u)(f,a)_{j_*}=\Ad(w)\big(\rho^u_1(f,a)\oplus\cdots\oplus\rho^u_k(f,a)\oplus\rho_{\mathrm{rest}}(f,a)\big).
	\end{equation}
	Then each $H_n(u)$ is a unital $*$-homomorphism, $u\mapsto H_n(u)(x)$ is continuous, $H_n(0)(f,a)=\psi_n(a)$ exactly, and $H_n(1)$ detects all the
	selected tests in the summand $j_*$, since $\rho_j^1\cong\rho_j$.
	
	The \emph{total changing dimension} $\bar d_n=\sum_j\dim\rho_j$ is finite and depends only on the source stage and the chosen tests. (It is not to be
	confused with the dimensions $d_n$ of Theorem~\ref{thm:realize}.) For a prescribed $\eta_n>0$ we also require
	\begin{equation}\label{eq:etasmall}
		\bar d_n/r_{n+1,j_*}\le\eta_n .
	\end{equation}
	Since $r_{n+1,j_*}=\sum_l(C_n)_{j_*,\pi_n(l)}r_{n,l}$, the requirements \eqref{eq:reserve} and \eqref{eq:etasmall} are lower bounds on the entries of $C_n$
	in terms of completed stage-$n$ data; they can be met by Theorem~\ref{thm:presentation}. In particular, the source endpoint multiplicities
	$(\Lambda_{n,0})_{p,l}$ entering $\mu_l(\rho_j)$ for an interval representation are already known when the bound is imposed.
	
	\subsection{Endpoint collars}\label{ss:collars}
	
	The homotopy $H_n$ changes the homogeneous component, and the interval part has to follow at its endpoints. We insert \emph{collars} at both ends,
	in which the endpoint values run through $H_n$. Fix $c\in(0,\frac12)$, for example $c=\frac14$, and define $\Xi_n(f,a)\in C([0,1],E'_{n+1})$ by
	\begin{equation}\label{eq:Psi}
		\Xi_n(f,a)(t)=
		\begin{cases}
			\beta_{n+1,0}\,q_{n+1}\,H_n(1-t/c)(f,a),&0\le t\le c,\\
			\Phi_n(f,a)\big((t-c)/(1-2c)\big),&c\le t\le1-c,\\
			\beta_{n+1,1}\,q_{n+1}\,H_n((t-1+c)/c)(f,a),&1-c\le t\le1.
		\end{cases}
	\end{equation}
	Set
	\begin{equation}\label{eq:xin}
		\xi_n=\big(\Xi_n,H_n(1)\big)\colon A_n\to A_{n+1},\qquad\xi_n(f,a)=\big(\Xi_n(f,a),H_n(1)(f,a)\big).
	\end{equation}
	
	\begin{figure}[ht]
		\centering
		\begin{tikzpicture}[scale=1.1,>={Stealth[length=2mm]}]
			\draw[thick] (0,0) -- (10,0);
			\foreach \x/\l in {0/0,2.5/c,7.5/{1-c},10/1}{\draw (\x,0.1) -- (\x,-0.1) node[below]{$\l$};}
			\fill[blue!10] (0,0.15) rectangle (2.5,0.9); \fill[red!8] (2.5,0.15) rectangle (7.5,0.9); \fill[blue!10] (7.5,0.15) rectangle (10,0.9);
			\node[font=\small] at (1.25,0.52) {$\beta_{n+1,0}q_{n+1}H_n(u)$};
			\node[font=\small] at (5,0.52) {$\Phi_n$, reparametrized};
			\node[font=\small] at (8.75,0.52) {$\beta_{n+1,1}q_{n+1}H_n(u)$};
			\draw[->] (0.1,1.2) -- node[above,font=\scriptsize]{$u$ from $1$ to $0$} (2.4,1.2);
			\draw[->] (7.6,1.2) -- node[above,font=\scriptsize]{$u$ from $0$ to $1$} (9.9,1.2);
			\node[font=\small,align=center] at (5,-1.0) {at $t=0,1$: $u=1$, matching the new homogeneous component $H_n(1)$;\\ at $t=c,1-c$: $u=0$, matching the exact endpoints of $\Phi_n$};
		\end{tikzpicture}
		\caption{The collar modification $\Xi_n$ of the interval part.}\label{fig:collar}
	\end{figure}
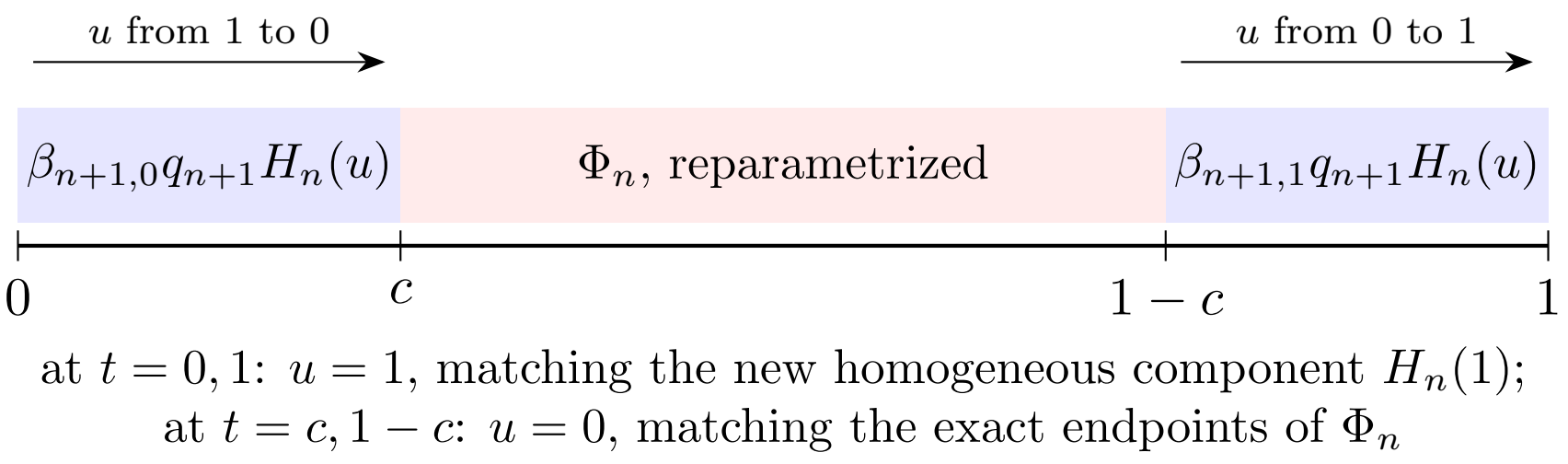
	
	\begin{lemma}\label{lem:xin}
		The map $\xi_n$ is a unital injective $*$-homomorphism $A_n\to A_{n+1}$, homotopic to $\varphi_n$ through unital $*$-homomorphisms. It maps diagonal,
		monomial and nonnegative monomial elements to elements of the same kind and satisfies $\EE_{n+1}\xi_n=\xi_n\EE_n$. In particular it maps $B_n$ into
		$B_{n+1}$, $\Nor(B_n)$ into $\Nor(B_{n+1})$ and $N_{+,n}$ into $N_{+,n+1}$, and it induces the same maps as $\varphi_n$ on $K_0$ and $K_1$.
	\end{lemma}
	
	\begin{proof}
		\emph{Well defined.} At $t=c$ the collar formula gives $\beta_{n+1,0}q_{n+1}H_n(0)(f,a)=\beta_{n+1,0}q_{n+1}\psi_n(a)$, and the middle formula gives
		$\Phi_n(f,a)(0)$; these agree by the exact endpoint identity of Section~\ref{ss:exactint}. The same holds at $t=1-c$ with $\beta_{n+1,1}$ and $\Phi_n(f,a)(1)$. So
		$\Xi_n(f,a)$ is continuous on $[0,1]$. At the endpoints, $\Xi_n(f,a)(0)=\beta_{n+1,0}q_{n+1}H_n(1)(f,a)$ and $\Xi_n(f,a)(1)=\beta_{n+1,1}q_{n+1}H_n(1)(f,a)$,
		which are exactly the endpoint conditions of $A_{n+1}$ for the homogeneous component $H_n(1)(f,a)$. So $\xi_n(f,a)\in A_{n+1}$.
		
		\emph{Homomorphism.} For each fixed $t$, $(f,a)\mapsto\Xi_n(f,a)(t)$ is a composite of unital $*$-homomorphisms ($H_n(u)$, $q_{n+1}$, $\beta_{n+1,i}$, or a fibre of
		$\Phi_n$), and so is $H_n(1)$. Operations in $A_{n+1}$ are pointwise, so $\xi_n$ is a unital $*$-homomorphism.
		
		\emph{Injective.} The distinguished summand of $H_n(1)$ is that of $\psi_n$ (as $j_*\ne\ib_{n+1}$), which is faithful on $F_n$ by Lemma~\ref{lem:psi}. So
		$\xi_n(f,a)=0$ implies $a=0$. Then the middle part of \eqref{eq:Psi} runs through all values $\Phi_n(f,0)(t')$, $t'\in[0,1]$, and the forward copies of
		$f_p(t')$ in \eqref{eq:Phi} give $f=0$ (every entry of $P_n$ is positive).
		
		\emph{Homotopy.} First, $\varphi_n$ is homotopic to the map $\varphi_n'$ with interval part $t\mapsto\Phi_n(f,a)(\varrho(t))$, where $\varrho\colon[0,1]\to[0,1]$
		is $0$ on $[0,c]$, linear on $[c,1-c]$ and $1$ on $[1-c,1]$: the maps with $\varrho_v=(1-v)\id+v\varrho$, $v\in[0,1]$, are homomorphisms into $A_{n+1}$ because
		$\varrho_v$ fixes $0$ and $1$. Next, for $u\in[0,1]$ consider the map with homogeneous component $H_n(u)$ and interval part given by \eqref{eq:Psi} with the
		collar parameters $1-t/c$ and $(t-1+c)/c$ replaced by $u(1-t/c)$ and $u(t-1+c)/c$; the middle part is unchanged. For every $u$ the joins match (the collar
		parameter is $0$ at $t=c,1-c$) and the endpoint conditions hold (the collar parameter is $u$ at $t=0,1$, matching the homogeneous component $H_n(u)$). At
		$u=0$ this is $\varphi_n'$ and at $u=1$ it is $\xi_n$. Homotopic homomorphisms induce the same maps on $K$-theory.
		
		\emph{Fibre matrices.} Every fibre matrix of $\xi_n(f,a)$ is, up to conjugation by a fixed permutation matrix, a block-diagonal sum of matrices obtained
		from fibre matrices of $(f,a)$ by reordering the canonical bases: in the homogeneous part this is Lemma~\ref{lem:psi} and the properties of $\rho_j^u$; in
		the middle part it is \eqref{eq:Phi}; in the collars it is the composite of $H_n(u)$ with $q_{n+1}$ and the labelled block homomorphism $\beta_{n+1,i}$
		(Lemma~\ref{lem:labelledprops}). By Lemma~\ref{lem:monomialprops}(d), diagonal, monomial and nonnegative monomial elements are preserved and diagonal
		compression commutes with $\xi_n$; as $\EE_n$ and $\EE_{n+1}$ are fibrewise diagonal compressions, $\EE_{n+1}\xi_n=\xi_n\EE_n$. Lemma~\ref{lem:normalizers}(a)
		gives the statement on normalizers.
	\end{proof}
	
	\begin{remark}
		In the collars, the target interval fibres near $t=0$ contain blocks $\rho^u_j$ which, for $u$ near $0$, are evaluations of the source near its
		basepoint or near the source interval endpoints. Thus the eigenvalue pattern of $\xi_n$ does not satisfy the restrictive conditions on connecting maps
		imposed in Li's construction \cite[\S\S4--5]{Li}, which are designed to compute the groupoid model of the limit directly. We do not need them:
		Proposition~\ref{prop:twist} requires only that nonnegative monomial elements be preserved, which is Lemma~\ref{lem:xin}.
	\end{remark}
	
	\subsection{Trace size of the modification}\label{ss:tracesize}
	
	\begin{definition}
		For a unital $C^*$-algebra $A$, $T(A)$ is the set of tracial states with the weak$^*$ topology, a compact convex set, and $\Aff T(A)$ is the real Banach
		space of continuous affine functions on it. For $x\in A_{sa}$ let $\hat x\in\Aff T(A)$, $\hat x(\tau)=\tau(x)$. A unital homomorphism $\phi\colon A\to A'$
		induces the continuous affine map $\phi^*\colon T(A')\to T(A)$, $\tau\mapsto\tau\circ\phi$, and the positive unital contraction
		\[
		\phi^\sharp\colon\Aff T(A)\to\Aff T(A'),\qquad\phi^\sharp(v)=v\circ\phi^* .
		\]
		On $T(A)$ we also use the \emph{dual norm} $\|\tau-\tau'\|=\sup\{|\tau(x)-\tau'(x)|:x\in A,\ \|x\|\le1\}$.
	\end{definition}
	
	\begin{lemma}\label{lem:affdense}
		Let $K$ be a compact convex subset of a locally convex space. A linear subspace of $\Aff(K)$ which contains the constants and separates the points of $K$
		is dense. In particular, for a unital $C^*$-algebra $A$ with $T(A)\ne\emptyset$, the functions $\hat x$, $x\in A_{sa}$, are dense in $\Aff T(A)$, and if
		$A=\overline{\bigcup_n A^{(n)}}$ for an increasing sequence of unital subalgebras, the $\hat x$ with $x\in\bigcup_nA^{(n)}_{sa}$ are dense.
	\end{lemma}
	
	\begin{proof}
		$\Aff(K)$ with the constant function $1$ is an order unit space, so every bounded linear functional on it is a difference of two positive ones, and the
		positive ones are nonnegative multiples of states; the states of $\Aff(K)$ are the point evaluations at points of $K$ \cite[Chapter~II]{Alf}. So a functional
		annihilating the subspace $L$ has the form $\alpha\,\ev_{k_1}-\beta\,\ev_{k_2}$ with $\alpha,\beta\ge0$; vanishing on $1$ gives $\alpha=\beta$, and vanishing on
		the separating subspace $L$ gives $k_1=k_2$ or $\alpha=0$. So it is zero, and $L$ is dense by the Hahn--Banach theorem. The functions $\hat x$ contain
		$\hat1=1$ and separate tracial states by definition; in the last case they separate because traces are continuous and the union is dense.
	\end{proof}
	
	\begin{lemma}[Traces of the building blocks]\label{lem:tracesblock}
		Every tracial state $\tau$ of $A_n$ has the form
		\[
		\tau(f,a)=\sum_q\int_{(0,1)}\tr_{s_{n,q}}(f_q(t))\,d\mu_q(t)+\int_{Z_n}\tr_{r_{n,\ib_n}}(a_{\ib_n}(z))\,d\nu(z)+\sum_{l\ne\ib_n}c_l\,\tr_{r_{n,l}}(a_l),
		\]
		with finite positive Borel measures $\mu_q$ on $(0,1)$ and $\nu$ on $Z_n$ and constants $c_l\ge0$, of total mass $1$; here $\tr_k$ is the normalized trace on
		$M_k$ (on a fibre of $V_{\lambda_n}$ it does not depend on the basis). In other words, every tracial state is an integral of normalized traces of irreducible
		representations.
	\end{lemma}
	
	\begin{proof}
		Let $I=C_0((0,1),E'_n)$. A positive tracial functional $\omega$ on $C_0(X,M_s)$ ($X$ locally compact) has the form $\omega(f)=\int\tr_s(f(t))\,d\mu(t)$: for
		$g\in C_0(X)_+$, traciality gives $\omega(g\otimes e_{ij})=\omega((\sqrt g\otimes e_{i1})(\sqrt g\otimes e_{1j}))=\omega((\sqrt g\otimes e_{1j})(\sqrt g\otimes
		e_{i1}))=\delta_{ij}\,\omega(g\otimes e_{11})$, so $\omega(f)=\omega_0(\Tr f)$ for the positive functional $\omega_0(g)=\omega(g\otimes e_{11})$ on $C_0(X)$, which
		is a measure by the Riesz representation theorem. Applying this to $\tau|_I$ gives measures $\mu_q$, and the formula $\tau_I(f,a)=\sum_q\int\tr(f_q(t))d\mu_q(t)$
		defines a positive tracial functional on $A_n$. If $(e_k)$ is an approximate unit of $I$ of the form $e_k=\chi_k\otimes1$ with $0\le\chi_k\uparrow1$ on $(0,1)$,
		then for $x\ge0$ in $A_n$, $\tau(x)\ge\tau(x^{1/2}e_kx^{1/2})=\tau(e_k^{1/2}xe_k^{1/2})\to\tau_I(x)$ by monotone convergence. So $\tau-\tau_I$ is a positive tracial
		functional, and it vanishes on $I$; hence it factors through $A_n/I\cong F_n$. On a matrix summand it is a multiple of $\tr$. On the covering algebra, let
		$\omega$ be a positive tracial functional and $\nu$ the measure $g\mapsto\omega(g\cdot1)$ on $Z_n$. Choose a finite cover of $Z_n$ by evenly covered open sets
		$U_j$ and a subordinate partition of unity $(\phi_j)$. Over $U_j$, $F_{\lambda_n}$ restricts to $C_0(U_j,M_r)$ in a local sheet trivialization, and the argument
		above gives $\omega(\phi_ja)=\int\tr(\phi_j(z)a(z))\,d\nu(z)$, since the measure for $\omega$ on $C_0(U_j,M_r)$ is $\nu|_{U_j}$ (evaluate on $g\cdot1$). Summing over $j$
		gives the formula. The total mass is $\tau(1)=1$.
	\end{proof}
	
	As a consequence, for a continuous affine function $v$ on $T(A_m)$ and a weak$^*$-continuous family of tracial states $\sigma_t$ indexed by a compact space
	with a probability measure $m$, the \emph{barycentre formula} $v\big(\int\sigma_t\,dm(t)\big)=\int v(\sigma_t)\,dm(t)$ holds: it holds for $v=\hat x$ by definition of
	the integral, and in general by Lemma~\ref{lem:affdense} and uniform approximation.
	
	\begin{lemma}[Comparison of block traces]\label{lem:blockcompare}
		Let $\pi,\pi'\colon A\to M_N$ be nondegenerate (for instance unital) $*$-homomorphisms such that, after conjugation by unitaries, $\pi=\pi_0\oplus\pi_1$ and
		$\pi'=\pi_0\oplus\pi_1'$ with $\pi_0$ of dimension $(1-w)N$ and $\pi_1,\pi'_1$ of dimension $wN$, $0\le w\le1$. Then $\tr_N\circ\pi$ and $\tr_N\circ\pi'$ are tracial states with
		$\|\tr_N\circ\pi-\tr_N\circ\pi'\|\le2w$ in the dual norm, and $|v(\tr_N\circ\pi)-v(\tr_N\circ\pi')|\le2w\|v\|$ for every $v\in\Aff T(A)$.
	\end{lemma}
	
	\begin{proof}
		$\tr_N\circ\pi=(1-w)\mu_0+w\mu$ and $\tr_N\circ\pi'=(1-w)\mu_0+w\mu'$ with $\mu_0=\tr\circ\pi_0$, $\mu=\tr\circ\pi_1$, $\mu'=\tr\circ\pi_1'$ tracial states (if $w\in\{0,1\}$
		one term is absent). The difference is $w(\mu-\mu')$, whose dual norm is at most $2w$; and $v$ is affine, so $v(\tr\circ\pi)-v(\tr\circ\pi')=w(v(\mu)-v(\mu'))$.
	\end{proof}
	
	\begin{lemma}\label{lem:tracesize}
		Under \eqref{eq:wsmall} and \eqref{eq:etasmall}:
		\begin{enumerate}[label=\textup{(\alph*)}]
			\item $\|\tau\circ\xi_n-\tau\circ\varphi_n\|\le2\eta_n+2\varepsilon_n$ in the dual norm of $A_n^*$, for every $\tau\in T(A_{n+1})$;
			\item $\|\xi_n^\sharp v-\varphi_n^\sharp v\|\le(2\eta_n+2\varepsilon_n)\|v\|$ for every $v\in\Aff T(A_n)$.
		\end{enumerate}
	\end{lemma}
	
	\begin{proof}
		\emph{Step 1: irreducible traces.} Let $\tau=\tr\circ\pi$ for an irreducible representation $\pi$ of $A_{n+1}$ (Lemma~\ref{lem:irreps}). We show that $\pi\circ\xi_n$ and
		$\pi\circ\varphi_n$ are connected by at most two comparisons as in Lemma~\ref{lem:blockcompare}, with $w\le\eta_n$ and $w\le\varepsilon_n$ respectively.
		\begin{itemize}
			\item \emph{Homogeneous evaluations and matrix summands other than $j_*$.} The homogeneous components of $\xi_n$ and $\varphi_n$ agree there, so
			$\pi\circ\xi_n=\pi\circ\varphi_n$.
			\item \emph{The summand $j_*$.} By \eqref{eq:Hn}, $\pi\circ\xi_n$ and $\pi\circ\varphi_n$ have the common part $\rho_{\mathrm{rest}}$ and differ on blocks of total
			dimension $\bar d_n$; so $w=\bar d_n/r_{n+1,j_*}\le\eta_n$ by \eqref{eq:etasmall}.
			\item \emph{Interval points $t\in[c,1-c]$ of a summand $q$.} $\pi\circ\xi_n$ is the fibre of $\Phi_n$ at $t'=(t-c)/(1-2c)$ and $\pi\circ\varphi_n$ is the fibre at $t$.
			Both contain the same constant part (the $R_n$-blocks of \eqref{eq:Phi}) and differ only in the variable part, of relative size $w=w_{n,q}/s_{n+1,q}\le\varepsilon_n$
			by \eqref{eq:wsmall}.
			\item \emph{Interval points in a collar,} say $t\in(0,c)$. Compare $\pi\circ\xi_n=\beta_{n+1,0}q_{n+1}H_n(u)(\cdot)_q$ first with
			$\beta_{n+1,0}q_{n+1}H_n(0)(\cdot)_q=\Phi_n(\cdot)_q(0)$. The arguments $q_{n+1}H_n(u)$ and $q_{n+1}H_n(0)$ differ only in the summand $j_*$, and there on a block
			of relative dimension at most $\eta_n$; the labelled block homomorphism $\beta_{n+1,0}(\cdot)_q$ consists of copies of the summands of $\widehat F_{n+1}$, so in the
			fibre the changed blocks have relative dimension at most $\eta_n$. Then compare $\Phi_n(\cdot)_q(0)$ with $\pi\circ\varphi_n=\Phi_n(\cdot)_q(t)$: they differ only in the
			variable part, relative size at most $\varepsilon_n$. The collar near $t=1$ is symmetric.
		\end{itemize}
		By Lemma~\ref{lem:blockcompare} and the triangle inequality,
		\[
		\|\tau\circ\xi_n-\tau\circ\varphi_n\|\le2\eta_n+2\varepsilon_n,\qquad|v(\tau\circ\xi_n)-v(\tau\circ\varphi_n)|\le(2\eta_n+2\varepsilon_n)\|v\|
		\]
		for irreducible traces $\tau$ (the intermediate traces in the collar case are tracial states of $A_n$).
		
		\emph{Step 2: all traces.} By Lemma~\ref{lem:tracesblock} (for $A_{n+1}$), a general $\tau$ is an integral $\int\tr\circ\pi_\omega\,dm(\omega)$ of irreducible traces,
		with $\omega\mapsto\tr\circ\pi_\omega$ weak$^*$-continuous on each of the finitely many pieces of the parameter space. For a contraction $x\in A_n$,
		$|\tau(\xi_nx)-\tau(\varphi_nx)|\le\int|\tr\pi_\omega(\xi_nx)-\tr\pi_\omega(\varphi_nx)|\,dm\le2\eta_n+2\varepsilon_n$, which proves (a). For (b), the barycentre formula gives
		$v(\tau\circ\xi_n)-v(\tau\circ\varphi_n)=\int\big(v(\tr\circ\pi_\omega\circ\xi_n)-v(\tr\circ\pi_\omega\circ\varphi_n)\big)\,dm(\omega)$, and Step~1 bounds the integrand.
	\end{proof}
	
	\begin{remark}
		(b) also follows from (a): by the Cuntz--Pedersen theorem \cite{CP}, $\Aff T(A_n)$ is isometrically the quotient of $(A_n)_{sa}$ by the closed span of
		the self-adjoint commutators, so every $v$ with $\|v\|<1$ is $\hat x$ for a self-adjoint contraction $x$, and $v(\tau\circ\xi_n)-v(\tau\circ\varphi_n)=\tau(\xi_nx)-\tau(\varphi_nx)$.
	\end{remark}
	
	\subsection{Simplicity}\label{ss:simplicity}
	
	\begin{proposition}\label{prop:simple}
		The tests and evaluation paths can be chosen so that $A_\infty=\lims(A_n,\xi_n)$ is simple.
	\end{proposition}
	
	\begin{proof}
		\emph{Choice of tests.} Each $A_m$ is separable; fix a countable norm-dense subset $\{b_{m,1},b_{m,2},\dots\}$ of its positive part, and enumerate as
		$c_{m,1},c_{m,2},\dots$ the nonzero elements among the \emph{rational cutdowns} $(b_{m,k}-\varrho)_+$, $\varrho\in\Q_{>0}$. At step $n$ (the construction of
		$\xi_n$), the tests are the elements
		\[
		\xi_{m,n}(c_{m,k}),\qquad1\le m\le n,\ 1\le k\le n ,
		\]
		of $A_n$, where $\xi_{m,n}=\xi_{n-1}\circ\cdots\circ\xi_m$ and $\xi_{n,n}=\id$. They are finitely many, and nonzero because the maps are injective; for each one we
		select an irreducible representation of $A_n$ in which it is nonzero. Thus every test $c_{m,k}$ is selected (through its image) at every step $n\ge\max\{m,k\}$.
		It is important to use tests from all earlier stages: the density argument below starts from an arbitrary stage.
		
		\emph{Step 1: detection.} If $y\in A_n$ is a selected test, then $\xi_n(y)$ has nonzero component in the matrix summand $j_*$ of $A_{n+1}$: by \eqref{eq:Hn},
		$H_n(1)(y)_{j_*}$ contains the block $\rho_j^1(y)\cong\rho_j(y)\ne0$.
		
		\emph{Step 2: spreading from an original matrix summand.} Let $z\in A_{n+1}$ be positive with nonzero component $z_j$ in an original matrix summand
		$j\in I_{n+1}\setminus\{\ib_{n+1}\}$. Then $z'=\xi_{n+1}(z)\in A_{n+2}$ is nonzero in the following irreducible representations of $A_{n+2}$:
		\begin{enumerate}[label=(\roman*)]
			\item every original matrix summand $l''\in I_{n+2}\setminus\{\ib_{n+2}\}$: if $l''$ is not the modified summand, $\psi_{n+1}(z)_{l''}$ contains
			$(M_{n+1})_{l'',j}=(C_{n+1})_{l'',j}\ge2$ copies of $z_j$; in the modified summand, $\rho_{\mathrm{rest}}$ contains at least one copy of $z_j$ by \eqref{eq:reserve};
			\item every homogeneous evaluation $\ev_w$, $w\in Z_{n+2}$: the distinguished target contains $(C_{n+1})_{\ib_{n+2},j}-1\ge1$ constant copies of $z_j$ (entry (i)
			of Section~\ref{ss:homog}); in particular $z'$ is nonzero at the basepoint $\theta_{n+2}$;
			\item every interior interval evaluation: in the middle part, \eqref{eq:Phi} contains $(R_{n+1})_{q,j}\ge1$ copies of $z_j$; in the collars, the fibre
			$\beta_{n+2,i}q_{n+2}H_{n+1}(u)(z)_q$ contains at least one copy of every summand of $\widehat F_{n+2}$ (as $\Lambda_{n+2,i}$ is strictly positive), in particular of an original
			matrix summand $l''\ne\ib_{n+2}$, whose value contains copies of $z_j$ by (i) (for the modified summand, through $\rho_{\mathrm{rest}}$, which does not depend on $u$);
			\item the extra matrix summand $x_j\in\mathcal X_{n+2}$, which receives $(M_{n+1})_{x_j,j}=1$ copy of $z_j$.
		\end{enumerate}
		The other extra summands $x_i$, $i\ne j$, receive only the source coordinates $l$ with $\pi_{n+1}(l)=i$, and may not see $z$.
		
		\emph{Step 3: one more step.} Let $z'$ be as in Step~2 and $z''=\xi_{n+2}(z')\in A_{n+3}$. Then $z''$ is nonzero in every irreducible representation of $A_{n+3}$:
		original matrix summands, homogeneous evaluations and interval evaluations as in Step~2 (using an original coordinate $l\ne\ib_{n+2}$ of $I_{n+2}$, where $z'_l\ne0$
		by Step~2(i)); and an extra summand $x_i$, $i\in I_{n+2}$, receives the source coordinate $i$ itself, where $z'$ is nonzero: by Step~2(i) if $i\ne\ib_{n+2}$, and by
		Step~2(ii) at the basepoint if $i=\ib_{n+2}$ (the block is $q_{n+2}(z')_{\ib_{n+2}}=z'_{\ib_{n+2}}(\theta_{n+2})$).
		
		\emph{Step 4: full images.} By Steps~1--3, if $y\in A_n$ is a selected test at step $n$, then $\xi_{n,n+3}(y)$ is a positive element of $A_{n+3}$ which is nonzero in
		every irreducible representation. It is \emph{full}: a proper closed ideal of a $C^*$-algebra is contained in the kernel of some irreducible representation
		\cite[\S5.4]{Mur}, so the closed ideal generated by $\xi_{n,n+3}(y)$ is $A_{n+3}$. Hence the closed ideal of $A_\infty$ generated by $\xi_{n,\infty}(y)$ contains
		$\xi_{n+3,\infty}(A_{n+3})\ni1$, i.e.\ it is $A_\infty$.
		
		\emph{Step 5: every nonzero ideal contains the image of a test.} Let $\mathcal J\neq0$ be a closed ideal of $A_\infty$, and $y\in\mathcal J$ positive with $\|y\|=1$. Since
		$\bigcup_m\xi_{m,\infty}(A_m)$ is dense, there are $m$ and a self-adjoint $x\in A_m$ with $\|\xi_{m,\infty}(x)-y^{1/2}\|<\frac1{20}$; then $x^2\ge0$ and
		$\|\xi_{m,\infty}(x^2)-y\|<\frac18$. Choose $b=b_{m,k}$ in the dense family with $\|b-x^2\|<\frac18$, so $\|\xi_{m,\infty}(b)-y\|<\frac14$. Let $\varpi\colon A_\infty\to A_\infty/\mathcal J$
		be the quotient map. Then $\|\varpi(\xi_{m,\infty}(b))\|<\frac14$, so the spectrum of the positive element $\varpi(\xi_{m,\infty}(b))$ lies in $[0,\frac14)$, and
		\[
		\varpi\big(\xi_{m,\infty}((b-\tfrac12)_+)\big)=\big(\varpi(\xi_{m,\infty}(b))-\tfrac12\big)_+=0 .
		\]
		Thus $\xi_{m,\infty}((b-\frac12)_+)\in\mathcal J$. It is nonzero, because $\|b\|=\|\xi_{m,\infty}(b)\|>\frac34$ (the maps are injective, hence isometric). So
		$(b-\frac12)_+$ is one of the tests $c_{m,k'}$, and it is selected at step $n=\max\{m,k'\}$, through its image $\xi_{m,n}(c_{m,k'})$. By Step~4, the ideal generated by
		$\xi_{m,\infty}(c_{m,k'})$ is $A_\infty$, and it is contained in $\mathcal J$. So $\mathcal J=A_\infty$.
	\end{proof}
	
	\begin{remark}
		In Step~5 a one-line quotient argument takes the place of a cutdown lemma (as in \cite[Lemma~2.2]{KR}). The point of using the rational cutdowns of a
		dense family, rather than the dense family itself, is that $(b-\frac12)_+$ lies in the ideal generated by $y$ when $b$ is close to $y$, whereas $b$ itself need not.
	\end{remark}
	%======================================================================
	\section{The unital case: traces and the pairing}\label{sec:traces}
	
	The trace simplex of the limit has to be $T=T(A)$, with pairing $\rho$. We first identify the traces for the unmodified maps $\varphi_n$, by comparing
	with an AF algebra, and then transfer the result to the modified maps $\xi_n$ by a perturbation argument. Put
	\[
	A_\infty^{(0)}=\lims(A_n,\varphi_n),\qquad\widehat F=\lims(\widehat F_n,\widehat\psi_n),\qquad A_\infty=\lims(A_n,\xi_n).
	\]
	By \eqref{eq:qpsi} and the exact endpoint identities, the maps $\bar q_n\colon A_n\to\widehat F_n$, $(f,a)\mapsto q_n(a)$, satisfy
	\begin{equation}\label{eq:qbar}
		\bar q_{n+1}\circ\varphi_n=\widehat\psi_n\circ\bar q_n ,
	\end{equation}
	so they induce a surjective unital homomorphism $\bar q\colon A_\infty^{(0)}\to\widehat F$. (The modified maps $\xi_n$ need not satisfy \eqref{eq:qbar}, because $H_n(1)$
	changes the summand $j_*$; this is why the comparison is made for the $\varphi_n$ first.)
	
	\subsection{Traces on inductive limits and on AF algebras}\label{ss:tracelimits}
	
	\begin{lemma}\label{lem:Tlimit}
		Let $A=\lims(A_n,\phi_n)$ be an inductive limit of unital $C^*$-algebras with unital connecting maps.
		\begin{enumerate}[label=\textup{(\alph*)}]
			\item The map $\tau\mapsto(\tau\circ\phi_{n,\infty})_n$ is an affine homeomorphism of $T(A)$ onto the inverse limit
			$\limp(T(A_n),\phi_n^*)=\{(t_n)\in\prod_nT(A_n):\ t_n=t_{n+1}\circ\phi_n\}$, with the product topology.
			\item For $v\in\Aff T(A_n)$, the function $v\circ\phi_{n,\infty}^*$ on $T(A)$ has norm $\lim_{m\to\infty}\|\phi_{n,m}^\sharp v\|$ (a decreasing limit).
			\item The functions $v\circ\phi^*_{n,\infty}$, $v\in\Aff T(A_n)$, $n\ge1$, are dense in $\Aff T(A)$.
		\end{enumerate}
	\end{lemma}
	
	\begin{proof}
		(a) The map is affine and continuous. It is injective because $\bigcup_n\phi_{n,\infty}(A_n)$ is dense. If $(t_n)$ is compatible, the formula
		$\tau(\phi_{n,\infty}(x))=t_n(x)$ is well defined on the dense $*$-subalgebra $\bigcup_n\phi_{n,\infty}(A_n)$; it is a tracial positive functional of norm $1$ there, and it
		extends by continuity to a tracial state of $A$. So the map is onto, and a continuous bijection from a compact space to a Hausdorff space is a homeomorphism.
		
		(b) $\phi^\sharp_{n,m}$ is contractive, so the norms decrease. Let $K_m=\phi_{n,m}^*(T(A_m))\subset T(A_n)$; these compact sets decrease, and $\|\phi_{n,m}^\sharp
		v\|=\sup_{K_m}|v|$. By (a), $\phi_{n,\infty}^*(T(A))\subset\bigcap_mK_m$; conversely, a point of $\bigcap_mK_m$ is the $n$th coordinate of a compatible sequence, by a
		compactness argument (the sets of compatible finite sequences $(t_n,\dots,t_m)$ with prescribed $t_n$ are nonempty compact and form an inverse system). So
		$\phi_{n,\infty}^*(T(A))=\bigcap_mK_m$, and the supremum of $|v|$ over a decreasing intersection of compact sets is the limit of the suprema.
		
		(c) These functions contain the functions $\hat x$, $x\in\bigcup_n\phi_{n,\infty}(A_n)_{sa}$, which are dense by Lemma~\ref{lem:affdense}.
	\end{proof}
	
	\begin{lemma}[Traces on AF algebras]\label{lem:AFtraces}
		For the unital AF algebra $\widehat F$, the map $\sigma\mapsto\sigma_*$ is an affine homeomorphism from $T(\widehat F)$ onto the state space of $(K_0(\widehat F),
		K_0(\widehat F)^+,[1])$. With Lemma~\ref{lem:delayed} and Theorem~\ref{thm:Fdimgroup}, $(K_0(\widehat F),K_0(\widehat F)^+,[1])\cong(F,F^+,\bar u)$, and hence
		$T(\widehat F)\cong S(F)\cong T$.
	\end{lemma}
	
	\begin{proof}
		For $\widehat F_n=\bigoplus_lM_{r_{n,l}}$, a tracial state is $\sum_l\lambda_l\tr_{r_{n,l}}$ with $\lambda_l\ge0$, $\sum\lambda_l=1$, and its $K_0$-map is the state
		$h\mapsto\sum_l\lambda_lh_l/r_{n,l}$ of $(\Z^{J_n},\Z^{J_n}_{\ge0},r_n)$; every state of this ordered group has this form. So $T(\widehat F_n)\cong S(K_0(\widehat F_n))$, naturally
		for the connecting maps. Pass to inverse limits, using Lemma~\ref{lem:Tlimit}(a) and the corresponding (elementary) fact for states of inductive limits of ordered
		groups. $K_0(\widehat F)=\lims(\Z^{J_n},M_n)$ with the positive cone and unit $r_n$, since $\widehat\psi_n$ has type $M_n$ (Lemma~\ref{lem:psi}); Lemma~\ref{lem:delayed}
		identifies this with $(F,F^+)$, and $r_n\mapsto\bar u$ by Section~\ref{ss:unit}. Theorem~\ref{thm:Fdimgroup} identifies $S(F)$ with $T$.
	\end{proof}
	
	\subsection{Comparison with the AF system}\label{ss:AFcompare}
	
	Set $V_n=\Aff T(A_n)$ and $W_n=\Aff T(\widehat F_n)$. Let
	\[
	\mathsf e_n\colon V_n\to W_n,\qquad(\mathsf e_nv)(\sigma)=v(\sigma\circ\bar q_n),
	\]
	be restriction to the traces which factor through $\bar q_n$, and write $\alpha_n=\varphi_n^\sharp$ and $b_n=\widehat\psi_n^\sharp$. By \eqref{eq:qbar},
	\begin{equation}\label{eq:quotrel}
		\mathsf e_{n+1}\circ\alpha_n=b_n\circ \mathsf e_n .
	\end{equation}
	
	\begin{lemma}\label{lem:section}
		$W_n\cong\R^{J_n}$, and $\mathsf e_n$ is a surjective positive unital contraction with a positive unital linear section of norm one.
	\end{lemma}
	
	\begin{proof}
		$T(\widehat F_n)$ is the simplex with vertices $\tr_l=\tr_{r_{n,l}}\circ(\text{projection to the }l\text{th summand})$, so an affine function on it is determined by its values at the
		vertices, and $W_n\cong\R^{J_n}$ with the supremum norm. Let $e_l$ be the unit of the $l$th summand of $\widehat F_n$ and $\tilde e_l$ the corresponding central projection
		of $F_n$ (so $q_n(\tilde e_l)=e_l$). The element
		\[
		c_l=\big(t\mapsto(1-t)\beta_{n,0}(e_l)+t\beta_{n,1}(e_l),\ \tilde e_l\big)
		\]
		belongs to $A_n$ (the endpoint conditions hold since $q_n(\tilde e_l)=e_l$), is positive (a convex combination of two projections at every $t$), satisfies $\bar
		q_n(c_l)=e_l$, and $\sum_lc_l=1$. Define $\sigma_n\colon W_n\to V_n$ by $\sigma_n(\lambda)(\tau)=\sum_l\lambda_l\tau(c_l)$. Then $\sigma_n$ is linear, positive, unital
		($\sum_l\tau(c_l)=1$), and contractive ($|\sum\lambda_l\tau(c_l)|\le\max|\lambda_l|\sum_l\tau(c_l)$). At a vertex, $(\mathsf e_n\sigma_n(\lambda))(\tr_{l'})=\sum_l\lambda_l\tr_{l'}(e_l)=\lambda_{l'}$,
		so $\mathsf e_n\sigma_n=\id$. Since $\mathsf e_n$ is contractive, $\sigma_n$ is isometric, and $\mathsf e_n$ is onto.
	\end{proof}
	
	\begin{lemma}\label{lem:almostquot}
		Under \eqref{eq:wsmall} and \eqref{eq:msmall},
		\begin{equation}\label{eq:almostquot}
			\|\alpha_nv\|\le\|\mathsf e_nv\|+\varepsilon_n\|v\|\qquad(v\in V_n).
		\end{equation}
	\end{lemma}
	
	\begin{proof}
		\emph{Normalized traces of covering entries.} For a basic homomorphism \eqref{eq:basic} with $\deg\pi=r$ and $\deg\nu=d$, normalized fibre traces satisfy
		\begin{equation}\label{eq:covtrace}
			\tr_{rd}\big(\Phi(a)(x')\big)=\frac1d\sum_{w\in\nu^{-1}(x')}\tr_r\big(a(g(w))\big),
		\end{equation}
		because $\Phi(a)(x')$ is block-diagonal with $d$ blocks of size $r$. In particular, a matrix-source covering entry (entries (i) and (ii) of Section~\ref{ss:homog})
		contributes exactly the normalized trace of its source matrix, regardless of the monodromy and the degree of the covering.
		
		\emph{Irreducible traces.} Let $\tau=\tr\circ\pi$ for an irreducible representation $\pi$ of $A_{n+1}$. We show that $\tau\circ\varphi_n=(1-w)\,\sigma\circ\bar q_n+w\mu$ for
		tracial states $\sigma\in T(\widehat F_n)$, $\mu\in T(A_n)$ and some $w\le\varepsilon_n$.
		\begin{itemize}
			\item Matrix targets $l'\ne\ib_{n+1}$: $\pi\circ\varphi_n=\psi_n(\cdot)_{l'}$ factors through $\bar q_n$; $w=0$.
			\item Homogeneous evaluation $\ev_z$ of $A_{n+1}$: by \eqref{eq:covtrace}, the entries (i), (ii) and the $m_n-2$ basepoint sheets of (iii) contribute traces of the form
			$\sigma'\circ\bar q_n$. Only the two sheets $\kappa_n$ and $h_n$ do not factor through $\bar q_n$; their blocks have total dimension $2r_{n,\ib_n}$, while
			$r_{n+1,\ib_{n+1}}\ge m_nr_{n,\ib_n}$. So $w\le2/m_n\le\varepsilon_n$ by \eqref{eq:msmall}.
			\item Interval evaluation of a summand $q$ at $t\in(0,1)$: in \eqref{eq:Phi} the constant part factors through $\bar q_n$, and the variable part has relative size
			$w_{n,q}/s_{n+1,q}\le\varepsilon_n$ by \eqref{eq:wsmall}.
		\end{itemize}
		Then $(\alpha_nv)(\tau)=v(\tau\circ\varphi_n)=(1-w)\,v(\sigma\circ\bar q_n)+w\,v(\mu)=(1-w)(\mathsf e_nv)(\sigma)+wv(\mu)$, so $|(\alpha_nv)(\tau)|\le\|\mathsf e_nv\|+\varepsilon_n\|v\|$.
		
		\emph{All traces.} By Lemma~\ref{lem:tracesblock} and the barycentre formula, $(\alpha_nv)(\tau)$ is an average of the values at irreducible traces, so the bound holds for
		every $\tau\in T(A_{n+1})$.
	\end{proof}
	
	In words: up to $\varepsilon_n$, the map $\varphi_n$ on traces factors through the finite-dimensional quotient $\widehat F_n$. If $\varepsilon_n\to0$, the limit should therefore
	have the same traces as the AF algebra $\widehat F$; this is the content of the next proposition.
	
	\begin{proposition}\label{prop:tracesAF}
		If $\varepsilon_n\to0$, then $\bar q^*\colon T(\widehat F)\to T(A_\infty^{(0)})$, $\sigma\mapsto\sigma\circ\bar q$, is an affine homeomorphism. Hence $T(A_\infty^{(0)})\cong T(\widehat
		F)\cong T$, and under this identification the pairing of $T(A_\infty^{(0)})$ with $K_0(A_\infty^{(0)})=G_0$ is $\rho$.
	\end{proposition}
	
	\begin{proof}
		\emph{Injective and continuous.} $\bar q$ is surjective, so $\bar q^*$ is injective; it is affine and weak$^*$-continuous. Its image $K'$ is a compact convex subset of
		$T(A_\infty^{(0)})$.
		
		\emph{The image carries all norms.} Let $v\in V_n$, and let $\tilde v=v\circ\varphi_{n,\infty}^*\in\Aff T(A_\infty^{(0)})$. By Lemma~\ref{lem:Tlimit}(b) applied to both inductive systems
		and by \eqref{eq:quotrel},
		\[
		\|\tilde v\|=\lim_m\|\alpha_{n,m}v\|,\qquad\sup_{K'}|\tilde v|=\sup_{\sigma\in T(\widehat F)}|v(\sigma\circ\bar q\circ\varphi_{n,\infty})|=\lim_m\|b_{n,m}\mathsf e_nv\|=\lim_m\|\mathsf e_m\alpha_{n,m}v\|,
		\]
		where $\alpha_{n,m}=\alpha_{m-1}\cdots\alpha_n$ and $b_{n,m}$ similarly, and we used $\bar q\circ\varphi_{n,\infty}=\widehat\psi_{n,\infty}\circ\bar q_n$. Lemma~\ref{lem:almostquot} at
		stage $m$ gives $\|\alpha_{n,m+1}v\|\le\|\mathsf e_m\alpha_{n,m}v\|+\varepsilon_m\|\alpha_{n,m}v\|\le\|\mathsf e_m\alpha_{n,m}v\|+\varepsilon_m\|v\|$. Letting $m\to\infty$ and using $\varepsilon_m\to0$,
		\[
		\|\tilde v\|\le\sup_{K'}|\tilde v| ,
		\]
		and the reverse inequality is trivial. By Lemma~\ref{lem:Tlimit}(c) the functions $\tilde v$ are dense in $\Aff T(A_\infty^{(0)})$, so $\sup_{T(A_\infty^{(0)})}|a|=\sup_{K'}|a|$ for every
		$a\in\Aff T(A_\infty^{(0)})$.
		
		\emph{Surjective.} Suppose $\tau_0\in T(A_\infty^{(0)})\setminus K'$. By the Hahn--Banach separation theorem in the weak$^*$ topology there is $x\in (A_\infty^{(0)})_{sa}$ with
		$\tau_0(x)>M:=\sup_{K'}\tau(x)$. Let $m_0=\inf_{K'}\tau(x)$ and $a=\hat x-\frac{M+m_0}2$. Then $\sup_{K'}|a|=\frac{M-m_0}2<a(\tau_0)$, contradicting the previous step. So $K'=T(A_\infty^{(0)})$,
		and $\bar q^*$ is a continuous affine bijection of compact spaces, hence an affine homeomorphism. Lemma~\ref{lem:AFtraces} gives $T(\widehat F)\cong T$.
		
		\emph{The pairing.} Let $x\in K_0(A_\infty^{(0)})$ be represented at stage $n$ by $x_n\in K_0(A_n)=G_n\subset\Z^{J_n}\oplus E_n$. For $\tau=\sigma\circ\bar q$,
		$\tau_*(x)=\sigma_*(\bar q_*(x))$, and $\bar q_{n*}(x_n)$ is the free part $h$ of $x_n=(h,e)$ (Lemma~\ref{lem:KFn}). In the limit, $K_0(A_\infty^{(0)})=G_0\to K_0(\widehat F)=F$ is the
		quotient map $G_0\subset H\to H/R=F$, because $j_n(h,e)\equiv\mu_n(\pi_nh)$ modulo $R$ (Lemma~\ref{lem:delayed}). If $\sigma$ corresponds to $\tau'\in T$, then
		$\sigma_*=s_{\tau'}$, and $s_{\tau'}(\bar g)=\theta(g,0)(\tau')=\rho(g)(\tau')$ for $g\in G_0$, by the definition of $\theta$ in Section~\ref{ss:enlarge}. So the pairing is $\rho$.
	\end{proof}
	
	\begin{remark}
		One can also argue with order-unit spaces: Lemma~\ref{lem:almostquot} shows that the map $\mathsf e$ induced by the $\mathsf e_n$ between the completed limits of $(V_n,\alpha_n)$ and
		$(W_n,b_n)$ is isometric, Lemma~\ref{lem:section} shows that it has dense range, and a surjective unital isometry of complete order-unit spaces is an order
		isomorphism; dualizing gives the proposition. The argument above avoids the completion.
	\end{remark}
	
	\subsection{Summable perturbations of trace maps}\label{ss:perturb}
	
	The modified maps $\xi_n$ do not satisfy \eqref{eq:qbar}, but by Lemma~\ref{lem:tracesize} they differ from the $\varphi_n$ on traces by a summable amount. The
	following lemma turns this into an identification of inverse limits.
	
	\subsubsection*{The topologies}
	Trace spaces carry two natural topologies: the weak$^*$ topology, in which they are compact, and the dual-norm topology, in which the perturbation estimates
	hold. They are very different: the traces of $C([0,1])$ given by Dirac measures at distinct points are at dual-norm distance $2$, so $T(C[0,1])$ is not compact in
	the dual norm. The lemma uses compactness for the weak$^*$ topology and estimates for the dual norm. Note that a weak$^*$-compact set of states is norm
	closed (the dual-norm topology is finer).
	
	\begin{lemma}\label{lem:perturb}
		Let $S_n$ be weak$^*$-compact convex sets of states of $C^*$-algebras $A_n$, and let $a_n,b_n\colon S_{n+1}\to S_n$ be affine, weak$^*$-continuous maps which are
		contractions for the dual norms. Suppose
		\[
		\sup_{t\in S_{n+1}}\|a_n(t)-b_n(t)\|\le\omega_n,\qquad\sum_n\omega_n<\infty .
		\]
		Then the inverse limits $\limp(S_n,a_n)$ and $\limp(S_n,b_n)$ are affinely homeomorphic. If $v_n\colon S_n\to\R$ are dual-norm-continuous affine functions with
		$v_n\circ a_n=v_n\circ b_n=v_{n+1}$ for all $n$, the homeomorphism preserves the functions induced by the $v_n$ (that is, $v_n((St)_n)=v_n(t_n)$ in the notation below).
	\end{lemma}
	
	\begin{figure}[ht]
		\centering
		\begin{tikzcd}[column sep=3.5em,row sep=2.4em]
			\cdots\arrow[r,"b_{m}"] & S_{m}\arrow[r,"b_{m-1}"]\arrow[d,dashed] & \cdots\arrow[r,"b_n"] & S_n & \text{(the given $b$-compatible $t$)}\\
			\cdots\arrow[r,"a_{m}"'] & S_{m}\arrow[r,"a_{m-1}"'] & \cdots\arrow[r,"a_n"'] & S_n & (St)_n=\lim_ma_{n,m}(t_m)
		\end{tikzcd}
		\caption{The map $S$ of Lemma~\ref{lem:perturb}: push $t_m$ down along the $a$-maps and let $m\to\infty$.}\label{fig:perturb}
	\end{figure}
	
	\begin{proof}
		Write $a_{n,m}=a_n\circ\cdots\circ a_{m-1}\colon S_m\to S_n$ ($a_{n,n}=\id$), and similarly $b_{n,m}$.
		
		\emph{Construction of $S$.} Let $t=(t_m)$ be $b$-compatible, i.e.\ $t_m=b_m(t_{m+1})$. For fixed $n$ consider the sequence $a_{n,m}(t_m)$, $m\ge n$. Since
		$t_m=b_m(t_{m+1})$ and $a_{n,m}$ is a dual-norm contraction,
		\[
		\begin{aligned}
			\|a_{n,m+1}(t_{m+1})-a_{n,m}(t_m)\|&=\|a_{n,m}(a_m(t_{m+1}))-a_{n,m}(b_m(t_{m+1}))\|\\
			&\le\|a_m(t_{m+1})-b_m(t_{m+1})\|\le\omega_m .
		\end{aligned}
		\]
		As $\sum\omega_m<\infty$, the sequence is Cauchy in the dual norm, uniformly in $t$; let $(St)_n$ be its limit. It lies in $S_n$, because $S_n$ is norm closed. The
		sequence $St$ is $a$-compatible: $a_n((St)_{n+1})=\lim_ma_n(a_{n+1,m}(t_m))=(St)_n$, since $a_n$ is norm continuous (a contraction).
		
		\emph{Properties of $S$.} Each coordinate $t\mapsto(St)_n$ is a uniform (in the dual norm) limit of the weak$^*$-continuous maps $t\mapsto a_{n,m}(t_m)$; a uniform norm limit
		of weak$^*$-continuous maps into a bounded set of states is weak$^*$-continuous (weak$^*$ convergence at a fixed element is controlled by norm convergence). So $S$ is
		continuous for the product of the weak$^*$ topologies, and it is affine. Moreover
		\[
		\|(St)_m-t_m\|\le\sum_{k\ge m}\omega_k ,
		\]
		by summing the increments from $a_{m,m}(t_m)=t_m$.
		
		\emph{The inverse map.} Interchanging the roles of $a$ and $b$ gives $R\colon\limp(S_n,a_n)\to\limp(S_n,b_n)$. For a $b$-compatible $t$ and $n\le m$,
		\[
		\|b_{n,m}((St)_m)-t_n\|=\|b_{n,m}((St)_m)-b_{n,m}(t_m)\|\le\|(St)_m-t_m\|\le\sum_{k\ge m}\omega_k\longrightarrow0 ,
		\]
		so $(RSt)_n=\lim_mb_{n,m}((St)_m)=t_n$, i.e.\ $RS=\id$; symmetrically $SR=\id$. So $S$ is an affine homeomorphism (a continuous bijection between compact Hausdorff
		spaces).
		
		\emph{Invariant functions.} If $v_n\circ a_n=v_n\circ b_n=v_{n+1}$, then $v_n\circ a_{n,m}=v_m=v_n\circ b_{n,m}$, and by norm continuity of $v_n$,
		\[
		v_n((St)_n)=\lim_mv_n(a_{n,m}(t_m))=\lim_mv_m(t_m)=\lim_mv_n(b_{n,m}(t_m))=v_n(t_n).\qedhere
		\]
	\end{proof}
	
	\subsubsection*{Application}
	Choose, from now on, positive sequences $\varepsilon_n,\eta_n$ with
	\begin{equation}\label{eq:summable}
		\sum_n(\varepsilon_n+\eta_n)<\infty;\qquad\text{for example }\varepsilon_n=\eta_n=2^{-n-6}.
	\end{equation}
	Apply Lemma~\ref{lem:perturb} with $S_n=T(A_n)$, $a_n=\varphi_n^*$ and $b_n=\xi_n^*$. These maps are affine, weak$^*$-continuous and dual-norm contractions, and
	Lemma~\ref{lem:tracesize}(a) gives the hypothesis with $\omega_n=2\eta_n+2\varepsilon_n$. By Lemma~\ref{lem:Tlimit}(a), $T(A_\infty^{(0)})=\limp(T(A_n),\varphi_n^*)$ and
	$T(A_\infty)=\limp(T(A_n),\xi_n^*)$, so with Proposition~\ref{prop:tracesAF},
	\begin{equation}\label{eq:TA}
		T(A_\infty)\cong T(A_\infty^{(0)})\cong T .
	\end{equation}
	
	\emph{The pairing is preserved.} The maps $\xi_n$ and $\varphi_n$ are homotopic (Lemma~\ref{lem:xin}), so they induce the same maps on $K$-theory and $K_0(A_\infty)=K_0(A_\infty^{(0)})=G_0$
	compatibly. Let $x\in K_0(A_\infty)$ be represented by a projection $p\in M_k(A_n)$ (differences are handled by linearity). For $m\ge n$ define $v_m(t)=(t\otimes\Tr_k)(\varphi_{n,m}(p))$ on
	$T(A_m)$, an affine function with $|v_m(t)-v_m(s)|\le k\|t-s\|$, hence dual-norm continuous. Then $v_m\circ\varphi_m^*=v_{m+1}$ by definition, and $v_m\circ\xi_m^*=v_{m+1}$ because $\xi_m(\varphi_{n,m}(p))$ and
	$\varphi_m(\varphi_{n,m}(p))$ are homotopic, hence unitarily equivalent, projections, which have the same trace. Apply the last assertion of Lemma~\ref{lem:perturb} to the tail
	systems from stage $n$ on, which have the same inverse limits: for $\tau\in T(A_\infty)$ with coordinates $t_m=\tau\circ\xi_{m,\infty}$ and the corresponding $\tau'\in T(A_\infty^{(0)})$ with
	coordinates $(St)_m$,
	\[
	\tau_*(x)=v_n(t_n)=v_n((St)_n)=\tau'_*(x).
	\]
	Hence the identification \eqref{eq:TA} carries the pairing of $A_\infty$ to the pairing of $A_\infty^{(0)}$, which is $\rho$ by Proposition~\ref{prop:tracesAF}.
	%======================================================================
	\section{The unital case: completing the induction and the proof}\label{sec:proof}
	
	The goal of this section is the unital case of Theorem~\ref{thm:main}:
	
	\begin{theorem}[Unital case]\label{thm:unital}
		Let $A$ be a unital, separable, simple, nuclear, stably finite, $\JS$-stable $C^*$-algebra satisfying the UCT. There is a second countable, locally compact
		Hausdorff, principal étale groupoid $\mathcal G$ with compact unit space such that $A\cong C^*(\mathcal G)$.
	\end{theorem}
	
	We apply the construction of Sections~\ref{sec:alg}--\ref{sec:traces} to the invariant of $A$, in Case (U).
	
	\subsection{The order of the choices}\label{ss:choices}
	
	The preceding sections impose many conditions: lower bounds on the entries of the connecting matrices $C_n$, on the common parts $Q_{n+1}$ of the
	endpoint matrices, and on the dimensions $d_n$ of the torsion stages. Some of these quantities depend on each other: for example, the multiplicity $m_n$
	must dominate the covering degrees $\cc_n$, while the torsion groups $E_{n+1}$, and hence the complexes $T_{n+1}$ and their covering degrees, depend on the
	choice of $C_n$. In this subsection we show that all choices can be made in a single successive construction, each choice depending only on quantities
	already fixed. This is the purpose of the delay of Section~\ref{sec:alg}.
	
	\subsubsection*{The data of a completed stage}
	We say that \emph{stage $n$ is completed} when the following data have been fixed (Table~\ref{tab:stage}).
	
	\begin{table}[ht]
		\centering
		\small
		\begin{tabular}{>{\raggedright\arraybackslash}p{0.32\textwidth}>{\raggedright\arraybackslash}p{0.55\textwidth}}
			\toprule
			datum & where it is defined\\
			\midrule
			$I_n$, $\ib_n$, the lifts $\ell_n$, the defects $t_{n-1}$ & Sections~\ref{ss:positive}, \ref{ss:defects}\\
			$E_n\supset t_{n-1}(\Z^{I_{n-1}})$, its decomposition $\Theta_n$ & Section~\ref{ss:defects}, (T1)\\
			$J_n=I_n\sqcup\mathcal X_n$, $\pi_n$, $\delta_n$, the unit $(r_n,\zeta_n)$ & Sections~\ref{ss:delayed}, \ref{ss:unit}\\
			$G_n$, $\Delta_n\colon\Z^{J_n}\to\Z^{l_n}$ & Section~\ref{ss:unit}\\
			the torsion stage $\mathsf S_n$: $d_n$, $X_n$, $T_n$, $a_n$, $\phi_n$, generator coverings, $\cc_n=\cc(\mathsf S_n)$ & Definition~\ref{def:stage}, Corollary~\ref{cor:bound}\\
			$K_n$, $Y_n$; $Z_n=T_n\vee Y_n\vee[0,1]$; a surjection $p_n\colon[0,1]\to Z_n$ with $p_n(0)=\theta_n$ & Sections~\ref{ss:K1}, \ref{ss:Zn}, \ref{ss:homog}\\
			the labelled covering $\lambda_n$ of $Z_n$ with $\deg\lambda_n=r_{n,\ib_n}$, $\tc(\lambda_n)=\zeta_n$ & Section~\ref{ss:homog}\\
			$F_n$, $\widehat F_n$, $q_n$ & \eqref{eq:Fn}\\
			$\Lambda_{n,0},\Lambda_{n,1}$ with labellings $\iota^q_{n,i}$; $s_n$; $E_n'$; $\beta_{n,0},\beta_{n,1}$ & Section~\ref{ss:interval}\\
			the building block $A_n$, $B_n$, $\EE_n$, $N_{+,n}$ & \eqref{eq:An}\\
			the earlier maps $\xi_1,\dots,\xi_{n-1}$ and the dense families of positive elements of $A_1,\dots,A_n$ & Section~\ref{ss:simplicity}\\
			\bottomrule
		\end{tabular}
		\caption{The data of a completed stage $n$.}\label{tab:stage}
	\end{table}
	
	\subsubsection*{The step from stage $n$ to stage $n+1$}
	Suppose stage $n$ is completed. Proceed as follows.
	\begin{enumerate}[label=\textup{(\arabic*)}]
		\item \emph{Tests.} Choose the finitely many tests of step $n$ (Proposition~\ref{prop:simple}), their detecting irreducible representations
		$\rho_1,\dots,\rho_k$ of $A_n$, and the evaluation paths $\rho^u_j$ (Section~\ref{ss:paths}). This fixes the multiplicities $\mu_l(\rho_j)$ and the changing dimension
		$\bar d_n$.
		\item \emph{The connecting matrix $C_n$.} Using Theorem~\ref{thm:presentation}, choose the next stage of the presentation so that every entry of $C_n$ is at
		least a bound $b_n$ which guarantees:
		\begin{enumerate}[label=(\alph*)]
			\item all entries of $C_n$ are at least $2$ (so $M_n\ge0$ and its original rows are strictly positive);
			\item $m_n=(C_n)_{\ib_{n+1},\ib_n}-1\ge\max\{\cc_n,3\}$ (covering degrees of the extra sources, and room for $\kappa_n$, $h_n$);
			\item $2/m_n\le\varepsilon_n$, i.e.\ \eqref{eq:msmall};
			\item $(C_n)_{j_*,\pi_n(l)}\ge\sum_j\mu_l(\rho_j)+1$ for all $l\in J_n$, i.e.\ \eqref{eq:reserve};
			\item $\bar d_n\le\eta_n\,r_{n+1,j_*}$, i.e.\ \eqref{eq:etasmall}; since $r_{n+1,j_*}\ge\min_{i,j}(C_n)_{j,i}\cdot\sum_lr_{n,l}$, a bound on the entries suffices.
		\end{enumerate}
		Every quantity in (a)--(e) belongs to the completed stage $n$ or to step (1). All of them are lower bounds which we impose on \emph{every} entry of $C_n$ (in (d)
		and (e) for every original target row in place of $j_*$); so they do not depend on which coordinates of $I_{n+1}$ are afterwards designated as $\ib_{n+1}$ and $j_*$.
		Now designate $\ib_{n+1}\in I_{n+1}$ and choose $j_*\in I_{n+1}\setminus\{\ib_{n+1}\}$ (possible as $|I_{n+1}|\ge2$). This defines $I_{n+1}$, $C_n$, $M_n$, $M_n'$, $J_{n+1}$,
		and the next unit vector $r_{n+1}=M_nr_n$; the map $\gt_n$ is defined once $E_{n+1}$ is known.
		\item \emph{Lifts, defect, torsion group.} Choose lifts $\ell_{n+1}$ of the new stage, compute $t_n=\ell_n-\ell_{n+1}C_n$, and let
		$E_{n+1}=\langle E_n,t_n(\Z^{I_n}),\mathfrak r_1,\dots,\mathfrak r_{n+1}\rangle$ (Section~\ref{ss:defects}), with a decomposition $\Theta_{n+1}$. Now $\delta_{n+1}$, $\zeta_{n+1}$, $G_{n+1}$,
		$\Delta_{n+1}$ and $L_n$ are determined (Section~\ref{ss:unit}). Choose the next finitely generated group $K_{n+1}$ of the presentation of $G_1$, with the map
		$\chi_n$.
		\item \emph{Topology.} Choose $d_{n+1}\ge d_n$ by Theorem~\ref{thm:realize}(c) and apply Theorem~\ref{thm:realize}(b) to $\mathsf S_n$, $(E_{n+1},\Theta_{n+1})$ and $i_n$. This
		gives the torsion stage $\mathsf S_{n+1}$ (in particular $T_{n+1}$, its generator coverings and $\cc_{n+1}$) and the based map $g_n\colon T_{n+1}\to T_n$. Form $Y_{n+1}$,
		$k_n$, $Z_{n+1}$, $\kappa_n$, $h_n=p_n\circ\rho_{n+1}$ and a surjection $p_{n+1}\colon[0,1]\to Z_{n+1}$ (Section~\ref{ss:homog}). The new covering degrees, i.e.\ $\cc_{n+1}$,
		need satisfy no bound from step (2): the corrections $t_n$ which they realize are paid for in the \emph{next} step, by $C_{n+1}$.
		\item \emph{The homogeneous map.} For each $i\in I_{n-1}$, choose the covering $\nu_{n,i}$ of $T_n$ with $\tc(\nu_{n,i})=\phi_n(\delta_n(x_i))$ and
		$\deg\nu_{n,i}\le\cc_n\le m_n$ (Corollary~\ref{cor:bound}; trivial of degree one if $\delta_n(x_i)=0$). Construct $\psi_n$ and the labelled covering $\lambda_{n+1}$
		(Definition~\ref{def:lambda}); by Lemma~\ref{lem:lambdaclass} its degree and class are $r_{n+1,\ib_{n+1}}$ and $\zeta_{n+1}$. This defines $F_{n+1}$, $\widehat F_{n+1}$, $q_{n+1}$ and
		$\widehat\psi_n$. Construct the homotopy $H_n$ in the summand $j_*$ (Section~\ref{ss:paths}).
		\item \emph{Endpoint matrices.} Choose $P_n,N_n$ strictly positive with $P_n-N_n=L_n$, and then the common part $Q_{n+1}$ so large that $R_n>0$ entrywise and
		\eqref{eq:wsmall} holds; put $\Lambda_{n+1,0}=\Delta_{n+1}^-+Q_{n+1}$, $\Lambda_{n+1,1}=\Delta_{n+1}^++Q_{n+1}$, $s_{n+1}=\Lambda_{n+1,0}r_{n+1}$.
		\item \emph{Labellings and the maps.} Choose the permutations $\sigma_{q,i}$ and hence the labelled block homomorphisms $\beta_{n+1,i}$ so that the endpoint identities
		of Section~\ref{ss:exactint} hold. This defines the building block $A_{n+1}$, the exact map $\varphi_n$, and the collar modification $\xi_n$ (Section~\ref{ss:collars}). Fix
		a countable dense family of positive elements of $A_{n+1}$.
	\end{enumerate}
	After step (7), stage $n+1$ is completed.
	
	\begin{figure}[ht]
		\centering
		\begin{tikzpicture}[>={Stealth[length=2mm]},node distance=0.55cm and 0.8cm,
			box/.style={draw,rounded corners,fill=blue!5,align=center,font=\small,minimum height=0.8cm}]
			\node[box,fill=gray!10] (S) {completed stage $n$\\ $\cc_n$, $A_n$, $\lambda_n$, $p_n$};
			\node[box,right=of S] (t) {(1) tests,\\ $\mu_l(\rho_j)$, $\bar d_n$};
			\node[box,right=of t] (C) {(2) $C_n$\\ bounds (a)--(e)};
			\node[box,below=of C] (E) {(3) $\ell_{n+1}$, $t_n$,\\ $E_{n+1}$, $\Delta_{n+1}$, $L_n$};
			\node[box,left=of E] (T) {(4) $T_{n+1}$, $g_n$,\\ $\cc_{n+1}$, $Z_{n+1}$};
			\node[box,left=of T] (psi) {(5) $\nu_{n,i}$, $\psi_n$,\\ $\lambda_{n+1}$, $H_n$};
			\node[box,below=of psi] (Q) {(6) $P_n,N_n,Q_{n+1}$,\\ $\Lambda_{n+1,i}$};
			\node[box,right=of Q] (phi) {(7) $\beta_{n+1,i}$, $A_{n+1}$,\\ $\varphi_n$, $\xi_n$};
			\node[box,right=of phi,fill=gray!10] (S1) {completed\\ stage $n+1$};
			\draw[->] (S) -- (t); \draw[->] (t) -- (C); \draw[->] (C) -- (E); \draw[->] (E) -- (T); \draw[->] (T) -- (psi);
			\draw[->] (psi) -- (Q); \draw[->] (Q) -- (phi); \draw[->] (phi) -- (S1);
			\draw[->,dashed,red!70!black] (S.south) -- node[left,font=\scriptsize]{$\cc_n\le m_n$} (psi.north);
		\end{tikzpicture}
		\caption{One step of the induction. Every arrow points from data to choices which depend on them; no choice depends on later data. The dashed arrow is the
			delayed use of the covering degrees: the coverings used in (5) live on the \emph{source} complex $T_n$, whose degree bound $\cc_n$ entered the choice of $C_n$ in
			(2).}\label{fig:induction}
	\end{figure}
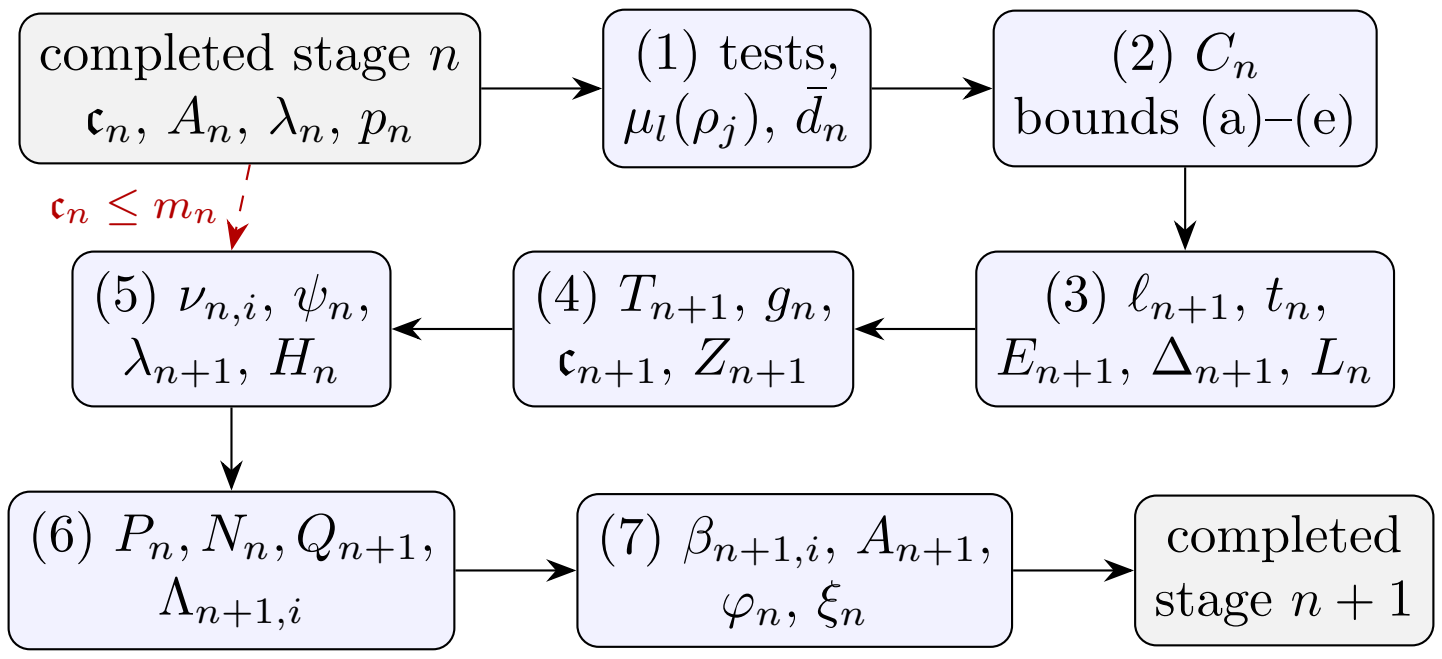
	
	\subsubsection*{The first stage}
	Choose the first stage of the presentation (Theorem~\ref{thm:presentation}) with $v_1$ strictly positive and distinguished coordinate at least $2$, lifts $\ell_1$, the
	defect $t_0=u-\ell_1(v_1)$, and $E_1\ni t_0$ (Section~\ref{ss:unit}). The unit is $(r_1,\zeta_1)=(v_1-h_{\ib_1}+x_0,0)$. Choose $d_1$ by Theorem~\ref{thm:realize}(c) (with
	$d_0=4$) and $\mathsf S_1$ by Theorem~\ref{thm:realize}(a); choose $K_1$ and $Y_1$; form $Z_1$ and $p_1$; let $\lambda_1$ be the trivial covering of $Z_1$ of degree $r_{1,\ib_1}$
	(its class $\zeta_1=0$ is correct), with any labelling; choose $\Lambda_{1,0},\Lambda_{1,1}$ by \eqref{eq:Bn} with standard labellings. This completes stage $1$.
	
	\begin{proposition}\label{prop:induction}
		The construction above can be carried out for all $n$. It produces building blocks $(A_n,B_n)$ and unital injective homomorphisms $\xi_n\colon A_n\to A_{n+1}$ such
		that all the requirements of Sections~\ref{sec:blocks}--\ref{sec:traces} hold, with $\varepsilon_n,\eta_n$ as in \eqref{eq:summable}.
	\end{proposition}
	
	\begin{proof}
		At each step, every choice is possible: in (2) by Theorem~\ref{thm:presentation}(b), which allows the entries of the next connecting matrix to exceed any bound
		determined by completed data; in (3) by Section~\ref{ss:defects}; in (4) by Theorem~\ref{thm:realize}, whose inputs ($\mathsf S_n$, $E_{n+1}$, $\Theta_{n+1}$, $i_n$, $d_{n+1}$) are
		known; in (5) by Corollary~\ref{cor:bound} and the bound $m_n\ge\cc_n$ imposed in (2); in (6) because $P_n,N_n$ depend on $L_n$, and $R_n$ and \eqref{eq:wsmall} are
		lower bounds on $Q_{n+1}$ given $M_n$, $P_n$, $N_n$, $\Lambda_{n,i}$, $s_n$ and $r_{n+1}$, all known; in (7) by Lemma~\ref{lem:labelledprops}(b). Every coordinate of each unit vector
		$r_n$ is positive (Section~\ref{ss:unit}), so every matrix algebra used has positive size. No requirement involves a quantity chosen later: the only quantity that could
		have caused a circularity, the covering degree realizing a correction $t_n$, is a degree on $T_{n+1}$, and it is bounded by $\cc_{n+1}$, which enters only the choice of
		$C_{n+1}$ at the next step. The interval target dimensions $s_{n+1}$ remain free until $L_n$, $P_n$, $N_n$ and the source endpoint matrices are known.
	\end{proof}
	
	\subsection{The invariant and the groupoid model of the first limit}\label{ss:firstlimit}
	
	Let $A_\infty=\lims(A_n,\xi_n)$ and $B=\lims(B_n,\xi_n|_{B_n})\subset A_\infty$.
	
	\begin{proposition}\label{prop:firstlimit}
		$A_\infty$ is a unital, simple, separable, nuclear, stably finite $C^*$-algebra satisfying the UCT, with
		\begin{equation}\label{eq:invA}
			K_0(A_\infty)=G_0,\qquad[1_{A_\infty}]=u,\qquad K_1(A_\infty)=G_1,\qquad(T(A_\infty),\rho_{A_\infty})=(T,\rho),
		\end{equation}
		where the identifications are compatible (the pairing of $A_\infty$ corresponds to $\rho$). Moreover there is a second countable, locally compact Hausdorff, principal étale
		groupoid $\bar{\mathcal G}$ with compact unit space such that
		\begin{equation}\label{eq:AGbar}
			(A_\infty,B)\cong\big(C^*_r(\bar{\mathcal G}),C(\bar{\mathcal G}^{(0)})\big).
		\end{equation}
	\end{proposition}
	
	\begin{proof}
		\emph{Simple, unital, separable.} Simplicity is Proposition~\ref{prop:simple}; the maps are unital, and each $A_n$ is separable.
		
		\emph{Nuclear and UCT.} Each $A_n$ is subhomogeneous: its irreducible representations have dimension at most $\max\{s_{n,q},r_{n,l}\}$ (Lemma~\ref{lem:irreps}). So $A_n$ is
		of type~I, hence nuclear. It belongs to the bootstrap class of Rosenberg and Schochet, which contains all type~I $C^*$-algebras and is closed under inductive limits
		\cite{RS}, \cite[\S22]{Bl}; inductive limits of nuclear $C^*$-algebras are nuclear. So $A_\infty$ is nuclear and satisfies the UCT \cite{RS}.
		
		\emph{$K$-theory.} $K$-theory commutes with inductive limits. By Lemma~\ref{lem:An}(d), $K_0(A_n)=G_n$ with unit $(r_n,\zeta_n)$ and $K_1(A_n)=K_n$, and by
		Lemmas~\ref{lem:phin} and~\ref{lem:xin} the maps $\xi_n$ induce $\gt_n|_{G_n}$ and $\chi_n$. By Section~\ref{ss:unit} and Lemma~\ref{lem:delayed},
		$\lims(G_n,\gt_n)=G_0$ with $(r_n,\zeta_n)\mapsto u$, and $\lims(K_n,\chi_n)=G_1$ by choice of the presentation.
		
		\emph{Traces and pairing.} \eqref{eq:TA} and the paragraph after it.
		
		\emph{Stably finite.} $T(A_\infty)\cong T$ is nonempty. Since $A_\infty$ is simple, every tracial state $\tau$ is faithful (its kernel $\{x:\tau(x^*x)=0\}$ is a closed two-sided ideal,
		not containing $1$), and so are the amplifications $\tau\otimes\tr_k$ on $M_k(A_\infty)$. A faithful tracial state rules out a proper isometry $v^*v=1\ne vv^*$ in $M_k(A_\infty)$, since
		$\tau(1-vv^*)=0$ would force $vv^*=1$. So $A_\infty$ is stably finite.
		
		\emph{The groupoid model.} We verify the hypotheses of Proposition~\ref{prop:twist} for $(A_n,B_n)$, the maps $\xi_n$ and the sets $S_n=N_{+,n}$:
		\begin{itemize}
			\item $(A_n,B_n)$ is a $C^*$-diagonal, $A_n$ unital and separable, with expectation $\EE_n$ (Lemma~\ref{lem:groupoid});
			\item $\xi_n$ is unital and injective (Lemma~\ref{lem:xin});
			\item (L1) $\xi_n(B_n)\subset B_{n+1}$, (L2) $\xi_n(\Nor(B_n))\subset\Nor(B_{n+1})$, (L3) $\EE_{n+1}\xi_n=\xi_n\EE_n$ (Lemma~\ref{lem:xin});
			\item (S1)--(S3) for $N_{+,n}$ (Lemma~\ref{lem:normalizers}(b)), and (S4) $\xi_n(N_{+,n})\subset N_{+,n+1}$ (Lemma~\ref{lem:xin}).
		\end{itemize}
		Proposition~\ref{prop:twist} gives \eqref{eq:AGbar}.
	\end{proof}
	
	We do not claim that $A_\infty$ itself is $\JS$-stable or has finite nuclear dimension; the complexes $T_n$ have unbounded dimension. This is repaired in the next step.
	
	\subsection{Recovering the order and applying classification}\label{ss:order}
	
	\begin{theorem}[{\cite[Corollary~1.9]{Li}}]\label{thm:LiZ}
		The Jiang--Su algebra $\JS$ has an untwisted groupoid model: there is a second countable, locally compact Hausdorff, principal étale groupoid $\mathcal G_\JS$ with compact
		unit space such that $\JS\cong C^*_r(\mathcal G_\JS)$.
	\end{theorem}
	
	(Deeley, Putnam and Strung gave another untwisted model of $\JS$, as the crossed product of a minimal homeomorphism of a point-like space restricted to an
	orbit-breaking subalgebra \cite{DPS18}.) Set
	\[
	\mathcal G=\bar{\mathcal G}\times\mathcal G_\JS .
	\]
	
	\begin{lemma}\label{lem:product}
		$\mathcal G$ is a second countable, locally compact Hausdorff, principal étale groupoid with unit space $\bar{\mathcal G}^{(0)}\times\mathcal G_\JS^{(0)}$, which is compact
		because both factors are, and
		\begin{equation}\label{eq:AZ}
			A_\infty\otimes\JS\cong C^*_r(\bar{\mathcal G})\otimes C^*_r(\mathcal G_\JS)\cong C^*_r(\mathcal G).
		\end{equation}
	\end{lemma}
	
	\begin{proof}
		The product of étale groupoids, with the product structure and topology, is étale (products of open bisections are open bisections), second countable,
		locally compact Hausdorff, and principal (isotropy of a product is the product of the isotropy groups). The map $C_c(\bar{\mathcal G})\odot C_c(\mathcal G_\JS)\to C_c(\mathcal G)$,
		$f\otimes g\mapsto((\gamma,\gamma')\mapsto f(\gamma)g(\gamma'))$, is a $*$-homomorphism with dense range (Stone--Weierstrass on compact pieces). For a unit $(x,x')$ the source
		fibre is $\bar{\mathcal G}_x\times(\mathcal G_\JS)_{x'}$, and the regular representation is $\lambda_x\otimes\lambda_{x'}$ on $\ell^2(\bar{\mathcal G}_x)\otimes\ell^2((\mathcal
		G_\JS)_{x'})$. The reduced norm is the supremum of these norms, which is the minimal tensor norm, because the families $(\lambda_x)$ and $(\lambda_{x'})$ are faithful
		families of representations and the minimal norm can be computed in any pair of faithful representations. So the map extends to an isomorphism
		$C^*_r(\bar{\mathcal G})\otimes_{\min}C^*_r(\mathcal G_\JS)\cong C^*_r(\mathcal G)$; both $C^*$-algebras are nuclear, so there is only one tensor norm.
	\end{proof}
	
	The proof uses neither compactness of the unit spaces nor any property of $\bar{\mathcal G}$ and $\mathcal G_\JS$ other than being second countable, locally compact Hausdorff,
	principal and étale: for any two such groupoids $\mathcal G',\mathcal G''$, the product $\mathcal G'\times\mathcal G''$ is of the same kind and
	$C^*_r(\mathcal G'\times\mathcal G'')\cong C^*_r(\mathcal G')\otimes_{\min}C^*_r(\mathcal G'')$. This general form is used in Part~III.
	
	In particular the model \eqref{eq:AZ} is untwisted, and $C(\bar{\mathcal G}^{(0)}\times\mathcal G_\JS^{(0)})=B\otimes C(\mathcal G_\JS^{(0)})$ is a $C^*$-diagonal of $A_\infty\otimes\JS$.
	
	\begin{lemma}\label{lem:AZinvariant}
		$A_\infty\otimes\JS$ is unital, simple, separable, nuclear, $\JS$-stable, stably finite and satisfies the UCT. The inclusion $a\mapsto a\otimes1$ induces isomorphisms
		$K_*(A_\infty)\cong K_*(A_\infty\otimes\JS)$ preserving the unit, and an affine homeomorphism $T(A_\infty\otimes\JS)\cong T(A_\infty)$, $\tau\mapsto\tau(\,\cdot\otimes1)$, compatible with the pairings. Hence
		\[
		K_0(A_\infty\otimes\JS)=G_0,\quad[1]=u,\quad K_1(A_\infty\otimes\JS)=G_1,\quad(T(A_\infty\otimes\JS),\rho_{A_\infty\otimes\JS})=(T,\rho).
		\]
	\end{lemma}
	
	\begin{proof}
		Tensor products of simple nuclear $C^*$-algebras are simple and nuclear; $A_\infty\otimes\JS$ is $\JS$-stable because $\JS\otimes\JS\cong\JS$ \cite{JS}; the UCT is preserved by tensor products with
		$\JS$, which is in the bootstrap class. The Jiang--Su algebra has the $K$-theory of $\C$, with $[1]=1$, and it is $KK$-equivalent to $\C$ via the unital inclusion
		\cite{JS}; by the Künneth theorem \cite{RS}, $a\mapsto a\otimes1$ induces isomorphisms on $K$-theory. For traces: let $\tau\in T(A_\infty\otimes\JS)$ and $a\in (A_\infty)_+$. Then $z\mapsto\tau(a\otimes
		z)$ is a positive tracial functional on $\JS$, so by uniqueness of the tracial state $\tau_\JS$ of $\JS$ \cite{JS} it equals $\tau(a\otimes1)\tau_\JS(z)$. By linearity
		$\tau(a\otimes z)=\tau(a\otimes1)\tau_\JS(z)$ for all $a,z$, i.e.\ $\tau=\tau|_{A_\infty}\otimes\tau_\JS$. Conversely $\tau_{A_\infty}\otimes\tau_\JS$ is a tracial state for every $\tau_{A_\infty}\in T(A_\infty)$. So
		restriction is an affine homeomorphism, and since projections of $M_k(A_\infty)$ map to $p\otimes1$, the pairings correspond. Stable finiteness follows as in
		Proposition~\ref{prop:firstlimit}.
	\end{proof}
	
	\subsubsection*{The positive cone}
	The Elliott invariant also contains the order on $K_0$. For the $C^*$-algebras at hand it is determined by the pairing.
	
	\begin{proposition}\label{prop:cone}
		Let $E$ be a unital, simple, separable, nuclear, stably finite, $\JS$-stable $C^*$-algebra. Then
		\begin{equation}\label{eq:cone}
			K_0(E)^+=\{0\}\cup\{x\in K_0(E):\ \tau_*(x)>0\text{ for all }\tau\in T(E)\}.
		\end{equation}
	\end{proposition}
	
	\begin{proof}
		If $p$ is a nonzero projection in $M_k(E)$, then $\tau_*[p]=(\tau\otimes\Tr_k)(p)>0$, because every tracial state of the simple $C^*$-algebra $M_k(E)$ is faithful. Conversely, let
		$x=[p]-[q]$ with $\tau_*(x)>0$ for all $\tau$, where $p,q\in M_k(E)$ are projections. The tracial states of $M_k(E)$ are the $\tau\otimes\tr_k=\frac1k(\tau\otimes\Tr_k)$,
		$\tau\in T(E)$, and the dimension functions $d(a)=\lim_n(\tau\otimes\tr_k)(a^{1/n})$ of $p$ and $q$ are $\frac1k\tau_*[p]$ and $\frac1k\tau_*[q]$; so $d(q)<d(p)$ for all tracial
		states of $M_k(E)$. Since $M_k(E)$ is exact, every bounded $2$-quasitrace is a trace \cite{Haa}, so this holds for all normalized quasitraces. By Rørdam's theorem, a simple,
		unital, exact, $\JS$-stable $C^*$-algebra has strict comparison of positive elements \cite[\S4]{Ro}: $d_\tau(a)<d_\tau(b)$ for all normalized quasitraces implies $a\precsim b$
		(Cuntz subequivalence). Applied to $M_k(E)$, which is again unital, simple, exact and $\JS$-stable, this gives $q\precsim p$, and for projections Cuntz subequivalence is Murray--von Neumann subequivalence: $q\sim q'\le p$. Then $x=[p-q']$ is the class of
		a projection, so $x\in K_0(E)^+$.
	\end{proof}
	
	\subsubsection*{Classification}
	\begin{theorem}[Unital classification {\cite[Theorem~A]{CGSTW}}; see also {\cite{GLN,EGLN,TWW}}]\label{thm:classification}
		Unital, simple, separable, nuclear, $\JS$-stable $C^*$-algebras satisfying the UCT are classified up to isomorphism by the invariant
		$\mathrm{KT}_u=(K_0,[1],K_1,T,\rho)$: two such $C^*$-algebras are isomorphic if and only if there are group isomorphisms of $K_0$ (preserving the unit class) and $K_1$ and an
		affine homeomorphism of trace simplices which are compatible with the pairings.
	\end{theorem}
	
	The invariant $\mathrm{KT}_u$ does not contain the order on $K_0$; for stably finite $C^*$-algebras the order is recovered from the pairing by Proposition~\ref{prop:cone}, so
	Theorem~\ref{thm:classification} is equivalent, for these $C^*$-algebras, to the classical formulation with the ordered Elliott invariant.
	
	\begin{proof}[Proof of Theorem~\ref{thm:unital}]
		Both $A$ and $A_\infty\otimes\JS$ are unital, simple, separable, nuclear, $\JS$-stable (for $A$ by hypothesis), and satisfy the UCT. By Lemma~\ref{lem:AZinvariant} and the definition of
		$(G_0,u,T,\rho,G_1)$ as the invariant of $A$, their invariants $\mathrm{KT}_u$ are isomorphic: the identity maps of $G_0=K_0(A)$ and $G_1=K_1(A)$, with $u=[1_A]$, and the
		affine homeomorphism $T(A_\infty\otimes\JS)\cong T=T(A)$, which is compatible with the pairings by construction. (By Proposition~\ref{prop:cone}, applied to $A$ and to $A_\infty\otimes\JS$, the
		positive cones correspond as well.) Theorem~\ref{thm:classification} gives $A\cong A_\infty\otimes\JS$, and by Lemma~\ref{lem:product},
		\[
		A\cong A_\infty\otimes\JS\cong C^*_r(\mathcal G)
		\]
		with $\mathcal G=\bar{\mathcal G}\times\mathcal G_\JS$ a second countable, locally compact Hausdorff, principal étale groupoid with compact unit space. Its unit space gives a
		$C^*$-diagonal of $A$ whose Weyl groupoid is $\mathcal G$ and whose Weyl twist is trivial (Theorem~\ref{thm:renault}). \qedhere
	\end{proof}
	
	\begin{remark}
		The passage to $A_\infty\otimes\JS$ is part of the construction: no assertion is made that the first limit $A_\infty$ already has finite nuclear dimension or is $\JS$-stable.
		Similarly, the finite refinements $T_n$ and their coverings have no dimension or degree bounds uniform in the invariant. These losses are compatible with the
		theorem, which asks for an untwisted groupoid model without a dimension constraint. Whether the groupoid can be chosen with finite dynamic asymptotic dimension, or
		with a finite-dimensional unit space, is not addressed here. A totally disconnected unit space, on the other hand, is impossible in general: a nonempty
		totally disconnected locally compact Hausdorff space $G^{(0)}$ contains a nonempty compact open subset $U$, and the characteristic function of $U$ is a nonzero
		projection in $C_0(G^{(0)})\subset C^*(G)$. Hence no nonzero stably projectionless $C^*$-algebra has such a model. Nor does $\JS$: for a principal model of $\JS$,
		a compact totally disconnected unit space with more than one point would give a projection in $\JS$ other than $0$ and $1$, and a one-point unit space
		would make the groupoid trivial, with $C^*$-algebra $\C$.
	\end{remark}
	%======================================================================
	\part*{Part III: The stably projectionless case}
	
	Throughout Part~III we are in Case (SP) of Section~\ref{sec:alg}: $u=(0,1)\notin G_0$, $\mathbf d_n=\Delta_nr_n\neq0$, the endpoint maps have zero blocks, and the
	building blocks $A_n$ are non-unital and stably projectionless (Lemma~\ref{lem:An}). The section numbers refer to the common Part~I for the spaces, the
	presentation and the building blocks.
	
	%======================================================================
	\section{Stably projectionless \texorpdfstring{$C^*$}{C*}-algebras: traces, scales and the reduction}\label{sec:SPprelim}
	
	Non-unital $C^*$-algebras require a more careful treatment of traces: traces may be unbounded, the tracial state space need not be compact, and the Elliott invariant
	contains an additional datum, the scale. In this section we recall the relevant notions (Section~\ref{ss:SPtraces}), the notion of continuous scale
	(Section~\ref{ss:cs}), and reduce the stably projectionless case of Theorem~\ref{thm:main} to two statements about the construction
	(Section~\ref{ss:reduction}).
	
	\subsection{Traces, dimension functions and the scale}\label{ss:SPtraces}
	
	\begin{definition}[Basic notions]
		Let $A$ be a $C^*$-algebra.
		\begin{enumerate}[label=(\alph*)]
			\item A $C^*$-subalgebra $B\subset A$ is \emph{hereditary} if $0\le a\le b$ with $b\in B$ and $a\in A$ implies $a\in B$; for $a\in A_+$, the hereditary subalgebra
			generated by $a$ is $\overline{aAa}$. A hereditary subalgebra is \emph{full} if it is not contained in a proper closed two-sided ideal; in a simple $C^*$-algebra every
			nonzero hereditary subalgebra is full.
			\item $A$ is \emph{$\sigma$-unital} if it has a countable approximate unit; equivalently, it has a \emph{strictly positive element} $e\in A_+$, i.e.\ one with
			$\overline{eAe}=A$, equivalently with $\varphi(e)>0$ for every state $\varphi$ \cite[3.10.5]{Ped}. If $e$ is strictly positive and $\|e\|\le1$, then $(e^{1/k})_k$ is an
			approximate unit. Separable $C^*$-algebras are $\sigma$-unital.
			\item The \emph{Pedersen ideal} $\Ped(A)$ is the minimal dense two-sided ideal of $A$; it is generated by the elements $g(a)$, $a\in A_+$, $g\in C_c((0,\infty))$
			\cite[5.6.1]{Ped}. If $A$ is unital, $\Ped(A)=A$.
			\item A \emph{lower semicontinuous (lsc) trace} is a map $\tau\colon A_+\to[0,\infty]$ which is additive, positively homogeneous, lower semicontinuous, and satisfies
			$\tau(x^*x)=\tau(xx^*)$. It is \emph{densely finite} if it is finite on $\Ped(A)_+$. We write $\Tt(A)$ for the cone of densely finite lsc traces, with the topology
			of pointwise convergence on $\Ped(A)$; $T(A)$ for the set of tracial states; and
			\[
			T_1(A)=\{\tau:\ \tau\text{ a bounded positive tracial functional},\ \|\tau\|\le1\},
			\]
			a weak$^*$-compact convex set. When $A$ is not unital, $T(A)$ need not be weak$^*$-closed in $T_1(A)$. The positive elements on which an lsc trace is finite span a
			two-sided ideal \cite[6.1.2]{Dix}.
			\item For $a,b\in A_+$, $a$ is \emph{Cuntz subequivalent} to $b$, written $a\precsim b$, if there are $x_k\in A$ with $x_kbx_k^*\to a$; and $a\sim b$ if $a\precsim b$
			and $b\precsim a$. The same definitions apply in $A\otimes\K$.
			\item For $\tau\in\Tt(A)$ and $a\in(A\otimes\K)_+$, the \emph{dimension function} is
			\[
			d_\tau(a)=\lim_{\varepsilon\to0}\tau(f_\varepsilon(a))=\lim_{n\to\infty}\tau(a^{1/n})\in[0,\infty],
			\]
			where $\tau$ also denotes $\tau\otimes\Tr$ and $f_\varepsilon$ is the continuous function which is $0$ on $[0,\varepsilon/2]$, $1$ on $[\varepsilon,\infty)$ and
			linear in between. It is invariant under Cuntz equivalence, monotone for $\precsim$, and additive on orthogonal elements. For a projection $p$, $d_\tau(p)=\tau(p)$.
			\item $A$ has \emph{strict comparison} (with respect to a set $S$ of traces) if $d_\tau(a)<d_\tau(b)$ for all $\tau\in S$ implies $a\precsim b$, for $a,b\in(A\otimes\K)_+$
			(with $b$ nonzero).
		\end{enumerate}
	\end{definition}
	
	For a strictly positive element $e\in A$, the function $\Sigma_A(\tau)=d_\tau(e)$ is independent of $e$ (all strictly positive elements are Cuntz equivalent); it is
	the \emph{scale} of $A$, a lower semicontinuous affine function $\Tt(A)\to[0,\infty]$, and $\Sigma_A(\tau)=\|\tau\|$ for bounded $\tau$. A bounded trace $\tau$ extends to
	the unitization by $\tilde\tau(1)=\|\tau\|$, and the pairing $\rho_A\colon K_0(A)\to\Aff(\Tt(A))$ is
	$\rho_A([p]-[q])(\tau)=(\tilde\tau\otimes\Tr)(p)-(\tilde\tau\otimes\Tr)(q)$ for projections $p,q$ over $\widetilde A$ with equal scalar parts; for unbounded traces see
	\cite[Definition~2.6]{GL4}. Following Gong and Lin \cite[2.7]{GL4}, the \emph{Elliott invariant} of a separable simple stably projectionless $C^*$-algebra $A$ with
	$\Tt(A)\neq\{0\}$ is
	\[
	\Ell(A)=\big((K_0(A),\Sigma(K_0(A)),\Tt(A),\Sigma_A,\rho_A),K_1(A)\big),\qquad\Sigma(K_0(A))=\{0\}
	\]
	(the second entry, the classes of projections in $A$, is trivial because $A$ has no nonzero projections).
	
	\begin{remark}\label{rem:normpreserving}
		If $\theta\colon X\to Y$ is a $*$-homomorphism with $\theta(e)$ strictly positive in $Y$ for some strictly positive $e\in X$, then for every bounded positive functional
		$\tau$ on $Y$ we have $\|\tau\circ\theta\|=\lim_k\tau(\theta(e)^{1/k})=\|\tau\|$, because $(\theta(e)^{1/k})_k$ is an approximate unit of $Y$. If $\theta$ maps one strictly
		positive element to a strictly positive element, it does so for all: if $e,e'$ are strictly positive in $X$, then
		$\overline{\theta(e')Y\theta(e')}\supset\theta(\overline{e'Xe'})=\theta(X)\ni\theta(e)$, hence it contains $\overline{\theta(e)Y\theta(e)}=Y$.
	\end{remark}
	
	\begin{theorem}[{\cite{Br}}]\label{thm:Brown}
		Let $A$ be a $\sigma$-unital $C^*$-algebra and $B\subset A$ a full $\sigma$-unital hereditary subalgebra. Then $B\otimes\K\cong A\otimes\K$.
	\end{theorem}
	
	We shall compare trace cones of $C^*$-algebras and of their hereditary subalgebras by restriction. The following fact is standard; we include a proof for completeness.
	
	\begin{lemma}\label{lem:hertraces}
		Let $X$ be a $\sigma$-unital $C^*$-algebra and $Y\subset X$ a $\sigma$-unital full hereditary $C^*$-subalgebra. Then restriction $\Tt(X)\to\Tt(Y)$, $\tau\mapsto\tau|_Y$, is a
		bijection, and a homeomorphism for the topologies of pointwise convergence on the Pedersen ideals. In particular, for $\sigma$-unital $Z$, every element of
		$\Tt(Z\otimes\K)$ has the form $\tau\otimes\Tr$ for a unique $\tau\in\Tt(Z)$.
	\end{lemma}
	
	\begin{proof}
		\emph{Injectivity and continuity.} An lsc trace $\tau$ is determined by its values on $\Ped(X)_+$, because $(x-\varepsilon)_+\in\Ped(X)$ increases to $x$ as
		$\varepsilon\searrow0$. Since $Y$ is full, the ideal generated by $\Ped(Y)$ is dense, hence contains the minimal dense ideal $\Ped(X)$ \cite[5.6.1]{Ped}. So every element of
		$\Ped(X)$ is a finite sum of elements $ubw$ with $b\in\Ped(Y)_+$ and $u,w\in\widetilde X$, and $\tau(ubw)=\tau(b^{1/2}wub^{1/2})$ by the trace property (polarize
		$\tau(z^*z)=\tau(zz^*)$, which holds on the domain of $\tau$); here $b^{1/2}wub^{1/2}$ is a linear combination of four elements of $\Ped(Y)_+$ which do not depend on $\tau$.
		Hence $\tau$ is determined by $\tau|_{\Ped(Y)}$, restriction is injective, and its inverse is continuous on its image. Restriction itself is continuous because
		$\Ped(Y)\subset\Ped(X)$: $\Ped(Y)$ is generated by the elements $g(y)$, $y\in Y_+$, $g\in C_c((0,\infty))$, and these lie in $\Ped(X)$.
		
		\emph{Stabilization.} Regard $Z$ as $Z\otimes e_{11}\subset Z\otimes\K$. For $\tau\in\Tt(Z)$, $\tau\otimes\Tr$ is an lsc trace on $Z\otimes\K$; the positive elements on which
		it is finite span a two-sided ideal containing $\Ped(Z)\otimes e_{11}$, and this ideal contains $\Ped(Z\otimes\K)$ by the first step (applied to the full hereditary
		subalgebra $Z\otimes e_{11}$). So $\tau\otimes\Tr\in\Tt(Z\otimes\K)$. Conversely let $\sigma\in\Tt(Z\otimes\K)$ and $\tau=\sigma|_{Z\otimes e_{11}}$. For $x=(x_{ij})\in M_k(Z)_+$,
		averaging over the unitaries $\diag(\omega^{a_1},\dots,\omega^{a_k})$, $\omega=e^{2\pi i/k}$, in the multiplier algebra (under which $\sigma$ is invariant) gives
		$\sigma(x)=\sum_i\sigma(x_{ii}\otimes e_{ii})$, and $\sigma(x_{ii}\otimes e_{ii})=\sigma(x_{ii}\otimes e_{11})$ by the trace property applied to $z=x_{ii}^{1/2}\otimes e_{1i}$.
		Thus $\sigma=\tau\otimes\Tr$ on $\bigcup_kM_k(Z)_+$. For general $x\in(Z\otimes\K)_+$ and $p_k=1\otimes(e_{11}+\dots+e_{kk})$, the elements $x^{1/2}p_kx^{1/2}$ increase to $x$ in
		norm and have the same trace as $p_kxp_k$ under any trace, so $\sigma(x)=(\tau\otimes\Tr)(x)$ by lower semicontinuity.
		
		\emph{Full corners.} Suppose $Y=pXp$ for a full projection $p$ in the multiplier algebra of $X$. By Brown's lemma \cite[Lemma~2.5]{Br} there is $v$ in the multiplier algebra
		of $X\otimes\K$ with $v^*v=1$ and $vv^*=p\otimes1$, so that $Y\otimes\K=v(X\otimes\K)v^*$. Let $\sigma\in\Tt(Y)$ and $\hat\sigma=\sigma\otimes\Tr$, and define
		$\hat\tau(x)=\hat\sigma(vxv^*)$ for $x\in(X\otimes\K)_+$. Since $x\mapsto vxv^*$ is an isomorphism of $X\otimes\K$ onto $Y\otimes\K$, $\hat\tau$ is lower semicontinuous and finite
		on $\Ped(X\otimes\K)$. It is a trace: for $z\in X\otimes\K$ put $w=vzv^*\in Y\otimes\K$; then $ww^*=vzz^*v^*$ and $w^*w=vz^*zv^*$. It extends $\hat\sigma$: for $y\in(Y\otimes\K)_+$ the
		element $u=vy^{1/2}$ lies in $(p\otimes1)(X\otimes\K)(p\otimes1)=Y\otimes\K$, and $\hat\tau(y)=\hat\sigma(uu^*)=\hat\sigma(u^*u)=\hat\sigma(y)$. By stabilization $\hat\tau=\tau\otimes\Tr$
		with $\tau=\hat\tau|_{X\otimes e_{11}}$, and then $\tau|_Y=\hat\sigma|_{Y\otimes e_{11}}=\sigma$.
		
		\emph{General case.} The linking algebra $L=\begin{pmatrix}Y&\overline{YX}\\ \overline{XY}&X\end{pmatrix}\subset M_2(X)$ is a $\sigma$-unital $C^*$-algebra (it is the hereditary
		subalgebra of $M_2(X)$ generated by $\diag(c,e)$ for strictly positive $c\in Y$, $e\in X$), and $Y=pLp$, $X=qLq$ for the full projections $p=e_{11}$, $q=e_{22}$ in its
		multiplier algebra. For $y\in Y_+$ the element $z=\begin{pmatrix}0&0\\ y^{1/2}&0\end{pmatrix}\in L$ satisfies $z^*z=y\oplus0$ and $zz^*=0\oplus y$, so every trace on $L$ takes the
		same value on $y\oplus0$ and on $0\oplus y$. Hence restriction $\Tt(X)\to\Tt(Y)$ is the composite of the inverse of the restriction $\Tt(L)\to\Tt(qLq)$ with the restriction
		$\Tt(L)\to\Tt(pLp)$, and both are bijective by the corner case.
	\end{proof}
	
	\subsection{Continuous scale}\label{ss:cs}
	
	\begin{definition}[{\cite{Lin1}}, {\cite[Definition~9.1]{GL1}}]\label{def:cs}
		A $\sigma$-unital non-unital $C^*$-algebra $A$ has \emph{continuous scale} if for an approximate unit $(e_n)$ with $e_{n+1}e_n=e_n$ and every $b\in A_+\setminus\{0\}$ there is
		$N$ with $e_m-e_n\precsim b$ for all $m>n\ge N$.
	\end{definition}
	
	Intuitively, the ``tail'' $e_m-e_n$ of the approximate unit becomes arbitrarily small in the Cuntz semigroup; in tracial terms, $\tau(e_n)\to\|\tau\|$ uniformly. The
	point of continuous scale is that the trace simplex is compact and all densely finite traces are bounded.
	
	\begin{proposition}\label{prop:cs}
		Let $A$ be a non-elementary $\sigma$-unital simple $C^*$-algebra.
		\begin{enumerate}[label=\textup{(\alph*)}]
			\item If $A$ is separable, then for every $a\in A_+\setminus\{0\}$ there is $e\in\overline{aAa}_+$ with $\|e\|=1$ such that $\overline{eAe}$ has continuous scale.
			\item If $A$ has continuous scale, then $A$ is algebraically simple; hence $\Ped(A)=A$, every densely finite lsc trace is bounded, and $\Tt(A)=\R_{\ge0}T(A)$. If moreover
			$T(A)\neq\emptyset$, then $T(A)$ is weak$^*$-compact.
		\end{enumerate}
	\end{proposition}
	
	\begin{proof}
		(a) is \cite[Remark~9.2]{GL1}, which rests on \cite[2.3]{Lin2}, applied to the non-elementary separable simple $C^*$-algebra $\overline{aAa}$. (b) Algebraic simplicity is
		\cite[3.3]{Lin1} (as used in the proof of \cite[Theorem~9.3]{GL1}); then $\Ped(A)$, a nonzero ideal, is $A$, so a densely finite lsc trace is finite everywhere, hence
		bounded (a finite positive tracial functional on a $C^*$-algebra is bounded). Compactness of $T(A)$ is \cite[Theorem~9.3]{GL1} (see also \cite[Theorem~5.3]{EGLN1}).
	\end{proof}
	
	For a separable simple stably projectionless $A$ with continuous scale and $T(A)\neq\emptyset$, Proposition~\ref{prop:cs}(b) shows that the invariant is determined by the
	quadruple $(K_0(A),T(A),\rho_A,K_1(A))$, with $\rho_A$ regarded as a map $K_0(A)\to\Aff(T(A))$ \cite[2.7]{GL4}: the cone is $\R_{\ge0}T(A)$, and the scale is the norm,
	i.e.\ the homogeneous function equal to $1$ on $T(A)$. For a strictly positive $e$ and $\tau\in T(A)$ we have $d_\tau(e)=\lim_k\tau(e^{1/k})=1$.
	
	\begin{theorem}[{\cite[Theorem~14.9]{GL4}}]\label{thm:GL}
		Let $A$ and $B$ be separable, simple, nuclear, stably projectionless, $\JS$-stable $C^*$-algebras which satisfy the UCT. Then $A\cong B$ if and only if
		$\Ell(A)\cong\Ell(B)$.
	\end{theorem}
	
	Stably projectionless $C^*$-algebras are not elementary, so for them $\JS$-stability is equivalent to finite nuclear dimension \cite{Ti,CE}; this is the setting of
	\cite{GL4}. Some regularity hypothesis is indispensable in a statement of this kind: $\C$ and $\JS$, or $\K$ and $\JS\otimes\K$, are separable, simple, nuclear and
	stably finite with isomorphic $K$-groups, trace cones, scales and pairings, but they are not isomorphic.
	
	Here an isomorphism of invariants consists of isomorphisms of the $K$-groups and an isomorphism of topological cones $\Tt(A)\cong\Tt(B)$ which is compatible with the
	scales and the pairings \cite[2.3, 2.7]{GL4}. For stably projectionless $C^*$-algebras with continuous scale this reduces to isomorphisms of $K_0,K_1$ and an affine
	homeomorphism of the compact trace simplices compatible with the pairings \cite[Corollary~14.7]{GL4}.
	
	\subsection{Reduction to continuous scale}\label{ss:reduction}
	
	Let $A$ be a classifiable stably projectionless $C^*$-algebra. Fix $a\in A_+\setminus\{0\}$ and $e\in\overline{aDa}_+$ as in Proposition~\ref{prop:cs}(a), and put
	\[
	A_c=\overline{eDe}.
	\]
	
	\begin{lemma}\label{lem:Dc}
		$A_c$ is separable, simple, nuclear, stably projectionless, $\JS$-stable, satisfies the UCT and has continuous scale. Its tracial state space $T:=T(A_c)$ is a nonempty
		metrizable Choquet simplex, $\Tt(A_c)=\R_{\ge0}T$, and there is an isomorphism $A\otimes\K\cong A_c\otimes\K$ carrying $A\otimes e_{11}$ onto a full hereditary subalgebra.
		Writing $G_0=K_0(A_c)\cong K_0(A)$, $G_1=K_1(A_c)\cong K_1(A)$ and $\rho=\rho_{A_c}\colon G_0\to\Aff(T)$, the invariant of $A_c$ is $(G_0,T,\rho,G_1)$, and $\rho$ satisfies
		\eqref{eq:P}: $\rho(G_0)\cap\Aff_{++}(T)=\emptyset$.
	\end{lemma}
	
	\begin{proof}
		Hereditary subalgebras of separable simple nuclear $C^*$-algebras are separable, simple and nuclear. $A_c$ is stably projectionless because $A_c\otimes\K\subset A\otimes\K$, and it is
		$\JS$-stable by Toms--Winter \cite[Corollary~3.1]{TW}. It is a full hereditary subalgebra of the simple $C^*$-algebra $A$, so Theorem~\ref{thm:Brown} gives
		$A_c\otimes\K\cong A\otimes\K$; hence $K_*(A_c)\cong K_*(A)$, and $A_c$ satisfies the UCT, which is invariant under $KK$-equivalence. Continuous scale holds by construction,
		and Proposition~\ref{prop:cs}(b) gives $\Ped(A_c)=A_c$ and $\Tt(A_c)=\R_{\ge0}T$.
		
		Next, $T\neq\emptyset$. The $C^*$-algebra $A$ is not purely infinite, since a simple purely infinite $C^*$-algebra has nonzero projections. Its Cuntz semigroup is almost unperforated by
		R{\o}rdam \cite[Theorem~4.5]{Ro}, and a simple exact $C^*$-algebra with almost unperforated Cuntz semigroup is purely infinite if and only if it is traceless
		\cite[Corollary~5.1]{Ro}, that is, has no nonzero densely finite lsc $2$-quasitrace. Such quasitraces are traces by exactness (Haagerup \cite{Haa}, Blanchard--Kirchberg
		\cite[Remark~2.29(i)]{BK}), and on a simple $C^*$-algebra an lsc trace which is finite at one nonzero positive element is finite on the Pedersen ideal. Hence
		$\Tt(A)\cong\Tt(A_c)\neq\{0\}$, and every nonzero element of $\Tt(A_c)$ is a positive multiple of a tracial state. By Proposition~\ref{prop:cs}(b) again, $T$ is
		weak$^*$-compact. It is a compact base of the cone of bounded traces of $A_c$, which is a lattice cone, so $T$ is a Choquet simplex \cite[\S2.1]{EGLN1}; it is metrizable
		because $A_c$ is separable. The isomorphism $A\otimes\K\cong A_c\otimes\K$ identifies $\Tt(A)\cong\Tt(A\otimes\K)\cong\Tt(A_c\otimes\K)\cong\R_{\ge0}T$ by restriction of traces
		(Lemma~\ref{lem:hertraces}), compatibly with the pairings.
		
		Condition \eqref{eq:P} is the necessity of weak unperforation of the pairing for $\JS$-stable stably projectionless $C^*$-algebras; see the comment after
		\cite[Theorem~1.3]{Li}, which refers to \cite[\S A.1]{EGLN2}, and \cite[Corollary~3.2]{GL3}. A direct appeal to strict comparison does not suffice, because the two
		projections representing a class of $K_0(A_c)$ live over the unitization, where strict comparison fails at the character. For completeness we give a self-contained
		proof in Appendix~\ref{app:P}; Theorem~\ref{thm:P} applies to $E=A_c$.
	\end{proof}
	
	\begin{remark}\label{rem:P}
		In Li's formulation \cite[Theorem~1.3]{Li}, \eqref{eq:P} is the condition that for every $g\in G_0$ there is $\tau\in T$ with $\rho(g)(\tau)=0$. The two are equivalent, since
		an affine function on a convex set taking both signs has a zero, and $\pm g\in G_0$. For compact $T$ there is even a common zero of all $\rho(g)$: a real combination of
		elements of $\rho(G_0)$ lying in the open cone $\Aff_{++}(T)$ could be perturbed to a rational one and then scaled to an element of $\rho(G_0)\cap\Aff_{++}(T)$, so
		$\mathrm{span}_\R\rho(G_0)$ misses $\Aff_{++}(T)$; by the Hahn--Banach theorem some state of the order-unit space $(\Aff(T),1)$ vanishes on $\rho(G_0)$, and the states of
		$\Aff(T)$ are the point evaluations \cite[Chapter~7]{Goo}.
	\end{remark}
	
	\subsubsection*{The plan}
	By Lemma~\ref{lem:Dc} and Theorem~\ref{thm:GL}, the stably projectionless case of Theorem~\ref{thm:main} follows from two statements, proved in
	Sections~\ref{sec:SPmaps}--\ref{sec:SPtraces} and in Section~\ref{sec:SPscale}, respectively.
	\begin{enumerate}[label=(\Roman*)]
		\item For every quadruple $(G_0,T,\rho,G_1)$ as in Theorem~\ref{thm:SPB} there is a second countable, locally compact Hausdorff, principal étale groupoid $G_c$ such that
		$C:=C^*(G_c)$ is a classifiable stably projectionless $C^*$-algebra with $K_*(C)\cong(G_0,G_1)$, $T(C)\cong T$ with pairing $\rho$, and with no unbounded densely finite lower
		semicontinuous traces: $\Tt(C)=\R_{\ge0}T(C)$. (Theorem~\ref{thm:SPB}.)
		\item For every lower semicontinuous affine function $\gamma\colon T\to(0,\infty]$ there is a positive element $a'$ of the canonical diagonal of $C\otimes\K$ with
		$d_{\tau\otimes\Tr}(a')=\gamma(\tau)$ for all $\tau\in T(C)\cong T$. (Proposition~\ref{prop:exactscale}.)
	\end{enumerate}
	The hereditary subalgebra generated by $a'$ is then the $C^*$-algebra of an open restriction of the groupoid of $C\otimes\K$ (Lemma~\ref{lem:restriction}), and for
	$\gamma$ the scale of $A$ it has the same invariant as $A$ (Section~\ref{ss:proofA}).
	
	\begin{figure}[ht]
		\centering
		\begin{tikzpicture}[>={Stealth[length=2mm]},node distance=0.6cm and 1.0cm,
			box/.style={draw,rounded corners,fill=blue!5,align=center,font=\small,minimum height=0.9cm}]
			\node[box,fill=gray!10] (D) {$A$ classifiable,\\ stably projectionless};
			\node[box,right=of D] (Dc) {$A_c\subset A$ with\\ continuous scale};
			\node[box,right=of Dc] (inv) {$(G_0,T,\rho,G_1)$\\ with \eqref{eq:P}};
			\node[box,below=of inv] (C) {$C=C^*(G_c)$\\ (Theorem~\ref{thm:SPB})};
			\node[box,left=of C] (a) {$a'\in C_0((G^s)^{(0)})$\\ with $d_\bullet(a')=\Sigma_A$};
			\node[box,left=of a,fill=gray!10] (D2) {$A\cong C^*(G^s|_{\{a'>0\}})$};
			\draw[->] (D) -- node[above,font=\scriptsize]{Prop.~\ref{prop:cs}} (Dc);
			\draw[->] (Dc) -- node[above,font=\scriptsize]{invariant} (inv);
			\draw[->] (inv) -- node[right,font=\scriptsize]{construction} (C);
			\draw[->] (C) -- node[above,font=\scriptsize]{scale} (a);
			\draw[->] (a) -- node[above,font=\scriptsize]{Thm.~\ref{thm:GL}} (D2);
		\end{tikzpicture}
		\caption{The proof in the stably projectionless case.}\label{fig:SPplan}
	\end{figure}
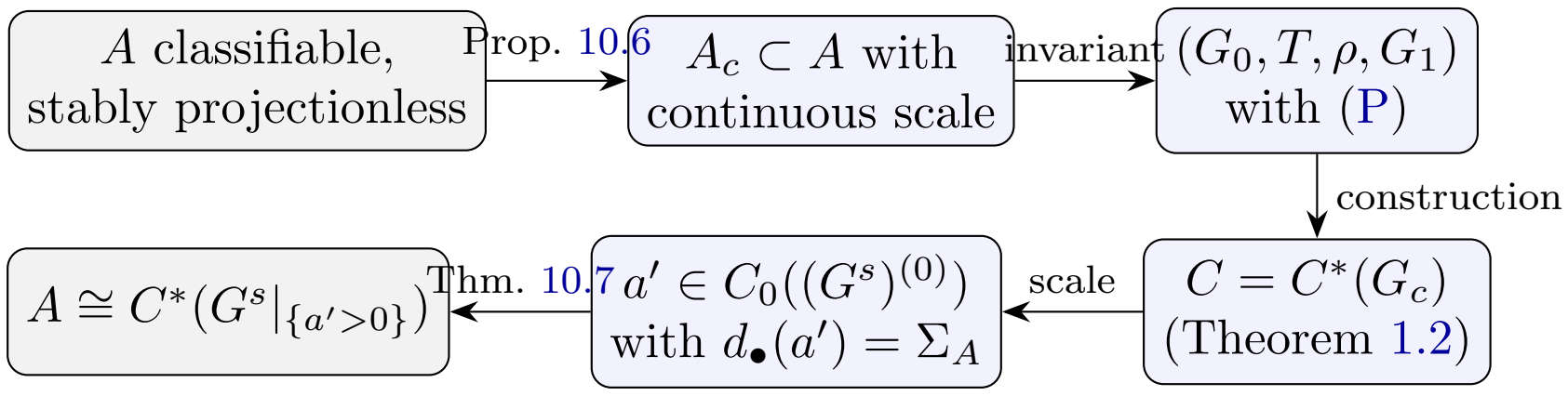
	
	From now on in Part~III we fix $(G_0,T,\rho,G_1)$ with $G_0,G_1$ countable abelian groups, $T$ a nonempty metrizable Choquet simplex and $\rho\colon G_0\to\Aff(T)$ a homomorphism
	satisfying \eqref{eq:P}, and we apply Sections~\ref{sec:alg}--\ref{sec:blocks} in Case (SP).
	%======================================================================
	\section{The stably projectionless case: connecting maps}\label{sec:SPmaps}
	
	In this section we construct connecting maps $\xi_n=(\Xi_n,\Psi_n)\colon A_n\to A_{n+1}$ for the non-unital building blocks of Case (SP), with homogeneous component
	$\Psi_n$ and interval component $\Xi_n$, such that the limit is simple, images of strictly positive elements are strictly positive, the $K$-theory is as prescribed,
	and diagonals, expectations and nonnegative monomial elements are preserved.
	
	\subsubsection*{The failure of the unital method}
	In the unital case (Section~\ref{sec:simple}) a selected irreducible representation was connected by a path of representations to a block sum of basepoint
	representations, and the path was inserted in place of reserved basepoint blocks. In Case (SP) some endpoint zero block is nonempty already at stage $1$, and
	from stage $2$ onward every endpoint zero block is nonempty (Remark~\ref{rem:spzeroblocks}). When the endpoint value $\beta_{n,i}q_n(\cdot)_p$ of an interval
	representation has a nonempty zero block, a path $u\mapsto f_p(us)$ ends at a value which is not a block sum of basepoint representations of any unital map, and
	cannot be put in place of basepoint blocks. Thus the unital endpoint-replacement argument does not apply uniformly to these building blocks. Instead we insert the detecting representations themselves into one
	matrix summand, with an explicit $K_0$-correction which is possible because $\mathbf d_n\neq0$, and we realize the resulting endpoint values inside the interval component by
	reparametrized copies whose parameters stay in the open interval. A pleasant by-product is that images of strictly positive elements remain strictly positive, which
	is needed for traces (Remark~\ref{rem:normpreserving}).
	
	Throughout we fix $\varepsilon_n=\eta_n=2^{-n-6}$.
	
	\subsection{Interval data}\label{ss:intervaldata}
	The interval components $\Xi_n$ have the following form. Besides copies of the interval functions and of basepoint blocks, they contain reparametrized copies whose
	endpoint values may be finitely many \emph{marked} interior parameters, and \emph{designated} copies of the distinguished basepoint block which move along paths in
	$Z_n$.
	
	\begin{definition}\label{def:intervaldata}
		Fix a stage $n$, a finite set $\mathfrak M\subset(0,1)$ of \emph{marked parameters} and, for each $t\in\mathfrak M$, a source summand $p(t)\in[l_n]$. For a target index $q$,
		\emph{interval data} consist of
		\begin{enumerate}[label=\textup{(\roman*)}]
			\item a finite family $\mathcal C_q$ of \emph{$f$-copies} $c$, each with a source summand $p_c\in[l_n]$ and a continuous function $\phi_c\colon[0,1]\to[0,1]$ with
			$\phi_c(t)\in(0,1)$ for $0<t<1$, endpoint values $\phi_c(0),\phi_c(1)\in\{0,1\}\cup\mathfrak M$, $(\phi_c(0),\phi_c(1))\notin\{(0,0),(1,1)\}$, and $p_c=p(t')$ whenever
			$\phi_c(i)=t'\in\mathfrak M$;
			\item a vector $\cv_q\in\Z_{\ge1}^{J_n}$ of \emph{constant copies} (copies of the basepoint blocks $q_n(\cdot)_l$), among which a subset $\mathcal D_q$ of copies of
			$q_n(\cdot)_{\ib_n}$ is \emph{designated}: each $c\in\mathcal D_q$ carries a continuous function $u_c\colon[0,1]\to[0,1]$ with $u_c(t)\in(0,1]$ for $0<t<1$ and
			$u_c(0),u_c(1)\in\{0,1\}$, and a path $z_c\colon[0,1]\to Z_n$ from $\theta_n$.
		\end{enumerate}
		The associated \emph{interval map} is
		\begin{equation}\label{eq:Xidef}
			\Xi(f,a)_q(t)=\diag\Big(\big(f_{p_c}(\phi_c(t))\big)_{c\in\mathcal C_q},\ \big(q_n(a)_l\big)_{\text{undesignated copies}},\ \big(a(z_c(u_c(t)))\big)_{c\in\mathcal D_q}\Big)\in
			M_{s_{n+1,q}},
		\end{equation}
		in a fixed arrangement of coordinates, where $a(z_c(\cdot))$ is taken in the labels transported along $z_c$ (sheet transport, Section~\ref{ss:paths} and
		Lemma~\ref{lem:Zn}(f)), and $s_{n+1,q}=\sum_{c\in\mathcal C_q}s_{n,p_c}+\sum_l\cv_{q,l}r_{n,l}$.
	\end{definition}
	
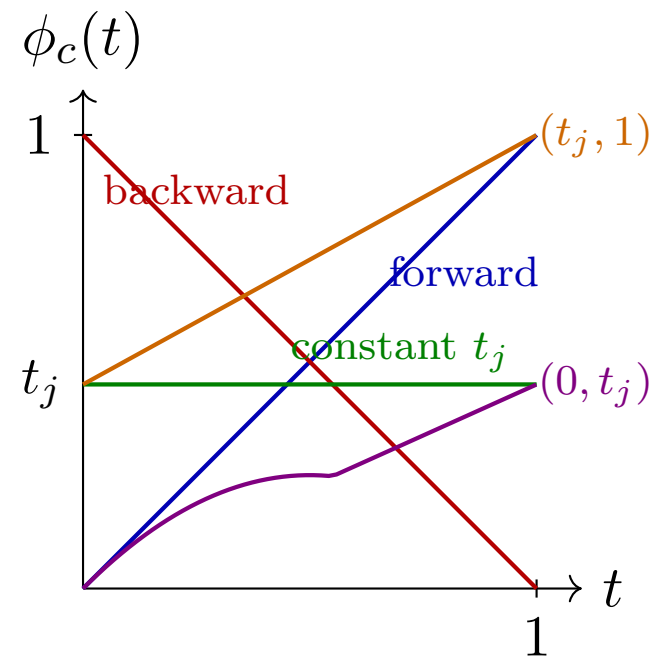
\begin{figure}[ht]
		\centering
		\begin{tikzpicture}[scale=3.0]
			\draw[->] (0,0) -- (1.1,0) node[right]{$t$}; \draw[->] (0,0) -- (0,1.1) node[above]{$\phi_c(t)$};
			\draw (1,0.02) -- (1,-0.02) node[below]{$1$}; \draw (0.02,1) -- (-0.02,1) node[left]{$1$};
			\draw[dashed,gray] (0,0.45) -- (1,0.45); \node[left,font=\small] at (0,0.45) {$t_j$};
			\draw[thick,blue!70!black] (0,0) -- (1,1); \node[blue!70!black,font=\scriptsize] at (0.84,0.7) {forward};
			\draw[thick,red!70!black] (0,1) -- (1,0); \node[red!70!black,font=\scriptsize] at (0.25,0.88) {backward};
			\draw[thick,green!50!black] (0,0.45) -- (1,0.45); \node[green!50!black,font=\scriptsize] at (0.7,0.52) {constant $t_j$};
			\draw[thick,orange!80!black,domain=0:1,samples=60] plot (\x,{max((1-\x)*0.45+\x*1,\x*(1-\x))});
			\node[orange!80!black,font=\scriptsize] at (1.13,1) {$(t_j,1)$};
			\draw[thick,violet,domain=0:1,samples=60] plot (\x,{max((1-\x)*0+\x*0.45,\x*(1-\x))});
			\node[violet,font=\scriptsize] at (1.13,0.45) {$(0,t_j)$};
		\end{tikzpicture}
		\caption{Kinds of $f$-copies: forward, backward, constant at a marked parameter $t_j$, and the two mixed kinds used in (D1). All stay in $(0,1)$ for interior $t$.}\label{fig:fcopies}
	\end{figure}

	An $f$-copy with $\phi_c(t)=t$ is called \emph{forward}, and one with $\phi_c(t)=1-t$ \emph{backward}. The excluded endpoint pairs $(0,0)$ and $(1,1)$ would produce a
	coordinate lying in the zero block at both ends, violating (Z) for the target; the requirement $\phi_c(t)\in(0,1)$ for interior $t$ guarantees that strictly positive
	elements are mapped to strictly positive elements (Figure~\ref{fig:fcopies}).
	
	\begin{lemma}\label{lem:intervaldata}
		For interval data as in Definition~\ref{def:intervaldata}:
		\begin{enumerate}[label=\textup{(\alph*)}]
			\item For $i\in\{0,1\}$ the fibre $\Xi(f,a)_q(i)$ is the direct sum of: the blocks $f_{p_c}(0)=\beta_{n,0}(q_na)_{p_c}$, with zero block $Z^{p_c}_{n,0}$, over the $c$ with
			$\phi_c(i)=0$; the blocks $f_{p_c}(1)=\beta_{n,1}(q_na)_{p_c}$, with zero block $Z^{p_c}_{n,1}$, over the $c$ with $\phi_c(i)=1$; the blocks $f_{p_c}(t')$ over the $c$ with
			$\phi_c(i)=t'\in\mathfrak M$; the basepoint blocks $q_n(a)_l$ of the undesignated copies and of the designated copies with $u_c(i)=0$; and the blocks $a(z_c(1))$ of the
			designated copies with $u_c(i)=1$. Its zero coordinates on $e_n$ form the set
			\[
			Z'^q_{i}:=\bigsqcup_{c:\,\phi_c(i)=0}Z^{p_c}_{n,0}\ \sqcup\ \bigsqcup_{c:\,\phi_c(i)=1}Z^{p_c}_{n,1}\subset[s_{n+1,q}] .
			\]
			\item $Z'^q_0\cap Z'^q_1=\emptyset$, and the complement of $Z'^q_0\cup Z'^q_1$ contains all coordinates of constant copies. Hence the sets $S'^q_i:=[s_{n+1,q}]\setminus
			Z'^q_i$ satisfy \textup{(Z)} and \textup{(C)}.
			\item For $0<t<1$ the fibre $\Xi(e_n)_q(t)$ is invertible.
			\item Every fibre matrix of $\Xi(f,a)$ is, after a fixed relabelling of coordinates, a direct sum of fibre matrices of $(f,a)$, taken in transported labels for the
			designated copies. Hence $\Xi$ maps diagonal, monomial and nonnegative monomial elements to elements of the same kind and commutes with fibrewise diagonal compression.
		\end{enumerate}
	\end{lemma}
	
	\begin{proof}
		(a) is read off from \eqref{eq:Xidef} and Lemma~\ref{lem:An}(c): the blocks $f_{p_c}(t')$ with $t'\in(0,1)$ and the blocks $a(z)$ are invertible on $e_n$, and the basepoint
		blocks are unital. (b) A coordinate of the copy $c$ lies in $Z'^q_0\cap Z'^q_1$ if and only if it lies in $Z^{p_c}_{n,\phi_c(0)}\cap Z^{p_c}_{n,\phi_c(1)}$, where we put
		$Z^{p}_{n,t'}:=\emptyset$ for $t'\in\mathfrak M$. This set is empty by (Z) at stage $n$ when $\{\phi_c(0),\phi_c(1)\}=\{0,1\}$, and trivially otherwise; the excluded kinds
		$(0,0)$ and $(1,1)$ are exactly those for which it could be nonempty. Coordinates of constant copies are never zero coordinates, and there is at least one constant copy
		since $\cv_q\ge1$. (c) For $0<t<1$ all blocks are invertible on $e_n$: $f_{p_c}(\phi_c(t))$ with $\phi_c(t)\in(0,1)$ by Lemma~\ref{lem:An}(c), $q_n(1)_l=1$ and $1(z)=1$. (d)
		follows from \eqref{eq:Xidef}, the properties of sheet transport (Section~\ref{ss:homog}) and Lemma~\ref{lem:labelledprops}.
	\end{proof}
	
	\begin{remark}\label{rem:spzeroblocks}
		If the interval data at stage $n$ contain at least one forward and one backward copy of every source summand (as they will, by \textup{(D3)} below), then by (a) each of
		$Z'^q_0$, $Z'^q_1$ contains a copy of $Z^{p}_{n,0}\sqcup Z^{p}_{n,1}$ for every $p$. Since some zero block at stage $n$ is nonempty, all zero blocks are nonempty at every stage
		$n\ge2$, in every summand.
	\end{remark}
	
	\subsection{The design}\label{ss:design}
	Fix a completed stage $n$ and finitely many nonzero positive elements of $A_n$ to be detected (they are chosen in Section~\ref{ss:SPsimplicity}). For each of them choose an
	irreducible representation of $A_n$ in which its image is nonzero (Lemma~\ref{lem:An}(b)). Matrix-type representations $(f,a)\mapsto a_l$ need no special treatment: the
	basepoint block $q_n(\cdot)_l$ occurs in all the places listed in Proposition~\ref{prop:xi}(d) below, which is all that the propagation argument in
	Proposition~\ref{prop:SPsimple} uses. Let $\varrho_1,\dots,\varrho_k$ be the remaining chosen representations. Each $\varrho_j$ is either of
	\begin{itemize}
		\item \emph{interval type}, $\varrho_j=\ev_{(p_j,t_j)}\colon(f,a)\mapsto f_{p_j}(t_j)$ with $t_j\in(0,1)$ and $\dim\varrho_j=s_{n,p_j}$, or of
		\item \emph{homogeneous type}, $\varrho_j\colon(f,a)\mapsto a(z_j(1))$, the evaluation at the endpoint of a chosen path $z_j\colon[0,1]\to Z_n$ with $z_j(0)=\theta_n$, taken in the
		sheet basis transported along $z_j$, with $\dim\varrho_j=r_{n,\ib_n}$.
	\end{itemize}
	Since a positive element which is nonzero at an interior point $(p,t)$ is nonzero at all nearby points, we may and do choose the parameters $t_j$ of the interval-type
	representations pairwise distinct. Set
	\begin{equation}\label{eq:spnu}
		\bb_j=\begin{cases}(\Lambda_{n,0})_{p_j,\cdot}\in\Z^{J_n},\quad z_j^{\mathrm{def}}:=|Z^{p_j}_{n,0}|,&\text{interval type},\\
			h_{\ib_n},\quad z_j^{\mathrm{def}}:=0,&\text{homogeneous type},\end{cases}
		\qquad\text{so that }\ \dim\varrho_j=\bb_j\cdot r_n+z_j^{\mathrm{def}} .
	\end{equation}
	Thus $\bb_j$ records the basepoint blocks to which $\varrho_j$ is homotopic (at $t=0$, or at $u=0$), and $z_j^{\mathrm{def}}$ the \emph{dimension defect}: the number of
	coordinates of $\varrho_j$ which are not accounted for by these blocks (the zero block at $t=0$).
	
	\begin{lemma}\label{lem:Krho}
		For $x=(h,e)\in G_n=K_0(A_n)$, the $K_0$-map of the finite-dimensional representation $\varrho_j$ is $(\varrho_j)_*(x)=\bb_j\cdot h$.
	\end{lemma}
	
	\begin{proof}
		An interval-type $\varrho_j=\ev_{(p,t)}$ is homotopic through the representations $\ev_{(p,t')}$, $t'\in(0,t]$, and $\ev_{(p,0)}=\beta_{n,0}(q_n\cdot)_p$ (continuity of
		$t'\mapsto f_p(t')$), whose $K_0$-map on $G_n\subset K_0(F_n)$ is $(h,e)\mapsto(\Lambda_{n,0})_{p,\cdot}h$ by Lemma~\ref{lem:labelledprops}(d). A homogeneous-type
		$\varrho_j$ is homotopic through the representations $a\mapsto a(z_j(u))$, $u\in[0,1]$ (sheet transport), to $a\mapsto a(\theta_n)=q_n(a)_{\ib_n}$, whose $K_0$-map is
		$(h,e)\mapsto h_{\ib_n}$. $K_0$-maps are homotopy invariant (also for non-unital maps).
	\end{proof}
	
	\paragraph{$K_0$-corrections.} Let $\mathfrak g_n$ be the greatest common divisor of the entries of $\mathbf d_n\neq0$ (see \eqref{eq:dnnonzero}) and fix $\mathbf a_n\in\Z^{l_n}$ with
	$\mathbf a_n\cdot \mathbf d_n=\mathfrak g_n$. For $\varrho_j$ of interval type put $\mu_j:=\mathfrak g_n$ and $y_j:=z_j^{\mathrm{def}}\mathbf a_n$; for homogeneous type put $\mu_j:=1$ and
	$y_j:=0$. Then
	\begin{equation}\label{eq:K0cond}
		y_j\cdot \mathbf d_n=\mu_j\,z_j^{\mathrm{def}}\qquad\text{for all }j,\qquad\text{and we write }y:=\sum_jy_j\in\Z^{l_n}.
	\end{equation}
	The point is the following. We want to replace, in the target matrix summand $j_*$, the basepoint blocks of type $\mu_j\bb_j$ by $\mu_j$ copies of $\varrho_j$. On $K_0$ this
	changes nothing (Lemma~\ref{lem:Krho}), but the dimension increases by $\mu_jz_j^{\mathrm{def}}$. We compensate by removing basepoint blocks of type $y\Delta_n$: on $K_0(A_n)=G_n$
	this changes nothing, because $\Delta_nh=0$ there, and it decreases the dimension by $(y\Delta_n)\cdot r_n=y\cdot \mathbf d_n=\sum_j\mu_jz_j^{\mathrm{def}}$. This is possible exactly
	because $\mathbf d_n\neq0$; the multiplicity $\mu_j=\mathfrak g_n$ makes the equation solvable in integers.
	
\begin{figure}[ht]
		\centering
		\begin{tikzpicture}[scale=0.62]
			\draw[thick] (0,0) rectangle (10,-10);
			\fill[orange!25] (0,0) rectangle (2.2,-2.2); \node[font=\scriptsize] at (1.1,-1.1) {$\varrho_1$};
			\fill[orange!25] (2.2,-2.2) rectangle (4.4,-4.4); \node[font=\scriptsize] at (3.3,-3.3) {$\varrho_1$};
			\fill[violet!20] (4.4,-4.4) rectangle (5.8,-5.8); \node[font=\scriptsize] at (5.1,-5.1) {$\varrho_2$};
			\fill[gray!20] (5.8,-5.8) rectangle (10,-10); \node[font=\small,align=center] at (7.9,-7.9) {$\beta_{\mathrm{rest}}(q_n\cdot)$\\ type $M'_{\mathrm{rest}}\ge1$};
			\draw[decorate,decoration={brace,amplitude=5pt}] (0,0.3) -- node[above=4pt,font=\scriptsize]{$\varrho_j^{\oplus\mu_j}$: total dim.\ $\le\eta_nr_{n+1,j_*}$} (5.8,0.3);
			\node[font=\small,align=left] at (15.0,-5) {the summand $j_*$ of $\Psi_n(f,a)$\\ (up to the permutation $w$)\\[4pt]
				$\varrho_1$: interval type, $\mu_1=\mathfrak g_n$\\ $\varrho_2$: homogeneous type, $\mu_2=1$\\[4pt]
				removed: basepoint blocks of type\\ $\sum_j\mu_j\bb_j+y\Delta_n$};
		\end{tikzpicture}
		\caption{The modified matrix summand $j_*$ (schematic).}\label{fig:Psidesign}
	\end{figure}
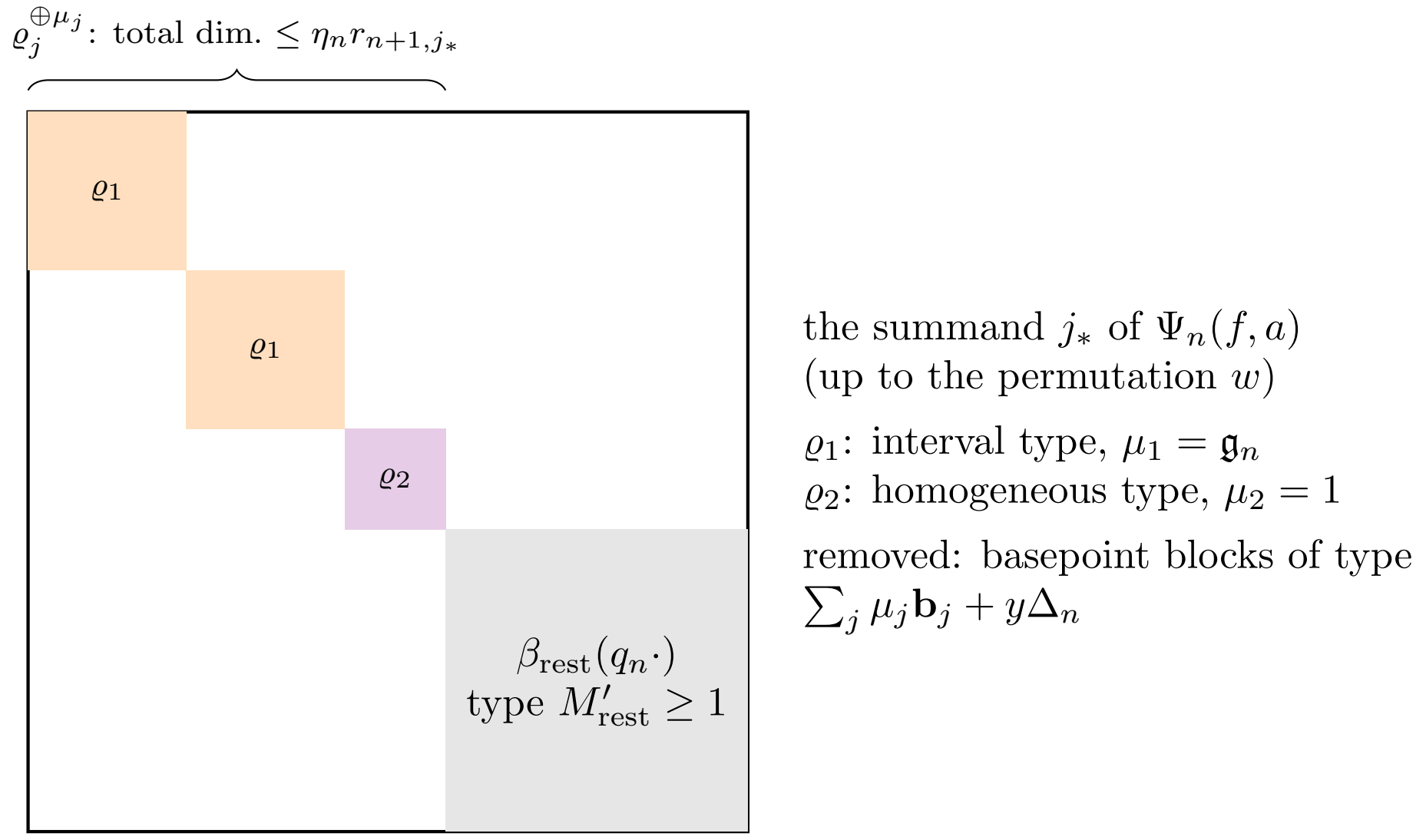

	\paragraph{The homogeneous component.} Choose an original target coordinate $j_*\in I_{n+1}\setminus\{\ib_{n+1}\}$ and require of $C_n$ that
	\begin{equation}\label{eq:Rone}
		\begin{gathered}
			M'_{\mathrm{rest}}:=(M_n)_{j_*,\cdot}-\sum_j\mu_j\bb_j-y\Delta_n\ \ge\ 1\quad\text{coordinatewise, and}\\
			\bar d'_n:=\sum_j\mu_j\dim\varrho_j+\|y\Delta_n\|_1\max_lr_{n,l}\ \le\ \eta_n\,r_{n+1,j_*}.
		\end{gathered}
	\end{equation}
	This holds once the entries $(C_n)_{j_*,\pi_n(l)}$ are large, because $\mu_j$, $\bb_j$, $y$, $\Delta_n$ and $\dim\varrho_j$ are determined at stage $n$ before $C_n$ is chosen
	(Section~\ref{ss:SPorder}). Define $\Psi_n\colon A_n\to F_{n+1}$ by $\Psi_n(f,a)_{l'}=\psi_n(a)_{l'}$ for $l'\neq j_*$ and
	\begin{equation}\label{eq:Psidef}
		\Psi_n(f,a)_{j_*}=\Ad(w)\Big(\bigoplus_j\varrho_j(f,a)^{\oplus\mu_j}\ \oplus\ \beta_{\mathrm{rest}}(q_na)\Big),
	\end{equation}
	where $\beta_{\mathrm{rest}}$ is a unital labelled block homomorphism of $\widehat F_n$ of type $M'_{\mathrm{rest}}$ and $w$ is a permutation matrix (Figure~\ref{fig:Psidesign}). The
	dimensions match: by \eqref{eq:spnu} and \eqref{eq:K0cond},
	\begin{align*}
		\sum_j\mu_j\dim\varrho_j+M'_{\mathrm{rest}}\cdot r_n&=\sum_j\mu_j(\bb_j\cdot r_n+z^{\mathrm{def}}_j)+(M_nr_n)_{j_*}-\sum_j\mu_j\bb_j\cdot r_n-y\cdot \mathbf d_n\\
		&=(M_nr_n)_{j_*}=r_{n+1,j_*} .
	\end{align*}
	
	\begin{lemma}\label{lem:Psi}
		$\Psi_n$ is a $*$-homomorphism with $\Psi_n(f,a)=0\Rightarrow a=0$, and $\Psi_n(e_n)$ is invertible in $F_{n+1}$. It maps diagonal, monomial and nonnegative monomial elements
		to elements of the same kind and commutes with diagonal compression. In all summands $l'\neq j_*$ we have $q_{n+1}\Psi_n=q_{n+1}\psi_n\pi^F_n=\widehat\psi_n\bar q_n$, where
		$\pi^F_n(f,a)=a$ is the quotient map $A_n\to F_n$, and on $K$-theory
		\[
		(\Psi_n)_*=(\psi_n\pi^F_n)_*\ \text{ on }K_0(A_n)=G_n\quad\text{and on }K_1(A_n)=K_n .
		\]
		Moreover $\Psi_n(\cdot)_{j_*}$ contains each $\varrho_j$ as a subrepresentation, and every basepoint block $q_n(\cdot)_l$ occurs in $\Psi_n(\cdot)_{j_*}$.
	\end{lemma}
	
	\begin{proof}
		The distinguished component of $\Psi_n$ is that of $\psi_n\pi^F_n$, which is injective on $F_n$ by Lemma~\ref{lem:psi}; so $\Psi_n(f,a)=0$ forces $a=0$. The preservation
		statements hold for $\psi_n$ by Lemma~\ref{lem:psi}, for the $\varrho_j$ by sheet transport, and for $\beta_{\mathrm{rest}}$ and $\Ad(w)$ by Lemma~\ref{lem:labelledprops}.
		$\Psi_n(e_n)=1$ in the summands $l'\neq j_*$, and in $j_*$ it is $\Ad(w)(\bigoplus_j\varrho_j(e_n)^{\oplus\mu_j}\oplus1)$, with $\varrho_j(e_n)$ invertible by
		Lemma~\ref{lem:An}(c). For $K_0$, let $x=(h,e)\in G_n$. The components $l'\neq j_*$ of $(\Psi_n)_*x$ and $(\psi_n\pi^F_n)_*x$ agree. In the matrix summand $j_*$ only ranks
		matter, and by Lemma~\ref{lem:Krho},
		\[
		\big((\Psi_n)_*x\big)_{j_*}=\sum_j\mu_j\,\bb_j\cdot h+M'_{\mathrm{rest}}\cdot h=(M_n)_{j_*,\cdot}h-(y\Delta_n)h=(M_nh)_{j_*}=\big((\psi_n\pi^F_n)_*x\big)_{j_*},
		\]
		since $\Delta_nh=0$. $K_1(A_n)\cong K_1(F_n)$ and $K_1(F_{n+1})$ live in the distinguished summands, where $\Psi_n=\psi_n\pi^F_n$. The last assertion holds because
		$M'_{\mathrm{rest}}\ge1$.
	\end{proof}
	
	\paragraph{The interval component.} Let $k^0_q:=(\Lambda_{n+1,0})_{q,j_*}$ and $k^1_q:=(\Lambda_{n+1,1})_{q,j_*}$, so that $k^1_q-k^0_q=(\Delta_{n+1})_{q,j_*}=:\dd_q$, where the
	strictly positive matrices $\Lambda_{n+1,i}$ with $\Lambda_{n+1,1}-\Lambda_{n+1,0}=\Delta_{n+1}$ are still to be fixed. The target endpoint map $\beta_{n+1,i}$ will contain $k^i_q$
	copies of the summand $j_*$, hence $\mu_jk^i_q$ copies of $\varrho_j$; the interval component must therefore contain exactly $\mu_jk^i_q$ blocks in the state $\varrho_j$ at $t=i$.
	We use interval data in the sense of Definition~\ref{def:intervaldata} with marked parameters $\mathfrak M=\{t_j:\ \varrho_j\text{ of interval type}\}$ and $p(t_j)=p_j$:
	\begin{enumerate}[label=\textup{(D\arabic*)}]
		\item For each interval-type $\varrho_j$: $\mu_j\min\{k^0_q,k^1_q\}$ copies of $f_{p_j}$ with $\phi_c\equiv t_j$; if $k^0_q>k^1_q$, in addition $\mu_j(k^0_q-k^1_q)$ copies with
		$(\phi_c(0),\phi_c(1))=(t_j,1)$; if $k^1_q>k^0_q$, in addition $\mu_j(k^1_q-k^0_q)$ copies with $(\phi_c(0),\phi_c(1))=(0,t_j)$. For the latter two kinds one may take
		$\phi_c(t)=\max\{(1-t)\phi_c(0)+t\phi_c(1),\,t(1-t)\}$. Thus exactly $\mu_jk^0_q$ copies are in the state $\varrho_j$ at $t=0$, and $\mu_jk^1_q$ at $t=1$.
		\item For each homogeneous-type $\varrho_j$: $\mu_j\max\{k^0_q,k^1_q\}$ designated copies of $q_n(\cdot)_{\ib_n}$ with path $z_j$, of which $\mu_j\min\{k^0_q,k^1_q\}$ have
		$(u_c(0),u_c(1))=(1,1)$ and the remaining $\mu_j|k^0_q-k^1_q|$ have $(u_c(0),u_c(1))=(1,0)$ if $k^0_q>k^1_q$, resp.\ $(0,1)$ if $k^1_q>k^0_q$; for instance
		$u_c(t)=\max\{(1-t)u_c(0)+tu_c(1),\,t(1-t)\}$. Thus exactly $\mu_jk^0_q$ (resp.\ $\mu_jk^1_q$) of them are in the state $\varrho_j$ at $t=0$ (resp.\ $t=1$).
		\item For each $p'\in[l_n]$: $\fw_q(p')\ge1$ forward copies and $\bw_q(p')\ge1$ backward copies of $f_{p'}$, with
		\begin{equation}\label{eq:Lprime}
			\fw_q(p')-\bw_q(p')=(L_n)_{q,p'}-\dd_q\,(y)_{p'}-\sum_{j\ \text{int},\ p_j=p'}\mu_j(k^0_q-k^1_q)_+ ,
		\end{equation}
		where $(y)_{p'}$ is the $p'$th coordinate of $y$ and the sum is over the interval-type $\varrho_j$ with $p_j=p'$.
		\item A vector $\cv_q\in\Z^{J_n}_{\ge1}$ of constant copies, containing the designated copies of (D2), determined by \eqref{eq:E0} below, with at least one undesignated copy of
		every $q_n(\cdot)_l$.
	\end{enumerate}
	Let $\Xi_n$ be the interval map of these data, with $s_{n+1,q}$ as in Definition~\ref{def:intervaldata}, and let $Z'^q_i$, $S'^q_i$ be as in Lemma~\ref{lem:intervaldata}. Write
	$\var_{n,q}=\sum_{p'}(\fw_q(p')+\bw_q(p'))s_{n,p'}$ for the total dimension of the forward and backward copies.
	
	The forward and backward copies play the role of $P_n$ and $N_n$ in the unital case (Section~\ref{sec:unitalmaps}): their difference \eqref{eq:Lprime} realizes the matrix
	$L_n$ on the quotient $\Z^{l_n}$, corrected for the $K_0$-correction $y$ and for the mixed copies of (D1).
	
	\begin{proposition}[endpoint bookkeeping]\label{prop:bookkeeping}
		For $i\in\{0,1\}$ let $N^i_j$ be the number of blocks in the state $\varrho_j$ at $t=i$, and let $\#^i_0(p')$, $\#^i_1(p')$ be the numbers of $f_{p'}$-copies with $\phi_c(i)=0$,
		resp.\ $\phi_c(i)=1$.
		\begin{enumerate}[label=\textup{(\alph*)}]
			\item (Type equations.) Put $x:=q_{n+1}\Psi_n(f,a)\in\widehat F_{n+1}$. The fibre $\Xi_n(f,a)_q(i)$ is, after a fixed permutation $w_{q,i}$, equal to
			$\beta^{\mathrm{std}}_{q,i}(x)$ for a labelled block homomorphism $\beta^{\mathrm{std}}_{q,i}$ of $\widehat F_{n+1}$ of type $((\Lambda_{n+1,i})_{q,l'})_{l'}$ into
			$M_{s_{n+1,q}}$, provided
			\begin{equation}\label{eq:M1}
				N^i_j=\mu_j(\Lambda_{n+1,i})_{q,j_*}\quad(\text{all }j),
			\end{equation}
			\begin{multline}\label{eq:E0}
				\sum_{p'}\Big[\#^i_0(p')(\Lambda_{n,0})_{p',\cdot}+\#^i_1(p')(\Lambda_{n,1})_{p',\cdot}\Big]+\cv_{q}-\sum_{j\ \text{hom}}N^i_jh_{\ib_n}\\
				=\sum_{l'\neq j_*}(\Lambda_{n+1,i})_{q,l'}(M_n)_{l',\cdot}+(\Lambda_{n+1,i})_{q,j_*}M'_{\mathrm{rest}} .
			\end{multline}
			In that case the support of $\beta^{\mathrm{std}}_{q,i}$ has cardinality $(\Lambda_{n+1,i}r_{n+1})_q=s_{n+1,q}-|Z'^q_i|$.
			\item With \textup{(D1)--(D3)}, equation \eqref{eq:M1} holds for $i=0,1$, and \eqref{eq:E0} for $i=1$ follows from \eqref{eq:E0} for $i=0$.
			\item Equation \eqref{eq:E0} for $i=0$ determines $\cv_q$. If $\Lambda_{n+1,i}$ is replaced by $\Lambda_{n+1,i}+\mathfrak e\,\mathbf 1$ (with $\mathbf 1$ the matrix all of whose
			entries are $1$), then for every $l\in J_n$ the number of undesignated constant copies of $q_n(\cdot)_l$ is $\mathfrak e$ times an integer $\ge2$, plus a term independent of
			$\mathfrak e$. Hence \textup{(D4)}, strict positivity of $\Lambda_{n+1,i}$ and any prescribed lower bound on the numbers of undesignated constant copies can be achieved by taking
			$\mathfrak e$ large.
		\end{enumerate}
		Consequently, with these choices and $\beta_{n+1,i}:=\Ad(w_{q,i})\circ\beta^{\mathrm{std}}_{q,i}$ on each summand $q$ (labelled block homomorphisms with supports $S'^q_i$, which
		satisfy \textup{(Z)} and \textup{(C)} by Lemma~\ref{lem:intervaldata}(b)), the exact endpoint identity
		\begin{equation}\label{eq:spendpoint}
			\Xi_n(f,a)(i)=\beta_{n+1,i}\,q_{n+1}\Psi_n(f,a),\qquad i=0,1,
		\end{equation}
		holds, and $\xi_n:=(\Xi_n,\Psi_n)$ is a $*$-homomorphism $A_n\to A_{n+1}$, where $A_{n+1}$ is formed with these $\beta_{n+1,i}$.
	\end{proposition}
	
	\begin{proof}
		(a) By Lemma~\ref{lem:intervaldata}(a), the fibre at $t=i$ is a direct sum of basepoint blocks $q_n(a)_l$, with the multiplicity vector on the left of \eqref{eq:E0}, of $N^i_j$
		blocks $\varrho_j(f,a)$, and of the zero coordinates $Z'^q_i$. On the other hand $\beta^{\mathrm{std}}_{q,i}(x)$ is a direct sum of $(\Lambda_{n+1,i})_{q,l'}$ copies of $x_{l'}$ for
		each $l'$, where for $l'\neq j_*$ the block $x_{l'}=(\widehat\psi_nq_na)_{l'}$ is a sum of basepoint blocks of type $(M_n)_{l',\cdot}$, and $x_{j_*}=\Psi_n(f,a)_{j_*}$ is, by
		\eqref{eq:Psidef}, a sum of $\mu_j$ blocks $\varrho_j(f,a)$ for each $j$ and of basepoint blocks of type $M'_{\mathrm{rest}}$. Two such direct sums are conjugate by a permutation
		matrix (respecting fixed internal labels of the blocks) as soon as the multiplicities of each kind of block agree, which is what \eqref{eq:M1} and \eqref{eq:E0} say. In that case
		the support of $\beta^{\mathrm{std}}_{q,i}$ has cardinality $\sum_{l'}(\Lambda_{n+1,i})_{q,l'}r_{n+1,l'}=(\Lambda_{n+1,i}r_{n+1})_q$, and this equals the number of nonzero coordinates of
		the fibre on $e_n$, namely $s_{n+1,q}-|Z'^q_i|$: by \eqref{eq:M1}, \eqref{eq:E0} and the dimension identity for the summand $j_*$, the total dimension of the $\varrho_j$-blocks and
		basepoint blocks of the fibre is $\sum_jN^i_j\dim\varrho_j+(\text{left side of \eqref{eq:E0}})\cdot r_n=\sum_{l'\neq j_*}(\Lambda_{n+1,i})_{q,l'}r_{n+1,l'}+(\Lambda_{n+1,i})_{q,j_*}
		\big(\sum_j\mu_j\dim\varrho_j+M'_{\mathrm{rest}}\cdot r_n\big)=(\Lambda_{n+1,i}r_{n+1})_q$.
		
		(b) Equation \eqref{eq:M1} holds by (D1) and (D2). Subtract \eqref{eq:E0} for $i=0$ from \eqref{eq:E0} for $i=1$. Since
		$\sum_{l'}(\Lambda_{n+1,1}-\Lambda_{n+1,0})_{q,l'}(M_n)_{l',\cdot}=(\Delta_{n+1}M_n)_{q,\cdot}=(L_n\Delta_n)_{q,\cdot}$ by \eqref{eq:quot}, and
		$(M_n)_{j_*,\cdot}-M'_{\mathrm{rest}}=\sum_j\mu_j\bb_j+y\Delta_n$, the difference of the right-hand sides is
		\begin{multline*}
			(L_n\Delta_n)_{q,\cdot}-\dd_q\Big(\sum_j\mu_j\bb_j+y\Delta_n\Big)\\
			=\sum_{p'}\big[(L_n)_{q,p'}-\dd_q(y)_{p'}\big](\Delta_n)_{p',\cdot}-\dd_q\sum_{j\ \text{int}}\mu_j(\Lambda_{n,0})_{p_j,\cdot}-\dd_q\sum_{j\ \text{hom}}\mu_jh_{\ib_n}.
		\end{multline*}
		On the left, the vectors $\cv_q$ cancel, and the designated terms contribute $-(N^1_j-N^0_j)h_{\ib_n}=-\dd_q\mu_jh_{\ib_n}$ for homogeneous $\varrho_j$, matching the last term. For
		the $f$-copies, (D1) and (D3) give
		\begin{align*}
			\#^0_0(p')&=\fw_q(p')+\sum_{j:\,p_j=p'}\mu_j(k^1_q-k^0_q)_+,& \#^0_1(p')&=\bw_q(p'),\\
			\#^1_1(p')&=\fw_q(p')+\sum_{j:\,p_j=p'}\mu_j(k^0_q-k^1_q)_+,& \#^1_0(p')&=\bw_q(p').
		\end{align*}
		Writing $(\Lambda_{n,1})_{p',\cdot}=(\Lambda_{n,0})_{p',\cdot}+(\Delta_n)_{p',\cdot}$, the remaining part of the difference of the left-hand sides is
		\[
		\sum_{p'}\Big[\big(\#^1_1(p')-\#^0_0(p')\big)(\Lambda_{n,0})_{p',\cdot}+\big(\#^1_1(p')-\bw_q(p')\big)(\Delta_n)_{p',\cdot}\Big].
		\]
		Here $\#^1_1(p')-\#^0_0(p')=\sum_{j:\,p_j=p'}\mu_j\big((k^0_q-k^1_q)_+-(k^1_q-k^0_q)_+\big)=-\dd_q\sum_{j\ \text{int},\,p_j=p'}\mu_j$, which matches the $\Lambda_{n,0}$-terms, and
		$\#^1_1(p')-\bw_q(p')=(L_n)_{q,p'}-\dd_q(y)_{p'}$ by \eqref{eq:Lprime}, which matches the $\Delta_n$-terms. Hence the two differences agree.
		
		(c) Solve \eqref{eq:E0} for $i=0$ for $\cv_q$. Replacing $\Lambda_{n+1,i}$ by $\Lambda_{n+1,i}+\mathfrak e\,\mathbf 1$ changes $k^0_q$ and $k^1_q$ by $\mathfrak e$ but not $\dd_q$,
		hence neither \eqref{eq:Lprime} nor the numbers $\fw_q,\bw_q$. The counts $\#^0_0(p')$ and $\#^0_1(p')$ are therefore independent of $\mathfrak e$, while $N^0_j=\mu_jk^0_q$ and the
		right-hand side grow linearly. Hence
		\[
		\cv_q=\mathfrak e\,v+(\text{fixed}),\qquad v=\sum_{l'\neq j_*}(M_n)_{l',\cdot}+M'_{\mathrm{rest}}+\sum_{j\ \text{hom}}\mu_jh_{\ib_n} ,
		\]
		where ``(fixed)'' denotes terms independent of $\mathfrak e$. Every column of $M_n$ has an entry $\ge1$ in an original row other than $j_*$ (as $|I_{n+1}|\ge2$), and
		$M'_{\mathrm{rest}}\ge1$, so
		\begin{equation}\label{eq:vbound}
			\sum_{l'\neq j_*}(M_n)_{l',l}+(M'_{\mathrm{rest}})_l\ \ge\ 2\qquad(l\in J_n).
		\end{equation}
		The designated copies of (D2) are copies of $q_n(\cdot)_{\ib_n}$ only, and their number is
		\[
		\textstyle\sum_{j\ \text{hom}}\mu_j\max\{k^0_q,k^1_q\}=\mathfrak e\sum_{j\ \text{hom}}\mu_j+(\text{fixed}).
		\]
		For $l\neq\ib_n$ all $\cv_{q,l}=\mathfrak e\,v_l+(\text{fixed})$ constant copies of $q_n(\cdot)_l$ are undesignated, and $v_l\ge2$ by \eqref{eq:vbound}. For $l=\ib_n$ the number of
		undesignated constant copies is $\mathfrak e\big(v_{\ib_n}-\sum_{j\ \text{hom}}\mu_j\big)+(\text{fixed})$, and $v_{\ib_n}-\sum_{j\ \text{hom}}\mu_j=\sum_{l'\neq
			j_*}(M_n)_{l',\ib_n}+(M'_{\mathrm{rest}})_{\ib_n}\ge2$ by \eqref{eq:vbound}.
		
		The final assertions follow from (a)--(c) and Lemma~\ref{lem:intervaldata}(b),(d).
	\end{proof}
	
	\begin{remark}\label{rem:detectionsize}
		The numbers of copies in (D1) and (D2) grow with the common part $\mathfrak e$ of $\Lambda_{n+1,i}$, but their total dimension in the fibre $q$ is at most
		\[
		\max\{k^0_q,k^1_q\}\sum_j\mu_j\dim\varrho_j\le\max_ik^i_q\,\eta_nr_{n+1,j_*}\le\eta_n\max_i(\Lambda_{n+1,i}r_{n+1})_q\le\eta_ns_{n+1,q}
		\]
		by \eqref{eq:Rone}, uniformly in $\mathfrak e$. If there are no detecting representations, $\Psi_n=\psi_n\pi^F_n$ and $\Xi_n$ consists of forward and backward copies and constant
		copies only.
	\end{remark}
	
	\begin{proposition}\label{prop:xi}
		The map $\xi_n=(\Xi_n,\Psi_n)\colon A_n\to A_{n+1}$ has the following properties.
		\begin{enumerate}[label=\textup{(\alph*)}]
			\item $\xi_n$ is injective, maps diagonal, monomial and nonnegative monomial elements to elements of the same kind, and $\EE_{n+1}\xi_n=\xi_n\EE_n$. Hence
			$\xi_n(B_n)\subset B_{n+1}$, $\xi_n(\Nor(B_n))\subset\Nor(B_{n+1})$ and $\xi_n(N_{+,n})\subset N_{+,n+1}$.
			\item $\xi_n(e_n)$ is strictly positive in $A_{n+1}$; its image in every irreducible representation of $A_{n+1}$ is invertible.
			\item $(\xi_n)_*=\gt_n|_{G_n}\colon K_0(A_n)=G_n\to K_0(A_{n+1})=G_{n+1}$, and $(\xi_n)_*=\chi_n$ on $K_1$.
			\item (Detection and spare copies.) $\Psi_n(\cdot)_{j_*}$ contains every $\varrho_j$ and every basepoint block of $A_n$. For every $l\in J_n$, the basepoint block $q_n(\cdot)_l$
			occurs in every irreducible representation of $A_{n+1}$ other than the extra matrix representations $x_i$ with $i\neq\pi_n(l)$; in particular it occurs in every interior
			interval fibre, at every point of the distinguished homogeneous summand, and in every original matrix summand.
			\item (Trace smallness.) Suppose $2/m_n\le\varepsilon_n$ and $\var_{n,q}\le\varepsilon_ns_{n+1,q}$ for all $q$. For every irreducible representation $\varpi$ of $A_{n+1}$ with
			normalized trace $\tau_\varpi$,
			\[
			\tau_\varpi\circ\xi_n=(1-\omega)\,\sigma'\circ\bar q_n+\omega\,\tau''\qquad\text{with }0\le\omega\le\varepsilon_n+\eta_n,\ \sigma'\in T(\widehat F_n),\ \tau''\in T(A_n).
			\]
			Moreover, $\|\sigma\circ\bar q_{n+1}\circ\xi_n-\sigma\circ\widehat\psi_n\circ\bar q_n\|\le2\eta_n$ for every $\sigma\in T_1(\widehat F_{n+1})$.
		\end{enumerate}
	\end{proposition}
	
	\begin{proof}
		(a) $\Psi_n$ detects $a$ (Lemma~\ref{lem:Psi}), and the forward copies in $\Xi_n$ run through all parameters of $f$; so $\xi_n$ is injective. The preservation statements are
		Lemma~\ref{lem:intervaldata}(d) and Lemma~\ref{lem:Psi}, and Lemma~\ref{lem:normalizers}(a) converts them into the statement on normalizers.
		
		(b) Interior interval fibres are invertible by Lemma~\ref{lem:intervaldata}(c), and homogeneous fibres by Lemma~\ref{lem:Psi}; apply Lemma~\ref{lem:fibrewise}.
		
		(c) The identifications $K_0(A_m)=G_m\subset K_0(F_m)$ of Lemma~\ref{lem:An}(d) are natural for the quotient maps, and the $F$-component of $\xi_n$ is $\Psi_n$ with
		$(\Psi_n)_*=(\psi_n\pi^F_n)_*$ (Lemma~\ref{lem:Psi}); by Lemma~\ref{lem:psi}, $(\psi_n)_*=\gt_n$ on $K_0$ and $\chi_n$ on $K_1$.
		
		(d) The first sentence is Lemma~\ref{lem:Psi}. Interior interval fibres contain undesignated constant copies of every $q_n(\cdot)_l$ by (D4). An original matrix summand
		$l'\notin\{\ib_{n+1},j_*\}$ contains $(M_n)_{l',l}=(C_n)_{l',\pi_n(l)}\ge2$ copies of $q_n(\cdot)_l$, and $j_*$ contains it by Lemma~\ref{lem:Psi}. In the distinguished summand, at
		every point of $Z_{n+1}$, the entries (i)--(iii) of $\psi_n$ contain $q_n(\cdot)_l$: for a matrix source $l$ the entry is, at every point, a direct sum of $\deg\nu_l\ge1$ copies of
		the source block, and for $l=\ib_n$ the $m_n-2\ge1$ constant sheets give the block $q_n(\cdot)_{\ib_n}$. The extra summand $x_i$ contains $q_n(\cdot)_l$ exactly when
		$\pi_n(l)=i$, by \eqref{eq:Mn}.
		
		(e) If $\varpi$ is a matrix summand $l'\neq j_*$, then $\tau_\varpi\circ\xi_n=\tau_\varpi\circ\psi_n\pi^F_n$ is a convex combination of normalized basepoint traces, and $\omega=0$.
		If $\varpi$ is the summand $j_*$, the blocks $\varrho_j^{\oplus\mu_j}$ have total dimension at most $\bar d'_n\le\eta_nr_{n+1,j_*}$ by \eqref{eq:Rone}, and the rest consists of
		basepoint blocks. If $\varpi$ is the evaluation at $z\in Z_{n+1}$ in the distinguished summand, then by the covering trace formula \eqref{eq:covtrace} every matrix-source
		entry of $\psi_n$ contributes exactly the normalized trace of its source block, and only the two sheets with the maps $\kappa_n,h_n$ in the distinguished-source entry are not
		basepoint evaluations; their share is at most $2/m_n\le\varepsilon_n$. If $\varpi$ is an interior interval evaluation $(q,t)$, the blocks which are not basepoint blocks are the
		forward and backward copies, of share at most $\varepsilon_n$, and the copies of (D1) and (D2), of share at most $\eta_n$ by Remark~\ref{rem:detectionsize}. In all cases
		$\omega\le\varepsilon_n+\eta_n$.
		
		For the last estimate, $\bar q_{n+1}\xi_n=q_{n+1}\Psi_n$ and $\widehat\psi_n\bar q_n=q_{n+1}\psi_n\pi^F_n$ differ only in the summand $j_*$. There the two representations
		$\Psi_n(\cdot)_{j_*}$ and $\psi_n\pi^F_n(\cdot)_{j_*}$ of $A_n$ have the same total dimension $r_{n+1,j_*}$ and share the basepoint blocks of type
		$\min\{(M_n)_{j_*,\cdot},M'_{\mathrm{rest}}\}$. Their unshared parts therefore have equal dimension; on the side of $\psi_n$ this is the dimension of the basepoint blocks of type
		$(M_n)_{j_*,\cdot}-\min\{(M_n)_{j_*,\cdot},M'_{\mathrm{rest}}\}\le\sum_j\mu_j\bb_j+|y\Delta_n|$, which is at most $\bar d'_n$. Hence the normalized traces of the two representations
		differ in norm by at most $2\bar d'_n/r_{n+1,j_*}\le2\eta_n$ (Lemma~\ref{lem:blockcompare}), and $\sigma$ puts weight at most $1$ on that summand.
	\end{proof}
	
	\subsection{The order of the choices}\label{ss:SPorder}
	The constructions above impose several size requirements. They can all be met in a single successive construction; this is the purpose of delaying the torsion
	corrections. The structure is the same as in the unital case (Section~\ref{ss:choices}, Table~\ref{tab:stage} and Figure~\ref{fig:induction}). Suppose stage $n$ has been
	completed. In particular we know the groups and lifts at that stage, the elements $\delta_n(x_i)$, the torsion stage $\mathsf S_n$ (the complex $T_n$, its generator coverings and
	the bound $\cc_n$), the covering $\lambda_n$, the $C^*$-algebra $A_n$ with both endpoint maps (hence the zero blocks and $\Lambda_{n,0}$), $\mathbf d_n$, the incoming maps, and the test elements
	chosen so far. Proceed as follows.
	\begin{enumerate}[label=(\roman*)]
		\item Choose the finitely many tests to be detected at this step (Section~\ref{ss:SPsimplicity}), the representations $\varrho_j$, the paths $z_j$, the multiplicities $\mu_j$
		and the corrections $y_j$. Their dimensions are now known.
		\item Use Theorem~\ref{thm:presentation} to choose the next matrix $C_n$ with sufficiently large entries: all entries at least $2$, $m_n\ge\max\{\cc_n,3\}$, $2/m_n\le\varepsilon_n$,
		and \eqref{eq:Rone} for every original target row in place of $j_*$. All these bounds use only completed data. Designate $\ib_{n+1}$ and choose $j_*\ne\ib_{n+1}$.
		\item Choose the next lifts $\ell_{n+1}$, compute $t_n$, and choose the finite group $E_{n+1}$ containing $E_n$, the values of $t_n$, and the next element of a fixed enumeration of
		$R$. Choose the next finitely generated group $K_{n+1}$ and the map $\chi_n$.
		\item Apply Theorem~\ref{thm:realize} to obtain $\mathsf S_{n+1}$ (in particular $T_{n+1}$, the generator coverings and $\cc_{n+1}$) and $g_n$; add $Y_{n+1}$ and the attached
		interval to obtain $Z_{n+1}$, and choose $p_{n+1}$. The new covering degrees need not satisfy any bound from step (ii): their corrections are delayed until the following step.
		\item Construct $\psi_n$, define $\lambda_{n+1}$ as the union of its image coverings (its class is the required $(r_{n+1},\zeta_{n+1})$), and construct $\Psi_n$ in the original
		matrix target $j_*$.
		\item Choose the quotient basis giving $\Delta_{n+1}$ and compute $L_n$. Choose $\fw_q,\bw_q\ge1$ with \eqref{eq:Lprime}, and then the common positive part $\mathfrak e$ of
		$\Lambda_{n+1,0},\Lambda_{n+1,1}$ so large that Proposition~\ref{prop:bookkeeping}(c) applies, that $\var_{n,q}\le\varepsilon_ns_{n+1,q}$ for all $q$ (note that $\var_{n,q}$ does not depend on
		$\mathfrak e$, while $s_{n+1,q}\ge\mathfrak e$), and that every $q_n(\cdot)_l$ has at least one undesignated constant copy. Determine $\cv_q$ by \eqref{eq:E0} and the permutations
		$w_{q,i}$, and form $\beta_{n+1,i}$, $A_{n+1}$ and $\xi_n$.
	\end{enumerate}
	The first stage is given by \eqref{eq:unit1}, with a trivial distinguished covering, endpoint matrices as in \eqref{eq:Bn}, and nested supports. Every coordinate of every $r_n$ is
	positive, so every matrix algebra used above has positive size. No requirement involves an unknown future quantity.
	
	\subsection{Simplicity, \texorpdfstring{$K$}{K}-theory and the trivial twist}\label{ss:SPsimplicity}
	Let $A_\infty=\lims(A_n,\xi_n)$, with the tests chosen as follows. At each stage $m$ fix a countable norm-dense family of positive elements of $A_m$, and let $\mathcal B_m$ consist of the
	nonzero cutdowns $(b-r)_+$, $r\in\Q_{>0}$, of its members; enumerate $\mathcal B_m$. At step $n$ the tests to be detected are the images in $A_n$ of the first $n$ elements of
	$\mathcal B_m$ for every $m\le n$; they are nonzero because the connecting maps are injective.
	
	\begin{proposition}\label{prop:SPsimple}
		$A_\infty$ is simple, separable and nuclear. For every $m$, $\xi_{m,\infty}(e_m)$ is strictly positive in $A_\infty$, and the maps $\xi_{m,k}$ ($k>m$) and $\xi_{m,\infty}$ send strictly positive
		elements to strictly positive elements.
	\end{proposition}
	
	\begin{proof}
		\emph{Strict positivity.} By Proposition~\ref{prop:xi}(b), induction and Remark~\ref{rem:normpreserving}, $\xi_{m,k}(e_m)$ is strictly positive in $A_k$ for all $k\ge m$. Put
		$x=\xi_{m,\infty}(e_m)$. For $k\ge m$, $\overline{xAx}\supset\xi_{k,\infty}\big(\overline{\xi_{m,k}(e_m)A_k\xi_{m,k}(e_m)}\big)=\xi_{k,\infty}(A_k)$, so $\overline{xAx}$ contains the dense
		subalgebra $\bigcup_k\xi_{k,\infty}(A_k)$ and equals $A_\infty$. The last assertion follows from Remark~\ref{rem:normpreserving}.
		
		\emph{Propagation.} Let $b\in A_n$ be a positive element whose image in $F_{n+1}$ is nonzero in the original matrix summand $j_*$ chosen at step $n$; by Lemma~\ref{lem:Psi} this
		holds for every test detected at step $n$ through some $\varrho_j$. By Proposition~\ref{prop:xi}(d) applied at step $n+1$ to $l=j_*$, the image of $b$ in $A_{n+2}$ is nonzero in
		every irreducible representation except possibly the extra matrix representations $x_i$, $i\neq j_*$, and nonzero in $x_{j_*}$. Applying Proposition~\ref{prop:xi}(d) at step $n+2$
		to the original coordinates $l\in I_{n+2}$, at which the image is nonzero, we find that the image of $b$ in $A_{n+3}$ is nonzero in every irreducible representation. The same
		argument, started one step later, applies to a test detected at step $n$ through a matrix representation $a\mapsto a_l$: by Proposition~\ref{prop:xi}(d) its image in $A_{n+1}$ is
		nonzero in all irreducible representations except possibly the extra representations $x_i$ with $i\neq\pi_n(l)$, and its image in $A_{n+2}$ is nonzero everywhere. A positive
		element of a $C^*$-algebra which is nonzero in every irreducible representation is full (a proper closed ideal is contained in the kernel of an irreducible representation
		\cite[\S5.4]{Mur}).
		
		\emph{Fullness in the limit.} If $b$ is full in $A_k$, the closed ideal of $A_\infty$ generated by $\xi_{k,\infty}(b)$ contains $\xi_{k,\infty}(A_k)\ni\xi_{k,\infty}(e_k)$, which is strictly
		positive in $A_\infty$, so the ideal is $A_\infty$.
		
		\emph{Simplicity.} Let $\mathcal J\ne0$ be a closed ideal and $y\in\mathcal J$ positive with $\|y\|=1$. As in the proof of Proposition~\ref{prop:simple} (Step~5), there are $m$ and a
		member $b$ of the dense family of stage $m$ with $\|\xi_{m,\infty}(b)-y\|<\frac14$; with the quotient map $\varpi\colon A_\infty\to A_\infty/\mathcal J$ we get
		$\varpi(\xi_{m,\infty}((b-\frac12)_+))=(\varpi(\xi_{m,\infty}(b))-\frac12)_+=0$, so $\xi_{m,\infty}((b-\frac12)_+)\in\mathcal J$, and it is nonzero because $\|b\|>\frac34$. The element
		$(b-\frac12)_+$ belongs to $\mathcal B_m$, so it is detected at some step, and by the propagation argument its image is full at a later stage. Hence $\mathcal J=A_\infty$. Separability and
		nuclearity pass to inductive limits of separable subhomogeneous $C^*$-algebras.
	\end{proof}
	
	It is essential that tests from all earlier stages are detected; a sequence of finer nets in the current stage alone would not give the density argument in the limit.
	
	\begin{corollary}\label{cor:SPK}
		$K_0(A_\infty)\cong G_0$ and $K_1(A_\infty)\cong G_1$. The $C^*$-algebra $A_\infty$ is stably projectionless, non-unital, and satisfies the UCT.
	\end{corollary}
	
	\begin{proof}
		By Proposition~\ref{prop:xi}(c) and continuity of $K$-theory, $K_0(A_\infty)=\lims(G_n,\gt_n)=G_0$ and $K_1(A_\infty)=\lims(K_n,\chi_n)=G_1$. A projection in $M_k(A_\infty)$ is unitarily equivalent to
		the image of a projection in some $M_k(A_n)$ (apply functional calculus to a nearby self-adjoint element of $M_k(A_n)$), and $A_n$ is stably projectionless by
		Lemma~\ref{lem:An}(e); in particular $A_\infty$ has no unit. Each $A_n$ is an extension of a homogeneous algebra over a finite CW complex by $C_0((0,1),E'_n)$, hence of type~I and in the
		bootstrap class, and so is the inductive limit; thus $A_\infty$ satisfies the UCT \cite{RS}.
	\end{proof}
	
	\begin{proposition}\label{prop:limittwist}
		Let $B=\lims B_n\subset A_\infty$. Then
		\begin{equation}\label{eq:AG}
			(A_\infty,B)\cong\big(C^*_r(\bar{\mathcal G}),C_0(\bar{\mathcal G}^{(0)})\big)
		\end{equation}
		for a second countable, locally compact Hausdorff, principal étale groupoid $\bar{\mathcal G}$, and $\bar{\mathcal G}^{(0)}$ is not compact.
	\end{proposition}
	
	\begin{proof}
		By Lemmas~\ref{lem:groupoid} and~\ref{lem:normalizers} and Proposition~\ref{prop:xi}(a), the pairs $(A_n,B_n)$, the maps $\xi_n$ and the sets $S_n=N_{+,n}$ satisfy the hypotheses
		of Proposition~\ref{prop:twist}, which does not require units. The unit space is not compact because $A_\infty$ is not unital.
	\end{proof}
	%======================================================================
	\section{The stably projectionless case: traces and the range theorem}\label{sec:SPtraces}
	
	Throughout, $A_\infty=\lims(A_n,\xi_n)$ is the $C^*$-algebra of Section~\ref{sec:SPmaps}, constructed with $\varepsilon_n=\eta_n=2^{-n-6}$. In particular the hypotheses $2/m_n\le\varepsilon_n$
	and $\var_{n,q}\le\varepsilon_ns_{n+1,q}$ of Proposition~\ref{prop:xi}(e) hold; they were imposed in steps (ii) and (vi) of Section~\ref{ss:SPorder}. In this section we
	identify the trace simplex of $A_\infty$ with $T$, compatibly with the pairing (Section~\ref{ss:SPAF}), show that all densely finite traces are bounded (Section~\ref{ss:SPbounded}),
	tensor with $\JS$ (Section~\ref{ss:SPAZ}) and prove continuous scale and Theorem~\ref{thm:SPB} (Section~\ref{ss:SPcs}).
	
	Compared with the unital case (Section~\ref{sec:traces}), two differences matter. The trace spaces are now the sets $T_1(\cdot)$ of tracial functionals of norm at most
	one, and one has to check that norms are preserved; this is where strict positivity of $\xi_n(e_n)$ is used. And the comparison with the AF algebra is done directly for
	the modified maps, by a quantitative estimate, which will also be used for the scale in Section~\ref{sec:SPscale}.
	
	\subsection{Trace spaces of non-unital inductive limits}
	For a $*$-homomorphism $\theta\colon X\to Y$ sending a strictly positive element to a strictly positive element, the pullback $\theta^*\tau=\tau\circ\theta$ is a norm-preserving
	(Remark~\ref{rem:normpreserving}), affine, weak$^*$-continuous map $T_1(Y)\to T_1(X)$, contractive for the dual norms and mapping $T(Y)$ into $T(X)$.
	
	\begin{lemma}\label{lem:invlim}
		Let $X=\lims(X_n,\theta_n)$ with injective $\theta_n$ sending strictly positive elements to strictly positive elements. Then $\tau\mapsto(\tau\circ\theta_{n,\infty})_n$ is a
		norm-preserving affine homeomorphism $T_1(X)\to\limp(T_1(X_n),\theta_n^*)$, and it restricts to a homeomorphism $T(X)\to\limp(T(X_n),\theta_n^*)$.
	\end{lemma}
	
	\begin{proof}
		The maps $\theta_{n,\infty}$ send strictly positive elements to strictly positive elements (the argument in the proof of Proposition~\ref{prop:SPsimple} is general), so restriction
		preserves norms and lands in the inverse limit. Conversely, a compatible family $(\tau_n)$ with $\|\tau_n\|\le1$ defines a tracial functional of norm at most $1$ on the dense
		$*$-subalgebra $\bigcup_n\theta_{n,\infty}(X_n)$, which extends uniquely to $X$; its norm is $\lim_n\|\tau_n\|=\|\tau_1\|$. The bijection is affine and weak$^*$-continuous, and its
		inverse is continuous because on the bounded set $T_1(X)$ pointwise convergence on a dense subalgebra implies weak$^*$-convergence. The norm-one parts correspond.
	\end{proof}
	
	\begin{lemma}\label{lem:extreme}
		Let $X$ be a separable subhomogeneous $C^*$-algebra. The extreme points of $T_1(X)$ are $0$ and the functionals $\tau_\pi=\frac1{\dim\pi}\Tr\circ\pi$, $\pi$ irreducible.
	\end{lemma}
	
	\begin{proof}
		An extreme point $\tau\neq0$ of $T_1(X)$ has norm one (otherwise $\tau=\|\tau\|\cdot(\tau/\|\tau\|)+(1-\|\tau\|)\cdot0$) and is an extreme point of $T(X)$. For a tracial state $\tau$,
		the cyclic vector of the GNS representation $\pi_\tau$ defines a faithful normal tracial state on $\pi_\tau(X)''$, and $\tau$ is extremal in $T(X)$ if and only if $\pi_\tau(X)''$ is a
		factor \cite[\S6.8]{Dix}. All irreducible representations of $X$ have dimension at most some $N$, so $X$ satisfies the standard polynomial identity of degree $2N$
		(Amitsur--Levitzki), and so does $\pi_\tau(X)''$, by Kaplansky density and strong$^*$-continuity of multiplication on bounded sets. A factor satisfying a polynomial identity is a
		matrix algebra $M_d$ with $d\le N$; hence $\pi_\tau$ is a multiple of an irreducible representation $\pi$ of dimension $d$, and $\tau=\tau_\pi$. Conversely each $\tau_\pi$ is extremal,
		since $\pi_{\tau_\pi}(X)''\cong M_d$ is a factor.
	\end{proof}
	
	\subsection{Comparison with the AF algebra}\label{ss:SPAF}
	Recall $\bar q_n\colon A_n\to\widehat F_n$, $\bar q_n(f,a)=q_n(a)$. Since $\bar q_n$ is surjective, the pullback $\bar q_n^{\,*}\colon T_1(\widehat F_n)\to T_1(A_n)$,
	$\sigma\mapsto\sigma\circ\bar q_n$, is isometric for the dual norms, affine and weak$^*$-continuous; since $\bar q_n(e_n)=1$, it preserves norms. The connecting maps do not
	satisfy $\bar q_{n+1}\xi_n=\widehat\psi_n\bar q_n$ exactly (the two sides differ in the summand $j_*$), so we work with two estimates.
	
	\begin{lemma}[almost quotient]\label{lem:SPalmostquot}
		For every $\tau\in T_1(A_{n+1})$ there is $\sigma\in T_1(\widehat F_n)$ with $\|\tau\circ\xi_n-\bar q_n^{\,*}\sigma\|\le2(\varepsilon_n+\eta_n)$.
	\end{lemma}
	
	\begin{proof}
		For $\tau=0$ take $\sigma=0$. For the normalized trace $\tau_\varpi$ of an irreducible representation $\varpi$, Proposition~\ref{prop:xi}(e) gives
		$\tau_\varpi\circ\xi_n=(1-\omega)\bar q_n^{\,*}\sigma'+\omega\tau''$ with $\omega\le\varepsilon_n+\eta_n$, so $\|\tau_\varpi\circ\xi_n-\bar q_n^{\,*}\sigma'\|\le2\omega$. For general $\tau$
		consider $f(\tau)=\operatorname{dist}\big(\tau\circ\xi_n,\bar q_n^{\,*}(T_1(\widehat F_n))\big)$ on $T_1(A_{n+1})$. It is convex, since $\xi_n^*$ is affine and $\bar
		q_n^{\,*}(T_1(\widehat F_n))$ is convex, and it is weak$^*$-lower semicontinuous: the dual norm is weak$^*$-lower semicontinuous, $(\tau,\sigma)\mapsto\tau\circ\xi_n-\bar
		q_n^{\,*}\sigma$ is weak$^*$-continuous, and an infimum over the compact set $T_1(\widehat F_n)$ of a jointly lower semicontinuous function is lower semicontinuous. It is at most
		$2(\varepsilon_n+\eta_n)$ at every extreme point of $T_1(A_{n+1})$ (Lemma~\ref{lem:extreme}), hence on their convex hull by convexity, hence on its closure $T_1(A_{n+1})$
		(Krein--Milman) by lower semicontinuity. The infimum is attained by compactness.
	\end{proof}
	
	\begin{lemma}[approximate compatibility]\label{lem:approxcompat}
		For every $\sigma\in T_1(\widehat F_{n+1})$, $\|\xi_n^*\bar q_{n+1}^{\,*}\sigma-\bar q_n^{\,*}\widehat\psi_n^*\sigma\|\le2\eta_n$.
	\end{lemma}
	
	\begin{proof}
		This is the last estimate of Proposition~\ref{prop:xi}(e), since $\xi_n^*\bar q_{n+1}^{\,*}\sigma=\sigma\circ\bar q_{n+1}\circ\xi_n$ and $\bar q_n^{\,*}\widehat\psi_n^*\sigma=\sigma\circ\widehat\psi_n\circ\bar q_n$.
	\end{proof}
	
	\begin{figure}[ht]
		\centering
		\begin{tikzcd}[column sep=4em,row sep=3em]
			\cdots\arrow[r] & T_1(A_{n+1})\arrow[r,"\xi_n^*"] & T_1(A_n)\arrow[r] & \cdots & T_1(A_\infty)\\
			\cdots\arrow[r] & T_1(\widehat F_{n+1})\arrow[r,"\widehat\psi_n^*"']\arrow[u,"\bar q_{n+1}^{\,*}"] & T_1(\widehat F_n)\arrow[r]\arrow[u,"\bar q_n^{\,*}"'] & \cdots & T_1(\widehat F)\arrow[u,"\Pi^{-1}"']
		\end{tikzcd}
		\caption{The squares commute up to $2\eta_n$ (Lemma~\ref{lem:approxcompat}); summability gives the homeomorphism $\Pi$ of the limits (Theorem~\ref{thm:SPtraces}).}\label{fig:SPtracesquares}
	\end{figure}
	
	\begin{theorem}\label{thm:SPtraces}
		There is an affine homeomorphism $\Pi\colon T_1(A_\infty)=\limp(T_1(A_n),\xi_n^*)\to T_1(\widehat F)=\limp(T_1(\widehat F_n),\widehat\psi_n^*)$, preserving norms, characterized by
		\begin{equation}\label{eq:quant}
			\|\tau_n-\bar q_n^{\,*}\sigma_n\|\le\vartheta_n:=\sum_{k\ge n}2\eta_k\qquad(n\ge1)\quad\text{for }\tau=(\tau_n),\ \Pi(\tau)=\sigma=(\sigma_n).
		\end{equation}
		In particular $T(A_\infty)\cong T(\widehat F)\cong T$ affinely and homeomorphically, and $T(A_\infty)$ is weak$^*$-compact. Under this identification the pairing $\rho_{A_\infty}\colon K_0(A_\infty)=G_0\to\Aff(T(A_\infty))$
		is $\rho$.
	\end{theorem}
	
	\begin{proof}
		\emph{Existence.} Let $\tau=(\tau_n)$ be compatible. For each $m$ choose $\sigma^{(m)}\in T_1(\widehat F_m)$ with $\|\tau_m-\bar q_m^{\,*}\sigma^{(m)}\|\le2(\varepsilon_m+\eta_m)$
		(Lemma~\ref{lem:SPalmostquot} applied to $\tau_{m+1}$, as $\tau_m=\xi_m^*\tau_{m+1}$), and put $\sigma^{(m)}_n=\widehat\psi^*_{n,m}\sigma^{(m)}$ for $n\le m$. Since all pullbacks are
		contractions and $\bar q_m^{\,*}$ is isometric,
		\begin{align*}
			\|\sigma^{(m)}-\widehat\psi_m^*\sigma^{(m+1)}\|&=\|\bar q_m^{\,*}\sigma^{(m)}-\bar q_m^{\,*}\widehat\psi_m^*\sigma^{(m+1)}\|\\
			&\le\|\bar q_m^{\,*}\sigma^{(m)}-\tau_m\|+\|\xi_m^*(\tau_{m+1}-\bar q_{m+1}^{\,*}\sigma^{(m+1)})\|\\
			&\qquad+\|\xi_m^*\bar q_{m+1}^{\,*}\sigma^{(m+1)}-\bar q_m^{\,*}\widehat\psi_m^*\sigma^{(m+1)}\|\\
			&\le2(\varepsilon_m+\eta_m)+2(\varepsilon_{m+1}+\eta_{m+1})+2\eta_m
		\end{align*}
		by Lemma~\ref{lem:approxcompat}. Hence $(\sigma^{(m)}_n)_{m\ge n}$ is norm-Cauchy for each $n$, with limit $\sigma_n\in T_1(\widehat F_n)$ (a norm-closed set), and
		$\widehat\psi_n^*\sigma_{n+1}=\sigma_n$ by continuity. Moreover, for $m>n$,
		\begin{align*}
			\|\tau_n-\bar q_n^{\,*}\sigma^{(m)}_n\|&\le\|\xi^*_{n,m}(\tau_m-\bar q_m^{\,*}\sigma^{(m)})\|+\sum_{k=n}^{m-1}\|\xi^*_{n,k}(\xi_k^*\bar q_{k+1}^{\,*}-\bar q_k^{\,*}\widehat\psi_k^*)\widehat\psi^*_{k+1,m}\sigma^{(m)}\|\\
			&\le2(\varepsilon_m+\eta_m)+\sum_{k=n}^{m-1}2\eta_k ,
		\end{align*}
		and letting $m\to\infty$ gives \eqref{eq:quant}.
		
		\emph{Uniqueness and affinity.} If $\sigma,\sigma'\in T_1(\widehat F)$ both satisfy \eqref{eq:quant} for the same $\tau$, then $\|\sigma_n-\sigma'_n\|\le2\vartheta_n$, and since
		$\sigma_n=\widehat\psi^*_{n,m}\sigma_m$, also $\|\sigma_n-\sigma'_n\|\le\|\sigma_m-\sigma'_m\|\le2\vartheta_m\to0$. So $\Pi$ is well defined, and it is affine because convex combinations of
		solutions of \eqref{eq:quant} are solutions. It preserves norms: $\|\sigma_n\|=\|\bar q_n^{\,*}\sigma_n\|$ and $\big|\|\tau_n\|-\|\bar q_n^{\,*}\sigma_n\|\big|\le\vartheta_n\to0$, while
		$\|\tau_n\|$ and $\|\sigma_n\|$ do not depend on $n$.
		
		\emph{Injectivity.} If $\Pi\tau=\Pi\tau'$, then $\|\tau_n-\tau'_n\|\le2\vartheta_n$ for all $n$, and $\tau_n=\xi^*_{n,m}\tau_m$ gives $\|\tau_n-\tau'_n\|\le2\vartheta_m\to0$.
		
		\emph{Surjectivity and continuity of the inverse.} Given $\sigma\in T_1(\widehat F)$, put $t^{(m)}_n:=\xi^*_{n,m}\bar q_m^{\,*}\sigma_m$ for $m\ge n$. By Lemma~\ref{lem:approxcompat},
		$\|t^{(m+1)}_n-t^{(m)}_n\|\le\|\xi^*_m\bar q_{m+1}^{\,*}\sigma_{m+1}-\bar q_m^{\,*}\widehat\psi_m^*\sigma_{m+1}\|\le2\eta_m$, so $t^{(m)}_n\to\tau_n$ in norm, uniformly in $\sigma$. The
		family $(\tau_n)$ is compatible, lies in $T_1(A_n)$ (a norm-closed set), and satisfies $\|\tau_n-\bar q_n^{\,*}\sigma_n\|\le\sum_{k\ge n}2\eta_k$, so $\Pi\tau=\sigma$. Each coordinate
		$\sigma\mapsto\tau_n$ is a uniform norm limit of weak$^*$-continuous maps, hence weak$^*$-continuous. Thus $\Pi^{-1}$ is a continuous affine bijection from the compact space
		$T_1(\widehat F)$ onto the Hausdorff space $\limp T_1(A_n)=T_1(A_\infty)$ (Lemma~\ref{lem:invlim}), hence a homeomorphism. Norm-one elements correspond, which gives $T(A_\infty)\cong T(\widehat
		F)\cong T$ by Lemma~\ref{lem:AFtraces}.
		
		\emph{Pairing.} Let $x\in K_0(A_\infty)$ be represented at stage $n_0$ by $x_{n_0}\in G_{n_0}$, that is, by projections $p,q\in M_k(\widetilde A_{n_0})$ with equal scalar parts, and let
		$x_n=(\xi_{n_0,n})_*x_{n_0}$, represented by $\xi_{n_0,n}(p)$, $\xi_{n_0,n}(q)$. By Proposition~\ref{prop:xi}(c), $(\xi_n)_*=\gt_n$ on $G_n$, so the free part $h_n$ of $x_n$ (its
		image in $K_0(\widehat F_n)=\Z^{J_n}$) is the image of $h_{n_0}$ under the connecting maps of $(\Z^{J_m},M_m)$; this $K$-theoretic input does not follow from the trace estimates. For
		$\tau\in T(A_\infty)$ with $\Pi\tau=\sigma$ and $n\ge n_0$, naturality gives $\rho_{A_\infty}(x)(\tau)=\rho_{A_n}(x_n)(\tau_n)=(\tau_n\otimes\Tr)(\xi_{n_0,n}(p-q))$, and
		$|(\tau_n\otimes\Tr)(z)-(\bar q_n^{\,*}\sigma_n\otimes\Tr)(z)|\le k\vartheta_n\|z\|$ for $z\in M_k(A_n)$, while $\|p-q\|\le1$. Further, $(\bar
		q_n^{\,*}\sigma_n\otimes\Tr)(\xi_{n_0,n}(p-q))$ is the pairing of $\sigma_n$ with the image of $x_n$ in $K_0(\widehat F_n)=\Z^{J_n}$, that is, with $h_n$, and this is the value
		$s(\bar x)$ of the state $s\in S(F,\bar u)$ corresponding to $\sigma$ at the image $\bar x\in F$ of $x\in G_0\subset H$, independently of $n$. Letting $n\to\infty$ gives
		$\rho_{A_\infty}(x)(\tau)=s(\bar x)$. Finally $s=s_{\tau'}$ for the point $\tau'\in T$ corresponding to $\sigma$, and $s_{\tau'}(\overline{(g,0)})=\theta(g,0)(\tau')=\rho(g)(\tau')$ by
		Theorem~\ref{thm:Fdimgroup}.
	\end{proof}
	
	\subsection{No unbounded densely finite traces}\label{ss:SPbounded}
	To use strict comparison with respect to tracial states, and to identify the cone $\Tt$, one needs to know that densely finite lower semicontinuous traces are bounded. Here
	the conditions (Z) and (C) on the endpoint maps enter.
	
	\begin{proposition}\label{prop:bounded}
		Every lower semicontinuous trace on $A_\infty$ which is finite on $\Ped(A_\infty)$ is bounded: $\Tt(A_\infty)=\R_{\ge0}\,T(A_\infty)$. The same holds for $A_\infty\otimes\JS$.
	\end{proposition}
	
	\begin{proof}
		Let $\tau\in\Tt(A_\infty)$ and fix $n$. Put $\tau_n=\tau\circ\xi_{n,\infty}$, an lsc trace on $A_n$. It is finite on $\Ped(A_n)$: the Pedersen ideal is the two-sided ideal generated by the
		elements $g(x)$, $x\in(A_n)_+$, $g\in C_c((0,\infty))$ \cite[5.6.1]{Ped}; moreover $\tau_n(g(x))=\tau(g(\xi_{n,\infty}x))<\infty$ because $g(\xi_{n,\infty}x)\in\Ped(A_\infty)$, and the positive
		elements on which an lsc trace is finite span a two-sided ideal \cite[6.1.2]{Dix}.
		
		On the ideal $\mathcal I_n=C_0((0,1),E'_n)=\bigoplus_qC_0((0,1),M_{s_{n,q}})$, an lsc trace which is finite on $\Ped(\mathcal I_n)=C_c((0,1),E'_n)$ has the form
		$\tau_n(z)=\sum_q\int_{(0,1)}\Tr(z_q(t))\,d\mu_q(t)$ for Radon measures $\mu_q$ on $(0,1)$, finite on compact subsets (compare the proof of Lemma~\ref{lem:tracesblock}). Choose
		$t_0\in(0,\frac12)$ and $\chi\in C_c((0,1])$ with $0\le\chi\le1$ and $\chi=1$ on $[t_0,1]$, and let $\chi'\in C_c((0,1])$ be $1$ on the support of $\chi$; then
		$\chi(e_n)\chi'(e_n)=\chi(e_n)$, so $\chi(e_n)\in\Ped(A_n)$ and $\tau_n(\chi(e_n))<\infty$. By Lemma~\ref{lem:An}(c) the fibre $e_n(t)$, $t\in(0,1)$, is diagonal with entry $1$ on the
		common support $S^q_{n,0}\cap S^q_{n,1}$, so $\Tr\big(\chi(e_n)_q(t)\big)\ge c_q:=|S^q_{n,0}\cap S^q_{n,1}|$, and $c_q>0$ by (C). For $z\in(A_n)_+$ and an approximate unit $(j_k)$ of
		$\mathcal I_n$ of the form $j_k=g_k\cdot1$ with $g_k\in C_c((0,1))$, $g_k\nearrow1$, monotonicity of $\tau_n$ and monotone convergence give
		$\tau_n(z)\ge\sup_k\tau_n(z^{1/2}j_kz^{1/2})=\sum_q\int\Tr(z_q(t))\,d\mu_q(t)$, with equality for $z\in(\mathcal I_n)_+$ by lower semicontinuity. Applied to $z=\chi(e_n)$ this yields
		$\sum_qc_q\,\mu_q((0,1))\le\tau_n(\chi(e_n))<\infty$: all $\mu_q$ are finite measures.
		
		Now $\|\tau_n\|=\sup_k\tau_n(e_n^{1/k})$. Write $e_n^{1/k}=(\chi h_k)(e_n)+((1-\chi)h_k)(e_n)$ with $h_k(x)=x^{1/k}$. The first summand is at most $\chi(e_n)$, so its trace is at most
		$\tau_n(\chi(e_n))$. For the second, the function $(1-\chi)h_k$ vanishes on $[t_0,1]$; so its $F_n$-component is $(1-\chi)(1)h_k(1)=0$, it lies in $\mathcal I_n$, and by
		Lemma~\ref{lem:An}(c) it is supported on the coordinates of $Z^q_{n,0}$ for $t<t_0$ and of $Z^q_{n,1}$ for $t>1-t_0$, with entries at most $1$. Its trace is therefore at most
		$\sum_q\big(|Z^q_{n,0}|\mu_q((0,t_0))+|Z^q_{n,1}|\mu_q((1-t_0,1))\big)<\infty$, uniformly in $k$. Hence $\tau_n$ is bounded.
		
		Finally, $\xi_{n,\infty}(e_n)$ is strictly positive in $A_\infty$ (Proposition~\ref{prop:SPsimple}), so $(\xi_{n,\infty}(e_n)^{1/k})_k$ is an approximate unit of $A_\infty$. For an lsc trace $\tau$ and
		an approximate unit $(u_k)$ we have $\tau(x)=\lim_k\tau(x^{1/2}u_kx^{1/2})\le\liminf_k\tau(u_k)$ for $0\le x\le1$, so $\|\tau\|\le\sup_k\tau(u_k)\le\|\tau\|$. Hence
		$\|\tau\|=\sup_k\tau_n(e_n^{1/k})=\|\tau_n\|<\infty$.
		
		For $A_\infty\otimes\JS$: $\JS$ is unital, so $A_\infty\otimes1\subset A_\infty\otimes\JS$ is a copy of $A_\infty$, $g(x)\otimes1=g(x\otimes1)\in\Ped(A_\infty\otimes\JS)$ for $x\in (A_\infty)_+$ and $g\in C_c((0,\infty))$, and
		$e\otimes1$ is strictly positive in $A_\infty\otimes\JS$ when $e$ is strictly positive in $A_\infty$. An lsc trace $\tau'$ on $A_\infty\otimes\JS$ which is finite on $\Ped(A_\infty\otimes\JS)$ therefore restricts to an
		lsc trace on $A_\infty\otimes1\cong A_\infty$ which is finite on $\Ped(A_\infty)$, hence bounded by the first part, and $\|\tau'\|=d_{\tau'}(e\otimes1)=\|\tau'|_{A_\infty\otimes1}\|<\infty$.
	\end{proof}
	
	\begin{remark}\label{rem:SPC}
		Condition (C) gives the uniform cutoff estimate used above; it is a convenient sufficient condition, not a necessary one. If (Z) holds and both endpoint supports are nonempty in
		every summand, then at every $t\in(0,1)$ some diagonal entry of $e_n(t)$ is at least $\frac12$: the entries are $1$, $1-t$ and $t$ on $S^q_{n,0}\cap S^q_{n,1}$, $Z^q_{n,1}$ and $Z^q_{n,0}$,
		and if (C) fails in the summand $q$, then $Z^q_{n,1}=S^q_{n,0}$ and $Z^q_{n,0}=S^q_{n,1}$ are both nonempty. The argument above, with $\chi=1$ on $[\frac12,1]$, then shows that every
		densely finite lsc trace is bounded. A small example in which (C) fails is
		\[
		\{(f,a)\in C([0,1],M_3)\oplus\C:\ f(0)=\diag(a,0,0),\ f(1)=\diag(0,a,a)\},
		\]
		with strictly positive element $(t\mapsto\diag(1-t,t,t),1)$. In more general pullbacks allowing an empty endpoint support, unbounded densely finite traces can occur: if
		$S^q_{n,0}$ were empty, then $(f,a)\mapsto\int_0^{1/2}\Tr f_q(t)\,t^{-2}\,dt$ would be such a trace. In our building blocks both endpoint supports are nonempty, since $\Lambda_{n,i}$
		and $r_n$ are strictly positive.
	\end{remark}
	
	\subsection{Tensoring with the Jiang--Su algebra}\label{ss:SPAZ}
	Let $\mathcal G_\JS$ be Li's untwisted model of $\JS$ (Theorem~\ref{thm:LiZ}). Then
	\[
	A_\infty\otimes\JS\cong C^*_r(\bar{\mathcal G}\times\mathcal G_\JS),
	\]
	using \eqref{eq:AG} and Lemma~\ref{lem:product} (whose proof does not use compactness of the unit spaces).
	
	\begin{proposition}\label{prop:SPAZ}
		$A_\infty\otimes\JS$ is separable, simple, nuclear, $\JS$-stable and stably projectionless, satisfies the UCT, and has nuclear dimension at most one. Moreover:
		\begin{enumerate}[label=\textup{(\alph*)}]
			\item every bounded trace on $A_\infty\otimes\JS$ has the form $\tau_{A_\infty}\otimes\tau_\JS$, and $\tau\mapsto\tau_{A_\infty}$ is an affine homeomorphism $T(A_\infty\otimes\JS)\cong T(A_\infty)\cong T$; under $K_0(A_\infty\otimes\JS)=K_0(A_\infty)=G_0$ the
			pairing is $\rho$, and $K_1(A_\infty\otimes\JS)=G_1$;
			\item $\Tt(A_\infty\otimes\JS)=\R_{\ge0}T(A_\infty\otimes\JS)$.
		\end{enumerate}
	\end{proposition}
	
	\begin{proof}
		Separability, simplicity (a minimal tensor product of simple $C^*$-algebras) and nuclearity are clear; $\JS$-stability holds since $\JS\otimes\JS\cong\JS$ \cite{JS}. The building blocks
		$A_n$ are of type~I, so $A_\infty$ belongs to the bootstrap class of Rosenberg and Schochet, which contains the type~I $C^*$-algebras and is closed under inductive limits
		\cite{RS}, \cite[\S22]{Bl}; so does $\JS$ \cite{JS}, and the bootstrap class is closed under tensor products \cite[\S22]{Bl}. Hence $A_\infty\otimes\JS$ satisfies the UCT. The
		K\"unneth theorem \cite{RS}, together with $K_*(\JS)=K_*(\C)$ \cite{JS}, gives $K_*(A_\infty\otimes\JS)=K_*(A_\infty)=(G_0,G_1)$. The nuclear dimension is at most one by Castillejos--Evington
		\cite{CE}, which applies to simple separable nuclear $\JS$-stable $C^*$-algebras, unital or not.
		
		(a) Let $\tau$ be a bounded trace on $A_\infty\otimes\JS$ and $a\in (A_\infty)_+$. Since $\JS$ is unital, $z\mapsto\tau(a\otimes z)$ is a bounded tracial positive functional on $\JS$, hence a multiple of
		the unique tracial state $\tau_\JS$: $\tau(a\otimes z)=\tau(a\otimes1)\tau_\JS(z)$. By linearity and density $\tau=\tau_{A_\infty}\otimes\tau_\JS$ with $\tau_{A_\infty}(a)=\tau(a\otimes1)$, and
		$\|\tau_{A_\infty}\|=\|\tau\|$ (use the approximate unit $e^{1/k}\otimes1$). Conversely $\tau_{A_\infty}\otimes\tau_\JS$ is a trace. The bijection is affine and weak$^*$-continuous in both directions, so
		$T(A_\infty\otimes\JS)\cong T(A_\infty)\cong T$ by Theorem~\ref{thm:SPtraces}. The pairing is unchanged because $x\mapsto x\otimes[1]$ is the K\"unneth isomorphism $K_0(A_\infty)\to K_0(A_\infty\otimes\JS)$ and
		$(\tau_{A_\infty}\otimes\tau_\JS)(p\otimes1)=\tau_{A_\infty}(p)$.
		
		Stable projectionlessness: a nonzero projection $p\in M_k(A_\infty\otimes\JS)$ would satisfy $\tau(p)>0$ for every $\tau\in T(A_\infty\otimes\JS)$, since traces on a simple $C^*$-algebra are faithful, and
		$[p]\in K_0(A_\infty\otimes\JS)=G_0$ with $\rho([p])(\tau)=\tau(p)$; thus $\rho([p])\in\Aff_{++}(T)$, contradicting \eqref{eq:P}.
		
		(b) is Proposition~\ref{prop:bounded}.
	\end{proof}
	
	\subsection{Continuous scale and the range theorem}\label{ss:SPcs}
	
	\begin{proposition}\label{prop:cscale}
		$A_\infty\otimes\JS$ has strict comparison of positive elements with respect to its tracial states, and it has continuous scale.
	\end{proposition}
	
	\begin{proof}
		\emph{Strict comparison.} The $C^*$-algebra $A_\infty\otimes\JS$ is simple, separable, exact and $\JS$-stable, so every lsc $2$-quasitrace on it is a trace (\cite{Haa}, \cite[Remark~2.29(i)]{BK}), and by
		Proposition~\ref{prop:SPAZ}(b) every densely finite lsc trace is a nonnegative multiple of an element of $T(A_\infty\otimes\JS)\neq\emptyset$. Hence \cite[Proposition~1.5]{CE}, which is based on
		R{\o}rdam's almost unperforation \cite[Theorem~4.5]{Ro}, gives: if $a,b\in M_m(A_\infty\otimes\JS)_+$ and $d_\tau(a)<d_\tau(b)$ for all $\tau\in T(A_\infty\otimes\JS)$, then $a\precsim b$.
		
		\emph{Continuous scale.} Let $e\in A_\infty\otimes\JS$ be strictly positive with $\|e\|=1$, and put $e_k=f_{2^{-k}}(e)$, so that $e_{k+1}e_k=e_k$. For $\tau\in T(A_\infty\otimes\JS)$ we have $\tau(e_k)\nearrow
		d_\tau(e)=\|\tau\|=1$, and by Dini's theorem on the compact space $T(A_\infty\otimes\JS)$ the convergence is uniform. For $m>k\ge2$, $d_\tau(e_m-e_k)\le\tau(e_{m+1}-e_{k-1})$, because
		$f_{2^{-m}}-f_{2^{-k}}$ is supported where $f_{2^{-m-1}}-f_{2^{-k+1}}=1$. Given $b\in(A_\infty\otimes\JS)_+\setminus\{0\}$, the lower semicontinuous function $\tau\mapsto d_\tau(b)$ is strictly positive on
		the compact set $T(A_\infty\otimes\JS)$, since traces are faithful, so $\beta:=\min_\tau d_\tau(b)>0$. Choosing $N$ with $\tau(e_j)>1-\beta$ for all $j\ge N$ and all $\tau\in T(A_\infty\otimes\JS)$, we get
		$d_\tau(e_m-e_k)<\beta\le d_\tau(b)$ for $m>k>N$ and all $\tau$, hence $e_m-e_k\precsim b$ by strict comparison.
	\end{proof}
	
	\begin{proof}[Proof of Theorem~\ref{thm:SPB}]
		Let $(G_0,T,\rho,G_1)$ be as in Theorem~\ref{thm:SPB}, construct $A_\infty$ as in Sections~\ref{sec:alg}--\ref{sec:SPmaps} in Case (SP), and put $G_c=\bar{\mathcal G}\times\mathcal G_\JS$ and
		$C=C^*_r(G_c)\cong A_\infty\otimes\JS$. By Propositions~\ref{prop:SPAZ} and~\ref{prop:cscale}, $C$ is a classifiable stably projectionless $C^*$-algebra with continuous scale, and its invariant is
		$(G_0,T,\rho,G_1)$. Every densely finite lsc trace on $C$ is bounded by Proposition~\ref{prop:SPAZ}(b); this also follows from continuous scale by Proposition~\ref{prop:cs}(b).
	\end{proof}
	
	In particular statement (I) of Section~\ref{ss:reduction} holds with $G_c=\bar{\mathcal G}\times\mathcal G_\JS$. Combined with \cite[Corollary~14.7]{GL4} (which classifies stably
	projectionless $C^*$-algebras with continuous scale by $(K_0,T,\rho,K_1)$), Theorem~\ref{thm:SPB} shows that the hereditary subalgebra $A_c$ of Section~\ref{ss:reduction} is isomorphic to
	$C^*(G_c)$ for suitable $G_c$; this already proves the stably projectionless case of Theorem~\ref{thm:main} for $C^*$-algebras with continuous scale.
	%======================================================================
	\section{The stably projectionless case: the scale and the proof}\label{sec:SPscale}
	
	The $C^*$-algebra $C$ of Theorem~\ref{thm:SPB} has continuous scale. A general classifiable stably projectionless $C^*$-algebra $A$ is stably isomorphic to such a $C^*$-algebra, and it is
	recovered as a hereditary subalgebra of $C\otimes\K$ determined by its scale. In this section we realize an arbitrary scale by an element $a'$ of the canonical diagonal of
	$C\otimes\K$. Then the hereditary subalgebra generated by $a'$ is again a groupoid $C^*$-algebra, of the restriction of the groupoid to the open set $\{a'>0\}$, and classification
	identifies it with $A$.
	
	\subsubsection*{Setting}
	Let $A_\infty=\lims(A_n,\xi_n)$ be the $C^*$-algebra of Sections~\ref{sec:SPmaps}--\ref{sec:SPtraces}, $B=\lims B_n=C_0(\bar{\mathcal G}^{(0)})$ its diagonal, and $C=A_\infty\otimes\JS=C^*_r(G_c)$ with
	$G_c=\bar{\mathcal G}\times\mathcal G_\JS$. Let $\mathcal P_\N$ be the pair groupoid of $\N$, so that $\K=C^*(\mathcal P_\N)$ with diagonal $c_0(\N)$, and
	\begin{equation}\label{eq:Gs}
		C\otimes\K\cong C^*_r(G^s),\qquad G^s:=\bar{\mathcal G}\times\mathcal G_\JS\times\mathcal P_\N,\qquad C_0((G^s)^{(0)})=B\otimes C(\mathcal G_\JS^{(0)})\otimes c_0(\N).
	\end{equation}
	We regard $C([k])$ as the algebra of diagonal matrices in $M_k$. Throughout we identify $T(A_\infty)$ with $T$ as follows: $\tau\in T(A_\infty)$ corresponds to $\Pi(\tau)\in T(\widehat F)$
	(Theorem~\ref{thm:SPtraces}), which corresponds to a state $s$ of $(F,\bar u)$ by Lemma~\ref{lem:AFtraces}, and $s=s_{\tau'}$ for a unique $\tau'\in T$ (Theorem~\ref{thm:Fdimgroup});
	we write $f(\tau):=f(\tau')$ for functions $f$ on $T$. Via Proposition~\ref{prop:SPAZ}(a) we also identify $T(A_\infty\otimes\JS)$ with $T$.
	
	\subsubsection*{The difficulty}
	Every lower semicontinuous affine function $\gamma\colon T\to(0,\infty]$ is an increasing limit of continuous affine functions, and one would like to realize these by
	diagonal elements $z_i$ in orthogonal blocks and add them up. The dimension function of a diagonal element $d$ is the measure of its open support, $d_\tau(d)=\mu_\tau(\{d>0\})$,
	which is only lower semicontinuous in $\tau$ in general. To run the induction, the partial sums must have \emph{continuous} dimension functions. We achieve this by choosing the
	supports so that their boundaries avoid a fixed countable set on which all the relevant measures are concentrated.
	
	\subsection{Traces at a finite stage are atomic on a fixed countable set}
	Let $\Omega_n=\mathcal G_n^{(0)}$ be the unit space of stage $n$ (Lemma~\ref{lem:groupoid}), so that $B_n=C_0(\Omega_n)$. Every bounded trace $\tau_n$ on $A_n$ restricts to a finite
	positive measure $\mu_{\tau_n}$ on $\Omega_n$ (Riesz representation), and for a positive diagonal element $d\in B_n\otimes C([k])=C_0(\Omega_n\times[k])$,
	\begin{equation}\label{eq:dmu}
		d_{\tau_n\otimes\Tr}(d)=(\mu_{\tau_n}\otimes\#)(\{d>0\}),
	\end{equation}
	where $\#$ is counting measure on $[k]$ (monotone convergence applied to $d^{1/j}\nearrow\chi_{\{d>0\}}$).
	
	\begin{lemma}\label{lem:atomic}
		There is a countable set $\mathcal S_n\subset\Omega_n$, depending only on the construction, such that for every $\tau\in T(A_\infty)$ the measure $\mu_{\tau_n}$, $\tau_n=\tau\circ\xi_{n,\infty}$, is
		concentrated on $\mathcal S_n$.
	\end{lemma}
	
	\begin{proof}
		By Theorem~\ref{thm:SPtraces}, for $M>n$ we have $\tau_n=\tau_M\circ\xi_{n,M}$ and $\|\tau_M-\sigma_M\circ\bar q_M\|\le\vartheta_M\to0$, so $\tau_n$ is the norm limit of the functionals
		$\sigma_M\circ\bar q_M\circ\xi_{n,M}$. Restricted to $B_M=C_0(\Omega_M)$, the functional $\sigma_M\circ\bar q_M$ is a positive combination of point masses at the finitely many units
		$w\in W^\theta_M$ (Lemma~\ref{lem:groupoid}). Since $\xi_{n,M}(B_n)\subset B_M$, each composite $\ev_w\circ\xi_{n,M}|_{B_n}$ is a character of $B_n$ or zero, that is, a point mass at a point
		of $\Omega_n$ or $0$. Thus the restriction of $\sigma_M\circ\bar q_M\circ\xi_{n,M}$ to $B_n$ is a measure supported on a finite set $\mathcal S_n^{(M)}\subset\Omega_n$ which depends only on
		the construction. Put $\mathcal S_n=\bigcup_{M>n}\mathcal S_n^{(M)}$. A norm limit of measures concentrated on $\mathcal S_n$ is concentrated on $\mathcal S_n$: if $\|\mu_k-\mu\|\to0$, then
		$\mu(\Omega_n\setminus\mathcal S_n)=\lim_k\mu_k(\Omega_n\setminus\mathcal S_n)=0$.
	\end{proof}
	
	\subsection{Approximate realization with continuous dimension functions}
	
	\begin{lemma}\label{lem:approxscale}
		Let $g\colon T\to(0,\infty)$ be continuous and affine, and let $\varepsilon>0$. There are $k\ge1$ and $z\in B\otimes C([k])\subset A_\infty\otimes M_k$ with $z\ge0$ and $\|z\|=1$, such that
		$\tau\mapsto d_\tau(z)$ is a \emph{continuous} affine function on $T(A_\infty)\cong T$ and $g-\varepsilon<d_\bullet(z)<g+\varepsilon$.
	\end{lemma}
	
\begin{figure}[ht]
		\centering
		\begin{tikzpicture}[scale=1.0]
			\draw[thick,fill=blue!4] (0,0) ellipse (4.2 and 1.7);
			\node at (-3.2,1.2) {$\Omega_n\times[k]$};
			\foreach \x/\y in {-1.8/0.2,0.6/0.5,2.2/-0.4}{\draw[thick,green!50!black,fill=green!15] (\x,\y) circle (0.55); \fill[green!50!black] (\x,\y) circle (2pt);}
			\foreach \x/\y in {-0.6/-0.6,1.4/0.9}{\fill[red!70!black] (\x,\y) circle (2pt);}
			\foreach \x/\y in {-2.6/-0.7,-1.25/0.75,-0.1/0.1,0.05/0.6,1.0/-0.2,2.75/-0.3,1.85/0.05,-1.0/-1.1,3.1/0.6,-2.3/0.5}{\draw[gray] (\x-0.05,\y-0.05) -- (\x+0.05,\y+0.05); \draw[gray] (\x-0.05,\y+0.05) -- (\x+0.05,\y-0.05);}
			\node[font=\small,align=left] at (7.6,0.3) {{\color{green!50!black}$\bullet$} points of $P_1$ (where $c=1$)\\ {\color{red!70!black}$\bullet$} points of $P_0$ (where $c=0$)\\
				{\color{gray}$\times$} the countable set $\mathcal S_n\times[k]$\\[3pt] $U$: union of the green balls;\\ $\partial U$ avoids $\mathcal S_n\times[k]$};
		\end{tikzpicture}
		\caption{The lift $d$ in the proof of Lemma~\ref{lem:approxscale}: $\{d>0\}=U$, and the boundary of $U$ carries no mass for any trace.}\label{fig:lift}
	\end{figure}
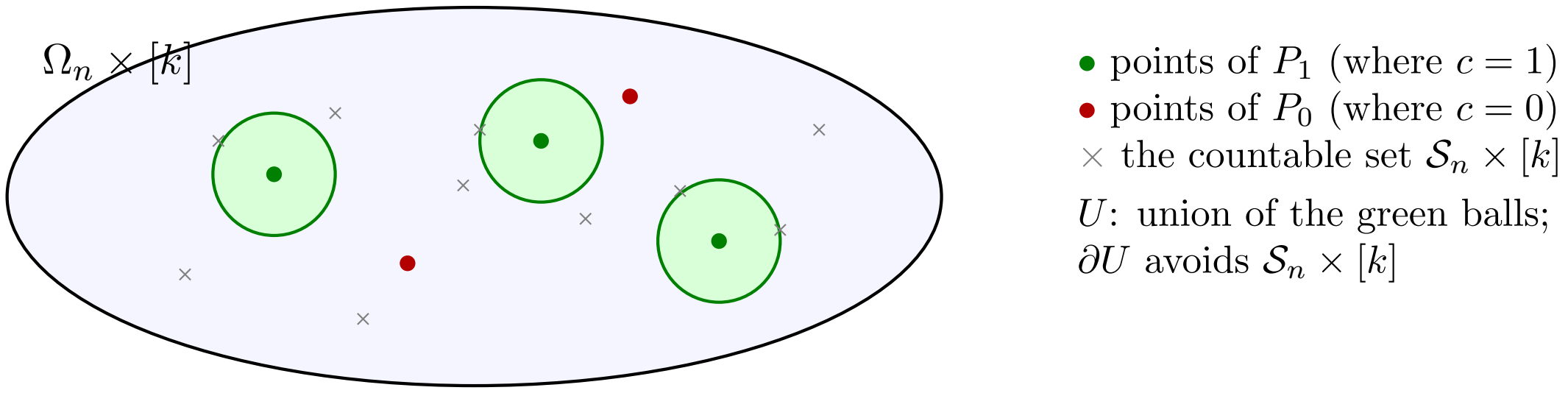

	\begin{proof}
		\emph{The target in the AF algebra.} Since $G'$ is dense in $\Aff(T)$, there is $f\in G'$ with $\|f-g\|<\min\{\varepsilon/2,\min_Tg\}$; in particular $f>0$. Then
		$(0,f)\in H^+\setminus\{0\}$, and its image $\bar f\in F^+\setminus\{0\}$ is represented at some stage by a nonzero vector $h\in(\Z_{\ge0})^{J_n}$, by Lemma~\ref{lem:delayed}. Since $\bar u$ is
		an order unit, $\bar f\le N\bar u$ for some $N$, and after passing to a later stage $h\le Nr_n$ coordinatewise. Passing to a still later stage (re-representing $h$ there does not
		change $N$), we fix $n$ with $N\vartheta_n<\varepsilon/2$; everything below is built at this stage. For $\tau\in T(A_\infty)$ let $\sigma=(\sigma_m)=\Pi(\tau)\in T(\widehat F)$; under the
		identifications just described, the state of $K_0(\widehat F_n)=\Z^{J_n}$ induced by $\sigma_n$ sends $h$ to the value of the corresponding state of $(F,\bar u)$ at $\bar f$, which is
		$f(\tau)$. Put $k=N$ and let $c$ be the projection in the diagonal of $\widehat F_n\otimes M_k$ with $\operatorname{rank}(c_l)=h_l$ in $M_{r_{n,l}}\otimes M_k$; then
		$d_{\sigma_n\otimes\Tr}(c)=\sum_lh_l\,\sigma_n(e^{\min}_l)=f(\tau)$, where $e^{\min}_l$ is a minimal projection of the $l$th summand.
		
		\emph{The lift.} Via $B_n=C_0(\Omega_n)$, the map $\bar q_n|_{B_n}\otimes\mathrm{id}$ is restriction to the finite closed subset $W^\theta_n\times[k]\subset\Omega_n\times[k]$
		(Lemma~\ref{lem:groupoid}), and $c$ is the indicator function of a subset $P_1\subset W^\theta_n\times[k]$; let $P_0=(W^\theta_n\times[k])\setminus P_1$. Choose a metric on the locally
		compact metrizable space $\Omega_n\times[k]$ and open balls $U=\bigcup_{p\in P_1}\mathrm{B}(p,r_p)$ with compact closures avoiding $P_0$, such that
		$\partial U\cap(\mathcal S_n\times[k])=\emptyset$ (Figure~\ref{fig:lift}). This is possible because $\mathcal S_n$ is countable and the spheres $\{x:\operatorname{dist}(x,p)=r\}$ for distinct $r$
		are disjoint, so all but countably many radii are admissible. Let $d\in C_0(\Omega_n\times[k])_+$ with $\{d>0\}=U$, $\|d\|=1$ and $d=1$ on $P_1$, for instance
		$d=\min\{1,\operatorname{dist}(\cdot,(\Omega_n\times[k])\setminus U)/r'\}$ for small $r'>0$. Then $(\bar q_n\otimes\mathrm{id})(d)=c$, and $d\in B_n\otimes C([k])$ is a positive diagonal element
		of norm one. Put $z=(\xi_{n,\infty}\otimes\mathrm{id})(d)\in B\otimes C([k])$.
		
		\emph{Continuity.} For $\tau\in T(A_\infty)$, $d_\tau(z)=\sup_j\tau(z^{1/j})=\sup_j(\tau_n\otimes\Tr)(d^{1/j})=d_{\tau_n\otimes\Tr}(d)=(\mu_{\tau_n}\otimes\#)(U)$ by \eqref{eq:dmu}. By
		Lemma~\ref{lem:atomic}, the measure $\mu_{\tau_n}\otimes\#$ gives no mass to $\partial U$, so $(\mu_{\tau_n}\otimes\#)(U)=(\mu_{\tau_n}\otimes\#)(\overline U)$. The map
		$\tau\mapsto\mu_{\tau_n}$ is weak$^*$-continuous; the measure of an open set is lower semicontinuous and that of a compact set is upper semicontinuous in the measure; hence
		$\tau\mapsto d_\tau(z)$ is continuous. It is affine because dimension functions are affine in the trace.
		
		\emph{The estimate.} We have $\|\tau_n-\sigma_n\circ\bar q_n\|\le\vartheta_n$ by \eqref{eq:quant}. For bounded positive functionals $\phi,\phi'$ on $A_n$ the sequences
		$(\phi\otimes\Tr)(d^{1/j})$ and $(\phi'\otimes\Tr)(d^{1/j})$ increase and converge, so $|d_{\phi\otimes\Tr}(d)-d_{\phi'\otimes\Tr}(d)|=\lim_j|((\phi-\phi')\otimes\Tr)(d^{1/j})|\le
		k\|\phi-\phi'\|$. Hence $|d_\tau(z)-d_{(\sigma_n\circ\bar q_n)\otimes\Tr}(d)|\le k\vartheta_n$. Moreover, for a functional factoring through a surjection,
		$d_{\sigma\circ\bar q}(d)=\lim_j\sigma(\bar q(d)^{1/j})=d_\sigma(\bar q(d))$, so $d_{(\sigma_n\circ\bar q_n)\otimes\Tr}(d)=d_{\sigma_n\otimes\Tr}(c)=f(\tau)$. As $k=N$ and
		$N\vartheta_n<\varepsilon/2$, we obtain $|d_\tau(z)-g(\tau)|\le k\vartheta_n+\|f-g\|<\varepsilon$.
	\end{proof}
	
	\begin{remark}\label{rem:SPLi}
		In Li's construction of the scale \cite[Lemma~7.3, Proposition~7.2]{Li}, an element $z$ is produced in a similar way, but the continuity of $d_\bullet(z)$ is not addressed. The
		inductive argument of \cite[Proposition~7.2]{Li} applies \cite[Lemma~7.3]{Li} to targets of the form $f_{i+1}-\sum_{h\le i}d_\bullet(z_h)$, which must be continuous and affine for the
		lemma to apply. In our setting continuity is secured by the choice of the lift $d$, whose boundary avoids the countable set $\mathcal S_n$.
	\end{remark}
	
	\begin{proposition}\label{prop:exactscale}
		For every lower semicontinuous affine function $\gamma\colon T\to(0,\infty]$ there is $a\in(B\otimes c_0(\N))_+\subset A_\infty\otimes\K$ with $d_\tau(a)=\gamma(\tau)$ for all $\tau\in T(A_\infty)\cong T$.
	\end{proposition}
	
	\begin{proof}
		\emph{The case $\gamma\equiv\infty$.} (This case occurs in Section~\ref{ss:proofA} when $A$ is stable.) Let $e:=\xi_{1,\infty}(e_1)$. Then $e\in B$, since $e_1\in B_1$
		(Lemma~\ref{lem:An}(c)) and the connecting maps preserve the diagonals (Proposition~\ref{prop:xi}(a)); $\|e\|\le1$; and $e$ is strictly positive in $A_\infty$ by
		Proposition~\ref{prop:SPsimple}. Hence $(e^{1/k})_k$ is an approximate unit of $A_\infty$, and $d_\tau(e)=\lim_k\tau(e^{1/k})=\|\tau\|=1$ for every $\tau\in T(A_\infty)$. Put
		\[
		a=\sum_{j=1}^\infty2^{-j}\,e\otimes e_{jj}\ \in\ (B\otimes c_0(\N))_+ ,
		\]
		a norm convergent sum. For $\varepsilon>0$ the element $f_\varepsilon(a)=\sum_jf_\varepsilon(2^{-j}e)\otimes e_{jj}$ has only finitely many nonzero terms, and $\tau(f_\varepsilon(2^{-j}e))$
		increases to $d_\tau(2^{-j}e)=d_\tau(e)=1$ as $\varepsilon\searrow0$. By additivity on orthogonal summands and monotone convergence, $d_\tau(a)=\sum_jd_\tau(e)=\infty=\gamma(\tau)$.
		
		\emph{An increasing sequence of continuous affine functions with gaps.} From now on $\gamma\not\equiv\infty$. Since $\gamma$ is lower semicontinuous with values in $(0,\infty]$ on the compact
		set $T$ and finite somewhere, $0<m_0:=\min_T\gamma<\infty$. A lower semicontinuous affine function $T\to(-\infty,\infty]$ is the pointwise supremum of its continuous affine minorants
		(the Hahn--Banach theorem applied to the closed convex epigraph; see \cite[Chapter~I, \S1]{Alf}). Hence $\gamma$ is the pointwise supremum of the continuous affine functions
		$g<\gamma$. If $\mathcal D\subset\Aff(T)$ is countable and dense, then also $\sup\{g\in\mathcal D:g<\gamma\}=\gamma$ (given a continuous affine $g<\gamma$ and $\varepsilon>0$, pick
		$g'\in\mathcal D$ with $\|g'-(g-2\varepsilon)\|<\varepsilon$); enumerate these functions as $(g_i)_{i\ge1}$. We construct continuous affine functions $f_i$ with
		\[
		f_1\ge m_0/2,\qquad f_i\ge g_i,\qquad f_i-f_{i-1}\ge\alpha_i>0\ \text{ on }T,\qquad f_i<\gamma\ \text{ on }T .
		\]
		Suppose $f_{i-1}$ has been constructed (put $f_0=m_0/4$). The function $v_i:=\max\{f_{i-1},g_i,m_0/2\}$ is continuous, convex and $v_i<\gamma$ pointwise, so $\gamma-v_i$ is lower
		semicontinuous with values in $(0,\infty]$, and $2\alpha_i:=\min\{1,\min_T(\gamma-v_i)\}>0$. Edwards' separation theorem on the Choquet simplex $T$ (an upper semicontinuous convex function
		below a lower semicontinuous concave one can be separated by a continuous affine function; see \cite[Chapter~II, \S3]{Alf} and Theorem~\ref{thm:edwards}), applied to the continuous
		convex function $v_i+\alpha_i$ and the lower semicontinuous concave real-valued function $\min\{\gamma,\max_Tv_i+2\}-\alpha_i$, which dominates $v_i+\alpha_i$ because $v_i+2\alpha_i\le\gamma$
		and $2\alpha_i\le1$, gives a continuous affine $f_i$ with $v_i+\alpha_i\le f_i\le\gamma-\alpha_i$. Then $f_i\ge g_i$, $f_i-f_{i-1}\ge\alpha_i$ and $f_i<\gamma$, and $f_i\nearrow\gamma$ pointwise
		because $g_i\le f_i<\gamma$ and $\sup_ig_i=\gamma$.
		
		\emph{Induction.} Choose $\upsilon_i\in(0,\alpha_{i+1}/2)$ with $\upsilon_i\to0$. Inductively choose $z_i\in B\otimes C([k_i])$ by Lemma~\ref{lem:approxscale}, applied to the target
		$f_i-S_{i-1}$ and to $\upsilon_i$, where $S_i:=\sum_{h\le i}d_\bullet(z_h)$ and $S_0=0$:
		\[
		\big|d_\bullet(z_i)-(f_i-S_{i-1})\big|<\upsilon_i .
		\]
		This is legitimate: $S_{i-1}$ is continuous and affine, being a finite sum of the continuous affine functions provided by Lemma~\ref{lem:approxscale}; for $i=1$ the target is
		$f_1\ge m_0/2>0$, and for $i\ge2$, $f_i-S_{i-1}=(f_i-f_{i-1})+(f_{i-1}-S_{i-1})>\alpha_i-\upsilon_{i-1}>0$, using $|S_{i-1}-f_{i-1}|<\upsilon_{i-1}<\alpha_i/2$. The choice gives $|S_i-f_i|<\upsilon_i$
		for all $i$, so $S_i\to\gamma$ pointwise, also where $\gamma=\infty$. Moreover the partial sums never overshoot: since $\upsilon_i<\alpha_{i+1}/2$, $f_{i+1}-f_i\ge\alpha_{i+1}$
		and $f_{i+1}<\gamma$,
		\[
		S_i<f_i+\upsilon_i<f_i+\alpha_{i+1}\le f_{i+1}<\gamma\qquad\text{on }T .
		\]
		
		\emph{The element.} Place the elements $2^{-i}z_i$ in pairwise orthogonal diagonal blocks of $B\otimes c_0(\N)$, and let $a=\sum_i2^{-i}z_i$ (norm convergent). For $\tau\in T(A_\infty)$ and
		$\varepsilon>0$, $f_\varepsilon(a)=\sum_if_\varepsilon(2^{-i}z_i)$ has only finitely many nonzero terms, so by monotone convergence
		$d_\tau(a)=\lim_{\varepsilon\to0}\sum_i\tau(f_\varepsilon(2^{-i}z_i))=\sum_id_\tau(z_i)=\lim_iS_i(\tau)=\gamma(\tau)$.
	\end{proof}
	
	\subsection{Restrictions of groupoids}
	
	\begin{lemma}\label{lem:restriction}
		Let $G$ be a second countable, locally compact Hausdorff étale groupoid and $a'\in C_0(G^{(0)})_+$, and put $U=\{a'>0\}$. Then $G|_U=\{\eta\in G:r(\eta),s(\eta)\in U\}$ is an open étale
		subgroupoid with unit space $U$; it is second countable, locally compact Hausdorff, and principal if $G$ is; and $\overline{a'C^*_r(G)a'}=C^*_r(G|_U)$, with diagonal $C_0(U)$.
	\end{lemma}
	
	\begin{proof}
		$G|_U=r^{-1}(U)\cap s^{-1}(U)$ is open in $G$, hence a second countable, locally compact Hausdorff étale groupoid with unit space $U$, principal if $G$ is. For $x\in G^{(0)}$ let
		$\lambda_x$ be the regular representation of $C_c(G)$ on $\ell^2(G_x)$, $G_x=s^{-1}(x)$, and similarly $\lambda^U_y$ for $G|_U$, $y\in U$. Let $f\in C_c(G|_U)\subset C_c(G)$. If
		$G_x\cap r^{-1}(U)=\emptyset$, then $\lambda_x(f)=0$. Otherwise pick $\eta_0\in G_x$ with $y:=r(\eta_0)\in U$; then $\eta\mapsto\eta\eta_0^{-1}$ is a bijection
		$G_x\cap r^{-1}(U)\to(G|_U)_y$, and $\lambda_x(f)$ vanishes on $\ell^2(G_x\setminus r^{-1}(U))$, maps into $\ell^2(G_x\cap r^{-1}(U))$, and is unitarily equivalent there to
		$\lambda^U_y(f)$. Hence $\|f\|_{C^*_r(G)}=\sup_x\|\lambda_x(f)\|=\sup_{y\in U}\|\lambda^U_y(f)\|=\|f\|_{C^*_r(G|_U)}$, and the inclusion extends to an isometric $*$-homomorphism
		$C^*_r(G|_U)\to C^*_r(G)$ with image $\overline{C_c(G|_U)}$. For $\varphi,\psi\in C_c(U)$ and $f\in C_c(G)$, $\varphi*f*\psi=(\varphi\circ r)f(\psi\circ s)\in C_c(G|_U)$; conversely, if
		$f\in C_c(G|_U)$ has support $K$, then $f=\varphi*f*\varphi$ for any $\varphi\in C_c(U)$ equal to $1$ on the compact set $r(K)\cup s(K)\subset U$. By density,
		$\overline{C_c(G|_U)}=\overline{C_0(U)C^*_r(G)C_0(U)}$, the hereditary subalgebra generated by $C_0(U)$, which equals $\overline{a'C^*_r(G)a'}$ because $a'$ is strictly positive in
		$C_0(U)$. The diagonal $C_0(U)$ of $C^*_r(G|_U)$ is $\overline{a'C_0(G^{(0)})a'}$.
	\end{proof}
	
	\subsection{Proof of the stably projectionless case}\label{ss:proofA}
	
	\begin{theorem}[Stably projectionless case]\label{thm:SPA}
		Let $A$ be a separable, simple, nuclear, stably projectionless, $\JS$-stable $C^*$-algebra satisfying the UCT. There is a second countable, locally compact Hausdorff, principal étale
		groupoid $\mathcal G$ such that $A\cong C^*(\mathcal G)$. The unit space $\mathcal G^{(0)}$ is not compact.
	\end{theorem}
	
	\begin{proof}
		Let $A_c$, $T=T(A_c)$ and $(G_0,\rho,G_1)$ be as in Lemma~\ref{lem:Dc}. Fix a strictly positive element $h_A\in A$ and an isomorphism $\Phi\colon A\otimes\K\to A_c\otimes\K$ as in
		Lemma~\ref{lem:Dc}, and define
		\begin{equation}\label{eq:gamma}
			\gamma(\tau):=d_{\tau\otimes\Tr}\big(\Phi(h_A\otimes e_{11})\big)\in(0,\infty]\qquad(\tau\in T),
		\end{equation}
		a lower semicontinuous affine function, strictly positive because $A_c\otimes\K$ is simple. Under the identification $\Tt(A)\cong\Tt(A\otimes\K)\cong\Tt(A_c\otimes\K)\cong\Tt(A_c)=\R_{\ge0}T$
		given by restriction of traces (Lemma~\ref{lem:hertraces}), $\gamma$ is the restriction to the base $T$ of the scale $\Sigma_A$, which is homogeneous of degree one; this identification also
		carries $\rho_A$ to $\rho$.
		
		Let $A_\infty$ be the $C^*$-algebra of Sections~\ref{sec:alg}--\ref{sec:SPmaps} built from $(G_0,T,\rho,G_1)$ in Case (SP), which is possible by \eqref{eq:P}; let $a$ be given by
		Proposition~\ref{prop:exactscale} for $\gamma$; and put
		\[
		a':=a\otimes1_\JS\in C_0((G^s)^{(0)})_+\subset C\otimes\K,\qquad U:=\{a'>0\}\subset(G^s)^{(0)},\qquad A':=\overline{a'(C\otimes\K)a'}.
		\]
		Since every bounded trace on $A_\infty\otimes\JS$ has the form $\tau\otimes\tau_\JS$ (Proposition~\ref{prop:SPAZ}(a)), $d_{(\tau\otimes\tau_\JS)\otimes\Tr}(a')=d_{\tau\otimes\Tr}(a)=\gamma(\tau)$ for all
		$\tau\in T(A_\infty)\cong T$.
		
		By Lemma~\ref{lem:restriction} and \eqref{eq:Gs}, $A'=C^*_r(G^s|_U)$, and $G^s|_U$ is second countable, locally compact Hausdorff, principal and étale. It remains to show that $A\cong A'$.
		
		\emph{Both $C^*$-algebras are classifiable.} $A'$ is a nonzero hereditary subalgebra of the simple $C^*$-algebra $C\otimes\K$, hence full. It is separable, simple, nuclear and stably
		projectionless (as $C$ is, Proposition~\ref{prop:SPAZ}), $\JS$-stable \cite[Corollary~3.1]{TW}, of finite nuclear dimension \cite{CE}, and it satisfies the UCT because
		$A'\otimes\K\cong C\otimes\K$ by Theorem~\ref{thm:Brown}. It is stably finite: an isometry $v\in M_m(\widetilde{A'})$ has a unitary scalar part, so $1-vv^*$ is a projection in $M_m(A')$,
		hence $0$. The $C^*$-algebra $A$ has the same properties, by hypothesis and the same argument. In particular both are separable, simple, nuclear, stably projectionless and
		$\JS$-stable and satisfy the UCT, so Theorem~\ref{thm:GL} applies.
		
		\emph{The invariants agree.} By Theorem~\ref{thm:Brown}, $K_*(A')\cong K_*(C)=(G_0,G_1)\cong K_*(A)$, compatibly with the identifications of Lemma~\ref{lem:Dc} and
		Corollary~\ref{cor:SPK}. Restriction of lower semicontinuous traces gives cone isomorphisms $\Tt(A')\cong\Tt(C\otimes\K)\cong\Tt(C)=\R_{\ge0}T(C)\cong\R_{\ge0}T$ (Proposition~\ref{prop:SPAZ}),
		carrying $\rho_{A'}$ to $\rho$ and the scale of $A'$ to the homogeneous extension of $\tau\mapsto d_{\tau\otimes\Tr}(a')=\gamma(\tau)$. On the side of $A$ we have $\Tt(A)\cong\R_{\ge0}T$,
		carrying $\rho_A$ to $\rho$ and $\Sigma_A$ to the homogeneous extension of $\gamma$. These identifications are isomorphisms of invariants in the sense of \cite[2.3, 2.7]{GL4}, which in
		addition requires the cone maps to be homeomorphisms for the topology of pointwise convergence on Pedersen ideals and to be compatible with the pairings:
		\begin{enumerate}[label=\textup{(\roman*)}]
			\item For $X\in\{C,A_c\}$ we have $\Tt(X)=\R_{\ge0}T(X)$ with $T(X)$ weak$^*$-compact, and on such a cone the topology of pointwise convergence on $\Ped(X)$ is the quotient topology of
			$\R_{\ge0}\times T(X)$. Indeed, the map $(t,\tau)\mapsto t\tau$ is continuous; the topology is Hausdorff, since lower semicontinuous traces which agree on $\Ped(X)$ agree; and if
			$t_\lambda\tau_\lambda$ converges on $\Ped(X)$, then $(t_\lambda)$ is bounded (otherwise a cluster point of $(\tau_\lambda)$ in $T(X)$ would vanish on the dense ideal $\Ped(X)$), so
			compactness of $[0,N]\times T(X)$ applies.
			\item The restriction maps used above ($A'\subset C\otimes\K$, $\Phi(A\otimes e_{11})\subset A_c\otimes\K$, and $X\otimes e_{11}\subset X\otimes\K$ for $X\in\{C,A,A_c\}$) are bijections and
			homeomorphisms by Lemma~\ref{lem:hertraces}.
			\item The pairing of \cite[Definition~2.6]{GL4} is computed in hereditary subalgebras on which the traces in question are bounded, independently of the choice; hence it is natural for
			the identifications in (ii), and on $C$ it is the pairing $\rho$ (Theorem~\ref{thm:SPtraces}, Proposition~\ref{prop:SPAZ}).
		\end{enumerate}
		Since $\Sigma(K_0)=\{0\}$ for both $C^*$-algebras, the invariants $((K_0,\Sigma(K_0),\Tt,\Sigma,\rho),K_1)$ of $A$ and $A'$ are isomorphic, and Theorem~\ref{thm:GL} gives $A\cong A'=C^*_r(G^s|_U)$.
		The unit space $U$ is not compact, since $A$ is not unital.
	\end{proof}
	
	Together with Theorem~\ref{thm:unital}, this proves Theorem~\ref{thm:main} for unital and for stably projectionless $C^*$-algebras; the remaining stably finite
	case is treated in Section~\ref{sec:proj}.
	%======================================================================
	\part*{Part IV: The remaining stably finite case}
	
	Parts~II and~III treat the classifiable stably finite $C^*$-algebras which are unital or stably projectionless. In this part we complete the stably finite
	case by treating non-unital $C^*$-algebras with projections in the stabilization.
	
	%======================================================================
	\section{Non-unital \texorpdfstring{$C^*$}{C*}-algebras with projections in the stabilization}\label{sec:proj}
	
	\subsection{The three stably finite cases}\label{ss:threecases}
	Recall that a $C^*$-algebra $A$ is \emph{stably finite} if the unital $C^*$-algebras $M_n(\widetilde A)$ are finite (every isometry is unitary) for all $n$. A
	classifiable stably finite $C^*$-algebra $A$ belongs to exactly one of the following classes:
	\begin{enumerate}[label=(\Roman*)]
		\item $A$ is unital (Part~II);
		\item $A$ is stably projectionless, i.e.\ $A\otimes\K$ has no nonzero projection (Part~III);
		\item $A$ is not unital, but $A\otimes\K$ contains a nonzero projection.
	\end{enumerate}
	In case~(III) the $C^*$-algebra $A$ itself need not contain a nonzero projection. For instance, let $a$ be a positive element of $\JS\otimes\K$ with
	$d_{\tau_\JS\otimes\Tr}(a)=\frac12$ (such elements exist: apply Proposition~\ref{prop:exactscaleU} below, with $\gamma\equiv\frac12$, to the $C^*$-algebra $A_\infty$ of Part~II built
	from the invariant of $\JS$, and map the resulting element to $A_\infty\otimes\JS\otimes\K$ by $x\otimes y\mapsto x\otimes1_\JS\otimes y$; here $A_\infty\otimes\JS\cong\JS$ by the proof
	of Theorem~\ref{thm:unital}), and let
	$A=\overline{a(\JS\otimes\K)a}$. A projection $q\in A$ satisfies $q\precsim a$, hence $\tau_\JS\otimes\Tr(q)=d_{\tau_\JS\otimes\Tr}(q)\le\frac12$; since
	$K_0(\JS)=\Z$ with $(\tau_\JS)_*$ the inclusion $\Z\subset\R$, this forces $q=0$. On the other hand $A\otimes\K\cong\JS\otimes\K$ by Brown's theorem.
	Another example in case~(III) is every stable $C^*$-algebra $E\otimes\K$ with $E$ unital, classifiable and stably finite.
	
	\begin{theorem}[Non-unital $C^*$-algebras with projections]\label{thm:proj}
		Let $A$ be a classifiable, stably finite, non-unital $C^*$-algebra such that $A\otimes\K$ contains a nonzero projection. There is a second countable, locally
		compact Hausdorff, principal étale groupoid $\mathcal G'$ with non-compact unit space such that $A\cong C^*(\mathcal G')$.
	\end{theorem}
	
	\subsubsection*{The idea}
	Let $p\in A\otimes\K$ be a nonzero projection. The corner $E_0=p(A\otimes\K)p$ is unital and classifiable, so it has an untwisted model $C^*(\mathcal G)$ by
	Part~II, and $A\otimes\K\cong E_0\otimes\K\cong C^*(\mathcal G\times\mathcal P_\N)$ by Brown's theorem. Inside this $C^*$-algebra, $A$ is the hereditary subalgebra
	generated by a single positive element $b$ (the image of $h_A\otimes e_{11}$ for a strictly positive $h_A\in A$). We shall find a positive element $a'$ of the canonical
	diagonal which generates an isomorphic hereditary subalgebra; that subalgebra is the $C^*$-algebra of the restriction of the groupoid to the open set
	$\{a'>0\}$ (Lemma~\ref{lem:restriction}).
	
\begin{figure}[ht]
		\centering
		\begin{tikzpicture}[>={Stealth[length=2mm]},node distance=0.9cm and 1.7cm,
			box/.style={draw,rounded corners,fill=blue!5,align=center,font=\footnotesize,minimum height=0.95cm,text width=3.25cm}]
			\node[box,fill=gray!10] (D) {$A$ non-unital,\\ $0\neq p\in A\otimes\K$};
			\node[box,right=of D] (E) {$E_0=p(A\otimes\K)p$,\\ unital, classifiable};
			\node[box,right=of E] (G) {$E_0\cong A_\infty\otimes\JS\cong C^*(\mathcal G)$\\ (Theorem~\ref{thm:unital})};
			\node[box,below=of G] (b) {$b\in A_\infty\otimes\JS\otimes\K$,\\ $\her(b)\cong A$};
			\node[box,left=of b] (a) {$a'\in C_0(\mathcal H^{(0)})$,\\ $d_\bullet(a')=d_\bullet(b)$};
			\node[box,left=of a,fill=gray!10] (D2) {$A\cong C^*(\mathcal H|_{\{a'>0\}})$};
			\draw[->] (D) -- node[above,font=\scriptsize]{corner} (E);
			\draw[->] (E) -- node[above,font=\scriptsize]{Part~II} (G);
			\draw[->] (G) -- node[right,font=\scriptsize]{Brown} (b);
			\draw[->] (b) -- node[above,font=\scriptsize]{Prop.~\ref{prop:exactscaleU}} (a);
			\draw[->] (a) -- node[above,font=\scriptsize]{Prop.~\ref{prop:soft}} (D2);
		\end{tikzpicture}
		\caption{The proof of Theorem~\ref{thm:proj}; $\mathcal H=\mathcal G\times\mathcal P_\N$.}\label{fig:projplan}
	\end{figure}
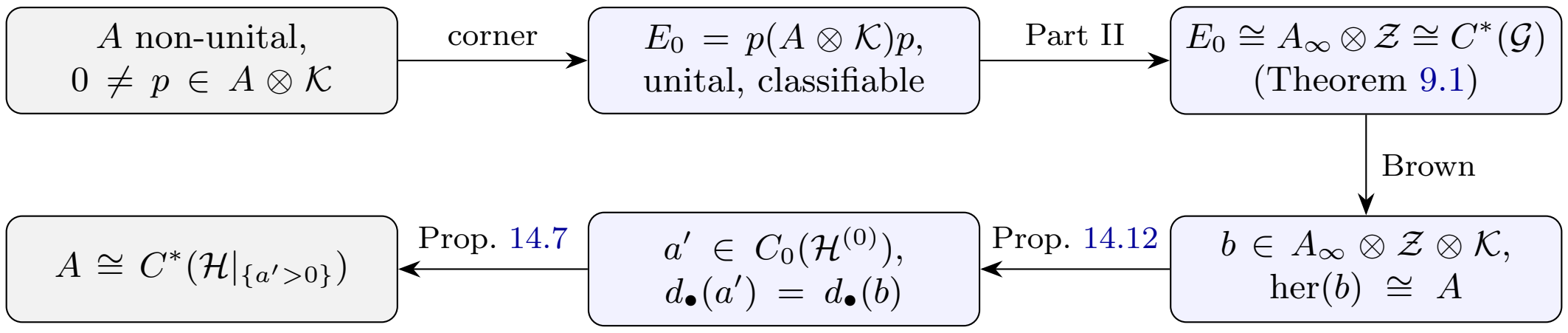

	Two hereditary subalgebras are recognized to be isomorphic as follows. In a $C^*$-algebra of stable rank one, Cuntz equivalent positive elements generate isomorphic
	hereditary subalgebras (Section~\ref{ss:Cuntziso}). Since $A$ is not unital, the spectrum of $b$ accumulates at $0$; such ``soft'' elements are determined
	up to Cuntz equivalence by their rank functions $\tau\mapsto d_\tau(b)$ (Section~\ref{ss:soft}). Finally, every lower semicontinuous affine rank function is
	realized in the diagonal of the unital model, by the method of Section~\ref{sec:SPscale} (Section~\ref{ss:realizeU}). In contrast with Part~III, no
	classification theorem for non-unital $C^*$-algebras is needed: the unital classification enters only through Part~II, and hereditary subalgebras are compared
	using R{\o}rdam's strict comparison and stable rank one theorems and the theorem of Coward, Elliott and Ivanescu. Figure~\ref{fig:projplan} shows the plan.
	
	\subsection{Unital hereditary subalgebras and the spectrum at zero}
	
	For a $C^*$-algebra $A$ and $b\in A_+$ we write $\her(b)=\overline{bAb}$ for the hereditary subalgebra generated by $b$.
	
	\begin{lemma}\label{lem:unitalher}
		Let $A$ be a $C^*$-algebra and $b\in A_+$. Then $\her(b)$ is unital if and only if $0$ is not an accumulation point of $\operatorname{sp}(b)$.
	\end{lemma}
	
	\begin{proof}
		Suppose $0$ is not an accumulation point of $\operatorname{sp}(b)$, say $\operatorname{sp}(b)\subset\{0\}\cup[\delta,\infty)$ with $\delta>0$. The functions
		$\chi=\chi_{(0,\infty)}$ and $h(t)=\chi(t)t^{-2}$ are continuous on $\operatorname{sp}(b)$ and vanish at $0$, so $q=\chi(b)$ and $h(b)$ lie in $C^*(b)\subset A$; $q$ is
		a projection with $qb=b$, and $q=b\,h(b)\,b$. Hence $\her(b)\subset qAq$ (as $b=qbq$) and $qAq=bh(b)bAbh(b)b\subset\her(b)$. Thus $\her(b)=qAq$, which is unital.
		
		Conversely, let $e$ be a unit of $\her(b)$. The element $b$ is strictly positive in $\her(b)$: if $\varphi$ is a state of $\her(b)$ with $\varphi(b)=0$, then for
		$x\in A$ we write $bxb=b^{1/2}\cdot(b^{1/2}xb)$ with both factors in $\her(b)$, and the Cauchy--Schwarz inequality gives $|\varphi(bxb)|^2\le\varphi(b)\,\|b^{1/2}xb\|^2=0$;
		so $\varphi=0$ on the dense subset $bAb$, which is absurd. In the unital $C^*$-algebra $\her(b)$
		a strictly positive element is invertible (if $0\in\operatorname{sp}_{\her(b)}(b)$, a state attaining the minimum of the spectrum vanishes at $b$). So
		$b\ge\delta e$ for some $\delta>0$, and $\operatorname{sp}(b)\subset\{0\}\cup[\delta,\|b\|]$, because $b$ and the projection $e$ commute and $b(1-e)=0$.
	\end{proof}
	
	\subsection{Cuntz equivalence and isomorphism of hereditary subalgebras}\label{ss:Cuntziso}
	
	We use the notions of Section~\ref{ss:SPtraces}: Cuntz subequivalence $\precsim$, Cuntz equivalence $\sim$, and dimension functions $d_\tau$. We also
	use three standard facts:
	\begin{enumerate}[label=(Cu\arabic*)]
		\item If $x,y\in A_+$ and $\|x-y\|<\varepsilon$, then $(x-\varepsilon)_+\precsim y$ \cite[Lemma~2.5(ii)]{KR00}.
		\item If $(x-\varepsilon)_+\precsim y$ for every $\varepsilon>0$, then $x\precsim y$ (immediate from the definition, since $(x-\varepsilon)_+\to x$); and $\precsim$ is
		transitive.
		\item If $z\in\her(x)_+$, then $z\precsim x$. (Put $h_k(t)=\min\{k,t^{-1/2}\}$; then $h_k(x)xh_k(x)=\min\{k^2x,1\}$ is an approximate unit of $\her(x)$, so
		$z^{1/2}h_k(x)\,x\,h_k(x)z^{1/2}\to z$.)
	\end{enumerate}
	
	\begin{lemma}\label{lem:polar}
		For every $x$ in a $C^*$-algebra $A$, the hereditary subalgebras $\her(x^*x)$ and $\her(xx^*)$ are isomorphic.
	\end{lemma}
	
	\begin{proof}
		Put $|x|=(x^*x)^{1/2}$. Then $\her(x^*x)=\overline{|x|A|x|}$, and $\her(xx^*)=\overline{xAx^*}$: indeed $xx^*Axx^*\subset xAx^*$, and $xyx^*=\lim_k
		(xx^*)^{1/k}xyx^*(xx^*)^{1/k}\in\her(xx^*)$ because $(xx^*)^{1/k}x\to x$. For a positive $c$ in the unitization,
		\[
		\|xcx^*\|=\|(c^{1/2}x^*)^*(c^{1/2}x^*)\|=\|c^{1/2}x^*xc^{1/2}\|=\||x|c|x|\|,
		\]
		and hence, for arbitrary $w$, $\|xwx^*\|^2=\|x(w|x|^2w^*)x^*\|=\||x|w|x|^2w^*|x|\|=\||x|w|x|\|^2$. Therefore $|x|w|x|\mapsto xwx^*$ is a well defined, linear,
		isometric map from $|x|A|x|$ onto $xAx^*$. It is $*$-preserving, and multiplicative:
		$(|x|w|x|)(|x|w'|x|)=|x|(w|x|^2w')|x|\mapsto x(w|x|^2w')x^*=(xwx^*)(xw'x^*)$. Its continuous extension is an isomorphism $\her(x^*x)\to\her(xx^*)$.
	\end{proof}
	
	\begin{theorem}[{\cite[Theorem~3]{CEI}}, {\cite[\S6.2]{ORT}}]\label{thm:ORT}
		Let $A$ be a separable $C^*$-algebra of stable rank one, and let $a,b\in A_+$ with $a\sim b$. Then there is $x\in A$ with $\her(a)=\her(x^*x)$ and
		$\her(b)=\her(xx^*)$.
	\end{theorem}
	
	This is \cite[\S6.2, (1)]{ORT} (``$a\approx b\Leftrightarrow a\sim_sb$'', where $\approx$ denotes Cuntz equivalence and $\sim_s$ is Blackadar's relation of
	\cite[Definition~2.1]{ORT}, i.e.\ the conclusion of the theorem); it rests on the
	theorem of Coward, Elliott and Ivanescu that over a $C^*$-algebra of stable rank one, Cuntz equivalent countably generated Hilbert modules are isomorphic
	\cite[Theorem~3]{CEI}.
	
	\begin{lemma}\label{lem:sr1}
		Let $E$ be a unital, simple, separable, stably finite, $\JS$-stable $C^*$-algebra. Then $E\otimes\K$ is separable and has stable rank one.
	\end{lemma}
	
	\begin{proof}
		For each $m$, $M_m(E)$ is unital, simple, finite (as $E$ is stably finite) and $\JS$-absorbing ($M_m(E)\otimes\JS\cong M_m(E\otimes\JS)\cong M_m(E)$), so it has
		stable rank one by R{\o}rdam's theorem \cite[Theorem~6.7]{Ro}. The $C^*$-algebra $E\otimes\K$ is the inductive limit of the $M_m(E)$ under the corner embeddings
		$x\mapsto x\oplus0$, and stable rank one passes to inductive limits \cite[Theorem~5.1]{Rie}.
	\end{proof}
	
	\begin{corollary}\label{cor:Cuntziso}
		Let $E$ be as in Lemma~\textup{\ref{lem:sr1}}. If $a,b\in(E\otimes\K)_+$ are Cuntz equivalent, then $\her(a)\cong\her(b)$.
	\end{corollary}
	
	\begin{proof}
		By Theorem~\ref{thm:ORT} and Lemma~\ref{lem:sr1}, $\her(a)=\her(x^*x)$ and $\her(b)=\her(xx^*)$ for some $x$, and these are isomorphic by Lemma~\ref{lem:polar}.
	\end{proof}
	
	\subsection{Soft elements are determined by their rank functions}\label{ss:soft}
	
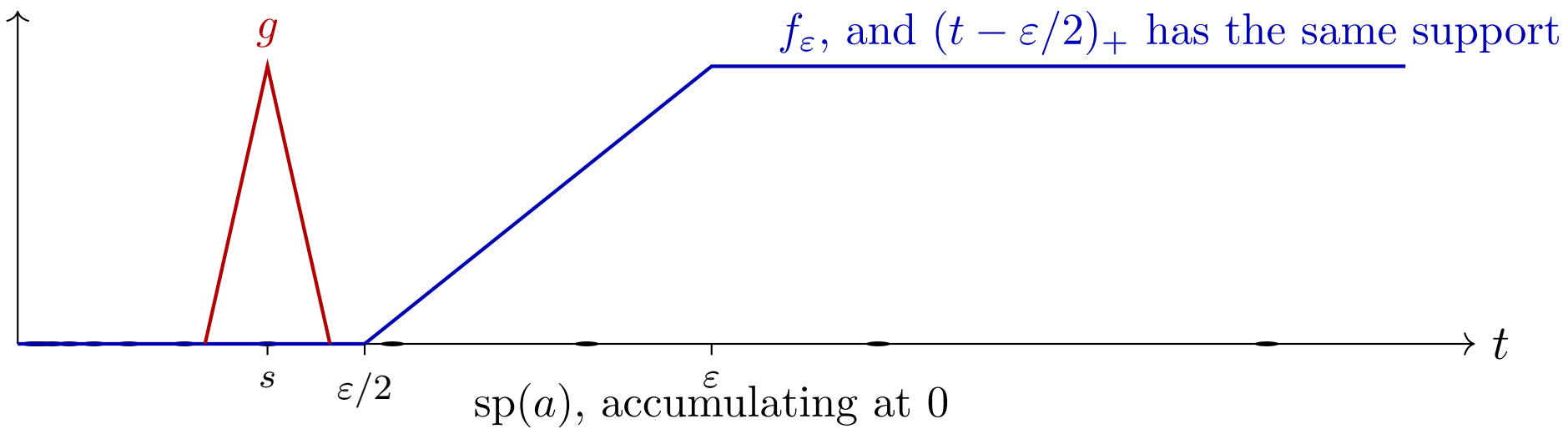
\begin{figure}[ht]
		\centering
		\begin{tikzpicture}[xscale=11,yscale=2.2]
			\draw[->] (0,0) -- (1.05,0) node[right]{$t$};
			\draw[->] (0,0) -- (0,1.2);
			\foreach \x in {0.9,0.62,0.41,0.27,0.18,0.12,0.08,0.055,0.037,0.025,0.017,0.011}{\fill (\x,0) circle (0.25pt);}
			\node[below,font=\small] at (0.5,-0.08) {$\operatorname{sp}(a)$, accumulating at $0$};
			\draw[thick,blue!70!black] (0,0) -- (0.25,0) -- (0.5,1) -- (1,1);
			\node[blue!70!black,font=\small] at (0.83,1.12) {$f_\varepsilon$, and $(t-\varepsilon/2)_+$ has the same support};
			\draw[thick,red!70!black] (0.135,0) -- (0.18,1) -- (0.225,0);
			\node[red!70!black,font=\small] at (0.18,1.12) {$g$};
			\draw[dashed] (0.25,0) -- (0.25,-0.04) node[below,font=\scriptsize]{$\varepsilon/2$};
			\draw[dashed] (0.5,0) -- (0.5,-0.04) node[below,font=\scriptsize]{$\varepsilon$};
			\draw[dashed] (0.18,0) -- (0.18,-0.04) node[below,font=\scriptsize]{$s$};
		\end{tikzpicture}
		\caption{The gap in the proof of Proposition~\ref{prop:soft}: $g$ is supported in $(0,\varepsilon/2)$ and $g(s)=1$ at a point $s\in\operatorname{sp}(a)$, so
			$c=g(a)$ is a nonzero element of $\her(a)$ orthogonal to $(a-\varepsilon/2)_+$.}\label{fig:gap}
	\end{figure}

	In a stabilized unital $C^*$-algebra, a projection $q$ and a non-projection $a$ may have the same rank function, $d_\tau(a)=\tau(q)$ for all $\tau$, without being
	Cuntz equivalent: by Lemma~\ref{lem:unitalher} and Corollary~\ref{cor:Cuntziso}, $\her(q)$ is unital, and $\her(a)$ is not if $\operatorname{sp}(a)$
	accumulates at $0$. The next proposition shows that this is the only obstruction. The point of the proof is a \emph{gap}: if $\operatorname{sp}(a)$ accumulates
	at $0$, then cutting $a$ down to $(a-\varepsilon)_+$ loses a nonzero piece $c$ of $\her(a)$, and so strictly decreases the rank at every trace
	(Figure~\ref{fig:gap}); strict comparison can then be applied.
	
	For a unital $C^*$-algebra $E$ and $\tau\in T(E)$ we write $d_\tau$ for $d_{\tau\otimes\Tr}$ on $(E\otimes\K)_+$, and $P_m=1_E\otimes(e_{11}+\dots+e_{mm})\in E\otimes\K$.
	Note that $P_m\in E\otimes\K$ because $E$ is unital, so $P_mxP_m\precsim x$ for $x\in(E\otimes\K)_+$ directly from the definition of $\precsim$; and for $z\in
	M_m(E)=P_m(E\otimes\K)P_m$ the function $\tau\mapsto(\tau\otimes\Tr)(z)=\sum_{i\le m}\tau(z_{ii})$ is continuous on $T(E)$.
	
	\begin{proposition}\label{prop:soft}
		Let $E$ be a unital, simple, separable, exact, $\JS$-stable $C^*$-algebra with $T(E)\neq\emptyset$. Let $a,b\in(E\otimes\K)_+$ be such that $0$ is an
		accumulation point of $\operatorname{sp}(a)$ and of $\operatorname{sp}(b)$, and
		\[
		d_\tau(a)=d_\tau(b)\qquad\text{for every }\tau\in T(E).
		\]
		Then $a\sim b$ and $\her(a)\cong\her(b)$.
	\end{proposition}
	
	\begin{proof}
		Put $\gamma(\tau)=d_\tau(a)=d_\tau(b)\in[0,\infty]$. By symmetry it suffices to show $a\precsim b$, and by (Cu2) it suffices to show $(a-\varepsilon)_+\precsim b$ for a
		fixed $\varepsilon>0$. The last assertion then follows from Corollary~\ref{cor:Cuntziso}: $E$ is stably finite, because a faithful tracial state rules out non-unitary
		isometries in $M_n(E)$ (as in the proof of Proposition~\ref{prop:firstlimit}).
		
		\emph{Step 1: strict comparison in matrix algebras.} If $x,y\in M_m(E)_+$ satisfy $d_\tau(x)<d_\tau(y)$ for all $\tau\in T(E)$, then $x\precsim y$. Indeed,
		$M_m(E)$ is unital, simple, exact and $\JS$-absorbing, its tracial states are the $\tau\otimes\tr_m=\frac1m(\tau\otimes\Tr)$, $\tau\in T(E)$, and R{\o}rdam's
		strict comparison theorem \cite[Corollary~4.6]{Ro} applies to $M_m(E)$.
		
		\emph{Step 2: the gap.} Since $0$ is an accumulation point of $\operatorname{sp}(a)$, there is $s\in\operatorname{sp}(a)$ with $0<s<\varepsilon/2$. Choose a
		continuous $g\colon[0,\infty)\to[0,1]$ with $g(s)=1$ and $g=0$ outside $(s/2,\varepsilon/2)$ (Figure~\ref{fig:gap}), and put $c=g(a)$. Then $c\neq0$, and $c$ is
		orthogonal to $(a-\varepsilon/2)_+$. The sum $(a-\varepsilon/2)_++c$ is a function of $a$ vanishing at $0$, hence lies in $\her(a)$ and is $\precsim a$ by (Cu3). By
		additivity of $d_\tau$ on orthogonal elements,
		\[
		d_\tau\big((a-\varepsilon/2)_+\big)+d_\tau(c)\le d_\tau(a)=\gamma(\tau)\qquad(\tau\in T(E)).
		\]
		Every $\tau\in T(E)$ is faithful ($E$ is simple), hence so is $\tau\otimes\Tr$ on $E\otimes\K$ (if $(\tau\otimes\Tr)(y)=0$ for $y\ge0$, then $\tau(y_{ii})=0$, so
		$y^{1/2}(1\otimes e_{ii})=0$ for all $i$); since $c\ne0$, $d_\tau(c)\ge(\tau\otimes\Tr)(c)>0$. Hence $d_\tau((a-\varepsilon/2)_+)<\gamma(\tau)$ whenever
		$\gamma(\tau)<\infty$.
		
		\emph{Step 3: passage to a matrix algebra.} Choose $m$ with $\|P_m(a-\varepsilon/2)_+P_m-(a-\varepsilon/2)_+\|<\varepsilon/2$, and put
		$x=P_m(a-\varepsilon/2)_+P_m\in M_m(E)_+$. By (Cu1), $(a-\varepsilon)_+=((a-\varepsilon/2)_+-\varepsilon/2)_+\precsim x$, and $x\precsim(a-\varepsilon/2)_+$. Hence
		$d_\tau(x)\le d_\tau((a-\varepsilon/2)_+)$, and $d_\tau(x)\le d_\tau(P_m)=m$. By Step~2,
		\begin{equation}\label{eq:xstrict}
			d_\tau(x)<\gamma(\tau)\qquad\text{for every }\tau\in T(E)
		\end{equation}
		(if $\gamma(\tau)=\infty$ this holds because $d_\tau(x)\le m$).
		
		\emph{Step 4: a uniform witness below $b$, by compactness.} Fix $\eta>0$ and put $u(\tau)=(\tau\otimes\Tr)(f_\eta(x))$, with $f_\eta$ as in
		Section~\ref{ss:SPtraces}. This is a continuous function on $T(E)$, since $f_\eta(x)\in M_m(E)$. Since $f_\delta((t-\eta)_+)\le f_\eta(t)$ for all $t\ge0$ and
		$\delta>0$, and $0\le f_\eta(x)\le1$ lies in $C^*(x)$,
		\[
		d_\tau\big((x-\eta)_+\big)\le u(\tau)\le d_\tau(x)<\gamma(\tau)\qquad(\tau\in T(E)).
		\]
		Let $y_k=P_kbP_k\in M_k(E)_+$. Then $y_k\precsim y_{k'}$ for $k\le k'$ (as $y_k=P_ky_{k'}P_k$), and $y_k\precsim b$. Fix $\tau_0\in T(E)$. For $\delta>0$ we have
		$\|y_k-b\|<\delta$ for large $k$, hence $(b-\delta)_+\precsim y_k$ by (Cu1); and $d_{\tau_0}((b-\delta)_+)\ge\tau_0(f_{2\delta}(b))$, because $f_{2\delta}(b)\in
		\her((b-\delta)_+)$ has norm at most one. Letting $\delta\to0$ gives $\liminf_kd_{\tau_0}(y_k)\ge d_{\tau_0}(b)=\gamma(\tau_0)>u(\tau_0)$. So there are $k$ with
		$d_{\tau_0}(y_k)>u(\tau_0)$ and then $\delta>0$ with $(\tau_0\otimes\Tr)(f_\delta(y_k))>u(\tau_0)$. Both sides of the last inequality are continuous in $\tau$, so it
		holds on a neighbourhood of $\tau_0$. By compactness of $T(E)$ there are finitely many pairs $(k_i,\delta_i)$ such that every $\tau$ satisfies
		$(\tau\otimes\Tr)(f_{\delta_i}(y_{k_i}))>u(\tau)$ for some $i$. Let $K\ge m$ be larger than all $k_i$ and $y=y_K$. For every $\tau$, choosing such an $i$,
		\[
		d_\tau(y)\ge d_\tau(y_{k_i})\ge(\tau\otimes\Tr)\big(f_{\delta_i}(y_{k_i})\big)>u(\tau)\ge d_\tau\big((x-\eta)_+\big).
		\]
		Both $(x-\eta)_+$ and $y$ lie in $M_K(E)_+$, so Step~1 gives $(x-\eta)_+\precsim y\precsim b$.
		
		\emph{Step 5: conclusion.} As $\eta>0$ was arbitrary, $x\precsim b$ by (Cu2), and $(a-\varepsilon)_+\precsim x\precsim b$ by Step~3.
	\end{proof}
	
	\begin{remark}
		The proposition is a special case of the description of the Cuntz semigroup of a simple $\JS$-stable $C^*$-algebra as the disjoint union of the Murray--von Neumann
		semigroup and the lower semicontinuous affine functions on the trace cone \cite{BPT,ERS}; the direct argument above uses only R{\o}rdam's strict comparison
		theorem.
	\end{remark}
	
	\subsection{Realizing a rank function in the diagonal of the unital model}\label{ss:realizeU}
	
	We return to the unital construction of Part~II: $A_\infty=\lims(A_n,\xi_n)$, $B=\lims B_n=C(\bar{\mathcal G}^{(0)})$ its diagonal (Proposition~\ref{prop:firstlimit}),
	$\Omega_n=\mathcal G_n^{(0)}$ the compact unit space of stage $n$ (Lemma~\ref{lem:groupoid}), $A_\infty^{(0)}=\lims(A_n,\varphi_n)$, the AF algebra $\widehat
	F=\lims(\widehat F_n,\widehat\psi_n)$, and the maps $\bar q_n\colon A_n\to\widehat F_n$. We identify $T(A_\infty)$ with $T$ by \eqref{eq:TA}. The method of
	Section~\ref{sec:SPscale} requires one quantitative input, the analogue of \eqref{eq:quant}: traces at a finite stage are uniformly close to traces factoring
	through $\bar q_n$. It follows from the perturbation argument of Section~\ref{ss:perturb}.
	
	\begin{lemma}\label{lem:quantU}
		Let $\tau\in T(A_\infty)$, let $\sigma\in T(\widehat F)$ be the trace which corresponds to $\tau$ under \eqref{eq:TA}, and put $\sigma_n=\sigma\circ\widehat\psi_{n,\infty}\in
		T(\widehat F_n)$. With $\omega_j=2\eta_j+2\varepsilon_j$ as in Section~\ref{ss:perturb},
		\begin{equation}\label{eq:quantU}
			\|\tau\circ\xi_{n,\infty}-\sigma_n\circ\bar q_n\|\le\vartheta'_n:=\sum_{j\ge n}\omega_j\qquad(n\ge1).
		\end{equation}
		Moreover, if $h\in\Z^{J_n}=K_0(\widehat F_n)$ has image $x\in F$, then $(\sigma_n)_*(h)=\theta(x)(\tau)$, where $\tau$ also denotes the corresponding point
		of $T$.
	\end{lemma}
	
	\begin{proof}
		The sequence $t=(\tau\circ\xi_{m,\infty})_m$ is compatible for the maps $\xi_m^*$. In Section~\ref{ss:perturb}, Lemma~\ref{lem:perturb} was applied with $a_m=\varphi_m^*$
		and $b_m=\xi_m^*$; the sequence $St$ is $\varphi^*$-compatible, so by Lemma~\ref{lem:Tlimit}(a) it is $(\tau'\circ\varphi_{m,\infty})_m$ for a unique $\tau'\in
		T(A_\infty^{(0)})$, and by Proposition~\ref{prop:tracesAF} $\tau'=\sigma\circ\bar q$ for a unique $\sigma\in T(\widehat F)$. This is how \eqref{eq:TA} identifies $\tau$ with
		$\sigma$. By \eqref{eq:qbar}, $\bar q\circ\varphi_{n,\infty}=\widehat\psi_{n,\infty}\circ\bar q_n$, so
		\[
		(St)_n=\sigma\circ\bar q\circ\varphi_{n,\infty}=\sigma\circ\widehat\psi_{n,\infty}\circ\bar q_n=\sigma_n\circ\bar q_n ,
		\]
		and the estimate $\|(St)_n-t_n\|\le\sum_{j\ge n}\omega_j$ from the proof of Lemma~\ref{lem:perturb} is \eqref{eq:quantU}. (We write $\vartheta'_n$ to distinguish
		these numbers from the $\vartheta_n$ of Part~III.) For the last assertion,
		$(\sigma_n)_*(h)=\sigma_*(\widehat\psi_{n,\infty*}h)=\sigma_*(x)$, because $K_0(\widehat F)=\lims(\Z^{J_n},M_n)$ is identified with $F$ via the maps induced by
		the $\pi_n$ (Lemmas~\ref{lem:AFtraces} and~\ref{lem:delayed}); and $\sigma_*=s_\tau$ under $T(\widehat F)\cong S(F)\cong T$ (Lemma~\ref{lem:AFtraces},
		Theorem~\ref{thm:Fdimgroup}), so $\sigma_*(x)=\theta(x)(\tau)$.
	\end{proof}
	
	Every tracial state $\tau_n$ of $A_n$ restricts to a probability measure $\mu_{\tau_n}$ on $\Omega_n$, and \eqref{eq:dmu} holds for diagonal elements of
	$B_n\otimes C([k])$.
	
	\begin{lemma}\label{lem:atomicU}
		There is a countable set $\mathcal S_n\subset\Omega_n$ such that for every $\tau\in T(A_\infty)$ the measure $\mu_{\tau_n}$, $\tau_n=\tau\circ\xi_{n,\infty}$, is concentrated on
		$\mathcal S_n$.
	\end{lemma}
	
	\begin{proof}
		This is the proof of Lemma~\ref{lem:atomic}, with \eqref{eq:quantU} in place of \eqref{eq:quant}. By \eqref{eq:quantU}, for $M>n$, $\tau_n=\tau_M\circ\xi_{n,M}$ is
		within $\vartheta'_M$ of $\sigma_M\circ\bar q_M\circ\xi_{n,M}$, and $\vartheta'_M\to0$. On $B_M=C(\Omega_M)$ the functional $\sigma_M\circ\bar q_M$ is a positive
		combination of point masses at the finitely many points of $W^\theta_M$, because $\bar q_M|_{B_M}$ is restriction to $W^\theta_M$ (Lemma~\ref{lem:groupoid}).
		Since $\xi_{n,M}$ is unital and $\xi_{n,M}(B_n)\subset B_M$ (Lemma~\ref{lem:xin}), each $\ev_w\circ\xi_{n,M}|_{B_n}$ is a character of $B_n$, i.e.\ a point mass at a
		point of $\Omega_n$. So $\sigma_M\circ\bar q_M\circ\xi_{n,M}$ restricts to a measure on a finite set $\mathcal S_n^{(M)}\subset\Omega_n$ determined by the
		construction; put $\mathcal S_n=\bigcup_{M>n}\mathcal S_n^{(M)}$. Norm limits of measures concentrated on $\mathcal S_n$ are concentrated on $\mathcal S_n$.
	\end{proof}
	
	\begin{lemma}\label{lem:approxscaleU}
		Let $g\in\Aff(T)$ with $g\gg0$, and let $\varepsilon>0$. There are $k\ge1$ and $z\in(B\otimes C([k]))_+\subset A_\infty\otimes M_k$ with $\|z\|=1$ such that
		$\tau\mapsto d_\tau(z)$ is a continuous affine function on $T(A_\infty)\cong T$ and $g-\varepsilon<d_\bullet(z)<g+\varepsilon$.
	\end{lemma}
	
	\begin{proof}
		We follow the proof of Lemma~\ref{lem:approxscale}, indicating the inputs from Part~II.
		
		\emph{The target.} Choose $f\in G'$ with $\|f-g\|<\min\{\varepsilon/2,\min_Tg\}$; then $f\gg0$, so $(0,f)\in H^+\setminus\{0\}$, and its image $\bar f\in F^+\setminus\{0\}$ is
		represented at some stage by a nonzero $h\in\Z_{\ge0}^{J_n}$ (Lemma~\ref{lem:delayed}). Since $\bar u$ is an order unit of $F$ (Theorem~\ref{thm:Fdimgroup}) and is
		represented by $r_n$ at stage $n$ (Section~\ref{ss:unit}), there is $N$ with $\bar f\le N\bar u$, and after passing to a later stage $h\le Nr_n$ coordinatewise.
		Passing to a still later stage, we may assume $N\vartheta'_n<\varepsilon/2$. Put $k=N$ and let $c$ be a projection in the diagonal of $\widehat F_n\otimes M_k$ with
		$\rank(c_l)=h_l\le kr_{n,l}$ in $M_{r_{n,l}}\otimes M_k$. By Lemma~\ref{lem:quantU}, $d_{\sigma_n\otimes\Tr}(c)=(\sigma_n)_*(h)=\theta(\bar f)(\tau)=f(\tau)$.
		
		\emph{The lift.} By Lemma~\ref{lem:groupoid}, $\bar q_n\otimes\id$ on $B_n\otimes C([k])=C(\Omega_n\times[k])$ is restriction to $W^\theta_n\times[k]$, and $c$ is the
		indicator function of a subset $P_1$ of this finite set; let $P_0$ be its complement. As in Lemma~\ref{lem:approxscale}, choose a finite union $U$ of open balls
		around the points of $P_1$, with closure disjoint from $P_0$, such that $\partial U\cap(\mathcal S_n\times[k])=\emptyset$ (Lemma~\ref{lem:atomicU}), and
		$d\in C(\Omega_n\times[k])_+$ with $\{d>0\}=U$, $\|d\|=1$ and $d=1$ on $P_1$ (Figure~\ref{fig:lift}). Then $(\bar q_n\otimes\id)(d)=c$. Put
		$z=(\xi_{n,\infty}\otimes\id)(d)$; it has norm one because the $\xi_n$ are injective (Lemma~\ref{lem:xin}).
		
		\emph{Continuity and estimate.} Exactly as in Lemma~\ref{lem:approxscale}: $d_\tau(z)=(\mu_{\tau_n}\otimes\#)(U)=(\mu_{\tau_n}\otimes\#)(\overline U)$ depends continuously
		on $\tau$, and $|d_\tau(z)-f(\tau)|\le k\|\tau_n-\sigma_n\circ\bar q_n\|\le k\vartheta'_n<\varepsilon/2$ by \eqref{eq:quantU}. Hence $|d_\tau(z)-g(\tau)|<\varepsilon$.
	\end{proof}
	
	\begin{proposition}\label{prop:exactscaleU}
		For every lower semicontinuous affine function $\gamma\colon T\to(0,\infty]$ there is $a\in(B\otimes c_0(\N))_+\subset A_\infty\otimes\K$ such that $d_\tau(a)=\gamma(\tau)$ for
		all $\tau\in T(A_\infty)\cong T$, and $0$ is an accumulation point of $\operatorname{sp}(a)$.
	\end{proposition}
	
	\begin{proof}
		If $\gamma\equiv\infty$, put $a=\sum_j2^{-j}\,1_{A_\infty}\otimes e_{jj}$; here $1_{A_\infty}\in B$, $d_\tau(a)=\sum_j\tau(1_{A_\infty})=\infty$, and $\operatorname{sp}(a)=\{0\}\cup\{2^{-j}\}$.
		Otherwise the proof of Proposition~\ref{prop:exactscale} applies word for word, with Lemma~\ref{lem:approxscaleU} in place of Lemma~\ref{lem:approxscale}: it uses only
		that $T$ is a metrizable Choquet simplex and that continuous affine targets can be approximated by dimension functions which are themselves continuous. It
		produces $a=\sum_i2^{-i}z_i$ with $z_i\ge0$, $\|z_i\|=1$, placed in pairwise orthogonal diagonal blocks. Then $2^{-i}=\|2^{-i}z_i\|\in\operatorname{sp}(a)$ for every
		$i$, so $0$ is an accumulation point of $\operatorname{sp}(a)$.
	\end{proof}
	
	\subsection{Proof of Theorem~\ref{thm:proj}}\label{ss:proofproj}
	
	\emph{The unital corner.} Let $p\in A\otimes\K$ be a nonzero projection and $E_0=p(A\otimes\K)p$. Then $E_0$ is unital, separable and simple (a hereditary
	subalgebra of the simple $C^*$-algebra $A\otimes\K$), nuclear (hereditary subalgebras of nuclear $C^*$-algebras are nuclear), and $\JS$-stable \cite[Corollary~3.1]{TW}. It is a
	full hereditary subalgebra of $A\otimes\K$, so $E_0\otimes\K\cong A\otimes\K$ by Theorem~\ref{thm:Brown}; in particular $E_0$ satisfies the UCT, which is invariant
	under stable isomorphism. It is stably finite. Indeed, $A\otimes\K$ contains no infinite projection: every projection of $A\otimes\K$ is equivalent to a
	projection $q\in M_n(A)$ for some $n$ \cite[Chapter~6]{RLL}, and if $q$ were equivalent to a subprojection $q'<q$, then $q'=qq'q$ and the implementing partial isometry
	$w=q'wq$ would lie in $M_n(A)$, and $w+(1-q)$ would be a non-unitary isometry in $M_n(\widetilde A)$, contradicting stable finiteness of $A$. The unit of
	$M_k(E_0)$, a corner of $A\otimes\K\otimes M_k\cong A\otimes\K$, is therefore a finite projection, so $M_k(E_0)$ is finite. So $E_0$ is a unital classifiable stably
	finite $C^*$-algebra.
	
	\emph{The unital model.} Apply the construction of Part~II to the invariant of $E_0$ (Case (U)). By the proof of Theorem~\ref{thm:unital}, $E_0\cong A_\infty\otimes\JS\cong
	C^*_r(\mathcal G)$ with $\mathcal G=\bar{\mathcal G}\times\mathcal G_\JS$, the isomorphism carrying $C(\mathcal G^{(0)})$ to $B\otimes C(\mathcal G_\JS^{(0)})$. Let
	$\mathcal P_\N$ be the pair groupoid of $\N$ and $\mathcal H=\mathcal G\times\mathcal P_\N$. By the general form of Lemma~\ref{lem:product},
	\[
	C^*_r(\mathcal H)\cong A_\infty\otimes\JS\otimes\K,\qquad C_0(\mathcal H^{(0)})=B\otimes C(\mathcal G_\JS^{(0)})\otimes c_0(\N),
	\]
	and $\mathcal H$ is a second countable, locally compact Hausdorff, principal étale groupoid. Fix an isomorphism $\Phi\colon A\otimes\K\to A_\infty\otimes\JS\otimes\K$ (which exists
	since $A\otimes\K\cong E_0\otimes\K\cong A_\infty\otimes\JS\otimes\K$), a strictly positive element $h_A\in A$, and put
	\[
	b=\Phi(h_A\otimes e_{11})\in(A_\infty\otimes\JS\otimes\K)_+ .
	\]
	Then $\her(b)=\Phi(\her(h_A\otimes e_{11}))=\Phi(A\otimes e_{11})\cong A$. Since $A$ is not unital, $0$ is an accumulation point of $\operatorname{sp}(b)$
	(Lemma~\ref{lem:unitalher}).
	
	\emph{The rank function.} By Lemma~\ref{lem:AZinvariant}, every tracial state of $A_\infty\otimes\JS$ has the form $\tau\otimes\tau_\JS$ with $\tau\in T(A_\infty)$. Define
	\[
	\gamma(\tau)=d_{(\tau\otimes\tau_\JS)\otimes\Tr}(b)\in(0,\infty]\qquad(\tau\in T(A_\infty)\cong T).
	\]
	It is affine, and lower semicontinuous as the supremum of the continuous functions $\tau\mapsto((\tau\otimes\tau_\JS)\otimes\Tr)(P_mf_\delta(b)P_m)$, $m\in\N$, $\delta>0$
	(with $P_m=1\otimes(e_{11}+\dots+e_{mm})$); and $\gamma>0$ because $b\neq0$ and the traces are faithful.
	
	\emph{The diagonal element.} Let $a\in(B\otimes c_0(\N))_+$ be given by Proposition~\ref{prop:exactscaleU} for $\gamma$, and let $a'$ be its image under the
	injective homomorphism $A_\infty\otimes\K\to A_\infty\otimes\JS\otimes\K$, $x\otimes y\mapsto x\otimes1_\JS\otimes y$. Then $a'\in C_0(\mathcal H^{(0)})_+$, $\operatorname{sp}(a')=
	\operatorname{sp}(a)$ accumulates at $0$, and, since functional calculus commutes with injective homomorphisms,
	\[
	d_{(\tau\otimes\tau_\JS)\otimes\Tr}(a')=d_{\tau\otimes\Tr}(a)=\gamma(\tau)=d_{(\tau\otimes\tau_\JS)\otimes\Tr}(b)\qquad(\tau\in T(A_\infty)).
	\]
	
	\emph{Comparison.} $A_\infty\otimes\JS$ is unital, simple, separable, nuclear (hence exact), $\JS$-stable and stably finite, with $T(A_\infty\otimes\JS)\cong T\neq\emptyset$
	(Lemma~\ref{lem:AZinvariant}). By Proposition~\ref{prop:soft}, $a'\sim b$ and $\her(a')\cong\her(b)\cong A$.
	
	\emph{The groupoid.} Let $U=\{a'>0\}\subset\mathcal H^{(0)}$ and $\mathcal G'=\mathcal H|_U$. By Lemma~\ref{lem:restriction}, $\mathcal G'$ is a second countable, locally
	compact Hausdorff, principal étale groupoid and $C^*_r(\mathcal G')=\her(a')\cong A$. If $U$ were compact, the indicator function of $U$ would be a unit of
	$C^*_r(\mathcal G')$; since $A$ is not unital, $U$ is not compact. \qed
	
	\begin{remark}[An untwisted range statement]\label{rem:projrange}
		The proof shows more. Let $(G_0,G_0^+,u,T,\rho,G_1)$ be the invariant of a unital classifiable stably finite $C^*$-algebra, and let $\gamma\colon T\to(0,\infty]$ be lower
		semicontinuous and affine. Then the hereditary subalgebra generated by the element $a'$ obtained from Proposition~\ref{prop:exactscaleU} is a non-unital classifiable
		$C^*$-algebra with $K$-theory $(G_0,G_1)$, trace cone $\R_{\ge0}T$ with pairing $\rho$, and scale $\gamma$, and it is the $C^*$-algebra of a principal étale groupoid. This is
		the untwisted analogue of Li's \cite[Corollary~1.6]{Li}. By Proposition~\ref{prop:soft}, the isomorphism class of such a hereditary subalgebra depends only on
		$\gamma$. It is never unital, and so it differs from the unital corners $q(A_\infty\otimes\JS\otimes\K)q$ even when $\gamma=\tau_*[q]$.
	\end{remark}
	
	\begin{remark}[The finite--infinite dichotomy]\label{rem:dichotomy}
		Every classifiable $C^*$-algebra $A$ is either stably finite or purely infinite. If $A$ is unital, this follows from R{\o}rdam's theorems: $A$ is purely
		infinite if it has no tracial state \cite[Theorem~4.5, Corollary~5.1]{Ro}, and stably finite otherwise (a faithful tracial state rules out non-unitary isometries
		in $M_n(A)$, as in the proof of Proposition~\ref{prop:firstlimit}). If $A$ is stably projectionless, it is stably finite: an isometry $v\in M_n(\widetilde A)$ has a
		unitary scalar part, so $1-vv^*$ is a projection in $M_n(A)$, hence $0$. Otherwise let $p\in A\otimes\K$ be a nonzero projection and $E_0=p(A\otimes\K)p$; as in the
		first paragraph of Section~\ref{ss:proofproj}, $E_0$ is unital, simple, separable, nuclear and $\JS$-stable, and $A\otimes\K\cong E_0\otimes\K$, so that $A$ is
		isomorphic to a hereditary subalgebra of $E_0\otimes\K$. If $E_0$ is purely infinite, so are $E_0\otimes\K$ and its nonzero hereditary subalgebras, among them $A$
		\cite{KR00}. If $E_0$ has a tracial state, it is stably finite, $E_0\otimes\K$ has stable rank one (Lemma~\ref{lem:sr1}), hence so do the matrix algebras over its
		unitization \cite{Rie}, and unital $C^*$-algebras of stable rank one are finite; $M_n(\widetilde A)$ is a unital subalgebra of $M_n((E_0\otimes\K)^\sim)$, hence finite.
	\end{remark}
	%======================================================================
	\appendix
	\section{Necessity of the pairing condition (P)}\label{app:P}
	
	In this appendix we give a self-contained proof that the pairing of a classifiable stably projectionless $C^*$-algebra with continuous scale contains no strictly positive
	function. The difficulty is that a class in $K_0$ is represented by two projections over the unitization, where strict comparison is not available at the character; we
	circumvent it by smearing the relative angle of the two projections with a diffuse element of $\JS$.
	
	\begin{theorem}\label{thm:P}
		Let $E$ be a simple, separable, exact, $\JS$-stable, stably projectionless $C^*$-algebra with $\Ped(E)=E$ and $T(E)$ nonempty and weak$^*$-compact. Then there is no $x\in K_0(E)$ with
		$\rho_E(x)(\tau)>0$ for all $\tau\in T(E)$.
	\end{theorem}
	
	We use the following facts.
	\begin{enumerate}[label=\textup{(F\arabic*)}]
		\item \emph{Strict comparison.} If $a,b\in M_m(E)_+$ and $d_\tau(a)<d_\tau(b)$ for all $\tau\in T(E)$, then $a\precsim b$. (Since $\Ped(E)=E$ every densely finite lsc trace is bounded;
		quasitraces are traces by exactness \cite{Haa,BK}; the statement is then \cite[Proposition~1.5]{CE}, which rests on R{\o}rdam's almost unperforation \cite[Theorem~4.5]{Ro}; see also
		\cite[Theorem~3.5]{GL1}.)
		\item \emph{Projections.} If $P,R$ are projections in a $C^*$-algebra $\mathcal M$ and $R\precsim P$, then $R$ is Murray--von Neumann equivalent to a subprojection of $P$. (Choose
		$r\in\mathcal M$ with $\|R-rPr^*\|<\frac12$. By \cite[Lemma~2.2]{KR} there is $d$ with $\frac12R=(R-\frac12)_+=d^*rPr^*d$; then $w=\sqrt2\,Pr^*d$ satisfies $w^*w=R$, and $Q:=ww^*$ is a
		projection in $P\mathcal MP$, so $Q\le P$.)
		\item \emph{Two projections.} The universal unital $C^*$-algebra $\mathcal U$ generated by two projections is $\{F\in C([0,\tfrac\pi2],M_2):F(0),F(\tfrac\pi2)\text{ diagonal}\}$, with
		$P_{\mathcal U}=\diag(1,0)$ and $R_{\mathcal U}(\theta)$ the projection onto $\C(\cos\theta,\sin\theta)$. Its irreducible representations are the two-dimensional ones at
		$\theta\in(0,\frac\pi2)$ and four characters $\chi_{11},\chi_{00}$ (at $\theta=0$, values $(P,R)=(1,1),(0,0)$) and $\chi_{10},\chi_{01}$ (at $\theta=\frac\pi2$, values $(1,0),(0,1)$). Every
		bounded tracial positive functional on $\mathcal U$ has the form $\omega(F)=\sum_\chi w_\chi\chi(F)+\int_{(0,\pi/2)}\Tr(F(\theta))\,d\mu(\theta)$ with $w_\chi\ge0$ and $\mu$ a finite measure
		(compare Lemma~\ref{lem:tracesblock}).
		\item \emph{Traces on $F\otimes\JS$}, $F$ unital, are $\sigma\otimes\tau_\JS$ with $\sigma\in T(F)$; $\JS\cong\JS^{\otimes\infty}$; $\JS$ is stably finite \cite{JS,TW}.
	\end{enumerate}
	
	\begin{proof}
		Suppose $\delta:=\min_{T(E)}\rho_E(x)>0$. Write $x=[p]-[r]$ with projections $p,r\in M_n(\widetilde E)$ whose scalar parts coincide (after conjugating by a scalar unitary):
		$\pi(p)=\pi(r)=r_0:=1_k\oplus0$. Then $(\tilde\tau\otimes\Tr)(p)-(\tilde\tau\otimes\Tr)(r)=\rho_E(x)(\tau)\ge\delta$ for $\tau\in T(E)$.
		
		\emph{Step 1.} Put $\mathcal M=M_n(\widetilde E\otimes\JS)$, $\mathcal J=M_n(E\otimes\JS)$ (an ideal with $\mathcal M/\mathcal J\cong M_n(\JS)$, quotient map $\pi_\JS$), $P=p\otimes1$,
		$R=r\otimes1$. Then $P-R\in\mathcal J$, $\pi_\JS(P)=\pi_\JS(R)=r_0\otimes1=:r_\JS$, and for $\omega=\tilde\tau\otimes\tau_\JS\otimes\Tr$ ($\tau\in T(E)$) we have
		$\omega(P)-\omega(R)=\rho_E(x)(\tau)\ge\delta$. Since $E\otimes\JS\cong E$, the ideal $\mathcal J$ is again a matrix algebra over a simple $C^*$-algebra, exact, $\JS$-stable, stably projectionless, with all
		densely finite lsc traces bounded and $T(E\otimes\JS)=\{\tau\otimes\tau_\JS\}\cong T(E)$.
		
		\emph{Step 2 (a diffuse element of $\JS$).} $\JS$ contains a unital copy of a dimension-drop algebra, hence of $C([0,1])$; let $h\in\JS$ correspond to the identity function. The law $\nu$ of
		$h$ under $\tau_\JS$ is not a point mass (the trace is faithful). In $\JS^{\otimes\infty}\cong\JS$ let $h_j$ be the copy of $h$ in the $j$th factor; the $h_j$ commute and their joint law
		under the (product) trace is $\nu^{\otimes\N}$. Put $z=\sum_j2^{-j}h_j$. Its law $\nu_z$ is that of $Y=\sum_j2^{-j}X_j$ with $X_j$ independent, $X_j\sim\nu$; it has no atoms. (If $a>0$
		were the largest atom of $Y$ and $s$ a point with $\mathbb P(Y=s)=a$, write $Y=(X_1+Y')/2$ with $Y'\sim Y$ independent of $X_1$; then
		$a=\mathbb P(Y=s)=\sum_{t}\nu(\{t\})\,\mathbb P(Y'=2s-t)\le a\sum_t\nu(\{t\})\le a$, the sum over the atoms $t$ of $\nu$. Hence $\nu$ is purely atomic and $\mathbb P(Y'=2s-t)=a$ for all atoms
		$t$; so the finite set $M$ of points of maximal mass satisfies $2s-t\in M$ for $s\in M$ and every atom $t$. Taking $s=\max M$ gives $s\le t$ for all atoms, taking $s=\min M$ gives
		$s\ge t$; so $\nu$ has a single atom, a contradiction.) Consequently $\kappa(\ell):=\sup\{\nu_z(I):I\text{ an interval of length }\ell\}\to0$ as $\ell\to0$.
		
		\emph{Step 3 (smearing the angle).} Let $D_0=C^*(1,p,r)\subset M_n(\widetilde E)$, a quotient of $\mathcal U$, and $\mathcal C=C^*(D_0\otimes1,1\otimes z)\subset\mathcal M$, a quotient of
		$\mathcal U\otimes C([0,1])$. Fix $0<\theta_1<\theta_1'<\theta_2'<\theta_2<\frac\pi2$ and $L>0$ with $\theta_2+L<\frac\pi2$ and $\theta_2'-\theta_1'>L$, and
		$\beta\in C([0,\frac\pi2],[0,1])$ with $\beta=1$ on $[\theta_1',\theta_2']$ and $\beta=0$ off $(\theta_1,\theta_2)$. With the rotation matrices
		$\mathrm{Rot}(\alpha)=\begin{pmatrix}\cos\alpha&-\sin\alpha\\ \sin\alpha&\cos\alpha\end{pmatrix}$, the function $V_{\mathcal U}(\theta,u)=\mathrm{Rot}(\beta(\theta)Lu)$ is a unitary of
		$\mathcal U\otimes C([0,1])$ (it is $1$ at $\theta\in\{0,\frac\pi2\}$); let $V\in\mathcal C$ be its image and $R'=VRV^*$. Since $\pi_\JS(P)=\pi_\JS(R)$ and $P_{\mathcal U}-R_{\mathcal U}$
		generates the ideal $\{F:F(0)=0\}$ of $\mathcal U$, the restriction of $\pi_\JS$ to $\mathcal C$ factors through the characters at $\theta=0$, where $V_{\mathcal U}=1$; hence
		$\pi_\JS(R')=r_\JS$ and $P-R'\in\mathcal J$. Moreover $\omega(R')=\omega(R)$ for every bounded trace $\omega$ of $\mathcal M$, $V$ being a unitary of $\mathcal M$. Conjugation by
		$\mathrm{Rot}(\alpha)$ maps $R_{\mathcal U}(\theta)$ to $R_{\mathcal U}(\theta+\alpha)$, while $P$ is not conjugated and $P_{\mathcal U}$ is constant; so the pair $(P,R')$ is the image
		of $(P_{\mathcal U},R_{\mathcal U})$ under the homomorphism $F\mapsto F\circ\Phi_L$, $\Phi_L(\theta,u)=\theta+\beta(\theta)Lu$ on the two-dimensional part, characters fixed. For
		$\omega=\tilde\tau\otimes\tau_\JS\otimes\Tr$, the restriction of $\omega$ to $\mathcal C$ is the product of $\omega|_{D_0\otimes1}$ and $\tau_\JS|_{C^*(z)}$. Hence, writing $(w_\chi,\mu)$ for
		the data (F3) of $\omega$ on $C^*(1,P,R)$, the data of $\omega$ on $C^*(1,P,R')$ are $(w_\chi,\mu')$ with $\mu'=(\Phi_L)_*(\mu\otimes\nu_z)$. In particular
		$w_{10}-w_{01}=\omega(P)-\omega(R')\ge\delta$ (two-dimensional representations and $\chi_{11},\chi_{00}$ contribute equally to $P$ and $R'$), and for every interval
		$W\subset(\theta_1'+L,\theta_2')$,
		\[
		\mu'(W)=\int\nu_z\{u:\theta+\beta(\theta)Lu\in W\}\,d\mu(\theta)\le\kappa(|W|/L)\,\mu\big((0,\tfrac\pi2)\big)\le n\,\kappa(|W|/L),
		\]
		uniformly in $\tau$ (only $\theta\in[\theta_1',\theta_2']$ can be moved into $W$, and then $\beta(\theta)=1$). Fix such a $W$ with $n\kappa(|W|/L)<\delta/2$.
		
		\emph{Step 4 (comparison).} Put $S=R'PR'$, $S'=PR'P$; with $x_0=R'P$, $S=x_0x_0^*$ and $S'=x_0^*x_0$. In a two-dimensional representation at angle $\theta'$, $S=\lambda R_{\mathcal U}$,
		$S'=\lambda P_{\mathcal U}$ with $\lambda=\cos^2\theta'$; at $\chi_{11}$, $S=S'=1$; elsewhere among the characters $S=S'=0$, with $R'=1$ at $\chi_{01}$ and $P=1$ at $\chi_{10}$. For
		$g\in C([0,1])$ let $g(S)$ be computed in the unital $C^*$-algebra $R'\mathcal MR'$ and $g(S')$ in $P\mathcal MP$. Let $[\lambda_a,\lambda_c]$ be an interval of $\lambda$-values corresponding to a
		subinterval of $W$, choose $\lambda_a<\lambda_b\le\lambda_b'<\lambda_c$ and continuous functions $\varphi_0,\varphi_1,\psi\colon[0,1]\to[0,1]$ with $\varphi_0=1$ on $[0,\lambda_a]$,
		$\varphi_0=0$ on $[\lambda_b,1]$; $\varphi_1=0$ on $[0,\lambda_b']$, $\varphi_1=1$ on $[\lambda_c,1]$; $\psi$ supported in $[\lambda_a,\lambda_c]$; and $\varphi_0+\varphi_1+\psi>0$ on $[0,1]$.
		Since $(\varphi_0+\varphi_1+\psi)(S)$ is invertible in $R'\mathcal MR'$, and $g(S)\sim g(S')$ (Cuntz equivalence) for $g\in C_0((0,1])$,
		\[
		\langle R'\rangle\le\langle\varphi_1(S)\rangle+\langle\psi(S)\rangle+\langle\varphi_0(S)\rangle=\langle\varphi_1(S')\rangle+\langle\psi(S')\rangle+\langle\varphi_0(S)\rangle .
		\]
		As $\pi_\JS(S)=\pi_\JS(S')=r_\JS$ and $\varphi_0(1)=\psi(1)=0$, the elements $\varphi_0(S)$, $\varphi_0(S')$, $\psi(S')$ lie in $\mathcal J$. By (F3), for $\omega=\tilde\tau\otimes\tau_\JS\otimes\Tr$,
		\[
		d_\omega(\varphi_0(S))=w_{01}+\mu'(\{\varphi_0>0\}),\quad d_\omega(\varphi_0(S'))=w_{10}+\mu'(\{\varphi_0>0\}),\quad d_\omega(\psi(S'))\le\mu'(W)<\delta/2,
		\]
		where $\{\varphi_0>0\}$ denotes the set of angles $\theta'$ with $\varphi_0(\cos^2\theta')>0$. Hence $d_\omega(\psi(S')\oplus\varphi_0(S))<w_{01}+\delta/2+\mu'(\{\varphi_0>0\})<d_\omega(\varphi_0(S'))$
		for every $\tau\in T(E)$. Up to positive scalar multiples, the restrictions of these traces $\omega$ to $\mathcal J=M_n(E\otimes\JS)$ are all the nonzero densely finite lower
		semicontinuous traces of $\mathcal J$: such a trace is $\sigma\otimes\Tr$ for a densely finite lower semicontinuous trace $\sigma$ on $E\otimes\JS$ (Lemma~\ref{lem:hertraces}),
		$\sigma$ is bounded since $\Ped(E\otimes\JS)=E\otimes\JS$, and the tracial states of $E\otimes\JS$ are the $\tau\otimes\tau_\JS$, $\tau\in T(E)$ (Step~1). Hence (F1), applied in
		$\mathcal J$ with $E\otimes\JS$ in place of $E$, gives $\psi(S')\oplus\varphi_0(S)\precsim\varphi_0(S')$. Since $\varphi_1\varphi_0=0$ and $0\le\varphi_1+\varphi_0\le1$, we have
		$0\le\varphi_1(S')+\varphi_0(S')\le P$, and
		\[
		\langle R'\rangle\le\langle\varphi_1(S')\rangle+\langle\varphi_0(S')\rangle=\langle\varphi_1(S')+\varphi_0(S')\rangle\le\langle P\rangle .
		\]
		
		\emph{Step 5.} By (F2), $R'\sim Q\le P$ for a projection $Q\in\mathcal M$. Then $\pi_\JS(Q)\le r_\JS$ and $\pi_\JS(Q)\sim\pi_\JS(R')=r_\JS$ in the finite $C^*$-algebra $M_n(\JS)$ (F4), so
		$\pi_\JS(Q)=r_\JS$ and $P-Q$ is a projection in $\mathcal J=M_n(E\otimes\JS)$. Bounded traces are invariant under Murray--von Neumann equivalence, so for
		$\omega=\tilde\tau\otimes\tau_\JS\otimes\Tr$
		\[
		\omega(P-Q)=\omega(P)-\omega(R')=\omega(P)-\omega(R)\ge\delta>0 .
		\]
		Hence $P-Q\ne0$, contradicting stable projectionlessness of $E\otimes\JS\cong E$.
	\end{proof}
	
	\begin{remark}
		Without the smearing in Step~3 one would need an angle window carrying uniformly small mass for all traces, which can fail: the mass of the two-dimensional representations may
		accumulate at angles approaching $\pi/2$ along a sequence of traces.
	\end{remark}
	%======================================================================

\end{document}